\documentclass[fleqn,10pt]{article}
\usepackage{latexsym, graphicx, epsfig, amsmath, amssymb,amsfonts}
\usepackage{natbib,amsthm,version}
\usepackage{amsbsy,bm,multirow}
\usepackage[titletoc,page]{appendix}
\usepackage[mathscr]{eucal}
\usepackage{mathtools}
\usepackage{color}
\usepackage[utf8]{inputenc}
\usepackage[english]{babel}
\usepackage{amsthm}
\usepackage[hidelinks]{hyperref}
\usepackage{url}
\usepackage{subfigure}
\usepackage[lined,boxed,linesnumbered,ruled]{algorithm2e}
\usepackage{float}
\usepackage{arydshln}
\usepackage{dsfont,bbm}

\usepackage[shortlabels]{enumitem}
\setlist{nolistsep}

\usepackage{natbib}
\newcommand{\mbs}[1]{\mathbf{#1}}
\newcommand{\mbb}[1]{\mathbb{#1}}

\newtheorem{theorem}{Theorem}[section]

\newtheorem{lemma}[theorem]{Lemma}
\newtheorem{remark}{Remark}[section]

\theoremstyle{definition}
\newtheorem{definition}{Definition}[section]
\newtheorem{assumption}{Assumption}[section]

\title{
  Learning the 
  Exact Time Integration Algorithm for Initial Value Problems
  by Randomized Neural Networks
} 
\author{
  Suchuan Dong\thanks{Author of correspondence. Email: sdong@purdue.edu}, \
  Naxian Ni
  \\
  Center for Computational and Applied Mathematics \\
  Department of Mathematics \\
  Purdue University, USA 
 } 

\date{(February 15, 2025)}

\begin{document}
\maketitle


\begin{abstract}

We present a method leveraging extreme learning machine (ELM) type
randomized neural networks (NNs) for learning the exact time integration
algorithm for initial value problems. The exact time integration
algorithm for non-autonomous (including autonomous) systems
can be represented by an algorithmic function in higher dimensions,
which satisfies an associated system of partial differential equations with
corresponding boundary conditions.
Our method learns the algorithmic function by solving this associated system
using ELM with a physics informed approach.
The trained ELM network serves as the learned algorithm and can be
used to solve the initial value problem with arbitrary initial data
or arbitrary step sizes from some domain.
When the right hand side of the non-autonomous system exhibits a
periodicity with respect to any of its arguments, while 
the solution itself to the problem is not periodic, we show that in this case
the algorithmic function is either periodic, or when it is not,
satisfies a well-defined relation for different periods.
This property can greatly simplify the network training for algorithm
learning  in many problems.
We consider explicit and implicit NN formulations, leading to
explicit and implicit time integration algorithms, and discuss
how to train the ELM network by the nonlinear least squares method.
Extensive numerical experiments with benchmark problems,
including non-stiff, stiff and chaotic systems, show that the learned NN
algorithm produces highly accurate solutions in long-time simulations,
with its time-marching errors decreasing nearly exponentially with increasing
degrees of freedom  in the neural
network. We compare extensively the computational performance (time-marching accuracy
vs.~time-marching cost) between the current NN algorithm and
the leading traditional time integration algorithms. The learned NN
algorithm is computationally competitive, markedly outperforming
the traditional algorithms in many problems.

\end{abstract}


\vspace{0.05cm}
Keywords: {\em
  extreme learning machine,
  randomized neural network,
  scientific machine learning,
  nonlinear least squares,
  exact time integration,
  exact flow map
}



\section{Introduction}
\label{sec_intro}


%
%


This work concerns the development of accurate and efficient 
methods for solving initial value problems (IVPs), by learning
their exact time integration algorithm with artificial neural networks.
Initial value problems with systems of
ordinary differential equations (ODEs)
are ubiquitous in modeling natural phenomena and
engineered systems. Problems as diverse as population dynamics,
the spread of diseases, chemical reactions, celestial mechanics, dynamical systems,
biological pattern formation, and the economic markets
can be modeled by such systems. They also arise, upon spatial discretization,
from time-dependent partial differential equations (PDEs)
describing the evolution of physical systems
based on first principles such as the conservation laws and
thermodynamic principles~\cite{GrootM1984,GurtinFA2010,Dong2014,Dong2018,YangD2020}.
The models arising from practical applications are usually nonlinear and can rarely
be solved by analytical methods. Numerical simulations are therefore
of fundamental importance and play a critical role
in understanding the dynamics of these systems.
%
Time marching algorithms for IVPs are traditionally based on Taylor series
expansions.
Research in the past decades has led to
accurate and robust methods with provable convergence
based on this approach~\cite{CellierK2006,HairerNW1993,HairerW1996}, which have become
the cornerstone in computational science and engineering.

The emergence of artificial neural networks~\cite{GoodfellowBC2016} and their application
in scientific computing, a.k.a.~scientific machine
learning (SciML)~\cite{Karniadakisetal2021,RaissiPK2019,SirignanoS2018,EY2018,TangWL2021,DongN2021,KrishnapriyanGZKM2021,DongY2022,WanW2020,Penwardenetal2023,QianZHD2023},
in recent years
have stimulated promising new approaches
for simulating dynamical systems~\cite{QinWX2019,RaissiPK2019,SirignanoS2018,Jietal2021,LiuKB2022,WangP2023,FlamantPS2020,Dellnitzetal2023,MoyaL2023}.
%
NN-based methods for such problems can be broadly classified into two types,
as pointed out in~\cite{Legaardetal2023}, direct solution models
and time-stepper models. Direct solution models seek the solution to the IVP
directly, in which the network input denotes the independent variable and
the output represents the solution. Representative techniques in
this category include the physics-informed neural network (PINN)
method~\cite{RaissiPK2019}, deep Galerkin method (DGM)~\cite{SirignanoS2018},
and their many extensions and variants; see~\cite{Karniadakisetal2021,Cuomo2022Scientific}
for a review of related techniques.
Because the trained NN corresponds to a specific input or initial
condition, the direct solution  model needs to be re-trained
if the initial data changes. 
The time-stepper models, on the other hand, follow an approach
analogous to traditional numerical solvers. Given the current state of
the system, the model attempts to compute the state at some point
in the future. Accordingly, the NNs for time-stepper models
are generally auto-regressive in nature, aiming to capture
the dynamics of the system. The trained time-stepper model can
handle different initial data or input, without the need for re-training.
Early works on time-stepper models have focused on the
data-driven modeling of dynamical systems, in which
the governing equations are
unknown and are approximated by the solution
trajectory data~\cite{QinWX2019,QinCJX2021,ChurchillX2022}.
A residual network (ResNet) structure and stacked ResNets with
a recurrent or recursive configuration
are proposed in~\cite{QinWX2019} to approximate the unknown governing
equations and learn a discrete flow map of the autonomous systems.
This technique is extended to non-autonomous systems in~\cite{QinCJX2021},
where the non-autonomous term is expanded in terms of a local set of bases,
and to cases with partially observed  state
variables~\cite{ChurchillX2022}; see recent related works
in~\cite{ChenW2023,LuT2024}.
A hierarchical time-stepper based on deep neural networks (DNN)
is developed in~\cite{LiuKB2022} to approximate the system flow map
 over a range of time scales,
in which multiple DNNs are trained corresponding to
a number of step sizes using the trajectory data;
see~\cite{HamidRNB2023} for an adaptive time-stepping scheme
with the hierarchical time-stepper.
%
In~\cite{Dellnitzetal2023} deep reinforcement learning has been
combined with classical numerical solvers to determine the step sizes in adaptive
time-stepping.
%
In~\cite{ZhangSK2022} NN  structures are designed  to
satisfy the symmetric degeneracy conditions in the GENERIC formalism
for simulating dynamical systems, and their universal approximation
property has been established.
%
We also refer the reader to~\cite{ChenRBD2019,Kimetal2021,FronkP2024a}
(among others) for neural ODEs
and to~\cite{Legaardetal2023} for a survey of other related techniques.


Given a non-autonomous system of ODEs, we seek the exact
time integration  algorithm for this system. 
An exact time integration (or time marching) algorithm refers to a discrete scheme
that can produce numerical solutions matching exactly  the true solutions
to the system for all step sizes within a range; see Section~\ref{sec_method}
for its definition. This type of algorithms is also
known as the exact finite difference schemes or
non-standard finite difference (NSFD) schemes~\cite{Mickens2021}.
The exact finite difference schemes for a number of problems
have been constructed explicitly in~\cite{Mickens2021}, by using denominator functions
for approximating discrete derivatives and by special treatments of the nonlinear
terms involved in the equations. However, 
constructing the exact discrete scheme for
an arbitrary given system of differential equations is virtually impossible, unfortunately,
because otherwise this will be tantamount to knowing the general forms of
the solution to the given system~\cite{Mickens2021}.

This motivates the current work for learning the exact time integration algorithms
through artificial neural networks, by leveraging the NN's universal approximation
power for function approximation.
%
Viewed from the dynamical systems perspective, the exact time integration
algorithm is closely related to the exact flow map (or evolution semigroup)~\cite{StuartH1996}
of the system. Given a system of differential equations, suppose
the system has a unique solution for any given initial condition.
Then the exact time integration algorithm  exists, and
we would like to learn this algorithm from the given system. The learned algorithm
then provides a method for solving this system, under any given initial condition
or step size.
The problem considered here is different from the so-called flow-map learning in
previous studies (see e.g.~\cite{QinWX2019,ChurchillX2023,LiuKB2022}, among others),
where the governing differential equations are unknown and a time-stepper model
is learned based on the solution trajectory data obtained from
either measurements or external numerical solvers.
In contrast, for the problems considered in this paper,
the system of differential equations is given and we seek
the exact time integration algorithm for this system. No other data
about the system is available,  and our method
does not rely on any external numerical solver.

%

Given a system of non-autonomous (including autonomous) ODEs,
the exact time integration algorithm for this system can be represented by a vector-valued
function in higher dimensions, which in the current paper will be referred to as
the algorithmic function of the given system.
The algorithmic function satisfies an associated system of
partial differential equations, together with corresponding boundary conditions.
We learn the exact time integration algorithm  by
solving the associated system of PDEs  for
the  algorithmic function with a physics informed approach,
leveraging a type of randomized feedforward neural networks
often known as extreme learning machines (ELMs)~\cite{DongL2021,DongY2022rm,WangD2024}.


Randomized NNs (or random-weight NNs)
are neural networks 
with a subset of the network parameters assigned to random values and fixed (non-trainable)
while only the rest of the network parameters are trained.
A simple strategy underlies randomized NNs. Since optimizing
the entire set of NN parameters is extremely difficult and  costly,
randomized NN attempts to randomly assign and fix a portion of the network parameters,
so that the ensuing network training problem can become simpler,
without severely compromising the network's achievable approximation
capability~\cite{DongY2022rm,NiD2023}.
ELM is a particular type of randomized feedforward NNs, in which all the hidden-layer
parameters are assigned to random values and fixed (non-trainable) while
only the output-layer parameters are trained. ELM was originally proposed
in~\cite{HuangZS2006} for linear classification and regression problems,
but has since found widespread applications  in many areas~\cite{Alabaetal2019}.
A close cousin to ELM is the 
random vector functional link (RVFL) networks~\cite{PaoPS1994}.
ELM type  randomized NNs, with a single hidden layer,
are known to be universal function
approximators~\cite{IgelnikP1995,HuangCS2006,Gonon2023}.


The use of ELM type randomized NNs  for scientific computing,
particularly for solving partial differential equations, occurs quite recently
and has flourished in the past few
years~\cite{DongL2021,Schiassietal2021,CalabroFS2021,WangD2024}.
We refer the reader to~\cite{PanghalK2020,DwivediS2020,DongL2021,LiuHWC2021,CalabroFS2021,ChenCEY2022,LiLX2023,QuanH2023,Calabroetal2023,SunDF2024} (among others) for ELM type methods
for linear PDEs, to~\cite{DongL2021,DongL2021bip,Schiassietal2021,FabianiCRS2021,DongY2022rm,NiD2023,DongW2023,WangD2024} (among others) for nonlinear PDEs,
and to a recent work~\cite{FabianiKSY2025} for neural operators.
%
Many studies (including those from our group) have shown that
ELM  can produce highly accurate solutions to linear and nonlinear
PDEs, exhibiting spectral-like accuracy~\cite{DongL2021,NiD2023},
with a competitive computational cost~\cite{DongY2022rm}.
The ELM errors decrease exponentially
with increasing number of degrees of freedom for smooth solutions, and
can reach the level near machine accuracy as the degrees of freedom
become large~\cite{DongY2022rm}. 
%
Several studies have looked into the computation of ODEs
based on randomized NNs~\cite{YangHL2018,LiuXWL2020,PanghalK2020,DongL2021,Schiassietal2022,FabianiGRS2023,FlorioSCF2023}.
These studies stem from the direct solution model, in which the NN
learns the solution corresponding to a specific  initial condition.
As a result, with different initial data,
the NN  has to be re-trained.
In contrast, with the method herein,
the trained NN  can solve the problem with arbitrary initial data from some
domain, with no need for re-training. This is a major difference between
these studies and the  current work.


In this paper we present a method based on ELM-type randomized NNs
for learning the exact time integration algorithm of non-autonomous and autonomous
systems. Our method learns the 
algorithmic function using ELM by solving its associated
system of partial differential equations  with a physics informed approach.
We explicitly incorporate an approximation of the solution field 
into the NN formulation, so that the ELM network effectively learns
the corresponding error function.
We consider  explicit and implicit NN formulations and discuss how to
train the ELM network by the nonlinear least squares (Gauss-Newton) method.
Accordingly, the trained NN
gives rise to explicit or implicit time integration algorithms for solving
the system.
%
In particular, we investigate the effect on the algorithm learning
when the right hand side (RHS) of the non-autonomous system
exhibits a periodicity with respect to time or to
any component of the solution variable.
We show that in this case, while the solution itself to the IVP is not periodic,
the algorithmic function  is either periodic, or when it is not, shows
a well-defined relation for different periods.
This means that, whenever the RHS of the non-autonomous system exhibits a
periodicity with respect to any of its arguments, the NN  only
needs to be trained on a domain covering one period along the corresponding
directions, and the trained NN can be used to solve the problem
on arbitrarily long time horizons (in case of temporal periodicity of RHS)
or for those corresponding solution components taking arbitrary values on
the real axis. These properties can greatly simplify the network training
and algorithm learning  for many problems.
%
In addition, we find that it is often crucial to incorporate domain decomposition
into the algorithm learning, especially for stiff problems.
In this case, a local ELM-type randomized NN learns the algorithmic function
on each sub-domain, and the local NNs  collectively represent
the algorithmic function over the entire  training domain.
What is most attractive lies in that different local NNs are not coupled,
due to the nature of the PDE system for the algorithmic function,
and thus they can be trained completely separately. 


We present extensive numerical experiments with several benchmark
problems (including non-stiff, stiff, and chaotic systems) to evaluate
the performance of the learned NN  algorithm.
We show that the learned algorithm produces highly accurate solutions
in long-time simulations. Its time-marching errors
decrease nearly exponentially as the number of
degrees of freedom (training collocation points, or
training parameters) in the NN  increases, while its time-marching
cost grows only quasi-linearly.
We compare extensively the computational performance (time-marching accuracy
vs.~time-marching cost) between the current NN algorithm
and the leading traditional time integration algorithms (with 
adaptive time-stepping and adaptive-order),
including DOP853 (Dormand Prince, order 8)~\cite{HairerNW1993},
RK45 (explicit Runge-Kutta, order 5)~\cite{DormandP1980},
and Radau (implicit Runge-Kutta Radau IIA, order 5)~\cite{HairerW1996}, among others.
The learned NN algorithm is computationally very competitive,
markedly outperforming the traditional algorithms in almost all problems.



The contributions of this paper lie in three aspects:
(i) the ELM-based method for learning the exact time integration algorithms
and the learned NN algorithms for accurately and efficiently solving
non-autonomous and autonomous systems;
(ii) the illumination of properties of the algorithmic function when the RHS
of non-autonomous (including autonomous) systems exhibits a periodicity
with respect to any of its arguments, which can be used to greatly
simplify the network training and algorithm learning; and
(iii) the comparison with traditional time integration algorithms
and  demonstration of the competitive and superior computational
performance of the learned NN algorithms. 


The method presented here has been implemented in Python
using the Tensorflow and Keras libraries. The numerical tests are performed
on a MAC computer (3.2GHz Quad-Core Intel Core i5 processor, 24GB memory)
in the authors' institution.


The rest of this paper is organized as follows.
In Section~\ref{sec_method} we discuss the properties of the algorithmic function,
its representation by ELM-type randomized NNs,
the network training by the nonlinear least squares method,
and how to use the learned algorithm
to solve non-autonomous and autonomous systems.
Section~\ref{sec_tests} consists of a comprehensive set of numerical tests
to assess the performance of the presented method with several benchmark problems.
We compare extensively the learned NN algorithms
with leading traditional time integration algorithms.
Section~\ref{sec_summary} provides further comments on several aspects of the current
method to  conclude the presentation.
The appendix includes proofs of several theorems from Section~\ref{sec_method}
concerning the properties of the algorithmic function.

\section{Learning the  Exact Time Integration Algorithm}
\label{sec_method}

\subsection{Exact Time Marching Scheme}


Let $n$ be a positive integer, and we
consider the initial value problem on $t\in(a,b)$,
\begin{subequations}\label{eq_1}
\begin{align}
  & \frac{dy}{dt} = f(y,t),  \label{eq_1a} \\
  & y(t_0) = y_0,  \label{eq_1b}
\end{align}
\end{subequations}
where $t$ denotes the time, $t_0\in(a,b)$ is
the initial time, 
$y: (a,b)\subset \mbb R\rightarrow \mbb R^n$
denotes the solution, $f: \mbb R^n\times\mbb R \rightarrow \mbb R^n$
is a prescribed function, and $y_0\in \mbb R^n$ denotes the initial data.

To numerically solve this system, a time integration algorithm
produces a sequence $\{ y_k\in\mbb R^n\ |\ \mathrm{integer}\ k\geqslant 0 \}$
as approximations to the solution $y(t)$  at
discrete instants, $\{y(t_k)\ |\ t_k=t_0+kh,\ k\geqslant 0 \}$,
where $h$ denotes the step size. Interestingly, there exist discrete schemes
that can produce  exact results to the system~\eqref{eq_1}.
The forms of such schemes for several problems
have been explicitly constructed  in~\cite{Mickens2021}.
Following~\cite{Mickens2021}, we define
an exact time integration (or time marching) algorithm as follows,

\begin{definition}
  A time marching algorithm is said to be exact if it produces results
  identical to the exact solution for arbitrary values of $h\in(0,h_{\max}]$
    ($h_{\max}>0$ denoting some constant),
    i.e.~$y_k=y(t_k)$ for $k\geqslant 0$.
\end{definition}


Suppose $f(y,t)$ is continuous and satisfies the Lipschitz condition
$\|f(y,t)-f(z,t) \|\leqslant \lambda \|y-z \|$ for some constant $\lambda$.
Then the system~\eqref{eq_1} has a unique solution within a neighborhood
of $t_0$~\cite{HairerNW1993}. We re-write this solution as,
\begin{equation}\label{eq_2}
  y(t) = \psi(y_0,t_0,t), \quad t\in[t_0-\delta, t_0+\delta]
\end{equation}
for some $\delta>0$, where
$  
  \psi(y_0,t_0,t_0) = y_0
$  
and the dependence of the solution on $y_0$ and $t_0$ has been made explicit.
The existence of the exact time marching scheme for~\eqref{eq_1}
is established in~\cite{StuartH1996,Mickens2021}, given by the following result
(see~pages 61-62 of~\cite{Mickens2021}).
\begin{theorem}(\cite{Mickens2021})\label{thm_1}
  The system~\eqref{eq_1} has an exact time integration scheme given by
  \begin{equation}\label{eq_4}
    y_{k+1} = \psi(y_k,t_k,t_k+h), \quad k\geqslant 0,
  \end{equation}
  where $\psi$ is given in~\eqref{eq_2} and $h$ is the step size.
\end{theorem}

This result indicates that by acquiring the function
$\psi(y_0,t_0,t)$, for $y_0\in\Omega\subset\mbb R^n$,
$t_0\in(a,b)\subset\mbb R$ and $t\in[t_0-\delta,t_0+\delta]$,
one can attain the exact time integration algorithm~\eqref{eq_4}
for solving~\eqref{eq_1}.
In the dynamical systems parlance, this function is often referred to
as the evolution map (or flow map, evolution semigroup)~\cite{StuartH1996}.
This is a high-dimensional vector-valued function, and
is unknown for an arbitrary given
$f(y,t)$. Attaining $\psi$ is in general even more challenging than
solving the original system~\eqref{eq_1}.
In the current work we use randomized NNs
to learn the function $\psi(y_0,t_0,t)$, and the trained network
provides an approximation to the exact time integration algorithm.
The trained NN can be used as
a time marching algorithm for solving~\eqref{eq_1},
with arbitrary initial data $(y_0,t_0)\in\Omega\times(a,b)$ and arbitrary
step size $h\in(0,h_{\max}]$.
It is crucial to note that it is only necessary to learn
$\psi(y_0,t_0,t)$ for $t$ within a neighborhood of $t_0$.

Specifically, we will learn the function $\psi(y_0,t_0,t)$
by ELM-type randomized neural networks~\cite{DongL2021,DongY2022rm,NiD2023,WangD2024},
leveraging ELM's universal approximation capability~\cite{IgelnikP1995,HuangCS2006},
high accuracy, faster training,
and effectiveness for function approximation in high
dimensions~\cite{DongL2021,WangD2024}.
It is necessary to distinguish non-autonomous and autonomous systems
when learning the exact time integration algorithm.
We will first discuss 
general non-autonomous systems, and then restrict our attention
to autonomous systems.

\subsection{Learning Exact Time Marching Algorithm for Non-Autonomous Systems}

%

Let us consider the general non-autonomous system~\eqref{eq_1} and
discuss how to learn the function $\psi$  by ELM.
Since it is only necessary to learn $\psi$ for $t$ within some neighborhood of $t_0$,
we define
\begin{equation}\label{eq_6}
  \xi = t-t_0, \quad y(t) = y(t_0+\xi) = Y(\xi),
\end{equation}
where $y(t)$ is the solution to~\eqref{eq_1}.
The system~\eqref{eq_1} can be reformulated  in terms of $Y(\xi)$ into,
\begin{subequations}\label{eq_7}
  \begin{align}
    & \frac{dY}{d\xi} = f(Y, t_0+\xi), \\
    & Y(0) = y_0.
  \end{align}
\end{subequations}
Since the solution to  system~\eqref{eq_7} depends on $y_0$, $t_0$ and $\xi$,
we re-write it as
\begin{equation}
  Y(\xi) = \psi(y_0,t_0,\xi),
\end{equation}
where  $\psi: \mbb R^n\times\mbb R\times \mbb R \rightarrow \mbb R^n$
is the function we are pursuing here
and provides the exact time marching algorithm.
Henceforth we will refer to $\psi(y_0,t_0,\xi)$ as the algorithmic
function for system~\eqref{eq_1}.

In light of~\eqref{eq_7}, the algorithmic function $\psi(y_0,t_0,\xi)$ is determined by the following
partial differential equations,
\begin{subequations}\label{eq_9}
  \begin{align}
    & \frac{\partial\psi}{\partial\xi} = f(\psi(y_0,t_0,\xi), t_0+\xi), \label{eq_9a} \\
    & \psi(y_0,t_0,0) = y_0. \label{eq_9b}
  \end{align}
\end{subequations}
Our objective  is to determine the function $\psi(y_0,t_0,\xi)$, for
$y_0\in\Omega\subset\mbb R^n$, $t_0\in[T_0,T_f]\subset\mbb R$,
and $\xi\in[0,h_{\max}]\subset\mbb R$
with prescribed $\Omega$, $T_0$, $T_f$ and $h_{\max}$, such that
the system~\eqref{eq_9} is satisfied.
$\psi(y_0,t_0,\xi)$ is an $n$-vector valued function of $(n+2)$ variables,
where $n$ denotes the dimension of system~\eqref{eq_1}.

\subsubsection{Property of Algorithmic Function $\psi(y_0,t_0,\xi)$ When $f(y,t)$
  Exhibits Some Periodicity}
\label{sec_a221}

We next discuss a useful property of the algorithmic
function $\psi(y_0,t_0,\xi)$
when the RHS of~\eqref{eq_1a}, $f(y,t)$,
exhibits a periodicity with respect to $t$ or  to
some components of $y$.
This property can be used to simplify the NN training
for learning  $\psi(y_0,t_0,\xi)$ when applicable.
In this subsection we assume that $(y,t)\in\mbb R^n\times\mbb R$ in
problem~\eqref{eq_1}
and $\xi\in[0,h_{\max}]$ for a prescribed constant $h_{\max}>0$.
We further make the following assumption:
\begin{assumption}\label{ass_1}
  The function $f(y,t)$, with $(y,t)\in\mbb R^n\times\mbb R$,
  is
  such that the problem~\eqref{eq_1}
  has a unique solution for all $t\in\mbb R$, with
  arbitrary $(y_0,t_0)\in\mbb R^n\times\mbb R$.
  
\end{assumption}
\noindent This assumption can be satisfied, e.g.~if $f(y,t)$ is globally
Lipschitz on $\mbb R^n\times\mbb R$ (see Theorem 2.1.3 in~\cite{StuartH1996}).
Under this assumption, it can be noted that problem~\eqref{eq_9} has a unique
solution $\psi(y_0,t_0,\xi)$ for all
$(y_0,t_0,\xi)\in\mbb R^n\times\mbb R\times[0,h_{\max}]$.

We look into the effect on $\psi(y_0,t_0,\xi)$ if $f(y,t)$
exhibits some periodicity, and consider the following two cases.
In the first case, 
$f(y,t)$ is periodic with respect to $t$, i.e.
\begin{equation}\label{eq_a9}
  f(y,t+T) = f(y,t), \quad \text{for all}\ (y,t)\in\mbb R^n\times\mbb R,
\end{equation}
where $T>0$ is the fundamental period.
In the second case, $f(y,t)=f(y_1,y_2,\dots,y_n,t)$
is periodic with respect to $y_i$ with the fundamental period $L_i>0$,
for some $1\leqslant i\leqslant n$, i.e.
\begin{equation}\label{eq_a10}
  f(y+L_i\mbs e_i,t) = f(y,t), \quad \text{for all}\ (y,t)\in\mbb R^n\times\mbb R,
\end{equation}
where $\mbs e_k=(0,\dots,0,1,0,\dots,0)\in\mbb R^n$ ($1\leqslant k\leqslant n$)
denote the standard basis vectors of $\mbb R^n$.

If $f(y,t)$ is periodic in $t$, the solution $y(t)$ to the initial value
problem~\eqref{eq_1} may or may not be a
periodic function. It depends on $y(t)$ for
$t\in[t_0,t_0+T]$, where $T$ is the period of $f(y,t)$ and
$t_0$ is the initial time.
The following lemma summarizes the result.
\begin{lemma}\label{lem_1}
  Suppose $f(y,t)$ is periodic in $t$ with the period $T$
  and $y(t)$ denotes the solution
  to problem~\eqref{eq_1}.
  If $y(t_0+T)=y(t_0)$, then $y(t)$ is periodic with a period $T$
  under the Assumption~\ref{ass_1},
  i.e.~$y(t+T)=y(t)$ for all $t\in\mbb R$.
\end{lemma}
\noindent A proof of this result is provided in the Appendix.
For a general periodic $f(y,t)$ with the period $T$, when
restricted to $t\in[t_0,t_0+T]$, the solution $y(t)$ to~\eqref{eq_1}
in general does not satisfy $y(t_0+T)=y(t_0)$, unless
$\int_{t_0}^{t_0+T}f(y(t),t)dt=0$.
This lemma therefore indicates that the solution to problem~\eqref{eq_1}
in general is not a periodic function, even though $f(y,t)$ is periodic in $t$.

The algorithmic function $\psi(y_0,t_0,\xi)$, on the other hand,
would be periodic with respect to $t_0$, if $f(y,t)$ is periodic in $t$.
This result is summarized in the following theorem.
\begin{theorem}\label{thm_a1}
  Suppose $f(y,t)$ is periodic in $t$ with the period $T$.
  Then $\psi(y_0,t_0,\xi)$ from the system~\eqref{eq_9} is
  periodic in $t_0$ with a period $T$ under the Assumption~\ref{ass_1},
  i.e.~$\psi(y_0,t_0+T,\xi)=\psi(y_0,t_0,\xi)$
  for all $(y_0,t_0,\xi)\in\mbb R^n\times\mbb R\times[0,h_{\max}]$.
\end{theorem}
\noindent A proof of this theorem is provided in the Appendix.
This result suggests that, if $f(y,t)$ is periodic in $t$,
we only need to learn the algorithmic function $\psi(y_0,t_0,\xi)$
over one period along the $t_0$ direction for non-autonomous systems.
The learned
algorithm can then be used to compute the solution $y(t)$
to problem~\eqref{eq_1} for all $t\in\mbb R$ (i.e.~arbitrarily long time horizons),
which, as noted above, is not periodic in general.

If $f(y,t)$ is periodic with respect to a component
of $y$, the algorithmic function $\psi(y_0,t_0,\xi)$
is not periodic with respect to the corresponding component of $y_0$.
Instead, along the direction of that component, $\psi$
satisfies a well-defined relation as given below.
\begin{theorem}\label{thm_a2}
  Suppose $f(y,t)$ is periodic in $y_i$ with the period $L_i$ (see~\eqref{eq_a10}),
  for some $1\leqslant i\leqslant n$.
  Then $\psi(y_0,t_0,\xi)$ from system~\eqref{eq_9} satisfies
  the following relation under the Assumption~\ref{ass_1},
  \begin{equation}\label{eq_a11}
    \psi(y_0+L_i\mbs e_i,t_0,\xi) = \psi(y_0,t_0,\xi) + L_i\mbs e_i,
    \quad \text{for all}\ (y_0,t_0,\xi)\in\mbb R^n\times \mbb R\times[0,h_{\max}].
  \end{equation}
\end{theorem}
\noindent A proof of this theorem is provided in the Appendix.
This result suggests that, if $f(y,t)$ is periodic with respect to some components
of $y$, we only need to learn $\psi(y_0,t_0,\xi)$ over one period
for the values of those  corresponding components of $y_0$.
The learned algorithm can then be used to compute the solution $y(t)$
to problem~\eqref{eq_1}, in which those components of $y(t)$
corresponding to the periodic directions of $f(y,t)$ can each take
arbitrary values on $\mbb R$.

When discussing the time integration using the learned
algorithmic function $\psi(y_0,t_0,\xi)$ later in Section~\ref{sec_224},
we will consider
how to exploit the above properties to simplify the computations.

\subsubsection{Representation of Algorithmic Function $\psi(y_0,t_0,\xi)$ by Randomized NNs}
\label{sec_221}

\begin{figure}
  \centerline{
    \includegraphics[width=2.8in]{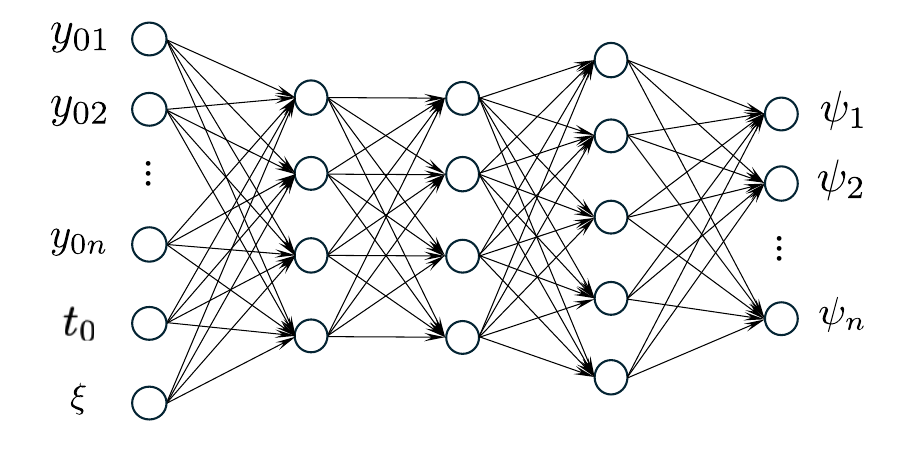}(a)
    \includegraphics[width=3in]{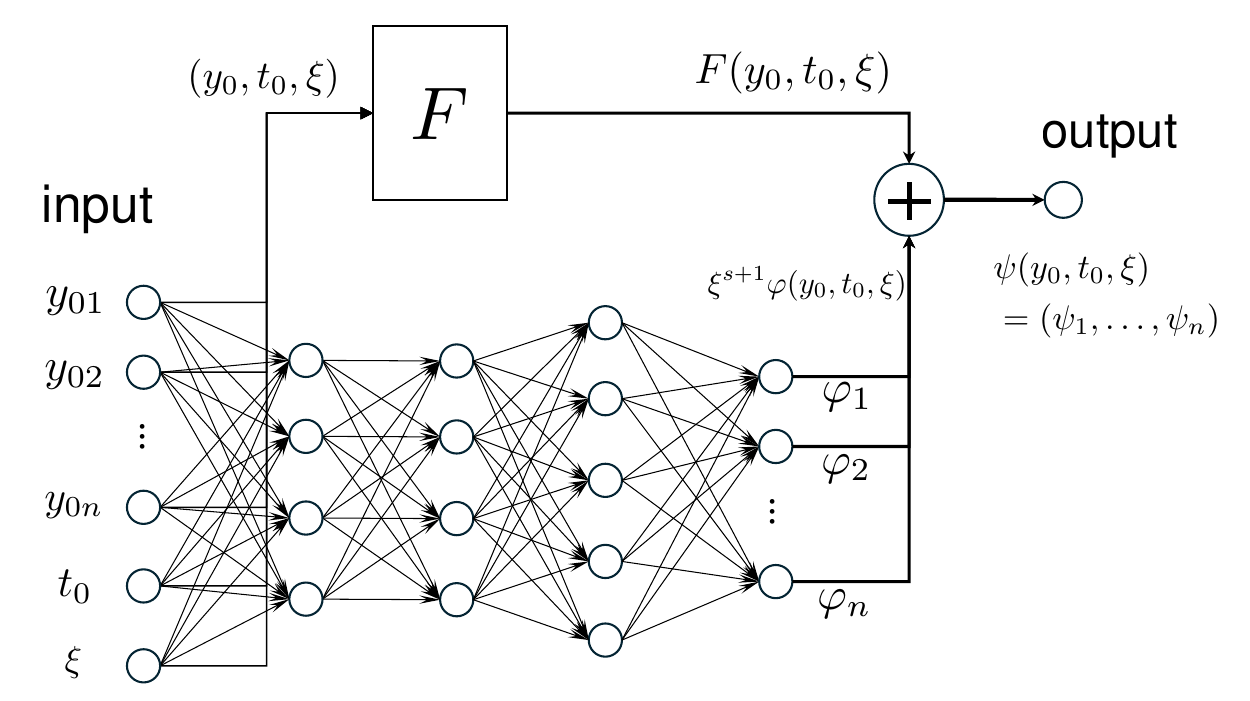}(b)
  }
  \caption{Non-autonomous system: (a) base NN architecture; (b)
    NN structure adopted in this paper.
    $F(y_0,t_0,\xi)$ is an approximation of $\psi(y_0,t_0,\xi)$ of $s$-th order.
    Three hidden layers are shown as an example.
  }
  \label{fg_1}
\end{figure}


We will learn $\psi(y_0,t_0,\xi)$ by solving the system~\eqref{eq_9}
using ELM, for $(y_0,t_0,\xi)$ over some prescribed domain
$\Omega\times[T_0,T_f]\times[0,h_{\max}]$,  with  a physics informed approach.
For this purpose, we represent $\psi(y_0,t_0,\xi)$ by an ELM-type randomized
feed-forward NN; see Figure~\ref{fg_1} for
the NN architecture for representing $\psi(y_0,t_0,\xi)$.
%
Figure~\ref{fg_1}(a) illustrates the base ELM  architecture.
Let $(L+1)$ denote the total number of layers in the network,
with $L\geqslant 2$ being an integer. We refer to the vector
\begin{equation}\label{eq_10}
  \mbs m = [m_0, m_1, \dots, m_L]
\end{equation}
as the architectural vector for this NN,
where $m_i$ ($0\leqslant i\leqslant L$) are positive integers
corresponding to the number of nodes in layer $i$.
Layer zero is the input layer, representing the input variables
$(y_0,t_0,\xi)$, with $m_0=n+2$. Layer
$L$ is the output layer, representing the function $\psi(y_0,t_0,\xi)$,
with $m_L=n$. The layers in between
are the hidden layers. From layer to layer the network logic
represents an affine transform, followed by a function composition
with an activation function
$\sigma: \mbb R\rightarrow \mbb R$~\cite{GoodfellowBC2016}.
The coefficients involved in the affine transforms (weights, biases) of different
layers are collectively referred to
as the network parameters, which may be adjustable (i.e.~trainable)
or fixed (non-trainable).

ELM distinguishes itself from
conventional feed-forward neural networks by the following additional
requirements~\cite{DongL2021,HuangZS2006}:
\begin{itemize}
\item The network parameters in all the hidden layers are pre-assigned to
  random values and fixed (non-trainable). Specifically, we
  set the hidden-layer parameters to uniform random values
  generated on the interval $[-R_m,R_m]$, where $R_m$ is a prescribed constant.
  Once assigned, these network parameters are fixed
  throughout the computation.
\item The network parameters in the output layer are trainable. They are
  the ELM training parameters.
\item The output layer should contain no activation (i.e.~$\sigma(x)=x$),
  with zero bias.
\end{itemize}
These ELM requirements  are imposed on the NNs  when representing
$\psi(y_0,t_0,\xi)$ in this paper.


In practice we find it beneficial to incorporate an approximation
of $\psi$ explicitly into the NN structure
when representing $\psi(y_0,t_0,\xi)$.
This is illustrated in Figure~\ref{fg_1}(b)
and is adopted in the current work.
The network structure in Figure~\ref{fg_1}(b)
implements the following representation
for $\psi(y_0,t_0,\xi)$,
\begin{equation}\label{eq_11}
  \psi(y_0,t_0,\xi) = F(y_0,t_0,\xi) + \xi^{s+1}\varphi(y_0,t_0,\xi),
\end{equation}
where $F(y_0,t_0,\xi)$ is a prescribed function, $s$ is a prescribed constant,
and $\varphi(y_0,t_0,\xi)=(\varphi_1,\dots,\varphi_n)\in\mbb R^n$
denotes a function to be determined and
is represented by an ELM-type randomized NN with its architectural
vector given by~\eqref{eq_10}.

In this paper we choose $F(y_0,t_0,\xi)$ to be an $s$-th order
approximation of $\psi(y_0,t_0,\xi)$ for some integer $s$. With this choice
the second term in~\eqref{eq_11}, $\xi^{s+1}\varphi$, effectively represents
the error of $F(y_0,t_0,\xi)$ in approximating $\psi(y_0,t_0,\xi)$.
Therefore, with the NN structure in Figure~\ref{fg_1}(b),
we are effectively learning the error function $\varphi(y_0,t_0,\xi)$,
and in turn the function $\psi(y_0,t_0,\xi)$, with ELM.
It is noted that by setting $F(y_0,t_0,\xi)=0$ and $s=-1$
the network in Figure~\ref{fg_1}(b) reduces to
the base NN structure of Figure~\ref{fg_1}(a).


We  consider the following representations of
$\psi(y_0,t_0,\xi)$ in this work:
\begin{subequations}\label{eq_12}
  \begin{align}
    & \label{eq_12a}
  \begin{array}{lll}
    s=0: & F(y_0,t_0,\xi) = y_0, & \psi(y_0,t_0,\xi) = y_0 + \xi\varphi(y_0,t_0,\xi);
  \end{array}  \\
  & \label{eq_12b}
  \begin{array}{lll}
    s=1: & F(y_0,t_0,\xi) = y_0+\xi f(y_0,t_0),
    & \psi(y_0,t_0,\xi) = y_0 + \xi f(y_0,t_0) + \xi^2\varphi(y_0,t_0,\xi).
  \end{array}
\end{align}
\end{subequations}
For some problems we also consider the following representation,
\begin{align}
  & \label{eq_12c}
  \begin{array}{ll}
    s = 2: & F_1(y_0,t_0,\xi) = f(y_0,t_0), \quad
    F_2(y_0,t_0,\xi) = f(y_0+\frac{\xi}{2} F_1,t_0+\frac{\xi}{2}), \\
    & F(y_0,t_0,\xi) = y_0 + \xi F_2(y_0,t_0,\xi), \quad
     \psi(y_0,t_0,\xi) = y_0 + \xi F_2(y_0,t_0,\xi) + \xi^3\varphi(y_0,t_0,\xi).
  \end{array}
\end{align}
Notice that these representations automatically satisfy the condition in~\eqref{eq_9b}.
The $F(y_0,t_0,\xi)$ in~\eqref{eq_12b} corresponds to a
forward Euler approximation of~\eqref{eq_9a}, and that in~\eqref{eq_12c}
a mid-point approximation of~\eqref{eq_9a}.
In these forms, $\varphi(y_0,t_0,\xi)$ corresponds to
the error function of $F(y_0,t_0,\xi)$ and is to be learned by ELM.

\begin{remark}\label{rem_1}
  In our implementation, the forms of $F(y_0,t_0,\xi)$
  and $\psi(y_0,t_0,\xi)$ in~\eqref{eq_12}-\eqref{eq_12c}
  are realized using a
  lambda layer from the Keras library,
  in which $\varphi(y_0,t_0,\xi)$ is implemented by a
  feedforward NN  with randomly-assigned
  and fixed (non-trainable) hidden-layer coefficients.
  We would like to comment that other representations of $\psi(y_0,t_0,\xi)$
  can be formulated, e.g.~by choosing $F(y_0,t_0,\xi)$ to be explicit
  Runge-Kutta approximations of~\eqref{eq_9a} of higher orders.

\end{remark}

\subsubsection{Neural Network Training to Learn Algorithmic Function $\psi(y_0,t_0,\xi)$}
\label{sec_222}


We next consider how to train the neural network to learn $\psi(y_0,t_0,\xi)$.
Our training is based on a physics-informed approach, 
employing the nonlinear least squares method~\cite{Bjorck1996}.

We employ the NN structure from Figure~\ref{fg_1}(b), corresponding to
the representations given in~\eqref{eq_12}-\eqref{eq_12c} for $\psi(y_0,t_0,\xi)$.
We will refer to the randomized feedforward NN in
the lower portion of Figure~\ref{fg_1}(b),
which represents $\varphi(y_0,t_0,\xi)$, as the $\varphi$-subnet.
The $\varphi$-subnet has an architecture characterized by the vector~\eqref{eq_10}.

Let $\phi(y_0,t_0,\xi)=(\phi_1(y_0,t_0,\xi),\phi_2(y_0,t_0,\xi),\dots,\phi_M(y_0,t_0,\xi))\in\mbb R^M$,
where $M=m_{L-1}$, denote the output fields of the last hidden layer of
the $\varphi$-subnet. Then the network logic of
the output layer of the $\varphi$-subnet leads to
the following relation,
\begin{equation}\label{eq_13}
  \varphi_i(y_0,t_0,\xi) = \sum_{j=1}^M\beta_{ij}\phi_j(y_0,t_0,\xi) = \bm\beta_i\cdot\phi(y_0,t_0,\xi),
  \quad 1\leqslant i\leqslant n,
\end{equation}
where $\varphi = (\varphi_1,\varphi_2,\dots,\varphi_n)$,
$\beta_{ij}$ ($1\leqslant i\leqslant n$, $1\leqslant j\leqslant M$)
are the weights in the $\varphi$-subnet's output layer,
and $\bm\beta_i = (\beta_{i1},\dots,\beta_{iM})\in\mbb R^M$.
Note that $\beta_{ij}$ are the training parameters of the $\varphi$-subnet,
and of the overall network for $\psi$.

The objective here is to train the parameters $\beta_{ij}$ so that $\psi(y_0,t_0,\xi)$
satisfies~\eqref{eq_9a} on $y_0\in\Omega\subset\mbb R^n$, $t_0\in[T_0,T_f]\subset\mbb R$,
and $\xi\in[0,h_{\max}]\subset\mbb R$, for prescribed $\Omega$, $T_0$, $T_f$ and $h_{\max}$.
Note that~\eqref{eq_9b}
is automatically satisfied by the NN formulation for~$\psi(y_0,t_0,\xi)$.

Define the residual function $R=(R_1,\dots,R_n)\in\mbb R^{n}$ of the problem,
\begin{equation}\label{eq_14}
  R(\bm\beta,y_0,t_0,\xi) = \frac{\partial\psi}{\partial \xi} - f(\psi,t_0+\xi),
\end{equation}
where $\psi$ is given by~\eqref{eq_11} and~\eqref{eq_12}-\eqref{eq_12c},
$\bm\beta = (\bm\beta_1,\bm\beta_2,\dots,\bm\beta_n)
= (\beta_{11}\dots,\beta_{1M}, \beta_{21},\dots,\beta_{nM})\in\mbb R^{nM}$,
and the dependence of $R$ on $\bm\beta$ has been made explicit.
We choose a set of $Q$ points,
$(y_0^{(i)},t_0^{(i)},\xi^{(i)})\in\Omega\times[T_0,T_f]\times[0,h_{\max}]$ ($1\leqslant i\leqslant Q$),
from a uniform random distribution and refer to them
as the collocation points hereafter.

Enforcing $R(\bm\beta,y_0,t_0,\xi)$ to be zero on these collocation points
gives rise to the following system,
\begin{equation}\label{eq_15}
  r^{(i)}(\bm\beta) = R(\bm\beta, y_0^{(i)},t_0^{(i)},\xi^{(i)})
  = \left.\frac{\partial\psi}{\partial \xi}\right|_{(y_0^{(i)},t_0^{(i)},\xi^{(i)})}
  -f(\psi(y_0^{(i)},t_0^{(i)},\xi^{(i)}),t_0^{(i)}+\xi^{(i)}) = 0,
  \quad 1\leqslant i\leqslant Q,
\end{equation}
where $r^{(i)}\in\mbb R^n$.
This is a system of nonlinear algebraic equations about $\bm\beta$,
consisting of $nQ$ equations with $nM$ unknowns.
Note that, for any given $\bm\beta$,
the term $\psi(y_0^{(i)},t_0^{(i)},\xi^{(i)})$ can be computed by
an evaluation of the neural network, and
the term $\left.\frac{\partial\psi}{\partial \xi}\right|_{(y_0^{(i)},t_0^{(i)},\xi^{(i)})}$
can be computed by a forward-mode automatic differentiation.

We seek a least squares solution  to the
algebraic system~\eqref{eq_15}, and determine $\bm\beta$
by the nonlinear least squares method (i.e.~Gauss-Newton
method)~\cite{Bjorck1996,Bjorck2015}.
Specifically, we compute $\bm\beta$ using the NLLSQ-perturb
(nonlinear least squares
with perturbations) algorithm developed in~\cite{DongL2021};
see also the Appendix A of~\cite{DongW2023} for a more detailed exposition
of NLLSQ-perturb. NLLSQ-perturb
employs the scipy implementation of the Gauss-Newton
method plus a trust-region strategy (scipy.optimize.least\_squares
routine) and additionally incorporates a perturbation scheme
to keep the method from being trapped to the worst local minima.
Upon attaining the least squares solution $\bm\beta$ by NLLSQ-perturb,
we set the output-layer parameters of the $\varphi$-subnet by this solution
to conclude the NN training.
For the NN training,
the input data to the network  consists of the set of collocation
points $(y_0^{(i)},t_0^{(i)},\xi^{(i)})\in \Omega\times[T_0,T_f]\times[0,h_{\max}]$
($1\leqslant i\leqslant Q$).

To elaborate on the training procedure,
the NLLSQ-perturb algorithm~\cite{DongL2021,DongW2023} requires two routines
for its input, one for computing the residual vector
$r(\bm\beta) = (\dots,r^{(i)}(\bm\beta),\dots)\in\mbb R^{nQ}$ and
the other for computing
the Jacobian matrix
$\frac{\partial r}{\partial\bm\beta}\in\mbb R^{nQ\times nM}$,
for an arbitrary given $\bm\beta$.
Computing the residual $r^{(i)}(\bm\beta)$ by~\eqref{eq_15} is straightforward,
noting that the terms $\psi(y_0^{(i)},t_0^{(i)},\xi^{(i)})$ and
$\left.\frac{\partial\psi}{\partial \xi}\right|_{(y_0^{(i)},t_0^{(i)},\xi^{(i)})}$ therein can be
attained by forward evaluations of the NN  or by automatic
differentiations, as discussed above.
Computing the Jacobian matrix is more involved.
We next discuss its computation and 
related implementation issues.


To facilitate computation of the Jacobian matrix, we note from~\eqref{eq_13}
and~\eqref{eq_11} that,
\begin{equation}\label{eq_16}
  \frac{\partial\varphi_i}{\partial\bm\beta_j} = \delta_{ij}\phi(y_0,t_0,\xi)
  \in\mbb R^{1\times M},
  \quad
  \frac{\partial\psi_i}{\partial\bm\beta_j} =
  \xi^{s+1}\frac{\partial\varphi_i}{\partial\bm\beta_j}
  = \delta_{ij}\xi^{s+1}\phi(y_0,t_0,\xi) \in\mbb R^{1\times M},
  \quad 1\leqslant i,j\leqslant n,
\end{equation}
where $\psi = (\psi_1,\dots,\psi_n)$ and
$\delta_{ij}$ denotes the Kronecker delta.
In light of~\eqref{eq_14},
for $1\leqslant i,j\leqslant n$,
\begin{equation}\label{eq_17}
  \frac{\partial R_i}{\partial\bm\beta_j} =
  \frac{\partial}{\partial\xi}\left(\frac{\partial\psi_i}{\partial\bm\beta_j} \right)
  -\sum_{k=1}^n\frac{\partial f_i}{\partial\psi_k}\frac{\partial\psi_k}{\partial\bm\beta_j}
  = (s+1)\xi^s\phi\delta_{ij} + \xi^{s+1}\frac{\partial\phi}{\partial\xi}\delta_{ij}
   - \xi^{s+1}\frac{\partial f_i}{\partial\psi_j}\phi \ \in\mbb R^{1\times M},
\end{equation}
where $f(\psi,t_0+\xi) = (f_1,\dots,f_n)$ and we have used~\eqref{eq_16}.
So the Jacobian matrix is given by
\begin{equation}\label{eq_18}
  \frac{\partial r}{\partial\bm\beta} = \begin{bmatrix}
    \frac{\partial r^{(1)}}{\partial\bm\beta_1} & \dots & \frac{\partial r^{(1)}}{\partial\bm\beta_n} \\
    \vdots & \ddots & \vdots \\
    \frac{\partial r^{(Q)}}{\partial\bm\beta_1} & \dots & \frac{\partial r^{(Q)}}{\partial\bm\beta_n}
  \end{bmatrix} \in\mbb R^{nQ\times nM},
  \quad
  \frac{\partial r^{(i)}}{\partial\bm\beta_j} =
  \left.\frac{\partial R}{\partial\bm\beta_j}\right|_{(y_0^{(i)},t_0^{(i)},\xi^{(i)})}=
  \begin{bmatrix}
    \left.\frac{\partial R_1}{\partial\bm\beta_j}\right|_{(y_0^{(i)},t_0^{(i)},\xi^{(i)})} \\
    \vdots \\ \left.\frac{\partial R_n}{\partial\bm\beta_j}\right|_{(y_0^{(i)},t_0^{(i)},\xi^{(i)})}
  \end{bmatrix} \in\mbb R^{n\times M},
\end{equation}
where $\frac{\partial R_i}{\partial\bm\beta_j}$ is given in~\eqref{eq_17}.

\begin{remark}\label{rem_2}
  The Jacobian matrix involves the terms like $\phi(y_0^{(i)},t_0^{(i)},\xi^{(i)})$ and
  $\left.\frac{\partial\phi}{\partial\xi}\right|_{(y_0^{(i)}),t_0^{(i)},\xi^{(i)})}$, which
  represent the output fields of the last hidden layer of
  the $\varphi$-subnet and their derivatives.
  To compute these terms,
  in our implementation we have created a Keras sub-model of the $\varphi$-subnet,
  which takes $(y_0,t_0,\xi)$ as its input
  and $\phi(y_0,t_0,\xi)$ (last hidden layer
  of the $\varphi$-subnet) as its output.
  The terms $\phi(y_0^{(i)},t_0^{(i)},\xi^{(i)})$ are then computed by a forward
  evaluation of this Keras sub-model,
  and the terms
  $\left.\frac{\partial\phi}{\partial\xi}\right|_{(y_0^{(i)},t_0^{(i)},\xi^{(i)})}$
  are computed by a forward-mode automatic differentiation with this sub-model.
  
\end{remark}

\begin{remark}\label{rem_23}
  When the problem~\eqref{eq_1} becomes more complicated (e.g.~being stiff),
  we find it necessary to incorporate a domain decomposition into the above
  algorithm for learning $\psi(y_0,t_0,\xi)$. The discussion below pertains to
  the use of domain decomposition in algorithm learning.
  Suppose the domain $\Omega\times[T_0,T_f]\subset\mbb R^{n+1}$ is
  partitioned into $N_e$ ($N_e\geqslant 1$)
  non-overlapping sub-domains,
  $\Omega\times[T_0,T_f]=\mathcal{D}_1\cup\mathcal{D}_2\cup\dots\cup\mathcal{D}_{N_e}$.
  On sub-domain $\mathcal{D}_i$ ($1\leqslant i\leqslant N_e$) we seek a
  function $\psi_i: \mathcal{D}_i\times[0,h_{\max}]\rightarrow\mbb R^n$ such that
  \begin{subequations}\label{eq_19}
    \begin{align}
      &
      \frac{\partial\psi_i}{\partial\xi} = f(\psi_i(y_0,t_0,\xi), t_0+\xi), \quad
      ((y_0,t_0),\xi)\in\mathcal{D}_i\times[0,h_{\max}];
      \\
      &
        \psi_i(y_0,t_0,0) = y_0, \quad\quad\quad\quad\quad\quad (y_0,t_0)\in\mathcal{D}_i.
  \end{align}
  \end{subequations}
  The algorithmic function
  $\psi(y_0,t_0,\xi)$ for $(y_0,t_0,\xi)\in\Omega\times[T_0,T_f]\times[0,h_{\max}]$
  is then given by
  \begin{equation}
    \psi(y_0,t_0,\xi) = \left\{
    \begin{array}{ll}
      \psi_1(y_0,\xi), & (y_0,t_0)\in\mathcal{D}_1; \\
      \psi_2(y_0,t_0,\xi), & (y_0,t_0)\in\mathcal{D}_2; \\
      \dots \\
      \psi_{N_e}(y_0,t_0,\xi), & (y_0,t_0)\in\mathcal{D}_{N_e}.
    \end{array}
    \right.
  \end{equation}
  %
  We use a local NN with an architecture given by Figure~\ref{fg_1}(b)
  to learn each $\psi_i(y_0,t_0,\xi)$, by solving the system~\eqref{eq_19}
  on $\mathcal D_i\times[0,h_{\max}]$, for $1\leqslant i\leqslant N_e$.
  This local learning problem on $\mathcal{D}_i\times[0,h_{\max}]$ is essentially
  the same as~\eqref{eq_9}.
  Hence the same algorithm as
  presented above can be used to
  learn $\psi_i$
  ($1\leqslant i\leqslant N_e$), noting that here the random collocation points
  shall be taken from $\mathcal{D}_i\times[0,h_{\max}]$. 
  We would like to emphasize that the learning problems
  on different sub-domains are not coupled, and they
  can be computed individually or in parallel.
  In our implementation we have created $N_e$ local ELM-type randomized NNs,
  implemented as Keras sub-models,
  with each for one sub-domain. The local NNs  each
  assumes a structure as given by Figure~\ref{fg_1}(b). They are
  trained individually  in an un-coupled fashion,
  using the procedure discussed in this section.
  The learned time-marching algorithm is represented by
  the collection of these local NNs.
  
\end{remark}

\begin{remark}\label{rem_24}
  The expressions~\eqref{eq_11} and~\eqref{eq_12}-\eqref{eq_12c} are explicit
  representations of $\psi(y_0,t_0,\xi)$.
  Given $(y_0,t_0,\xi)$, $\psi(y_0,t_0,\xi)$ can be explicitly computed
  by a forward NN evaluation.
  It is also possible to
  learn $\psi(y_0,t_0,\xi)$ by adopting an
  implicit representation, which would entail solving
  algebraic systems  during training and also during time marching,
  in addition to forward NN evaluations.
  Therefore, implicit representations of $\psi(y_0,t_0,\xi)$
  lead to implicit time integration algorithms.
  %
  There are several ways to represent $\psi(y_0,t_0,\xi)$ implicitly.
  One simple idea is to again adopt the
  representation~\eqref{eq_11} and~\eqref{eq_12}-\eqref{eq_12c}
  with the same NN as in Figure~\ref{fg_1}(b),
  but restrict $\xi$ to the domain $\xi\in[-h_{\max},0]$ during NN training.
  This amounts to an algorithm
  that marches backward in time, giving rise to an implicit time
  integration scheme. The procedure
  discussed in this section, with $\xi\in[-h_{\max},0]$,
  can train the NN for learning this implicit scheme.

  We next discuss another implicit representation
  of $\psi(y_0,t_0,\xi)$, by leveraging implicit Runge-Kutta type
  approximations. Let us consider
  \begin{equation}\label{eq_21}
    \psi(y_0,t_0,\xi) = G(y_0,t_0,\xi) + \xi^{s+1}\varphi(y_0,t_0,\xi),
  \end{equation}
  where $G(y_0,t_0,\xi)$ is an $s$-th order implicit 
  approximation of $\psi(y_0,t_0,\xi)$ in equation~\eqref{eq_9a}.
  We consider the following specific forms,
  \begin{subequations}\label{eq_22}
    \begin{align}
      & \label{eq_22b}
      \begin{array}{lll}
        s = 1: & K(y_0,t_0,\xi) = f(y_0+\xi K(y_0,t_0,\xi), t_0+\xi),\\
        & G(y_0,t_0,\xi) = y_0 + \xi K(y_0,t_0,\xi), \quad
         \psi(y_0,t_0,\xi) = y_0+\xi K(y_0,t_0,\xi) + \xi^2\varphi(y_0,t_0,\xi);
      \end{array} \\
      & \label{eq_22c}
      \begin{array}{lll}
        s = 2: & K_1(y_0,t_0,\xi) = f(y_0+\gamma\xi K_1(y_0,t_0,\xi), t_0+\gamma\xi),\\
        & K_2(y_0,t_0,\xi) = f(y_0+(1-\gamma)\xi K_1(y_0,t_0,\xi) + \gamma\xi K_2(y_0,t_0,\xi),t_0+\xi), \\
        & G(y_0,t_0,\xi) = y_0+(1-\gamma)\xi K_1(y_0,t_0,\xi) + \gamma\xi K_2(y_0,t_0,\xi), \\
        & \psi(y_0,t_0,\xi) = y_0+(1-\gamma)\xi K_1(y_0,t_0,\xi) + \gamma\xi K_2(y_0,t_0,\xi)
        + \xi^3\varphi(y_0,t_0,\xi);
      \end{array}
    \end{align} 
  \end{subequations}
  where $\gamma = 1-\frac{\sqrt{2}}{2}$.
  The form~\eqref{eq_22b}
  utilizes the first-order backward Euler approximation
  in the implicit Runge-Kutta form, and~\eqref{eq_22c} uses the 2nd-order
  diagonally implicit Runge-Kutta (DIRK) approximation~\cite{HairerW1996,KenedyC2016}.
  These representations automatically satisfy~\eqref{eq_9b}.
  %
  In~\eqref{eq_21}, $\varphi(y_0,t_0,\xi)$ is to be determined, and  is represented by
  an ELM-type randomized NN.
  This leads to the relation~\eqref{eq_13}, in which $\beta_{ij}$
  are the trainable parameters
  and $\phi(y_0,t_0,\xi)\in\mbb R^M$ are the output fields of
  the last hidden layer of the NN.
  To train this NN,
  we seek  $\beta_{ij}$
  such that the expression~\eqref{eq_21} for $\psi(y_0,t_0,\xi)$
  satisfies~\eqref{eq_9a}, in the least squares sense,
  on  $Q$
  random collocation points $(y_0^{(i)},t_0^{(i)},\xi^{(i)})\in\Omega\times[T_0,T_f]\times[0,h_{\max}]$
  ($1\leqslant i\leqslant Q$). We can similarly compute $\beta_{ij}$
  by the NLLSQ-perturb algorithm. In this case the residual function
  is defined by~\eqref{eq_14} and the residuals on
  the collocation points are given by~\eqref{eq_15}, where $\psi(y_0,t_0,\xi)$
  is given by~\eqref{eq_21}.

  Let us use the form~\eqref{eq_22b}
  to illustrate how to compute the residual
  vector $r(\bm\beta)\in\mbb R^{nQ}$ and the Jacobian matrix
  $\frac{\partial r}{\partial\bm\beta}\in\mbb R^{nQ\times nM}$
  for use by NLLSQ-perturb during training.
  The representation~\eqref{eq_22c} requires a procedure similar
  to what follows.
  With~\eqref{eq_22b}
  the residual~\eqref{eq_15} is reduced to
  \begin{equation}\label{eq_23}
    \begin{split}
    r^{(i)}(\bm\beta) =& \ K(y_0^{(i)},t_0^{(i)},\xi^{(i)})
    + \xi^{(i)}\left.\frac{\partial K}{\partial\xi} \right|_{(y_0^{(i)},t_0^{(i)},\xi^{(i)})}
    + 2\xi^{(i)}\varphi(y_0^{(i)},t_0^{(i)},\xi^{(i)})
    + (\xi^{(i)})^2\left.\frac{\partial\varphi}{\partial\xi} \right|_{(y_0^{(i)},t_0^{(i)},\xi^{(i)})}\\
    &- f(\psi(y_0^{(i)},t_0^{(i)},\xi^{(i)}),t_0^{(i)}+\xi^{(i)}),
    \quad 1\leqslant i\leqslant Q.
    \end{split}
  \end{equation}
  Here $K(y_0^{(i)},t_0^{(i)},\xi^{(i)})\in\mbb R^n$ is computed by solving
  the nonlinear equation (see the first equation of~\eqref{eq_22b})
  \begin{equation}\label{eq_24}
    K(y_0^{(i)},t_0^{(i)},\xi^{(i)}) = f(y_0^{(i)} + \xi^{(i)}K,t_0^{(i)}+\xi^{(i)}).
  \end{equation}
  In our implementation we solve this equation using the
  routine ``root'' in the scipy library (scipy.optimize.root).
  $\left.\frac{\partial K}{\partial\xi} \right|_{(y_0^{(i)},t_0^{(i)},\xi^{(i)})}\in\mbb R^n$
  in~\eqref{eq_23} is computed by solving the linear algebraic system
  \begin{equation}\label{eq_25}
    \left(\mbs I - \xi^{(i)}\left.\frac{\partial f}{\partial G} \right|_{(y_0^{(i)},t_0^{(i)},\xi^{(i)})} \right)
    \left.\frac{\partial K}{\partial\xi} \right|_{(y_0^{(i)},t_0^{(i)},\xi^{(i)})}
    = \left.\frac{\partial f}{\partial G} \right|_{(y_0^{(i)},t_0^{(i)},\xi^{(i)})}
    K(y_0^{(i)},t_0^{(i)},\xi^{(i)})
    + \left.\frac{\partial f}{\partial t}\right|_{(y_0^{(i)},t_0^{(i)},\xi^{(i)})}, 
  \end{equation}
  where $\mbs I\in \mbb R^{n\times n}$ is the identity matrix,
  $G = y_0 + \xi K(y_0,t_0,\xi)$, $t = t_0+\xi$, and
  $f = f(G,t)=f(y_0+\xi K, t_0+\xi)$.
  This linear system results from the differentiation of the
  first equation in~\eqref{eq_22b} with respect to $\xi$. We solve this system
  using the routine ``solve'' from the scipy library
  (scipy.linalg.solve).
  The terms $\varphi(y_0^{(i)},t_0^{(i)},\xi^{(i)})$ and 
  $\left.\frac{\partial\varphi}{\partial\xi} \right|_{(y_0^{(i)},t_0^{(i)},\xi^{(i)})}$
  in~\eqref{eq_23} are computed by a forward evaluation of the NN
  and by a forward-mode automatic differentiation.
  %
  The Jacobian matrix $\frac{\partial r}{\partial\bm\beta}$ is given by~\eqref{eq_18},
  in which
  \begin{equation}\label{eq_26}
    \frac{\partial R_i}{\partial\bm\beta_j} =
    2\xi\phi(y_0,t_0,\xi)\delta_{ij} + \xi^2\frac{\partial\phi}{\partial\xi}\delta_{ij}
    -\xi^2\frac{\partial f_i}{\partial\psi_j}\phi(y_0,t_0,\xi),
  \end{equation}
  to be evaluated on the collocation points.
  In summary, given $(y_0^{(i)},t_0^{(i)},\xi^{(i)})$ ($1\leqslant i\leqslant Q$)
  and an arbitrary $\bm\beta$, we take the following steps to compute
  $r(\bm\beta)$ and $\frac{\partial r}{\partial\bm\beta}$:
  \begin{enumerate}[(i), nosep]
  \item Compute $\varphi(y_0^{(i)},t_0^{(i)},\xi^{(i)})$,
    $\left.\frac{\partial\varphi}{\partial\xi}\right|_{(y_0^{(i)},t_0^{(i)},\xi^{(i)})}$,
    $\phi(y_0^{(i)},t_0^{(i)},\xi^{(i)})$ and
    $\left.\frac{\partial\phi}{\partial\xi}\right|_{(y_0^{(i)},t_0^{(i)},\xi^{(i)})}$
    by forward valuations of the neural network and by automatic
    differentiations;

  \item Solve equation~\eqref{eq_24} for $K(y_0^{(i)},t_0^{(i)},\xi^{(i)})$;
  \item Solve the linear system~\eqref{eq_25} for
    $\left.\frac{\partial K}{\partial\xi} \right|_{(y_0^{(i)},t_0^{(i)},\xi^{(i)})}$;
  \item Compute $\psi(y_0^{(i)},t_0^{(i)},\xi^{(i)})$ by the second equation in~\eqref{eq_22b};
  \item Compute $r(\bm\beta)$ by~\eqref{eq_23};
  \item Compute $\frac{\partial r}{\partial\bm\beta}$ by~\eqref{eq_18}
    and~\eqref{eq_26}.
    
  \end{enumerate}

\end{remark}

\begin{remark}\label{rem_c26}
  Hereafter, we refer to the learned NN algorithms employing the explicit
  formulations~\eqref{eq_12a},
  \eqref{eq_12b}, and \eqref{eq_12c} as ``NN-Exp-S0'',
  ``NN-Exp-S1'', and ``NN-Exp-S2'', respectively.
  We refer to those based on the implicit formulations~\eqref{eq_22b}
  and~\eqref{eq_22c} as ``NN-Imp-S1'' and ``NN-Imp-S2'', respectively.
\end{remark}


\subsubsection{Time Marching Based on Learned Algorithmic Function $\psi(y_0,t_0,\xi)$}
\label{sec_224}

%
%
%

With appropriately chosen domain $\Omega\times[T_0,T_f]\times[0,h_{\max}]$
for network training,
the trained ELM network can be used as a time marching algorithm
for solving problem~\eqref{eq_1}, with arbitrary initial data
$(y_0,t_0)\in\Omega\times[T_0,T_f]$ and step size $h\in(0,h_{\max})$.

Let $\psi_{\bm\beta}(y_0,t_0,\xi)$ denote the learned algorithmic function
represented by the trained neural network.
Suppose the initial time and data
are $(t_0,y_0)$, $h$ is the step size, and $y_k$ is the approximation to
the solution at $t_k=t_0+kh$ ($k\geqslant 0$).
Given $(y_k,t_k)$, we compute the  solution at $t_{k+1}=t_k+h$
by
\begin{equation}\label{eq_27}
  y_{k+1} = \psi_{\bm\beta}(y_k,t_k,h), \quad k\geqslant 0.
\end{equation}
It can be noted that only forward NN evaluations 
are needed for solving IVPs if the learned algorithm is based on
an explicit representation of $\psi(y_0,t_0,\xi)$.
On the other hand, if the algorithm
is based on an implicit representation (see Remark~\ref{rem_24}),
solving an algebraic system
is required during time marching, apart from the forward NN evaluations.
We will elaborate on the implicit case in a remark (Remark~\ref{rem_27})
below.


  Since one is not able to  learn the algorithmic function $\psi(y_0,t_0,\xi)$ perfectly,
  due to practical constraints (such as space, time, and computational resource),
  it should be emphasized that the learned algorithm $\psi_{\bm\beta}(y_0,t_0,\xi)$
  is but an approximation of
  the exact time integration algorithm $\psi(y_0,t_0,\xi)$.
  Nonetheless, we observe that the learned algorithms are
  highly competitive, in terms of their accuracy and time marching
  cost, compared with the leading traditional time integration algorithms.
  This point will be demonstrated
  with numerical experiments in Section~\ref{sec_tests}.

\begin{remark}\label{rem_25}
  Choosing an appropriate domain $\Omega\times[T_0,T_f]\times[0,h_{\max}]$, from
  which the collocation points are drawn, for training
  the NN  is important to the accuracy of the learned
  algorithm $\psi_{\bm\beta}(y_0,t_0,\xi)$.
  In general the domain $\Omega$ and $[T_0,T_f]$ should be
  sufficiently large so that $y_k\in\Omega$ and $t_k\in[T_0,T_f]$
  for all $0\leqslant k\leqslant N$ ($N$ denoting the number of time steps one plans to perform),
  and the parameter $h_{\max}$ should be sufficiently large so that
  the step size $h$ of interest falls in $(0,h_{\max})$.
  On the other hand, an overly large $\Omega$, $(T_f-T_0)$, or $h_{\max}$
  can increase the difficulty in the network training.
  In this regard, some knowledge about the system to be simulated can be helpful
  to the choice of $\Omega$, $h_{\max}$, and $[T_0,T_f]$.
  In the absence of any knowledge about the system, preliminary simulations
  are useful for making a choice about $\Omega$, $h_{\max}$ and $[T_0,T_f]$.
  
\end{remark}

\begin{remark}\label{rem_027}
  Suppose one intends to solve the non-autonomous
  system~\eqref{eq_1} for $0\leqslant t\leqslant t_f$,
  with $t_0=0$ in~\eqref{eq_1b} and $t_f$ being large (long time integration).
  When training the NN, the domain $[T_0,T_f]$
  should in general be chosen to be at least $[0,t_f]$.
  Since $t_f$ is large, using a single NN  to accurately
  learn $\psi(y_0,t_0,\xi)$ 
  on $(y_0,t_0,\xi)\in\Omega\times[0,t_f]\times[0,h_{\max}]$
  can become very challenging.
  To alleviate this problem, one can incorporate a domain decomposition
  of $[0,t_f]$ and employ local NNs
  to learn $\psi$ on the sub-domains (see Remark~\ref{rem_23}).
  Let us suppose $[0,t_f]$ is partitioned into $m$ ($m\geqslant 1$) sub-domains with
  $0=T_0<T_1<\cdots<T_m=t_f$.
  One only
  needs to train $m$ local NNs, each for a sub-domain
  $t_0\in[T_{i-1},T_i]$ ($1\leqslant i\leqslant m$),
  individually in an un-coupled fashion.
  By using a moderate size for the sub-domains,
  training the local NNs would become significantly
  easier. By incorporating domain decomposition and local neural
  networks, one can effectively learn the time marching algorithm
  for solving non-autonomous systems over long time horizons.

\end{remark}

\begin{remark}\label{rem_26}
  In the presence of domain decomposition (see Remark~\ref{rem_23}),
  the algorithm~\eqref{eq_27} needs to be modified accordingly for time
  integration. Let $\psi_{\bm\beta i}(y_0,t_0,\xi)$ 
  denote the learned $\psi_i(y_0,t_0,\xi)$ on
  the sub-domain $\mathcal{D}_i\times [0,h_{\max}]$
  for $1\leqslant i\leqslant N_e$. Then given
  $(y_k,t_k)$ we approximate the solution at 
  $t_{k+1}=t_k+h$ by
  \begin{equation}\label{eq_28}
    \begin{array}{l}
      (i)\ \text{determine}\ s\ (1\leqslant s\leqslant N_e)\
      \text{such that}\ (y_k,t_k)\in\mathcal{D}_s; \\
      (ii)\ y_{k+1} = \psi_{\bm\beta s}(y_k,t_k,h).
    \end{array}
  \end{equation}
  In the event $(y_k,t_k)$ falls on the shared
  boundary of two or more sub-domains, one can choose
  the $\psi_{\bm\beta}$ corresponding to any of these sub-domains
  for time integration. In our implementation, we have
  used the sub-domain with the lowest ID 
  for time integration.
  To reduce the influence of different choices of sub-domains
  in such cases, during NN training it is preferable, after
  the domain $\Omega\times[T_0,T_f]$ is partitioned into disjoint
  sub-domains $\mathcal{D}_i$,
  to enlarge each sub-domain slightly (e.g.~by a few percent)
  along the directions of domain decomposition.
  In this way, $\psi_i(y_0,t_0,\xi)$ will be learned/trained on
  a domain slightly larger than $\mathcal{D}_i$.
  To be more specific, suppose the sub-domain $\mathcal{D}_i$ has a dimension $[a,b]$
  along a direction with domain decomposition,
  and let $r\geqslant 0$ denote the enlargement
  factor. Then along this direction the enlarged sub-domain for NN training
  will have a dimension $\left[a-(b-a)\frac{r}{2}, b+(b-a)\frac{r}{2}\right]$.
  Note that the use of enlarged sub-domains is 
  for network training only. During time marching,
  the algorithm~\eqref{eq_28} will still choose the algorithmic
  function based on the disjoint partitions $\mathcal{D}_i$.

\end{remark}

\begin{remark}\label{rem_29}
  If $f(y,t)$ in~\eqref{eq_1}
  is periodic in $t$ with the period $T$, as shown by Theorem~\ref{thm_a1},
  it would be sufficient to learn $\psi(y_0,t_0,\xi)$
  on $(y_0,t_0,\xi)\in\Omega\times[T_0,T_f]\times[0,h_{\max}]$,
  with $[T_0,T_f]$ covering  one period of $f(y,t)$,
  and the learned algorithm $\psi_{\bm\beta}(y_0,t_0,\xi)$ can be used to solve
  system~\eqref{eq_1} for all $t\in\mbb R$ (arbitrarily long
  time horizons). This property can significantly simplify
  the algorithm learning and the long-time integration of
  non-autonomous systems with a  periodic RHS.
  %
  In this case,
  time integration using the learned algorithm
  $\psi_{\bm\beta}(y_0,t_0,\xi)$, which is trained on
  $(y_0,t_0,\xi)\in\Omega\times[0,T]\times[0,h_{\max}]$, takes the following
  form. Given ($y_k,t_k$) with $t_k\in\mbb R$, we compute
  ($y_{k+1},t_{k+1}$) by,
  \begin{equation}\label{eq_46}
    y_{k+1} = \psi_{\bm\beta}(y_k,t_k^*,h), \quad
    t_k^* = \mathrm{mod}(t_k,T), \quad
    t_{k+1} = t_k+h,
  \end{equation}
  where $\mathrm{mod}$ denotes the modulo operation.
  
\end{remark}

\begin{remark}\label{rem_210}
  If $f(y,t)$ in~\eqref{eq_1} is periodic with respect to
  one or more components of $y\in\mbb R^n$,
  when learning $\psi(y_0,t_0,\xi)$ on
  $(y_0,t_0,\xi)\in\Omega\times[T_0,T_f]\times[0,h_{\max}]$,
  it would be sufficient to choose $\Omega$ to cover
  one period of $f(y,t)$ along those directions, as shown by Theorem~\ref{thm_a2}.
  The resultant algorithm can be used to solve the system~\eqref{eq_1}
  with those components of $y_0\in\mbb R^n$
  in the periodic directions taking arbitrary values.

  Define a constant vector $\mbs L=(L_1,\dots,L_n)\in\mbb R^n$
  with $L_i\geqslant 0$ ($1\leqslant i\leqslant n$) and assume $\mbs L\neq 0$.
  Suppose
  \begin{equation}\label{eq_48}
    f(y+L_i\mbs e_i, t) = f(y,t), \quad
    \text{for all}\ (y,t)\in\mbb R^n\times\mbb R, \quad
    1\leqslant i\leqslant n.
  \end{equation}
  If $L_i>0$ ($1\leqslant i\leqslant n$),
  equation~\eqref{eq_48} indicates that
  $f(y,t)$ is periodic with respect to $y_i$
  with a period $L_i$.
  If $L_i=0$, equation~\eqref{eq_48} is reduced to an identity and
  we use this to denote that $f(y,t)$ is
  not periodic with respect to $y_i$.
  We will refer to $\mbs L$ as the periodicity vector of $f(y,t)$ hereafter.

  Given the periodicity vector $\mbs L=(L_1,\dots,L_n)$ for $f(y,t)$,
  we suppose the NN has been trained to learn
  $\psi(y_0,t_0,\xi)$ on $\Omega\times[T_0,T_f]\times[0,h_{\max}]$,
  where $\Omega$ covers one period of $f(y,t)$ along those
  periodic directions.  Solving
  system~\eqref{eq_1} based on the learned
  algorithm $\psi_{\bm\beta}(y_0,t_0,\xi)$ takes the following form.
  Given step size $h$ and ($y_k,t_k$), where $y_k=(y_{k1},\dots,y_{kn})$,
  \begin{enumerate}[(i),nosep] 
  \item for $i=1,\dots,n$:
    \begin{enumerate}[,nosep]
    \item if $L_i>0$, then set $y_{ki}^*=\mathrm{mod}(y_{ki},L_i)$
      and $q_i=\left\lfloor \frac{y_{ki}}{L_i}\right\rfloor$; \\
      else set $y_{ki}^*=y_{ki}$ and $q_i=0$.
    \end{enumerate}
    
  \item
    $y_{k+1} = \psi_{\bm\beta}(y_k^*,t_k,h) + \mathrm{mult}(q,\mbs L)$, where
    $y_k^*=(y_{k1}^*,\dots,y_{kn}^*)$,
    $q=(q_1,\dots,q_n)\in\mbb Z^n$, and $\mathrm{mult}(\cdot,\cdot)$ denotes
    the element-wise multiplication of two vectors and returns the resultant
    vector.

  \end{enumerate}
  In the above steps $\lfloor\cdot\rfloor$ denotes the floor function.
  %

\end{remark}

\begin{remark}\label{rem_27}
  Implicit
  representations of $\psi(y_0,t_0,\xi)$ (see Remark~\ref{rem_24}) lead to
  implicit time integration algorithms
  for solving~\eqref{eq_1}.
  For an implicit NN algorithm $\psi_{\bm\beta}$
  trained on $(y_0,t_0,\xi)\in\Omega\times[T_0,T_f]\times[-h_{\max},0]$,
  given $(y_k,t_k)$
  and step size $h\in(0,h_{\max})$,
  the solution $y_{k+1}$ at $t_{k+1}=t_k+h$ is determined by
  \begin{equation}\label{eq_29}
    y_k = \psi_{\bm\beta}(y_{k+1},t_{k+1},-h).
  \end{equation}
  This is a nonlinear system of algebraic equations involving 
  the NN function $\psi_{\bm\beta}$, which needs to be solved for $y_{k+1}$.
  This system can be solved using the scipy routine ``root''
  (scipy.optimize.root) in the implementation,
  in which the Jacobian matrix can be computed
  by automatic differentiation.
  For NN algorithms based on the implicit
  representations~\eqref{eq_22b} and~\eqref{eq_22c}, time marching takes the following forms:
  given $(y_k,t_k)$ and $h$,
    \begin{align}
      & K = f(y_k+hK, t_k+h), \label{eq_30a} \quad
      y_{k+1} = y_0 + hK + h^2\varphi_{\bm\beta}(y_k,t_k,h);
    \end{align}
  and 
  \begin{subequations}
    \begin{align}
      & K_1 = f(y_k + \gamma h K_1, t_k+\gamma h), \label{eq_31a} \quad
      K_2 = f(y_k + (1-\gamma)hK_1 + \gamma hK_2, t_k+h), 
      \\
      & y_{k+1} = y_k + (1-\gamma)hK_1 + \gamma hK_2
      + h^3\varphi_{\bm\beta}(y_k,t_k,h). \label{eq_31c}
    \end{align}
  \end{subequations}
  Here $\varphi_{\bm\beta}$ denotes the
  function $\varphi(y_0,t_0,\xi)$ in~\eqref{eq_22b} and~\eqref{eq_22c}
  learned by the neural network.
  Nonlinear equations need to be solved for
  computing $K$, $K_1$,
  and $K_2$ during time marching. They can be solved
  based on the scipy routine ``root'' in the implementation.
  
\end{remark}


\subsection{Learning Exact Time Marching Algorithm for Autonomous Systems}

If system~\eqref{eq_1} is autonomous, the algorithm presented in
the previous subsection will be simplified.
We briefly discuss this case here.
Consider an autonomous system in~\eqref{eq_1}, i.e.
\begin{align}
  & \frac{dy}{dt} = f(y), \label{eq_5}
\end{align}
where $f: \mbb R^n\rightarrow \mbb R^n$ is a prescribed function,
together with the initial condition~\eqref{eq_1b}.

By again introducing the transformation~\eqref{eq_6}, we re-write the
system consisting of~\eqref{eq_5} and~\eqref{eq_1b} into
\begin{subequations}\label{eq_32}
  \begin{align}
    & \frac{dY}{d\xi} = f(Y), \\
    & Y(0) = y_0.
  \end{align}
\end{subequations}
The solution to problem~\eqref{eq_32} depends only on $y_0$ and
$\xi$, and we re-write it as
$  
  Y(\xi) = \psi(y_0,\xi)
$  
to make the dependence on $y_0$ and $\xi$ explicit.
$\psi(y_0,\xi)$ is the algorithmic function for the system consisting of~\eqref{eq_5}
and~\eqref{eq_1b}. This function is determined by,
\begin{subequations}\label{eq_34}
  \begin{align}
    &
    \frac{\partial\psi}{\partial\xi} = f(\psi(y_0,\xi)),
    \label{eq_34a} \\
    & \psi(y_0,0) = y_0. \label{eq_34b}
  \end{align}
\end{subequations}
We learn the function $\psi(y_0,\xi)$, for $(y_0,\xi)\in\Omega\times[0,h_{\max}]$,
by solving the system~\eqref{eq_34} with ELM,
for prescribed $\Omega$ and $h_{\max}$.
The trained NN can be used to solve the problem consisting of~\eqref{eq_5}
and~\eqref{eq_1b}.


\begin{figure}
  \centerline{
    \includegraphics[width=3in]{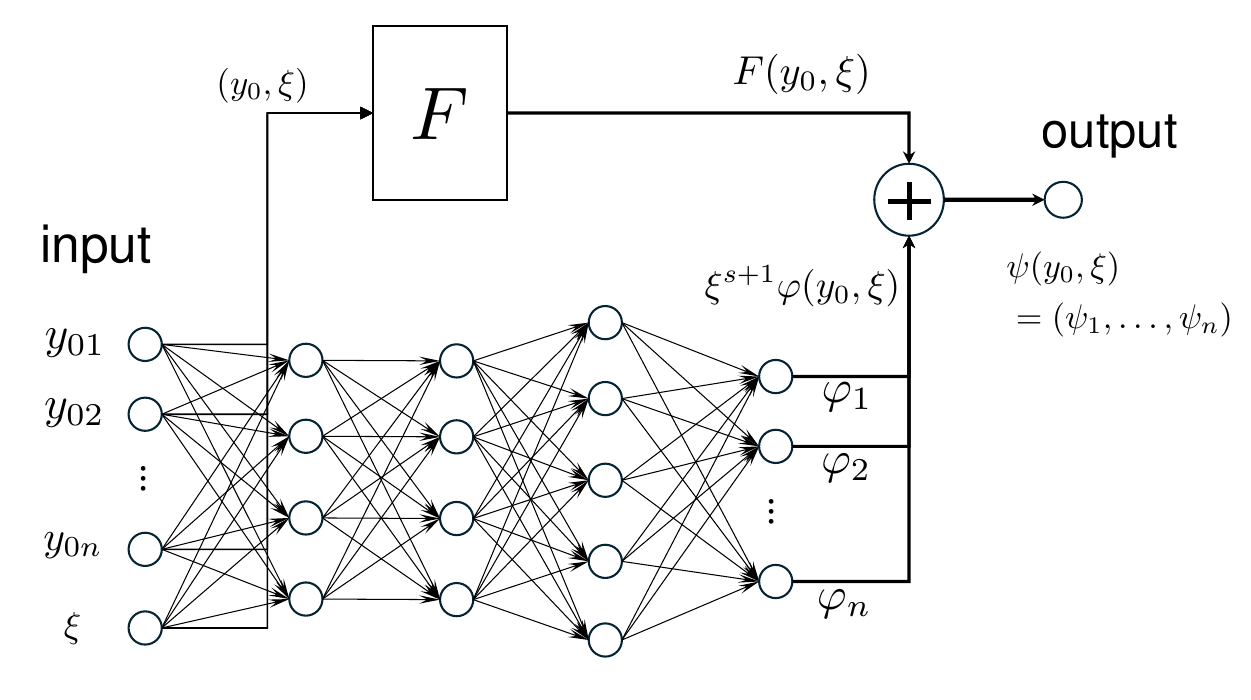}
  }
  \caption{Autonomous system: neural network structure,
    with three hidden layers shown
    as an example.
  }
  \label{fg_2}
\end{figure}

Figure~\ref{fg_2} illustrates the NN structure 
employed to learn $\psi(y_0,\xi)$,
which implements the following representation,
\begin{equation}\label{eq_35}
  \psi(y_0,\xi) = F(y_0,\xi) + \xi^{s+1}\varphi(y_0,\xi),
\end{equation}
where $F(y_0,\xi)\in\mbb R^n$ is a prescribed
$s$-th order approximation of $\psi(y_0,\xi)$ (for some integer $s$),
and $\varphi(y_0,\xi)\in\mbb R^n$
is an unknown function represented by the ELM-type
randomized NN.
The NN input nodes  represent
$(y_0,\xi)\in\mbb R^n\times\mbb R$, and the output
nodes represent $\psi(y_0,\xi)\in\mbb R^n$.
We again refer to the portion of the network
representing $\varphi(y_0,\xi)$ as the
$\varphi$-subset.
The $\varphi$-subset has an architecture characterized by~\eqref{eq_10},
with $m_0=n+1$ and $m_L=n$ and its hidden-layer coefficients
randomly assigned and fixed.

Analogous to~\eqref{eq_12} and~\eqref{eq_12c}, we consider the following
representations for $\psi(y_0,\xi)$:
\begin{subequations}\label{eq_36}
  \begin{align}
    & \label{eq_36a}
  \begin{array}{lll}
    s=0: & F(y_0,\xi) = y_0, & \psi(y_0,\xi) = y_0 + \xi\varphi(y_0,\xi);
  \end{array}  \\
  & \label{eq_36b}
  \begin{array}{lll}
    s=1: & F(y_0,\xi) = y_0+\xi f(y_0),
    & \psi(y_0,\xi) = y_0 + \xi f(y_0) + \xi^2\varphi(y_0,\xi);
  \end{array}
\end{align}
\end{subequations}
and for some problems also
\begin{align}
  & \label{eq_36c}
  \begin{array}{ll}
    s = 2: & F_1(y_0,\xi) = f(y_0), \quad F_2(y_0,\xi) = f(y_0+\frac{\xi}{2} F_1), \\
    & F(y_0,\xi) = y_0 + \xi F_2(y_0,\xi), \quad
    \psi(y_0,\xi) = y_0 + \xi F_2(y_0,\xi) + \xi^3\varphi(y_0,\xi).
  \end{array}
\end{align}
These forms automatically satisfy~\eqref{eq_34b}.
Therefore, only~\eqref{eq_34a} needs to be considered for NN training.


We train the NN such that $\psi(y_0,\xi)$
satisfies~\eqref{eq_34a}
for $y_0\in\Omega$ and $\xi\in[0,h_{\max}]$.
The procedure from Section~\ref{sec_222}  can be
adapted to train the NN here  for autonomous
systems.
The logic of the $\varphi$-subset's output layer gives rise to
a relation analogous to~\eqref{eq_13},
\begin{equation}\label{eq_37}
  \varphi_i(y_0,\xi) = \bm\beta_i\cdot\phi(y_0,\xi), \quad
  1\leqslant i\leqslant n,
\end{equation}
where $\varphi(y_0,\xi)=(\varphi_1,\dots,\varphi_n)\in\mbb R^n$,
$\phi(y_0,\xi)=(\phi_1,\dots,\phi_M)\in\mbb R^M$,
$M=m_{L-1}$ denoting the number of nodes in the last hidden layer
of the $\varphi$-subset, and
$\bm\beta_i = (\beta_{i1},\dots,\beta_{iM})\in\mbb R^M$
for $1\leqslant i\leqslant n$.
Note that $\varphi(y_0,\xi)$ is the output field of
the $\varphi$-subnet and $\phi(y_0,\xi)$ denotes
the output fields of the last hidden layer of the
$\varphi$-subnet. $\beta_{ij}$ ($1\leqslant i\leqslant n$,
$1\leqslant j\leqslant M$) are the training parameters of the network.


We determine the training parameters
$
\bm\beta = (\bm\beta_1,\dots,\bm\beta_n)
=(\beta_{11},\dots,\beta_{1M},\beta_{21},\dots,\beta_{nM})
\in\mbb R^{nM}
$
by the nonlinear least squares method.
Let $(y_0^{(i)},\xi^{(i)})$ ($1\leqslant i\leqslant Q$)
denote $Q$ random collocation points on
$\Omega\times[0,h_{\max}]$ drawn from
a uniform distribution. 
Enforcing~\eqref{eq_34a} on these collocation points
leads to the following nonlinear algebraic system about $\bm\beta$,
\begin{equation}\label{eq_38}
  r^{(i)}(\bm\beta) =
  \left.\frac{\partial\psi}{\partial\xi}\right|_{(y_0^{(i)},\xi^{(i)})}
  - f(\psi(y_0^{(i)},\xi^{(i)})) = 0,
  \quad 1\leqslant i\leqslant Q,
\end{equation}
in which $\psi$ is given by~\eqref{eq_36}-\eqref{eq_36c} and~\eqref{eq_37}
and the dependence of the residual $r^{(i)}\in\mbb R^n$ on $\bm\beta$
is made explicit.
This is a system of $nQ$ nonlinear algebraic equations
about $nM$ unknowns. We seek a least squares
solution to this system and compute $\bm\beta$
by the NLLSQ-perturb algorithm~\cite{DongL2021,DongW2023}.
NLLSQ-perturb requires the computation of the residual vector
and the Jacobian matrix for arbitrary given $\bm\beta$,
as noted previously.
Computing the residual vector
$r(\bm\beta) = (r^{(1)},\dots,r^{(Q)})\in\mbb R^{nQ}$
for NLLSQ-perturb is straightforward in light of~\eqref{eq_38},
noting that $\psi(y_0^{(i)},\xi^{(i)})$ and
$\left.\frac{\partial\psi}{\partial\xi}\right|_{(y_0^{(i)},\xi^{(i)})}$
therein can be obtained by forward NN evaluations
and by automatic differentiation.
%
The Jacobian matrix is given by
\begin{equation}
  \frac{\partial r}{\partial\bm\beta} = \begin{bmatrix}
    \frac{\partial r^{(1)}}{\partial\bm\beta_1} & \dots & \frac{\partial r^{(1)}}{\partial\bm\beta_n} \\
    \vdots & \ddots & \vdots \\
    \frac{\partial r^{(Q)}}{\partial\bm\beta_1} & \dots & \frac{\partial r^{(Q)}}{\partial\bm\beta_n}
  \end{bmatrix} \in\mbb R^{nQ\times nM},
  \quad
  \frac{\partial r^{(i)}}{\partial\bm\beta_j} =
  \begin{bmatrix}
    \frac{\partial r_1^{(i)}}{\partial\bm\beta_j} \\
    \vdots \\ \frac{\partial r_n^{(i)}}{\partial\bm\beta_j}
  \end{bmatrix} \in\mbb R^{n\times M},
\end{equation}
where
\begin{equation}
  \begin{split}
  \frac{\partial r^{(i)}_k}{\partial\beta_j}
  =&\ (s+1)(\xi^{(i)})^s\phi(y_0^{(i)},\xi^{(i)})\delta_{kj}
  + (\xi^{(i)})^{s+1}\left.\frac{\partial\phi}{\partial\xi} \right|_{(y_0^{(i)},\xi^{(i)})} \delta_{kj} \\
  & - (\xi^{(i)})^{s+1} \left.\frac{\partial f_k}{\partial\psi_j} \right|_{(y_0^{(i)},\xi^{(i)})} \phi(y_0^{(i)},\xi^{(i)})
  \ \ \in\mbb R^{1\times M},
  \quad 1\leqslant i\leqslant Q, \ 1\leqslant k, j\leqslant n.
  \end{split}
\end{equation}
%
The terms $\phi(y_0^{(i)},\xi^{(i)})$ and
$\left.\frac{\partial\phi}{\partial\xi} \right|_{(y_0^{(i)},\xi^{(i)})}$
involved in the Jacobian matrix can be computed
by forward evaluations of a network sub-model  and
by automatic differentiations (see Remark~\ref{rem_2}).


The trained NN  contains the learned algorithm for
solving the system consisting of~\eqref{eq_5} and~\eqref{eq_1b},
with any initial condition $y_0\in\Omega$, initial time $t_0\in\mbb R$, and
step size $h\in(0,h_{\max})$.
Suppose $(y_k,t_k)$ provide the solution and the time
at step $k$ ($k\geqslant 0$),
and that $\psi_{\bm\beta}(y_0,\xi)$ denotes the learned algorithmic function. 
The solution at the new time step is given by ($h$ denoting the step size),
\begin{equation}\label{eq_41}
  y_{k+1} = \psi_{\bm\beta}(y_k,h),
  \quad t_{k+1} = t_k+h,
\end{equation}
which involves the forward evaluation of the neural network.

\begin{remark}
  The discussion on implicit representations of $\psi$
  for non-autonomous systems in Remarks~\ref{rem_24} and~\ref{rem_27}
  can be adapted to autonomous
  systems for learning and using $\psi(y_0,\xi)$.
  The equations~\eqref{eq_21}
  and~\eqref{eq_22} and other related expressions
  need to be modified accordingly to exclude the $t_0$ effects.

\end{remark}


\section{Computational Examples}
\label{sec_tests}

In this section we evaluate the performance of the learned time integration algorithms 
from Section~\ref{sec_method} by simulating the dynamics of
several non-autonomous and autonomous systems.
Some of these systems are chaotic, or can become
stiff for a range of problem parameters.
In particular, we compare the learned NN algorithms with
the leading traditional time integration algorithms
as implemented in the Scipy library (scipy.integrate.solve\_ivp routine,
with methods ``DOP853'', ``RK45'', ``RK23'', ``Radau'' and ``BDF'')
in terms of the solution accuracy and the time marching cost.
All the scipy methods are adaptive in the step size and in the integration order,
with ``DOP853'', ``RK45'', and ``RK23'' being explicit and
``Radau'' and ``BDF'' being implicit schemes.
With the learned NN algorithms, we employ a constant time step size
for the majority of
simulations, and a quasi-adaptive time step for
certain stiff cases.

We define the maximum error $e_{\max}$ and the root-mean-squares (rms)
error $e_{\text{rms}}$ of a solution
$y(t)=(y_1(t),y_2(t),\dots,y_n(t))\in\mbb R^n$ for $t\in[t_0,t_f]$ as follows,
  \begin{equation}
    e_{\max} = \max\{\ \max\{\ |y_i^c(t_j) - y_i^{ex}(t_j)|\  \}_{j=1}^m  \ \}_{i=1}^n, \quad
    e_{\text{rms}} = \sqrt{\frac{1}{mn}\sum_{i=1}^n\sum_{j=1}^m|y_i^c(t_j) - y_i^{ex}(t_j)|^2 }.
  \end{equation}
Here $y^c$ is the numerical solution obtained by the NN algorithms (or
the scipy methods), $y^{ex}$ is the exact solution (if available)
or a reference solution, and $t_j\in[t_0,t_f]$ ($1\leqslant j\leqslant m$)
denotes the time instants resulting from the
time marching algorithm corresponding to a constant or a quasi-adaptive
time step. When the exact solution is unavailable, the reference solution
is computed using the scipy method ``DOP853'' (for non-stiff problems)
or ``Radau'' (for stiff problems) with absolute tolerance $10^{-16}$
and relative tolerance $10^{-13}$.  


For the learned NN algorithms, the time-marching cost lies primarily in
the forward evaluation of the trained NN (see~\eqref{eq_27})
at every time step. Therefore, efficient NN evaluation
is critical to the performance of the current NN algorithms.
As mentioned before, our implementation 
is based on the Tensorflow and Keras libraries.
We observe that the built-in NN evaluation methods from the Keras library, such as
the direct call against a Keras model (or the ``predict()'' method),
even in the graph mode,
induces a significant overhead, slowing down the NN time marching.

\begin{table}[tb]
  \centering
  \begin{tabular}{lcc}
    \hline
    $M$ & current-NN-evaluation (seconds) & keras-NN-evaluation (seconds) \\
    $400$ & $0.00138$ & $0.0258$ \\
    $600$ & $0.00148$ & $0.0270$ \\
    $800$ & $0.00155$ & $0.0271$ \\
    \hline
  \end{tabular}
  \caption{Comparison of the time-marching cost (wall time) of
    NN-Exp-S1 employing (i) the current customized NN evaluation method, and
    (ii) the Keras built-in NN evaluation method (in graph mode),
    for the test problem from Section~\ref{sec_lin_model} ($\lambda=100$).
    The tests correspond to those in Figure~\ref{fg_5}.
    $M$ is the number of training parameters.
    The built-in evaluation method calls the Keras model 
    directly with the input data.
  }
  \label{tab_1}
\end{table}

To reduce the NN evaluation overhead, we have employed a customized
evaluation method for the trained NN  using routines from
the numpy library. Here is the main idea for the customized evaluation.
After the neural network is trained,
we extract the weight and bias coefficients, as well as the activation functions,
for all layers from the trained Keras model. Then we form a routine in which
the extracted data are employed to implement the NN logic using
plain numpy functions. This routine is used as the customized NN evaluation
method for time marching.

%
We would like to mention
another point in our implementation, regarding the computation of $f(y,t)$
for given $(y,t)$ ($y\in\mbb R^n$) in the customized NN evaluation method, which is needed
for computing the $F(y_0,t_0,\xi)$ term in~\eqref{eq_11}.
We have implemented a routine for computing $f(y,t)$, used
specifically for time marching, in which the input data $(y,t)$ is
a vector with shape $(n+1,)$ (in the Python notation). Note that the computation of
$f(y,t)$ is also required during NN training, where the input data is a matrix
with shape ($N_c,n+1$) that represents $(y,t)$ on $N_c$ collocation points.
In the version specifically for time marching,
the simpler vectorial input data of $(y,t)$ allows
more efficient implementation for computing $f(y,t)$.

The learned algorithms
employing the above customized NN evaluation method is significantly faster
than those based on the built-in evaluation methods from Keras.
This point is demonstrated by Table~\ref{tab_1}, which shows the time-marching time
of the NN-Exp-S1 algorithm (see Remark~\ref{rem_c26})
for the test problem of Section~\ref{sec_lin_model}
employing the customized NN evaluation method and the Keras built-in method.
The algorithm using the customized NN evaluation is faster by more
than an order of magnitude.


When implementing the $\varphi$-subnet for representing $\varphi(y_0,t_0,\xi)$,
we have added a normalization between the input layer
(representing $(y_0,t_0,\xi)$) and the first hidden layer. The normalization implements
an affine transform for each component of $(y_0,t_0,\xi)$, transforming
the input $(y_0,t_0)$ data
from $[a_1,b_1]\times\dots\times[a_n,b_n]\times[T_1\times T_2]\subset\mbb R^n\times\mbb R$ to
the standard domain $[-1,1]^n\times[-1,1]$ and
the input $\xi$ data from $[0,h_{\max}]$ to the interval
$[0,\frac{1}{\delta_m}]$. Here the prescribed constant $\delta_m>0$ will be referred
to as the $\xi$-domain map factor, with $\delta_m=1$ by default.
For some stiff problems, we find it more favorable to employ a smaller $\delta_m$
with the NN algorithms. These $\delta_m$ values will be provided when discussing
the specific test problems. The use of $\delta_m$ essentially stipulates
that the normalization should map the interval $[0,\delta_m h_{\max}]$
to the standard domain $[0,1]$ for $\xi$.
In the simulations  we will
employ $\delta_m=1$, unless otherwise specified.
The normalization procedure discussed here
is implemented using a ``lambda'' layer from the Keras library.
When domain decomposition is present, the above normalization is implemented
for the local NN  on each sub-domain.
In the simulations we employ no domain decomposition by default,
unless otherwise specified. 

When the training domain is decomposed along some direction,
in Remark~\ref{rem_26} we have discussed enlarging the sub-domains
by a factor $r$ along this direction during network training.
In the following numerical tests
we employ an enlargement factor $r=0$ (i.e.~no enlargement) by default,
unless otherwise specified.

\subsection{A Linear System: Non-Stiff and Stiff Cases}
\label{sec_lin_model}

\begin{table}[tb]
  \centering
  \begin{tabular}{l|l|l}
    \hline
   $\lambda=100$ & domain: $(y_{0},t_0,\xi)\in [-1.1,1.1]\times[-0.05,1.05]\times[0,h_{\max}]$
    & NN: $[3, M, 1]$ ($M$ varied) \\
    & no domain decomposition & activation function: Gaussian \\
   & $r$: $0.0$  & $\delta_m$: $1$  \\
   & $Q$: varied & $R_m$: to be specified \\
    & $\Delta t$: $0.02$ (for time marching) & time: $t\in[0,1]$  \\
    \hline
    $\lambda=10^6$ & domain: $(y_{0},t_0,\xi)\in [-1.1,1.1]\times[-0.05,1.05]\times[0,h_{\max}]$
    & NN: $[3, M, 1]$ ($M$ varied) \\
    & sub-domains: 2 along $t_0$ (uniform) & activation function: Gaussian \\
    & $r$: $0.05$ along $t_0$ & $\delta_m$: $0.02$ or $0.01$  \\
    & $Q$: $1000$ or $900$ (random) & $R_m$: to be specified \\
    & $\Delta t$: $0.02$ (for time marching) &  time: $t\in[0,1]$ \\
    \hline
  \end{tabular}
  \caption{NN simulation parameters for the linear system
    (Section~\ref{sec_lin_model}).
    Values for some parameters are specified in the text.
  }
  \label{tab_a2}
\end{table}

In this test we consider a linear non-autonomous  problem,
\begin{subequations}\label{eq_53}
  \begin{align}
    &
    \frac{dy}{dt} = -\lambda[y - \cos(\pi t)], \quad t\in[0,1], \\
    &
    y(t_0) = y_0,
  \end{align}
\end{subequations}
where $\lambda>0$ is a constant, $y(t)$ is the unknown
to be computed, $t_0=0$, and $y_0=0$ is the initial condition.
This problem becomes very stiff with large $\lambda$ values.
It has the following exact solution,
\begin{equation}\label{eq_54}
  y(t) = \left[y_0-\frac{\lambda^2}{\lambda^2+\pi^2}\cos(\pi t_0)
     - \frac{\lambda\pi}{\lambda^2+\pi^2}\sin(\pi t_0)
     \right] e^{-\lambda(t-t_0)}
  + \frac{\lambda^2}{\lambda^2+\pi^2}\cos(\pi t)
  + \frac{\lambda\pi}{\lambda^2+\pi^2}\sin(\pi t).
\end{equation}

\begin{figure}
  \centerline{
    \includegraphics[width=2.0in]{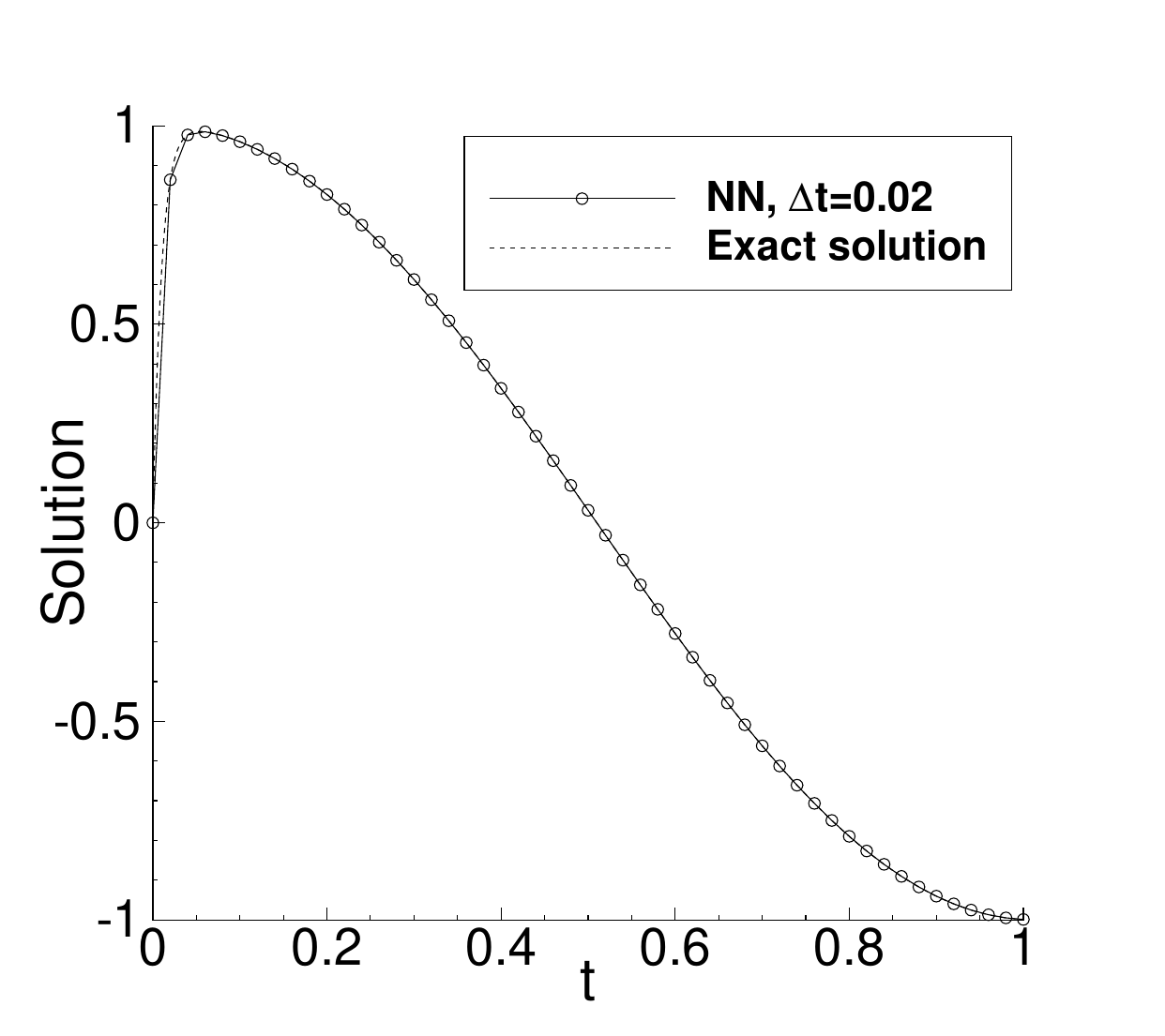}(a)
    \includegraphics[width=2.0in]{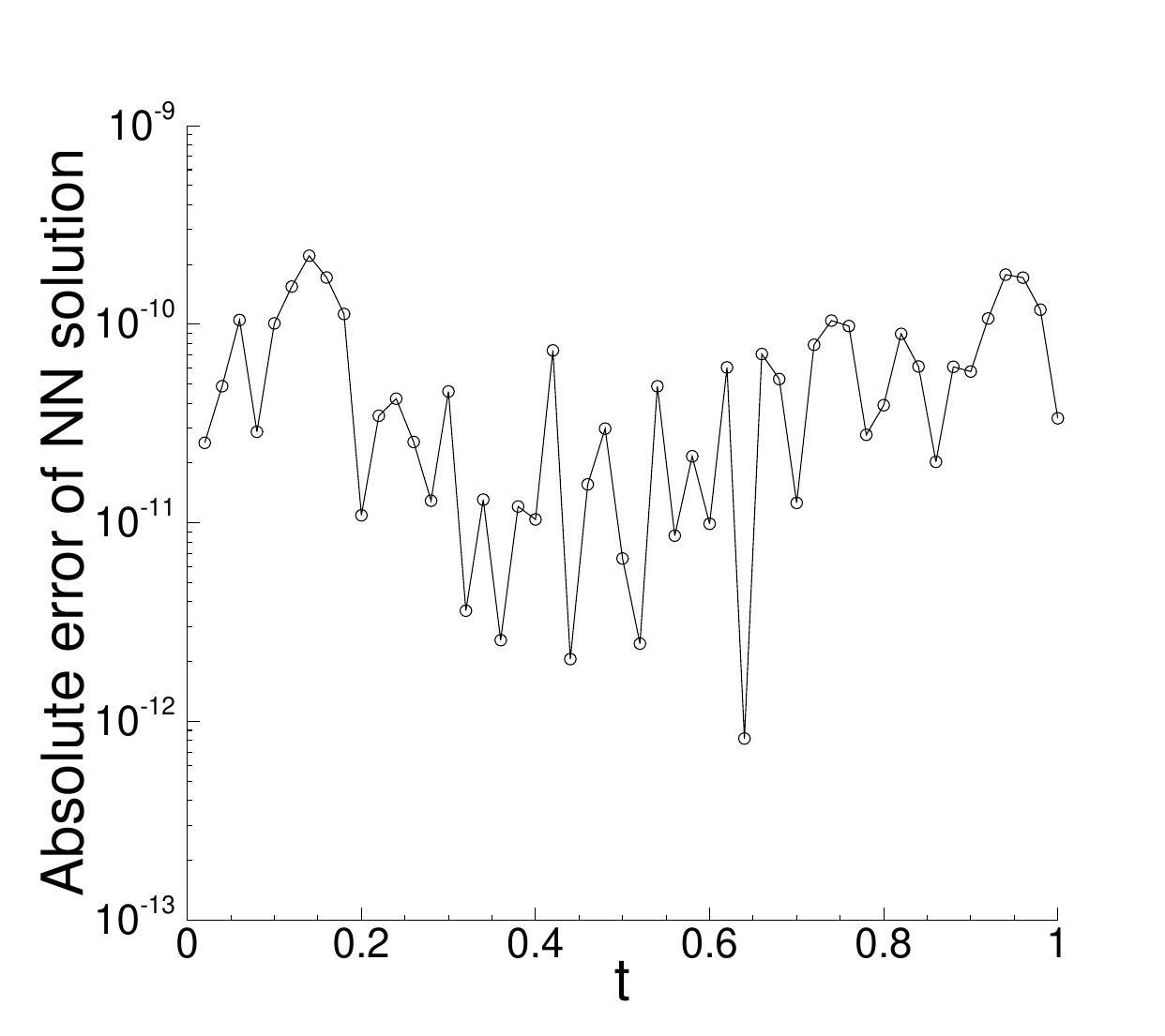}(b)
    \includegraphics[width=2.0in]{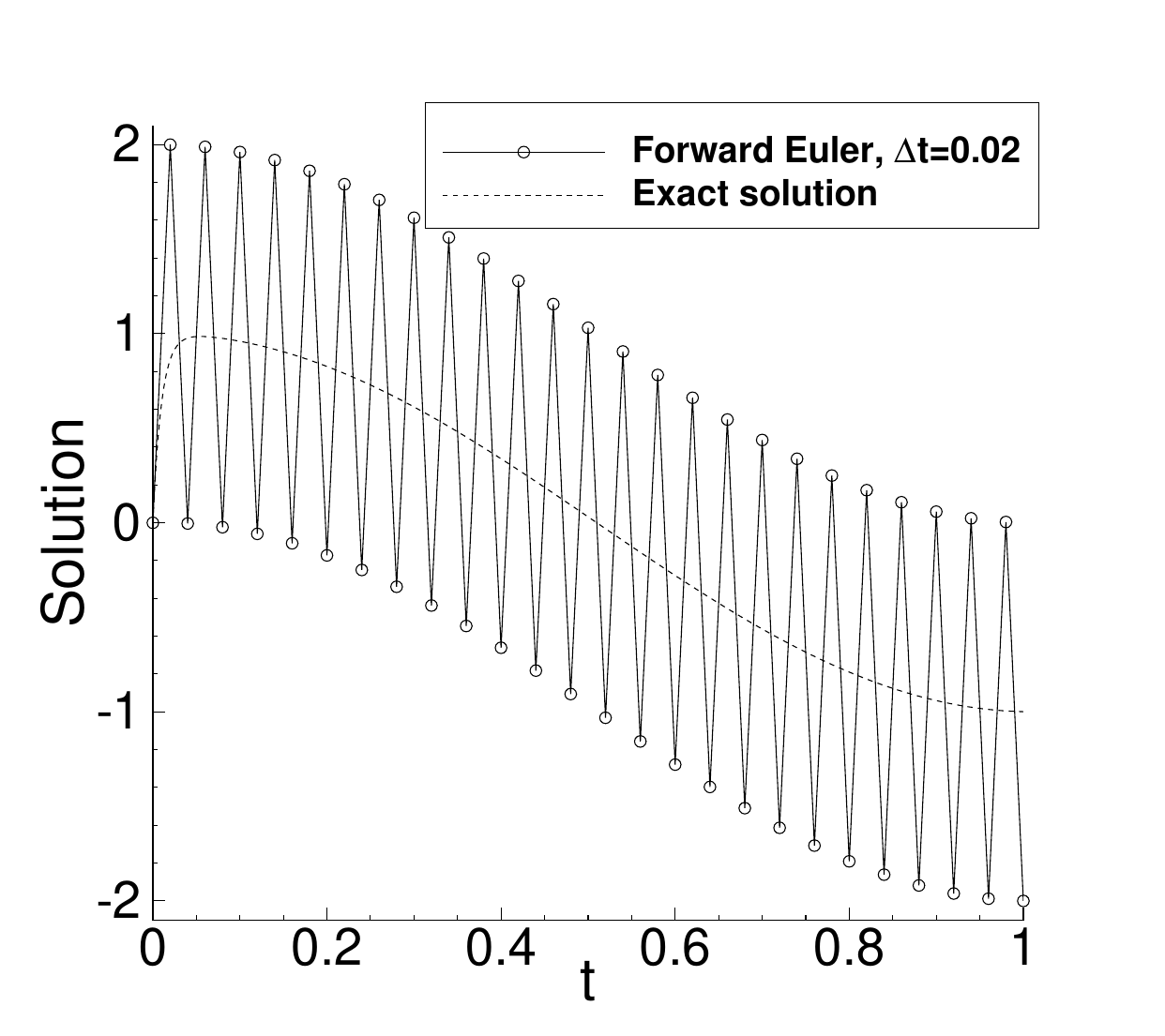}(c)
  }
  \caption{Linear model ($\lambda=100$):
    (a) Comparison of $y(t)$ between the NN-Exp-S1 solution and the exact solution.
    (b) Absolute error history of  the NN-Exp-S1 solution.
    (c) History of $y(t)$ obtained by the forward Euler method
    (i.e.~the $F(y_0,t_0,\xi)$ component in NN-Exp-S1, see~\eqref{eq_12b}).
    Training domain: $h_{\max}=0.03$,
    NN: [3, 800, 1], 
    $Q=2500$, and $R_m=0.5$. See Table~\ref{tab_a2} for the other parameter values.
  }
  \label{fg_3}
\end{figure}


We use an ELM network with architecture $[3,M,1]$
and the Gaussian activation function
$\sigma(x)=e^{-x^2}$ for the $\varphi$-subset
 to learn the  algorithmic function $\psi(y_0,t_0,\xi)$,
where $M$ is varied. The hidden-layer coefficients of the $\varphi$-subnet are
assigned to random values generated on the interval $[-R_m,R_m]$ drawn from
a uniform distribution, where the constant $R_m$ is specified below.
The NN is trained on a domain
$(y_0,t_0,\xi)\in[-1.1,1.1]\times[-0.05,1.05]\times[0,h_{\max}]$, where $h_{\max}$
is specified below. We employ $Q$ random collocation points,
where $Q$ is  varied,
uniformly drawn from the domain to train the NN  as discussed in
Section~\ref{sec_method}.
We employ the trained NN
to solve the problem~\eqref{eq_53}  for $t\in[0,1]$,
using a step size $\Delta t=0.02$ and the initial condition $(y_0,t_0)=(0,0)$.
The maximum and rms errors ($e_{\max}$, $e_{\text{rms}}$) of the numerical
solution against the exact solution~\eqref{eq_54} are then computed and
the time marching cost (wall time) of the NN algorithms is recorded for analysis.
Table~\ref{tab_a2} summarizes the simulation parameters related to the NN algorithms.
The same parameters will appear in the subsequent
test problems.


We first consider a non-stiff case, with $\lambda=100$ in~\eqref{eq_53}, and
illustrate the characteristics of the learned NN algorithms. 
Figures~\ref{fg_3}(a,b) provide an overview of the NN solution, comparing
the solution histories from the NN-Exp-S1 algorithm (see Remark~\ref{rem_c26})
and the exact solution~\eqref{eq_54} and showing the absolute-error history
of NN-Exp-S1. The parameter values for NN-Exp-S1 are listed in the figure caption
or in Table~\ref{tab_a2}.
The NN solution is highly accurate, with a maximum error on the order of $10^{-10}$
over $t\in[0,1]$.
Notice that the NN-Exp-S1 algorithm consists of two components (see~\eqref{eq_12b}),
$F(y_0,t_0,\xi)$ and $\varphi(y_0,t_0,\xi)$, with the former being the forward
Euler formula and the latter a neural network correction.
We find that the $\varphi(y_0,t_0,\xi)$ component is critical to
the NN-Exp-S1 accuracy. Without this component, the solution becomes very poor.
Figure~\ref{fg_3}(c) shows the history of $y(t)$ obtained by
the forward Euler method ($\Delta t=0.02$), i.e.~the $F(y_0,t_0,\xi)$ component
solely. The numerical solution is highly oscillatory, with no accuracy at all. 

\begin{figure}
  \centerline{
    \includegraphics[width=2in]{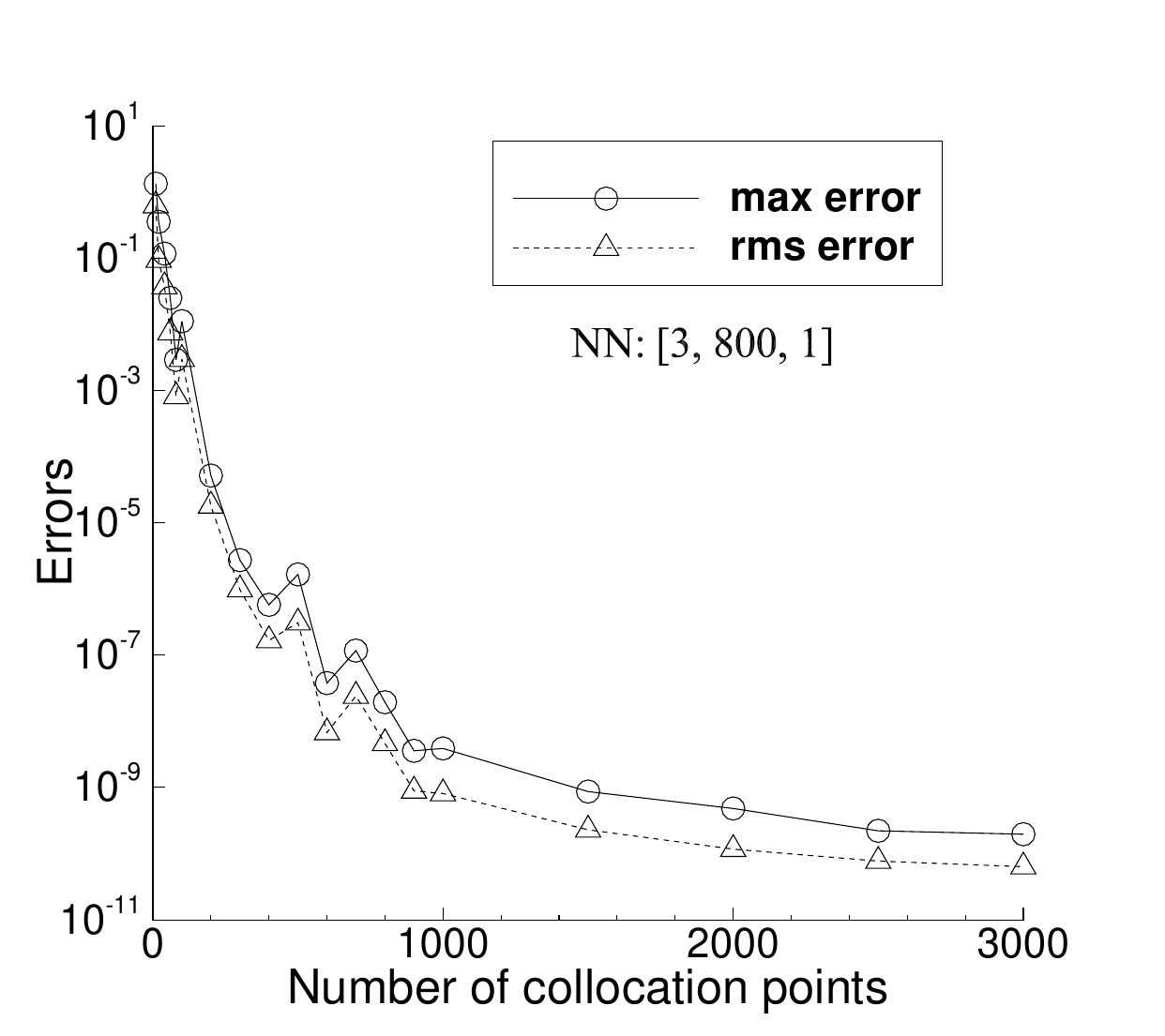}(a)
    \includegraphics[width=2in]{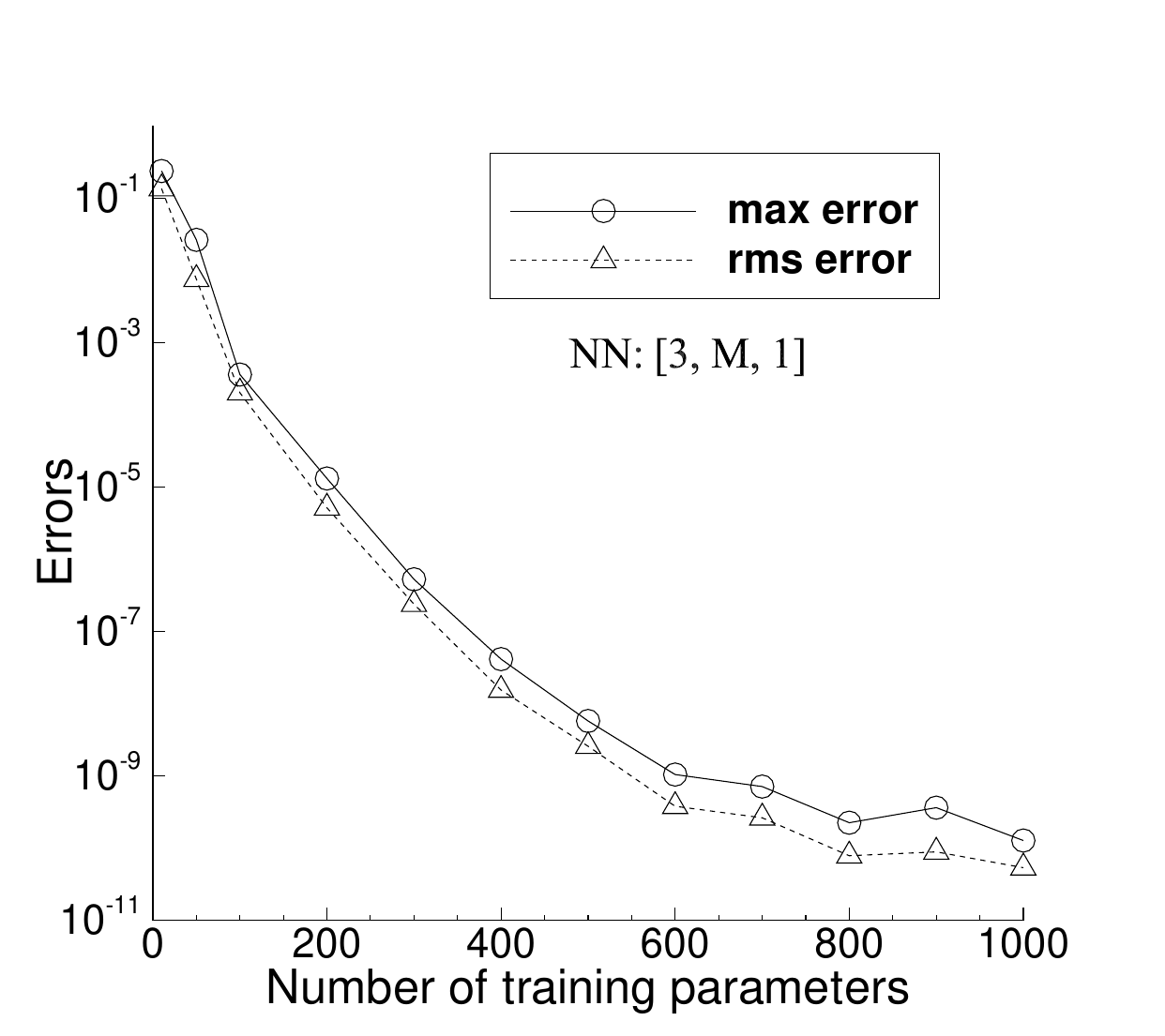}(b)
    \includegraphics[width=2in]{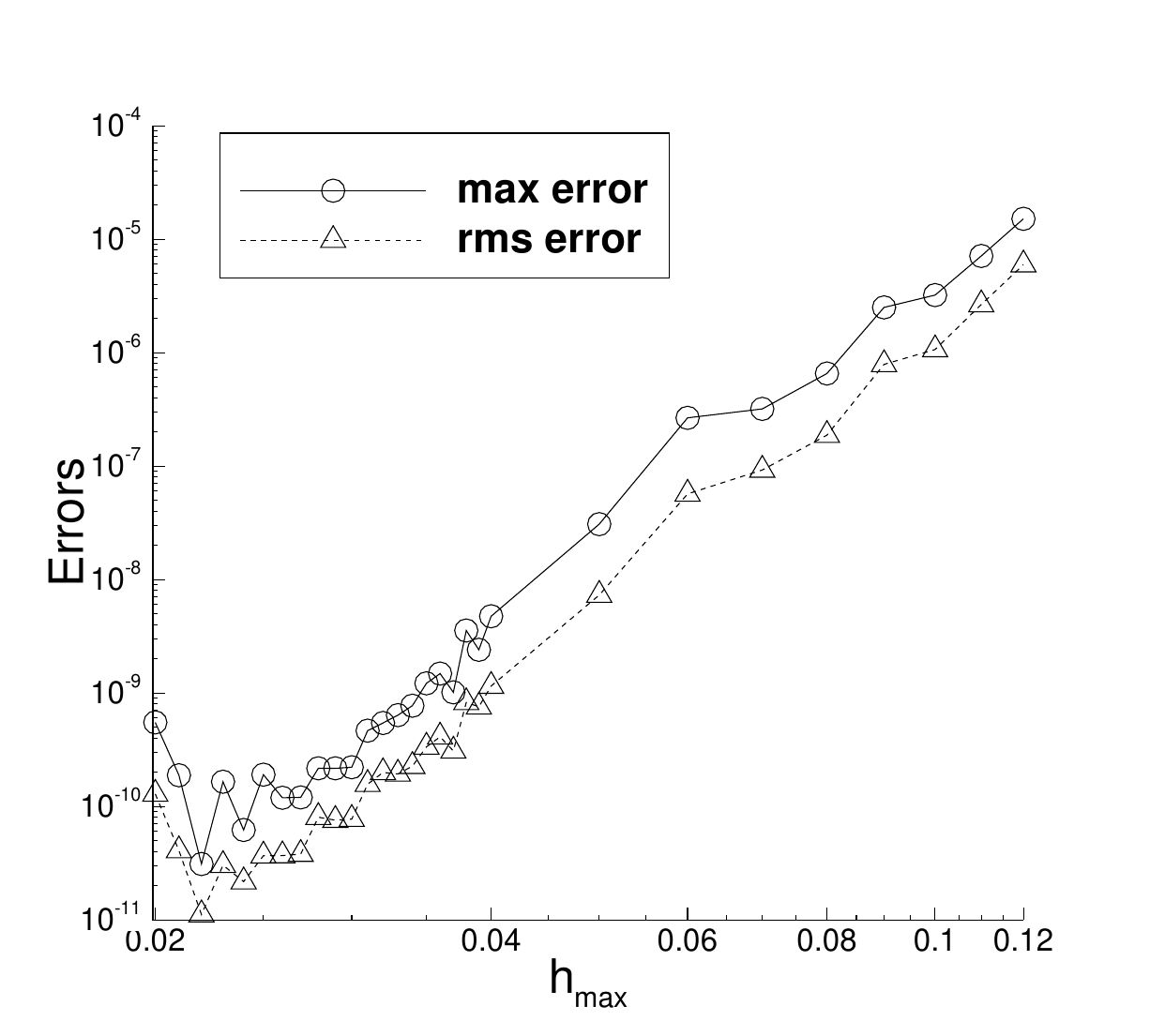}(c)
  }
  \caption{Linear model ($\lambda=100$):
    Time marching errors ($e_{\max}$, $e_{\text{rms}}$) versus (a) the number of training
    collocation points ($Q$), (b) the number of training parameters ($M$)
    in NN, and (c) the domain size $h_{\max}$ along $\xi$, obtained
    with the NN-Exp-S1 algorithm.
    $h_{\max}=0.03$ in (a,b) and is varied in (c);
    $M=800$ in (a,c) and is varied in (b);
    $Q=2500$ in (b,c) and is varied in (a);
    $R_m=0.5$; Other simulation parameters are given in Table~\ref{tab_a2}.
  }
  \label{fg_4}
\end{figure}

Figure~\ref{fg_4} illustrates the effects of several simulation parameters
on the attained NN algorithm.
It shows the maximum and rms time-marching errors of NN-Exp-S1
as a function of the number of collocation points ($Q$), the number of training parameters ($M$),
and the domain size $h_{\max}$ along $\xi$,
used for training the algorithm. The errors are
obtained on the points corresponding to a time step $\Delta t=0.02$ for $t\in[0,1]$.
The parameter values are specified in the figure caption or Table~\ref{tab_a2}.
Each plot represents a group of tests, in which one parameter is varied 
while the other parameters are fixed when training the NN-Exp-S1 network,
and the trained algorithm is employed in time marching to obtain the errors.
For example, in Figure~\ref{fg_4}(a) the number of training collocation points
$Q$ is varied systematically when training NN-Exp-S1, while the other parameters are fixed.
%
We observe that the NN errors decrease 
nearly exponentially with increasing number of collocation points ($Q$) or
 training parameters ($M$). 
 Using a smaller domain along $\xi$ (i.e.~smaller $h_{\max}$)
 generally leads to improved accuracy in
the attained NN algorithm.
Numerical experiments suggest that the accuracy of the NN algorithm seems to
be determined at the time of training. Once the network is trained,
the accuracy of the resultant NN algorithm will be fixed.
Varying the step size $\Delta t\in[0,h_{\max}]$ in time marching
using the trained NN algorithm appears to have little influence on the
time-marching error. This characteristic is somewhat different from traditional
numerical algorithms, whose accuracy generally
improves with decreasing time step size.

\begin{figure}
  \centerline{
    \includegraphics[width=2in]{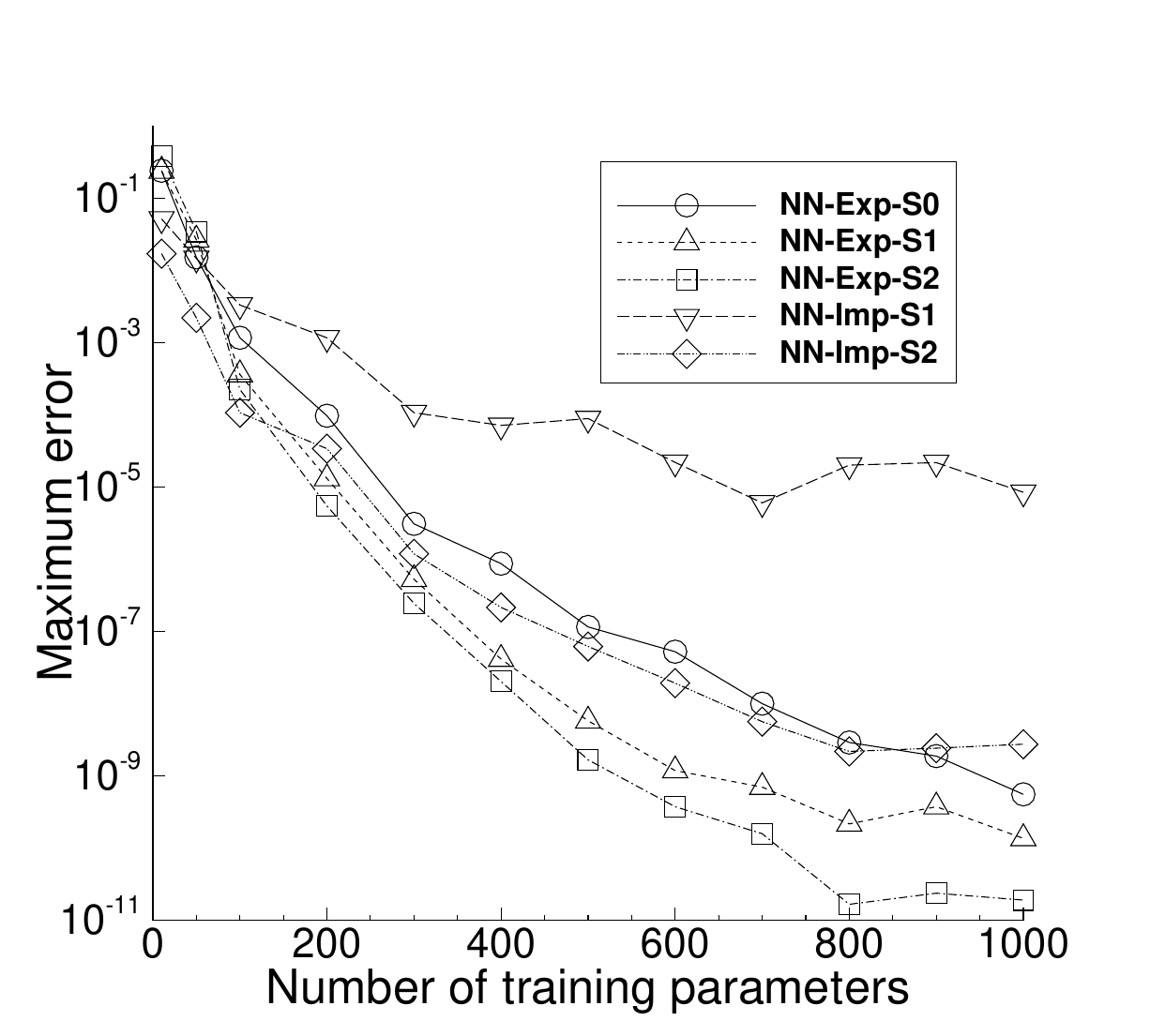}(a)
    \includegraphics[width=2in]{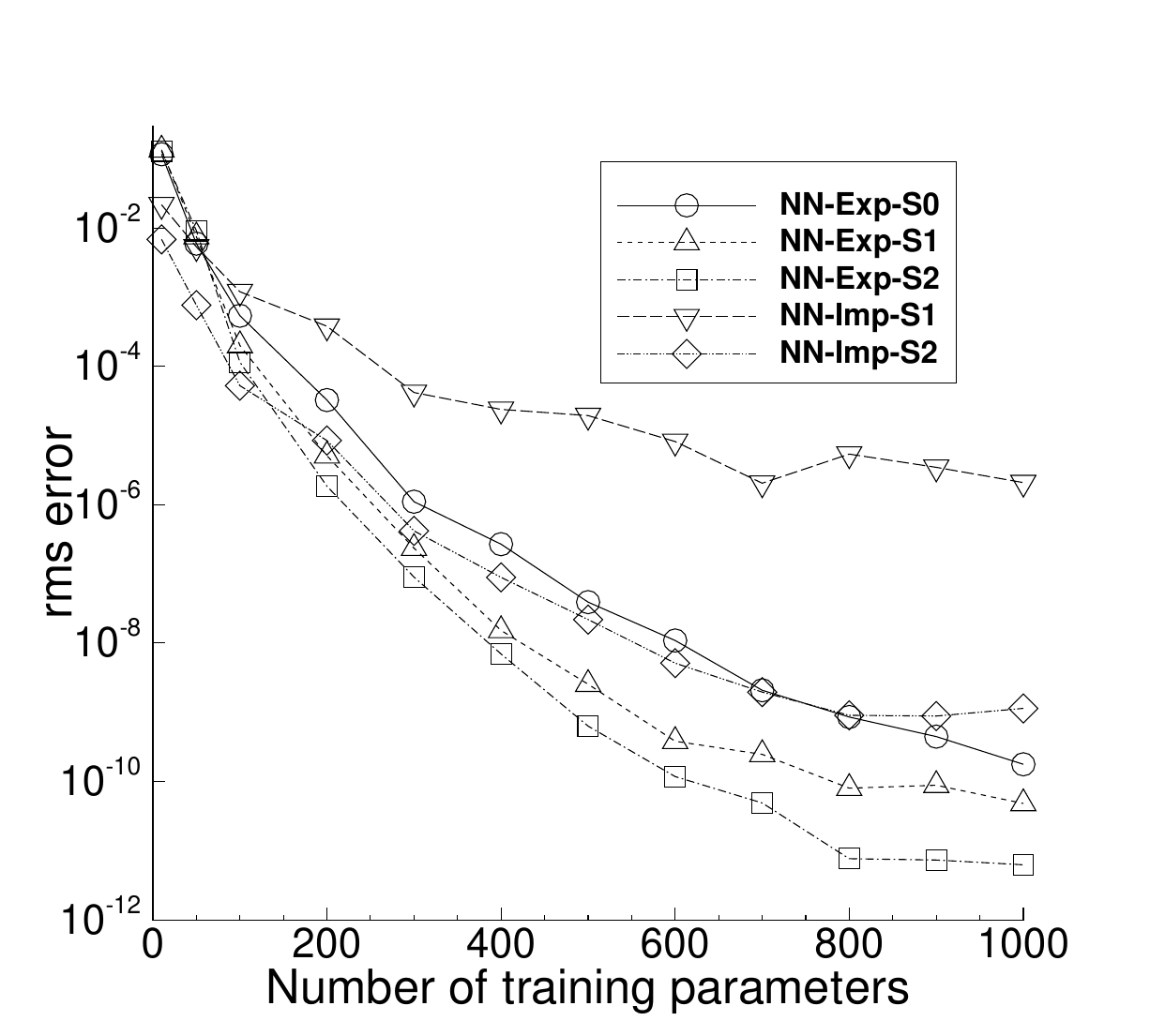}(b)
    \includegraphics[width=2in]{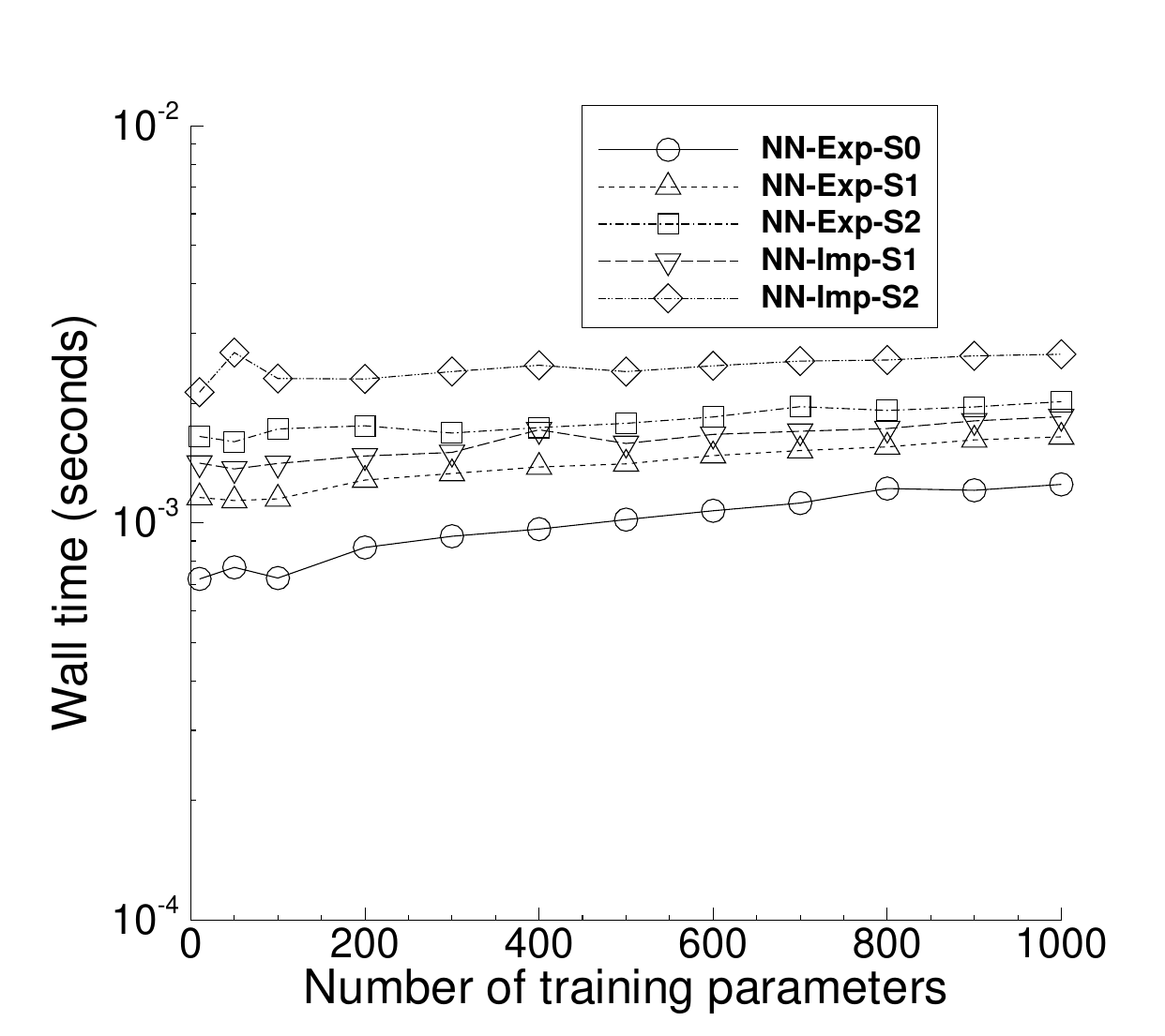}(c)
  }
  \caption{Linear model ($\lambda=100$):
    Comparison of  (a) the maximum and (b) the rms time-marching errors,
    and (c) the time-marching cost (wall time)
     versus the number of training parameters ($M$)
     among the explicit and implicit NN algorithms.
     $h_{\max}=0.03$, and $M$ is varied in the tests.
    NN-Exp-S0: $R_m=0.06$, $Q=2500$;
    NN-Exp-S1: $R_m=0.5$, $Q=2500$;
    NN-Exp-S2: $R_m=0.6$, $Q=3500$;
    NN-Imp-S1: $R_m=0.5$, $Q=2500$;
    NN-Imp-S2: $R_m=0.5$, $Q=3500$.
    Other simulation parameters are given in Table~\ref{tab_a2}.
  }
  \label{fg_5}
\end{figure}

A comparison of the accuracy and the time-marching cost (wall time)
of different NN algorithms with explicit (see~\eqref{eq_12}-\eqref{eq_12c})
and implicit (see~\eqref{eq_22})  formulations
is provided in Figure~\ref{fg_5}. These plots show the
maximum and rms errors and the time-marching time
as a function of the number of
training parameters ($M$) in the NN obtained by
NN-Exp-S0, NN-Exp-S1, NN-Exp-S2, NN-Imp-S1,
and NN-Imp-S2 (see Remark~\ref{rem_c26}).
The simulation parameters for this group of tests are specified in
the figure caption or in Table~\ref{tab_a2}.
Among the explicit NN algorithms, NN-Exp-S2 tends to be more
accurate than NN-Exp-S1, which in turn is generally more accurate
than NN-Exp-S0. In terms of time-marching cost 
NN-Exp-S0 is the fastest, followed by NN-Exp-S1 and then NN-Exp-S2.
The implicit NN algorithms 
are observed to be less accurate than the explicit ones.
For example, NN-Imp-S1 is the least accurate among these algorithms, with error levels
several orders of magnitude larger than those of the other algorithms,
and NN-Imp-S2 exhibits an accuracy comparable to NN-Exp-S0 for this problem.
The time marching cost of the implicit NN algorithms is notably higher than 
that of the explicit ones.
Since this test problem is linear, here the implicit NN algorithms do not actually
entail the solution of nonlinear algebraic equations during time marching.
For nonlinear problems in subsequent subsections
the implicit NN algorithms would be considerably slower compared with
the explicit ones, due to the need for solving nonlinear algebraic
systems during time integration.
Overall the implicit NN algorithms are not as competitive as the explicit ones.


\begin{figure}
  \centerline{
    \includegraphics[width=2in]{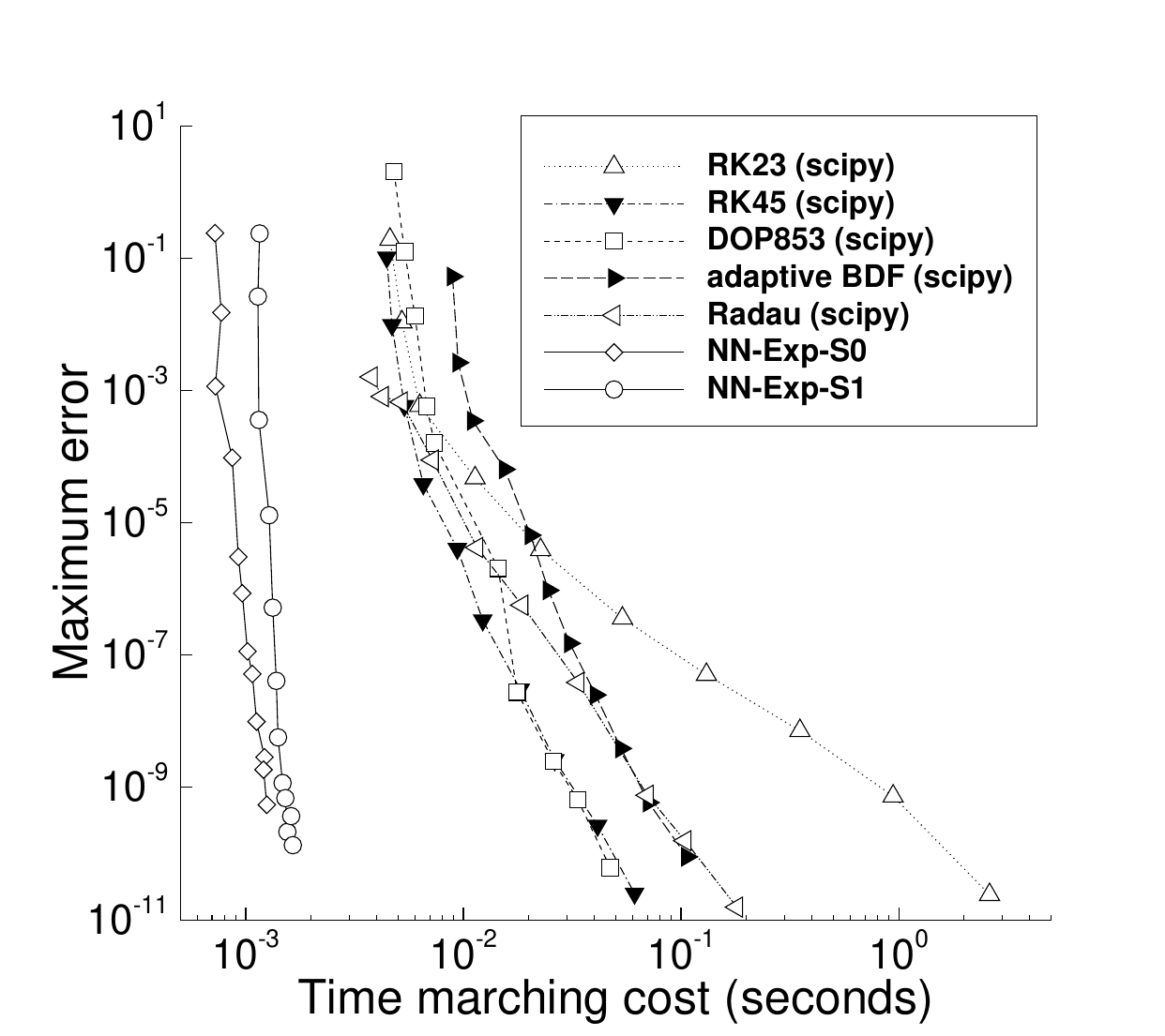}(a)
    \includegraphics[width=2in]{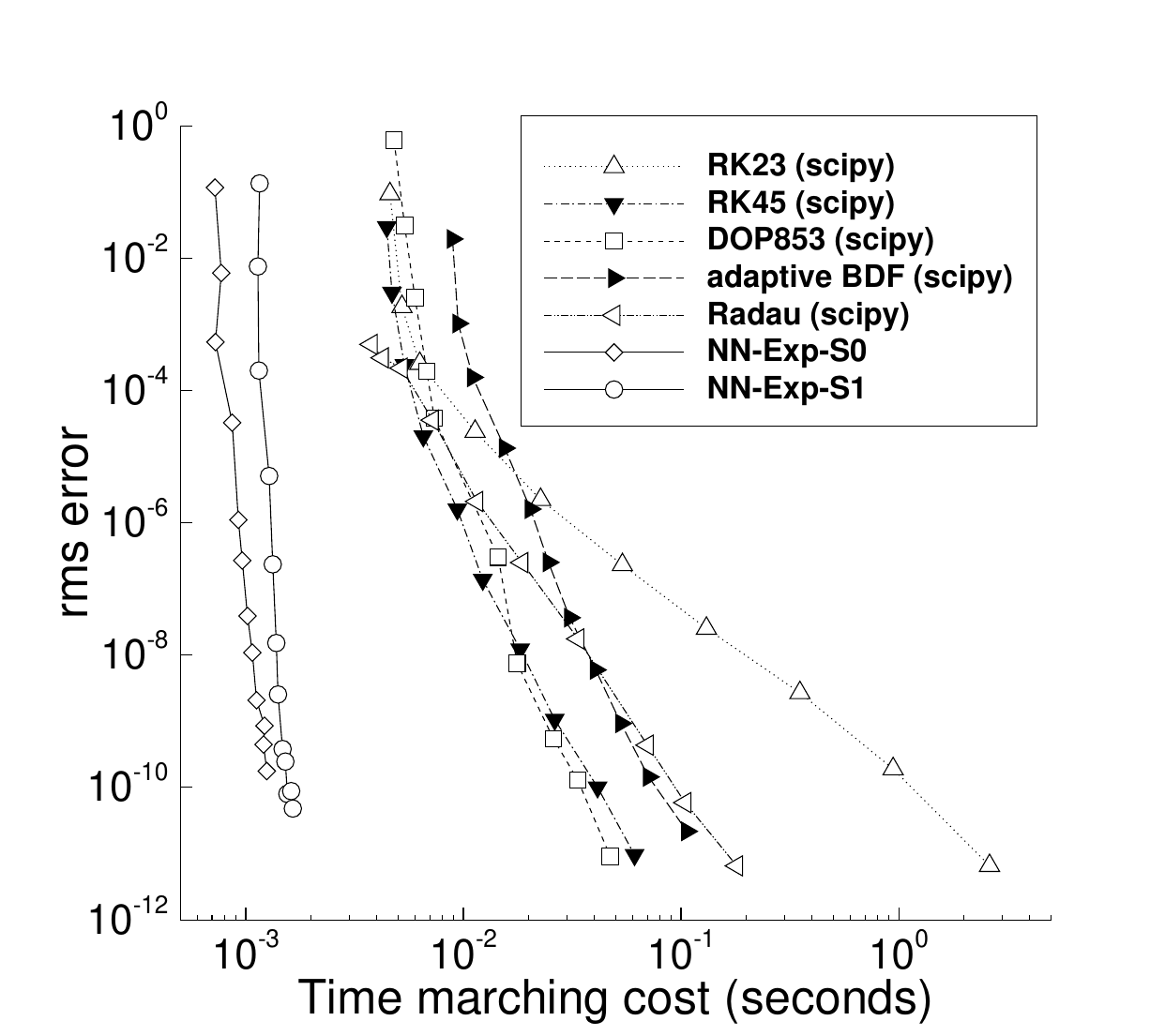}(b)
  }
  \caption{Linear model ($\lambda=100$):
    Comparison of
    (a) the maximum and (b) the rms time-marching errors 
    versus the time-marching cost (wall time)
    between the current NN algorithms (NN-Exp-S0 and NN-Exp-S1)
    and the scipy methods.
    The settings and parameters for NN-Exp-S0 and NN-Exp-S1
    follow those of Figure~\ref{fg_5}, with
    data points corresponding to different $M$ in the NN architecture.
    Scipy methods: adaptive time step, adaptive order, absolute tolerance $10^{-16}$,
    data points corresponding to different relative tolerance values,
    dense output on points corresponding to $\Delta t=0.02$ for $t\in[0,1]$.
  }
  \label{fg_6}
\end{figure}

Figure~\ref{fg_6} shows a comparison of the computational performance (accuracy
versus cost)
between the NN algorithms and the traditional time integration algorithms from
the scipy library.
It depicts the maximum and rms time-marching errors
as a function of the time-marching cost (wall time) obtained
by the current NN-Exp-S0 and NN-Exp-S1 algorithms
and several scipy methods, including ``RK23''
(explicit Runge-Kutta method or order 3(2))~\cite{BogackiS1989},
``RK45'' (explicit Runge-Kutta method of order 5(4))~\cite{DormandP1980},
``DOP853'' (explicit Runge-Kutta
method of order 8)~\cite{HairerNW1993},
``Radau'' (implicit Runge-Kutta method of Radau IIA family
of order 5)~\cite{HairerW1996}, and ``BDF'' (implicit multi-step
variable-order (orders 1 to 5)
based on backward differentiation formulas enhanced with NDF modification)~\cite{ByrneH1975,ShampineR1997}.
The data for NN-Exp-S0 and NN-Exp-S1 correspond to those
in Figure~\ref{fg_5}.
The scipy results are obtained by employing
the routine ``scipy.integrate.solve\_ivp()'' from the scipy library
with different methods.
Different data points for the scipy methods correspond to
different relative tolerance values, while the absolute tolerance is
fixed at $10^{-16}$. Since the scipy methods are adaptive in time step and
order, we have used the dense output option to attain their solutions on
points corresponding to a time step size $\Delta t=0.02$ for $t\in[0,1]$
to compute the errors.
Among the scipy methods, DOP853 and RK45 show the best performance,
followed by Radau and BDF, and then by RK23.
The NN-Exp-S0 algorithm appears to perform slightly better than NN-Exp-S1
for this problem. 
Both NN-Exp-S0 and NN-Exp-S1 significantly outperform the scipy methods.
Here by ``outperform'' we refer to the ability to achieve superior accuracy
under the same time-marching cost or achieve the same accuracy under a lower
time-marching cost.

\begin{figure}
  \centerline{
    \includegraphics[width=2in]{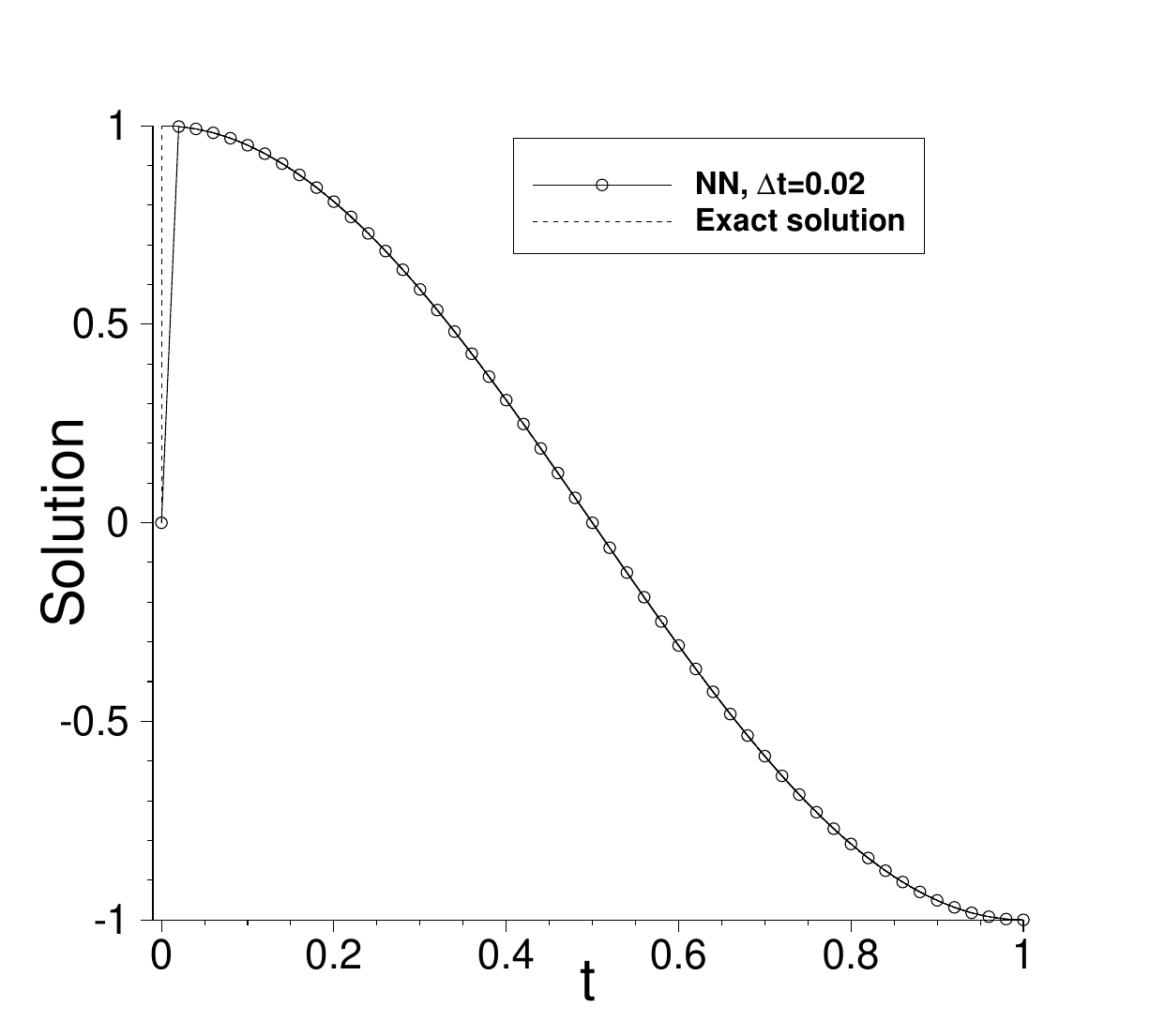}(a)
    \includegraphics[width=2in]{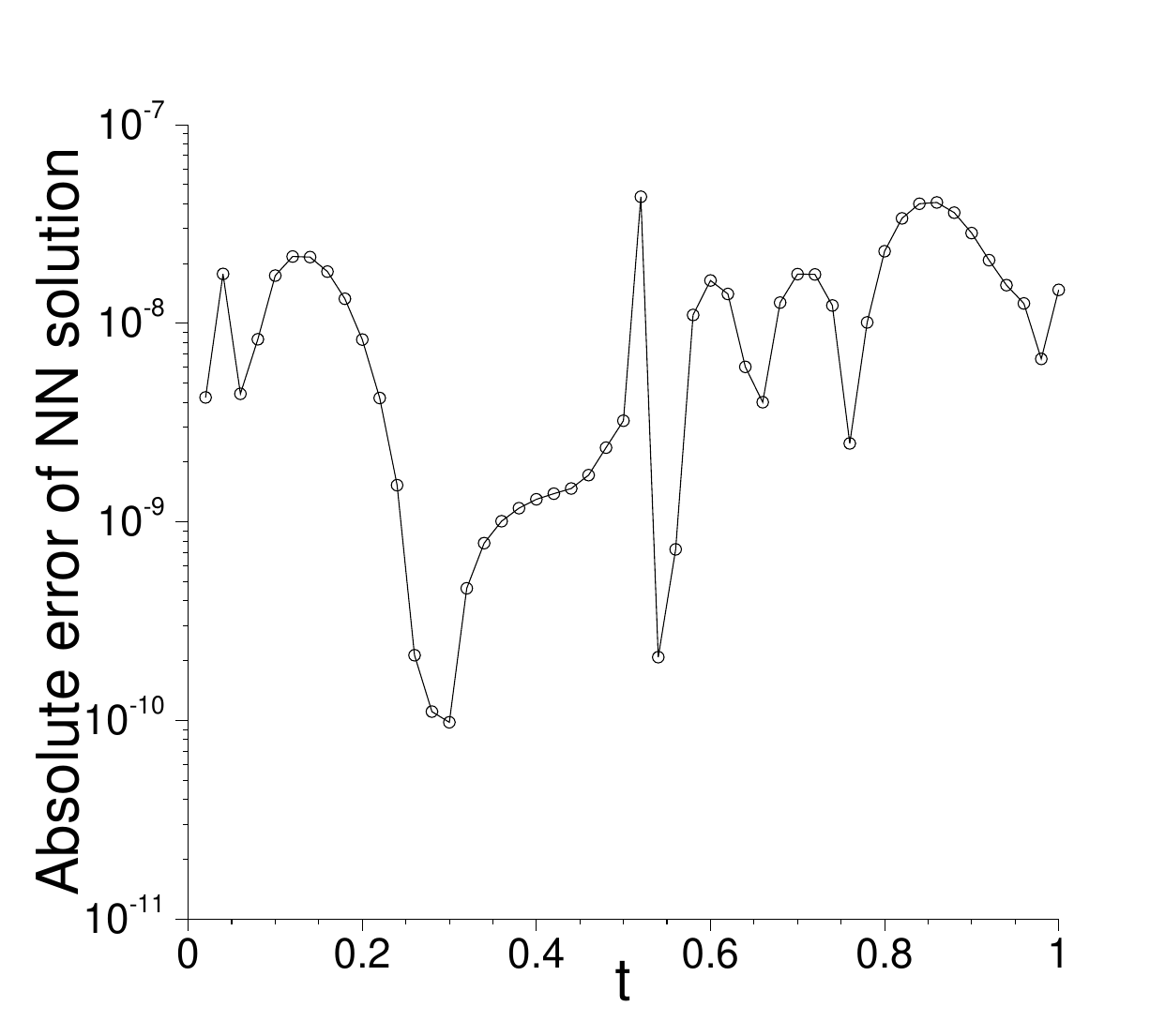}(b)
  }
  \caption{Linear model ($\lambda=10^6$, stiff):
    (a) Comparison between the NN-Exp-S0 solution and the exact solution.
    (b) Absolute-error history of the NN-Exp-S0 solution.
    $h_{\max}=0.025$,
    2 uniform sub-domains along $t_0$
    with enlargement factor $r=0.05$; 
    $M=1000$, $Q=1000$, $R_m=0.4$, and $\delta_m=0.02$.
    See Table~\ref{tab_a2} for the other parameter values.
  }
  \label{fg_7}
\end{figure}

We next consider a stiff case, with $\lambda=10^6$ in problem~\eqref{eq_53}.
Among the explicit NN algorithms in~\eqref{eq_12}-\eqref{eq_12c}, we observe that
NN-Exp-S0 works well for the
stiff problem, while NN-Exp-S1
works not as well as NN-Exp-S0 and NN-Exp-S2 fails to work (NN training fails
to converge, with large loss values).
The implicit NN algorithms from~\eqref{eq_22}
also work well for this stiff problem, but is not as competitive
as the explicit NN-Exp-S0 algorithm in performance.

Figure~\ref{fg_7} illustrates the characteristics of
the NN-Exp-S0 solution for this stiff case. 
Figure~\ref{fg_7}(a) compares the NN-Exp-S0 solution with
the exact solution, and Figure~\ref{fg_7}(b) shows the absolute error
of the NN-Exp-S0 solution.
The training domain, the NN structure, and the other parameter values are provided
in the figure caption or in Table~\ref{tab_a2}. In particular, two uniform sub-domains are
employed along the $t_0$ direction, with an enlargement factor $r=0.05$
for network training (see Remark~\ref{rem_26}), and
we have employed a $\xi$-domain map factor $\delta_m=0.02$
when normalizing the input data on each sub-domain
(see discussions at the beginning of Section~\ref{sec_tests}).
A time step size $\Delta t=0.02$ is used in time marching for
$t\in[0,1]$. The NN-Exp-S0 method is highly accurate,
with a maximum error on the order of $10^{-8}$ over the domain.

\begin{figure}
  \centerline{
    \includegraphics[width=2in]{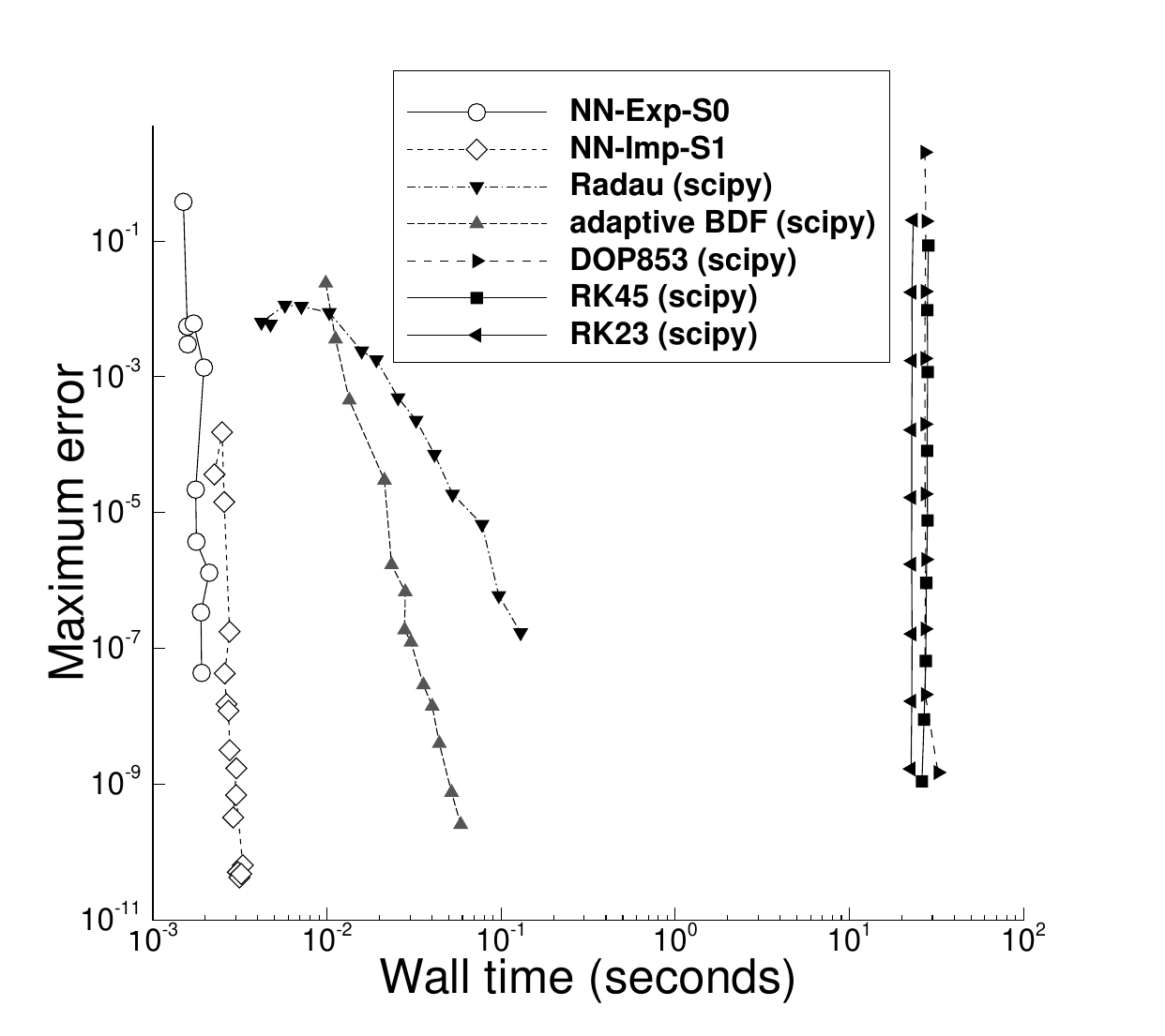}(a)
    \includegraphics[width=2in]{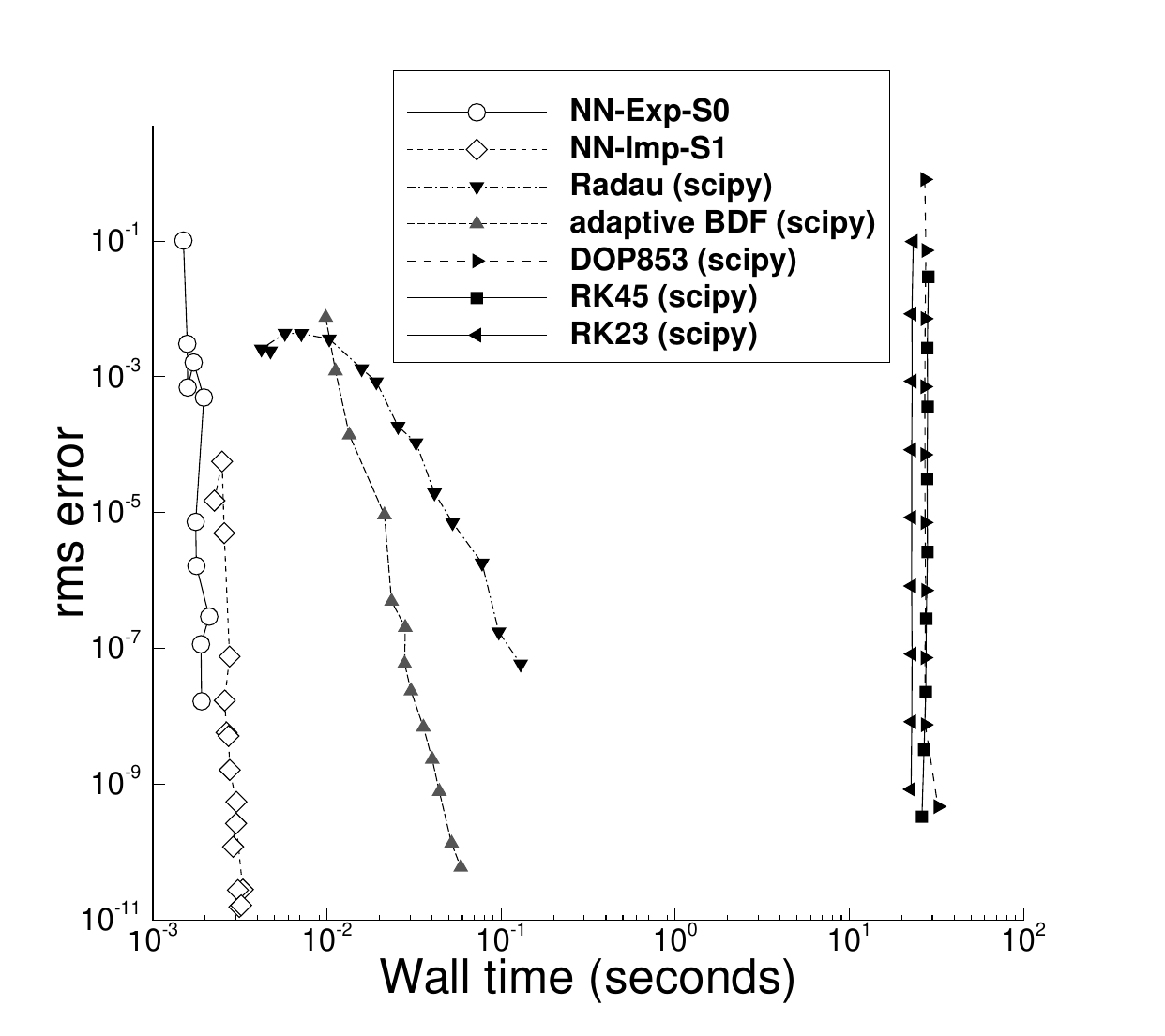}(b)
  }
  \caption{Linear model ($\lambda=10^6$, stiff):
    Comparison of (a) the maximum and (b) the rms time-marching errors
    versus the time-marching cost (wall time) 
    between the NN algorithms (NN-Exp-S0, NN-Imp-S1)
    and the  scipy methods.
    NN-Exp-S0: $h_{\max}=0.025$, $Q=1000$, $R_m=0.4$, $\delta_m=0.02$;
    NN-Imp-S1: $h_{\max}=0.024$, $Q=900$, $R_m=0.35$, $\delta_m=0.01$;
    The other parameter values are provided in Table~\ref{tab_a2}.
    Scipy methods: absolute tolerance $10^{-16}$, 
    data points corresponding to different relative tolerance values,
    dense output on points corresponding to $\Delta t=0.02$ for $t\in[0,1]$.
  }
  \label{fg_8}
\end{figure}

In Figure~\ref{fg_8} we compare the maximum and rms solution errors
versus the time-marching cost (wall time) obtained
by the NN algorithms (NN-Exp-S0 and NN-Imp-S1)
and the scipy methods for this stiff problem. 
The parameter values employed in
these tests are provided in the figure caption or in Table~\ref{tab_a2}.
The explicit scipy methods (DOP853, RK45, and RK23)
are not competitive for this stiff problem,
as expected, inducing a large time marching cost.
In contrast, the implicit scipy methods (Radau and BDF) perform considerably better.
The NN-Exp-S0 algorithm is more competitive than NN-Imp-S1, but
appears slightly less accurate than the latter for this problem.
Both NN-Exp-S0 and NN-Imp-S1 significantly outperform
the implicit and explicit scipy methods for this stiff problem.



\subsection{Pendulum Problem}

We test the learned NN algorithms using
the nonlinear pendulum problem in this section.
We study an autonomous system (free pendulum) first,
followed by a non-autonomous system (forced pendulum).

\subsubsection{Free Pendulum Oscillation}
\label{sec_321}

\begin{table}[tb]
  \centering
  \begin{tabular}{l|l}
    \hline
    domain: $(y_{01},y_{02},\xi)\in [-2,2]\times[-4,4]\times[0,0.25]$
    & sub-domains: 3, along $y_{01}$, uniform\\
    NN ($\varphi$-subnet) architecture: $[3, M, 2]$ ($M$ varied) & $r$: $0.1$ \\
    activation function: Gaussian & $\delta_m$: $1.0$  \\
    $Q$: $900$, $1000$ or $1500$ (random) & $R_m$: to be specified  \\
    $\Delta t$: $0.2$ (for time marching) & time: $t\in[0,200]$ (for time marching) \\
    \hline
  \end{tabular}
  \caption{NN simulation parameters for the free pendulum problem (Section~\ref{sec_321}).
  }
  \label{tab_a3}
\end{table}

\begin{figure}
  \centerline{
    \includegraphics[width=2in]{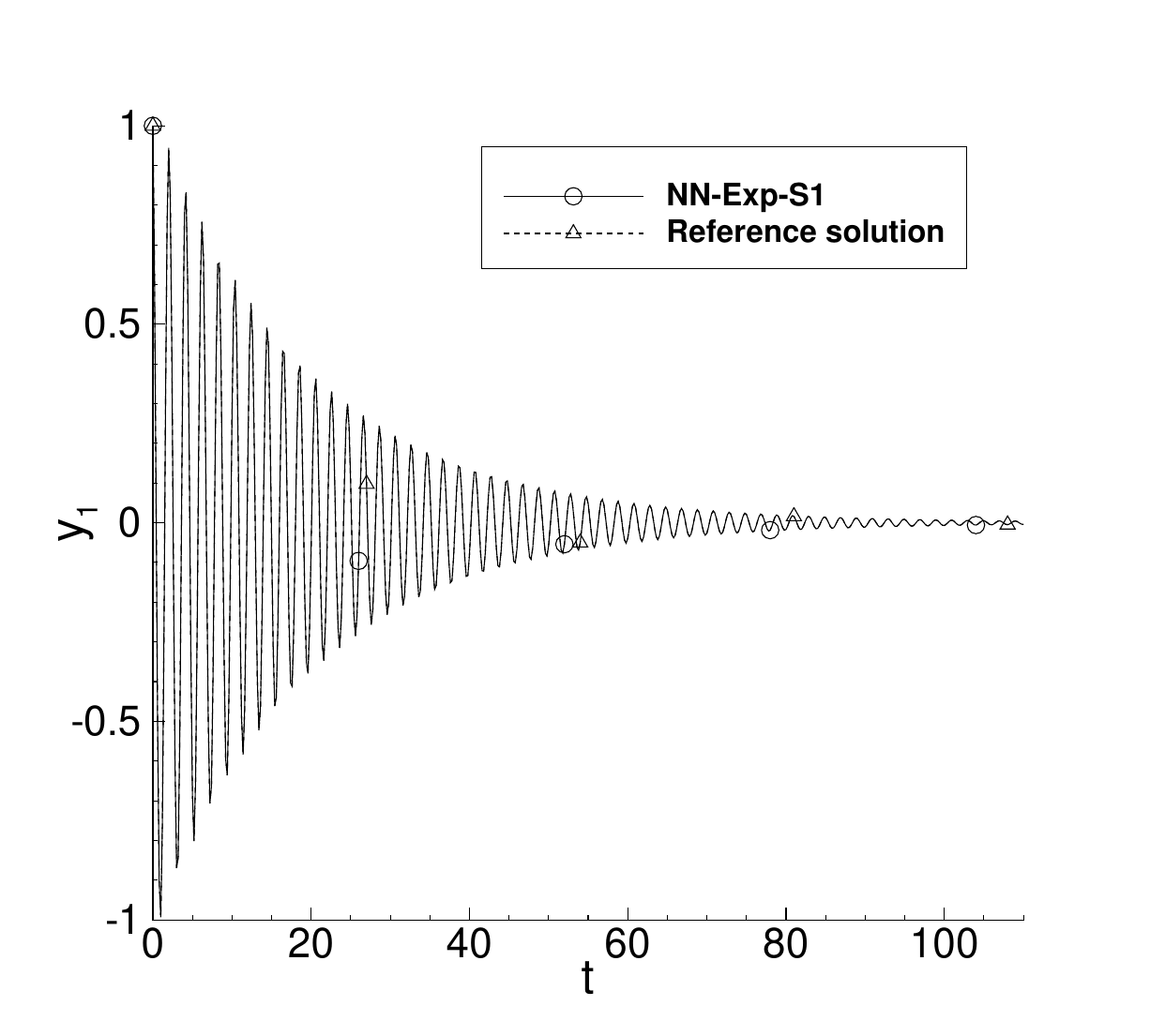}(a)
    \includegraphics[width=2in]{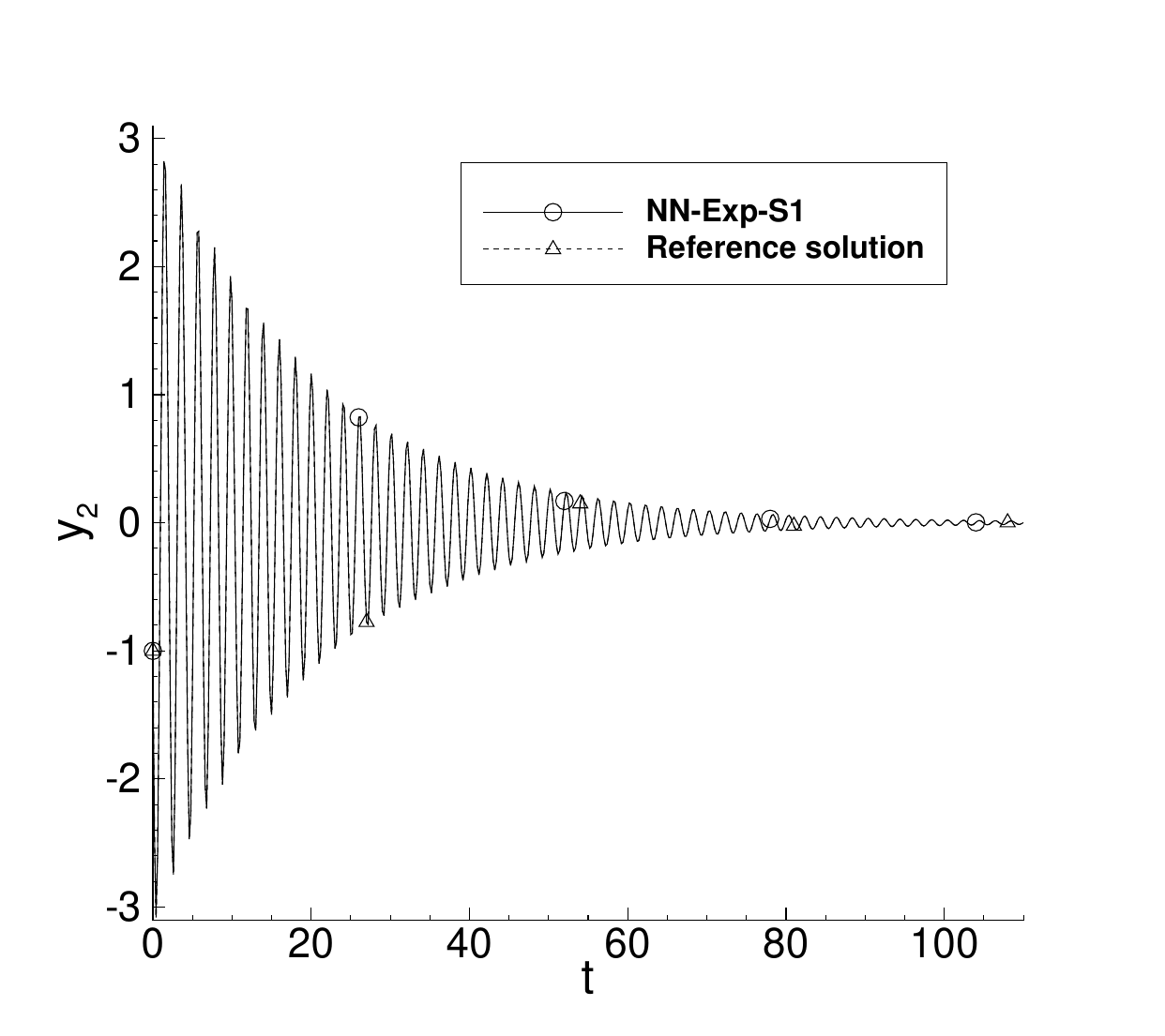}(b)
    \includegraphics[width=2in]{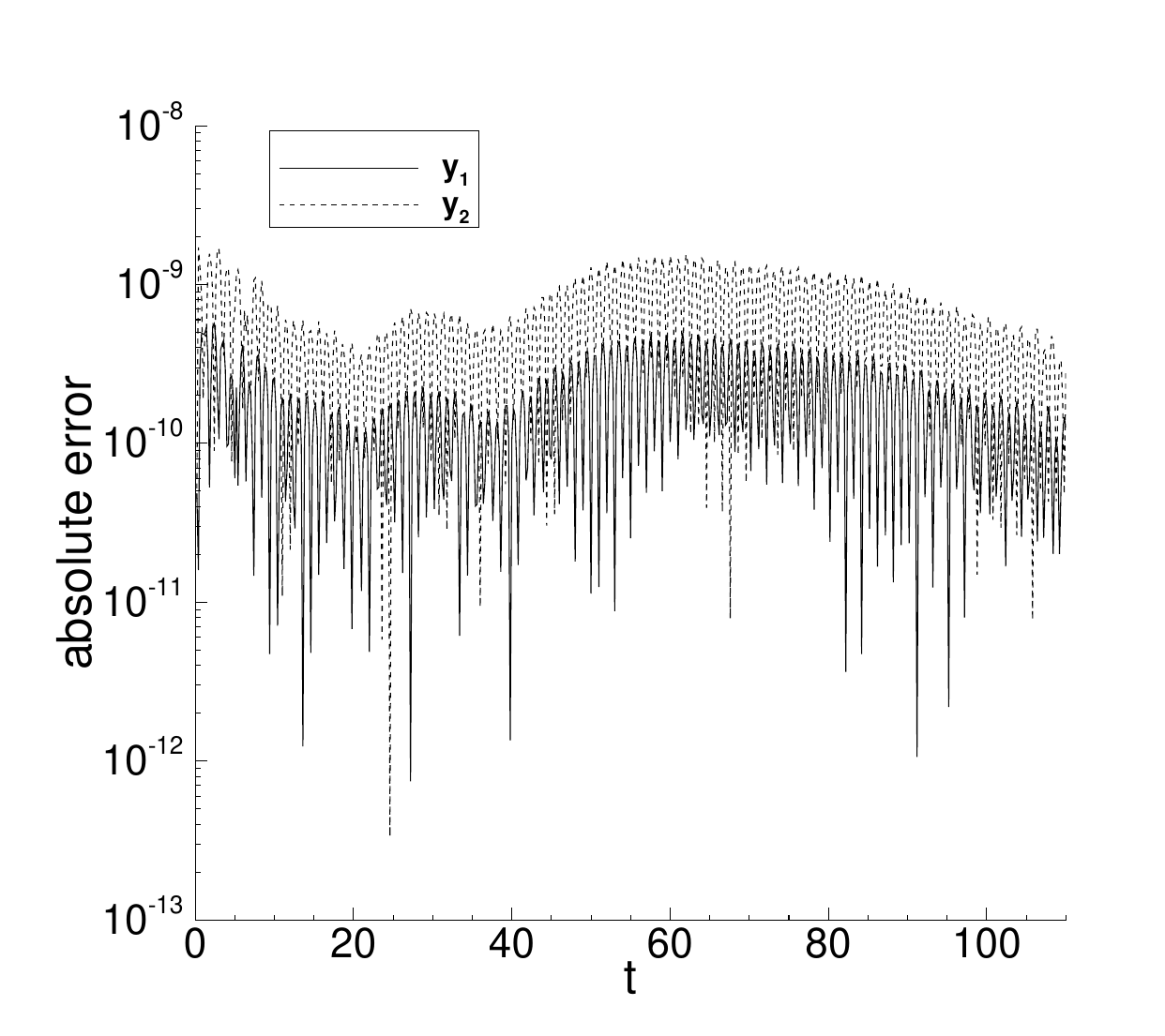}(c)
  }
  \caption{Free pendulum: Comparison 
    of (a) $y_1(t)$ and (b) $y_2(t)$ between the NN-Exp-S1 solution and the reference
    solution. (c) Absolute-error histories of the NN-Exp-S1 solution for
    $y_1(t)$ and $y_2(t)$.
    NN-Exp-S1: $M=800$, $Q=1000$, $R_m=0.6$; See Table~\ref{tab_a3} for the other
    parameter values.
  }
  \label{fg_9}
\end{figure}

Consider the initial value problem with the free pendulum equation,
\begin{subequations}\label{eq_56}
  \begin{align}
    &
    \frac{dy_1}{dt} = y_2, \quad
    \frac{dy_2}{dt} = -\alpha y_2 - \beta \sin y_1, \\
    &
    y_1(t_0) = y_{01}, \quad y_2(t_0) = y_{02},
  \end{align}
\end{subequations}
where $\alpha$ and $\beta$ are prescribed positive constants,
$y(t) = (y_1(t), y_2(t))\in\mbb R^2$ are the unknowns to be
computed, $t_0$ is the initial time, and
$y_0=(y_{10},y_{20})\in\mbb R^2$ are the initial conditions.
Physically $y_1(t)$ represents the pendulum angle, and $y_2(t)$
is the angular velocity.
In the following tests we employ $(\alpha,\beta)=(0.1,9.8)$, $t_0=0$, and
the initial condition $y_0=(y_{01},y_{02})=(1.0,-1.0)$, unless otherwise specified.

We employ an ELM network to learn the algorithmic function $\psi(y_{01},y_{02},\xi)$.
Table~\ref{tab_a3} summarizes the simulation parameters related to the NN algorithm.
In particular, we employ three uniform sub-domains along the $y_{01}$ direction,
with an enlargement factor $r=0.1$. The parameters $R_m$, $\Delta t$ and $\delta_m$
have the same meanings as in Section~\ref{sec_lin_model}, and $Q$ and $M$ refer to
the number of random collocation points and the hidden-layer width
of the $\varphi$-subnet on each sub-domain.
After the NN is trained, we employ the learned algorithm to solve
the problem~\eqref{eq_56} for $t\in[0,200]$ with a step size
$\Delta t=0.2$ and the initial condition $(y_{01},y_{02})=(1,-1)$.
The errors of the NN solution against a reference solution and the NN time-marching
time are recorded for analysis.
The reference solution is obtained by the using the scipy DOP853 method
with a sufficiently small tolerance (absolute tolerance $10^{-16}$, relative tolerance $10^{-13}$).

Figure~\ref{fg_9} provides 
an overview of the solution characteristics obtained by the NN algorithm.
It compares the NN-Exp-S1 solution for $y_1(t)$ and $y_2(t)$
with their reference solutions, and shows the absolute errors of the
NN-Exp-S1 solutions. The parameter values are given in
the figure caption or in Table~\ref{tab_a3}.
The NN solution is highly accurate, with a maximum error on the order of $10^{-9}$
for the time range.

\begin{figure}
  \centerline{
    \includegraphics[width=2in]{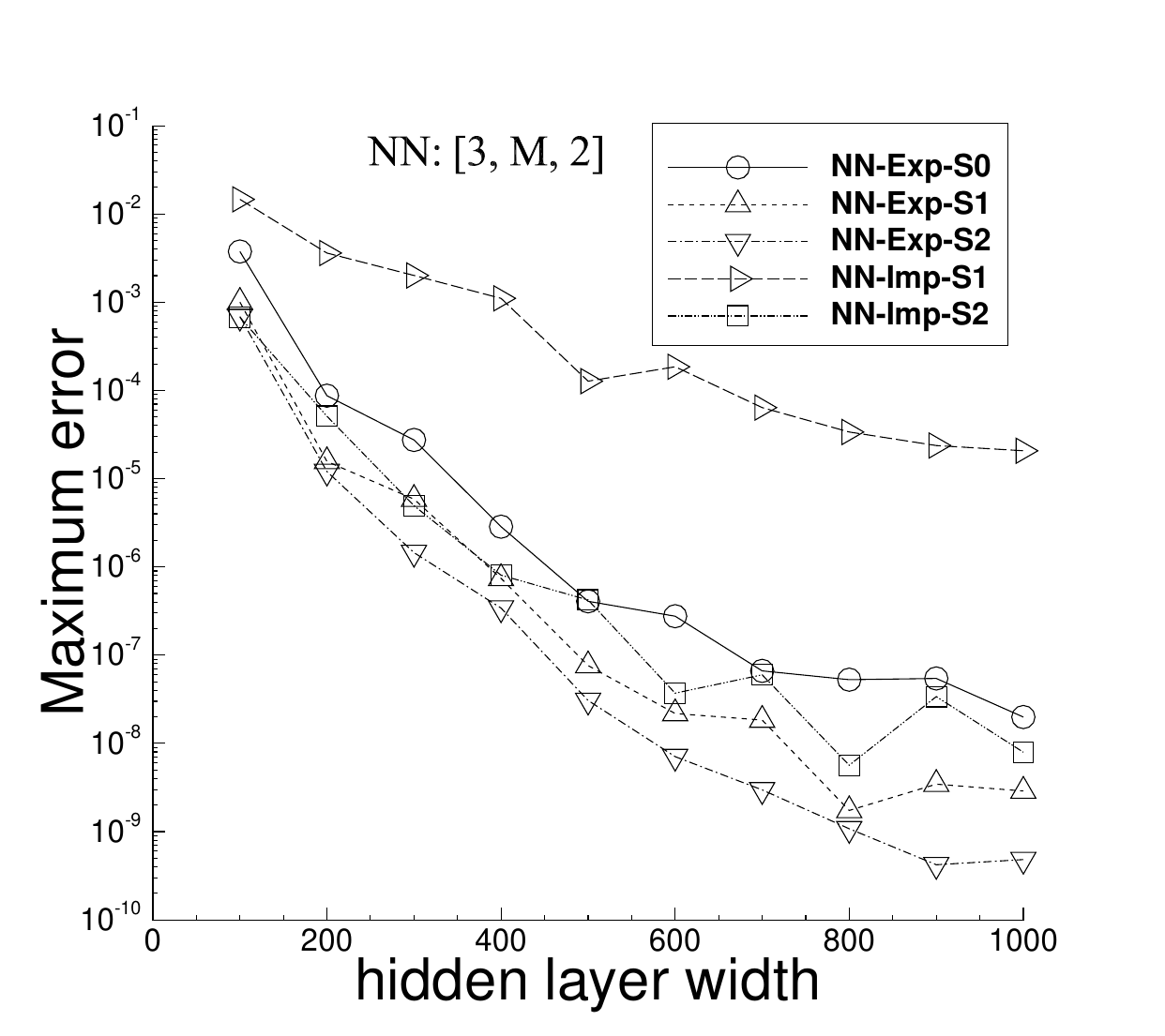}(a)
    \includegraphics[width=2in]{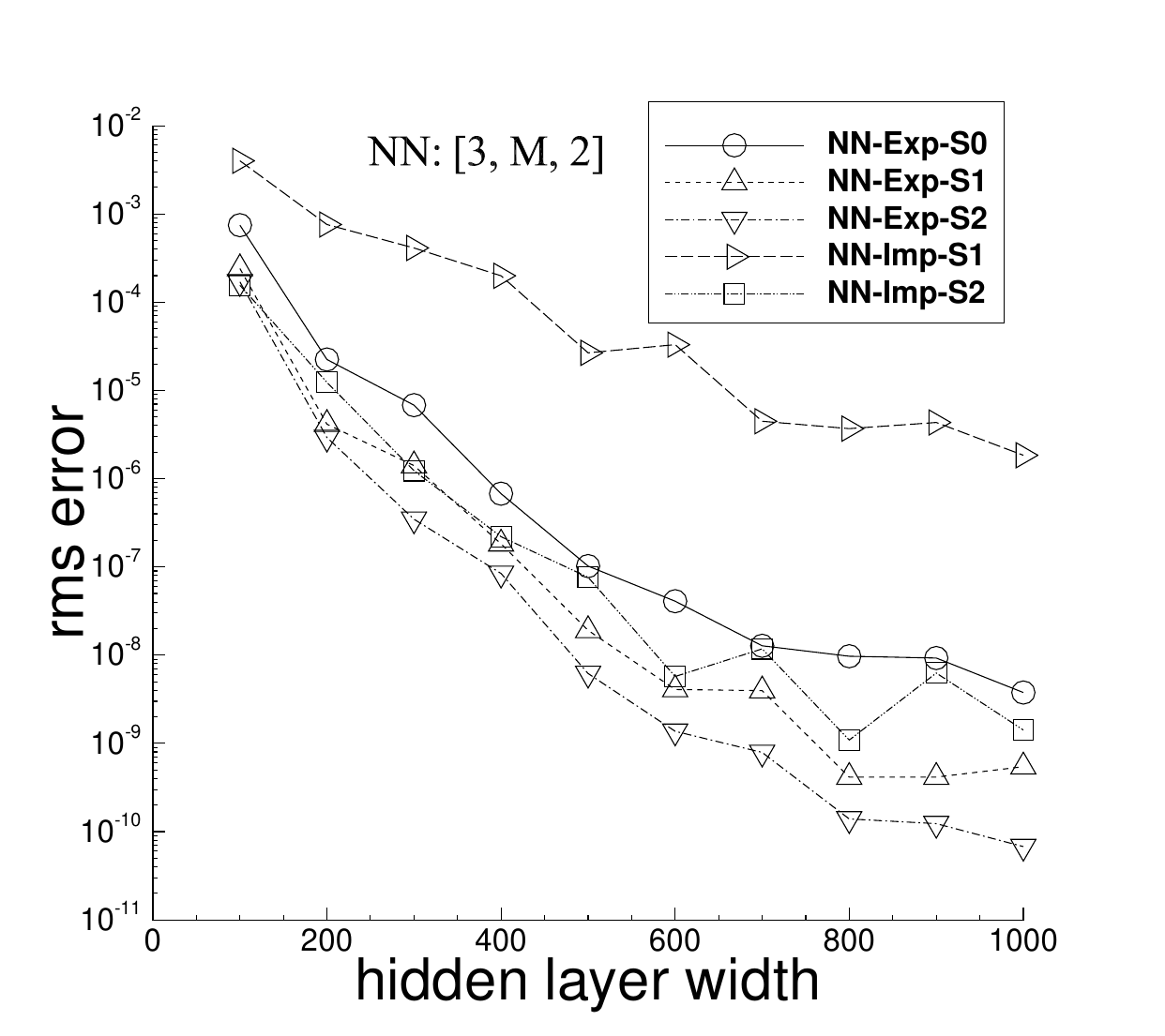}(b)
    \includegraphics[width=2in]{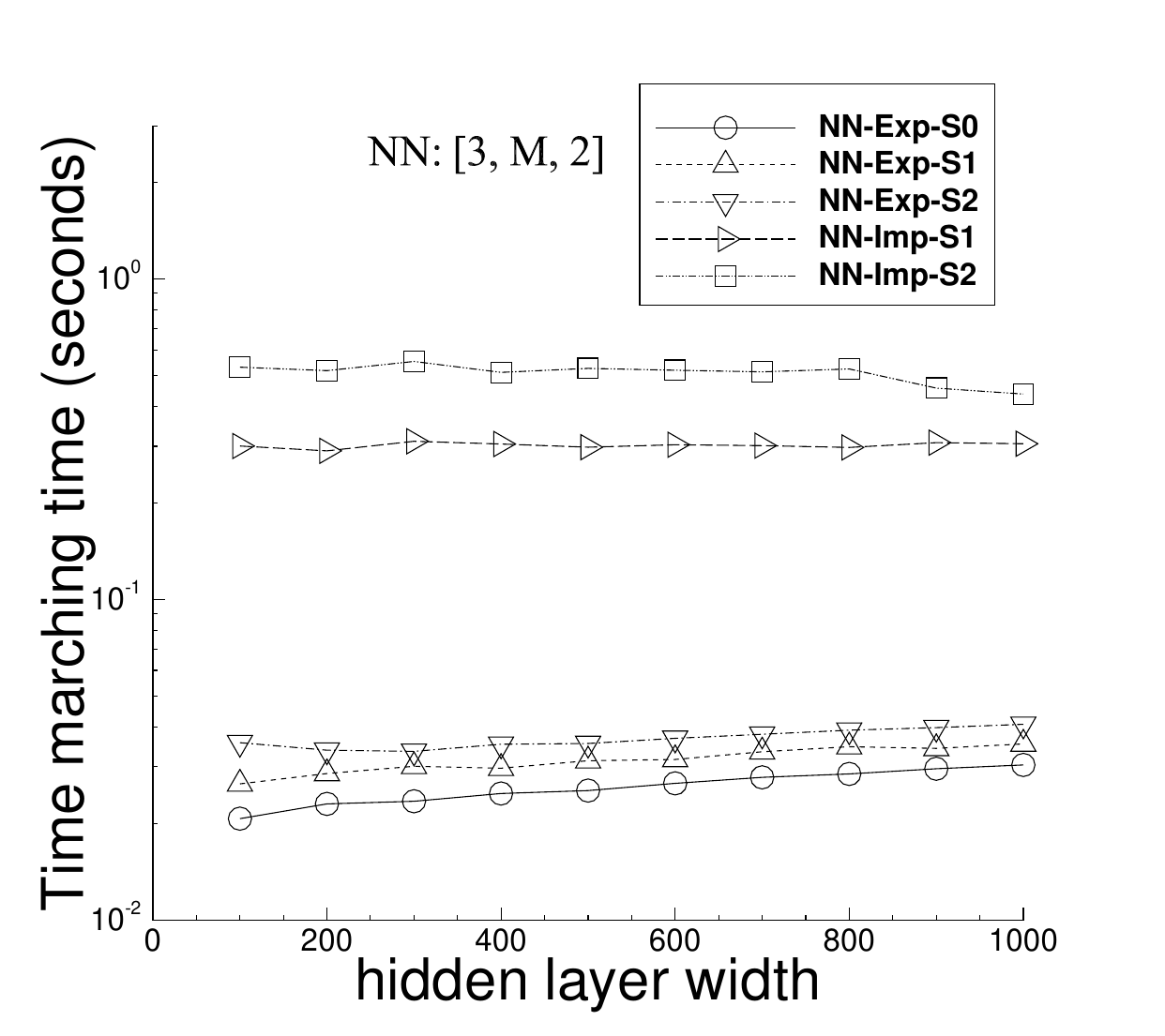}(c)
  }
  \caption{Free pendulum: Comparison of (a) the maximum and (b) the rms time-marching errors,
    and (c) the time-marching cost (wall time) versus the hidden-layer width ($M$)
    in ELM network for different NN algorithms.
    NN-Exp-S0: $R_m=0.6$, $Q=1000$.
    NN-Exp-S1: $R_m=0.6$, $Q=1000$.
    NN-Exp-S2: $R_m=0.65$, $Q=900$.
    NN-imp-S1: $R_m=0.85$, $Q=1500$.
    NN-Imp-S2: $R_m=0.85$, $Q=1000$.
    Other parameter values are given in Table~\ref{tab_a3}.
  }
  \label{fg_10}
\end{figure}

Figure~\ref{fg_10} is a comparison of the accuracy and the time-marching cost
of different NN algorithms for solving this problem. 
Here we plot the maximum and rms time-marching errors ($e_{\max}$, $e_{\text{rms}}$),
as well as the time-marching time,
versus the hidden-layer width ($M$) of the $\varphi$-subnet for
the NN algorithms with explicit (NN-Exp-S0, NN-Exp-S1, and NN-Exp-S2) and
implicit (NN-Imp-S1, and NN-Imp-S2) formulations.
The parameter values are listed in the figure caption or in Table~\ref{tab_a3}.
We can make several observations.
First, the NN solution errors decrease nearly exponentially
with increasing number of hidden-layer nodes, while the time-marching
cost only grows quasi-linearly.
Second, the explicit NN algorithms tend to be more accurate than the implicit ones.
Third, among the explicit NN algorithms the accuracy generally increases
from NN-Exp-S0 to NN-Exp-S1, and to NN-Exp-S2. Between the implicit
NN algorithms, NN-Imp-S2 is significantly more accurate than NN-Imp-S1.
Fourth, in terms of time-marching cost the explicit NN algorithms
are much faster than the implicit ones, by an order of magnitude
for this problem.
Among the explicit NN algorithms, NN-Exp-S0 is the fastest, followed by NN-Exp-S1 and
NN-Exp-S2. Of the implicit ones, NN-Imp-S1 is notably faster than NN-Imp-S2.
Overall, the explicit NN algorithms are more competitive than the implicit ones.
This is similar to what has been observed for the test problem in Section~\ref{sec_lin_model}.

The observation that implicit NN algorithms are less competitive than the explicit ones
appears to be a common characteristic in the test problems we have considered, including
those in the subsequent subsections. Implicit NN algorithms
are computationally much more expensive, and
generally not as accurate as the explicit ones.
For this reason, the computation results using the implicit NN algorithms
will not be included for the simulations in the following subsections.

\begin{figure}
  \centerline{
    \includegraphics[width=2in]{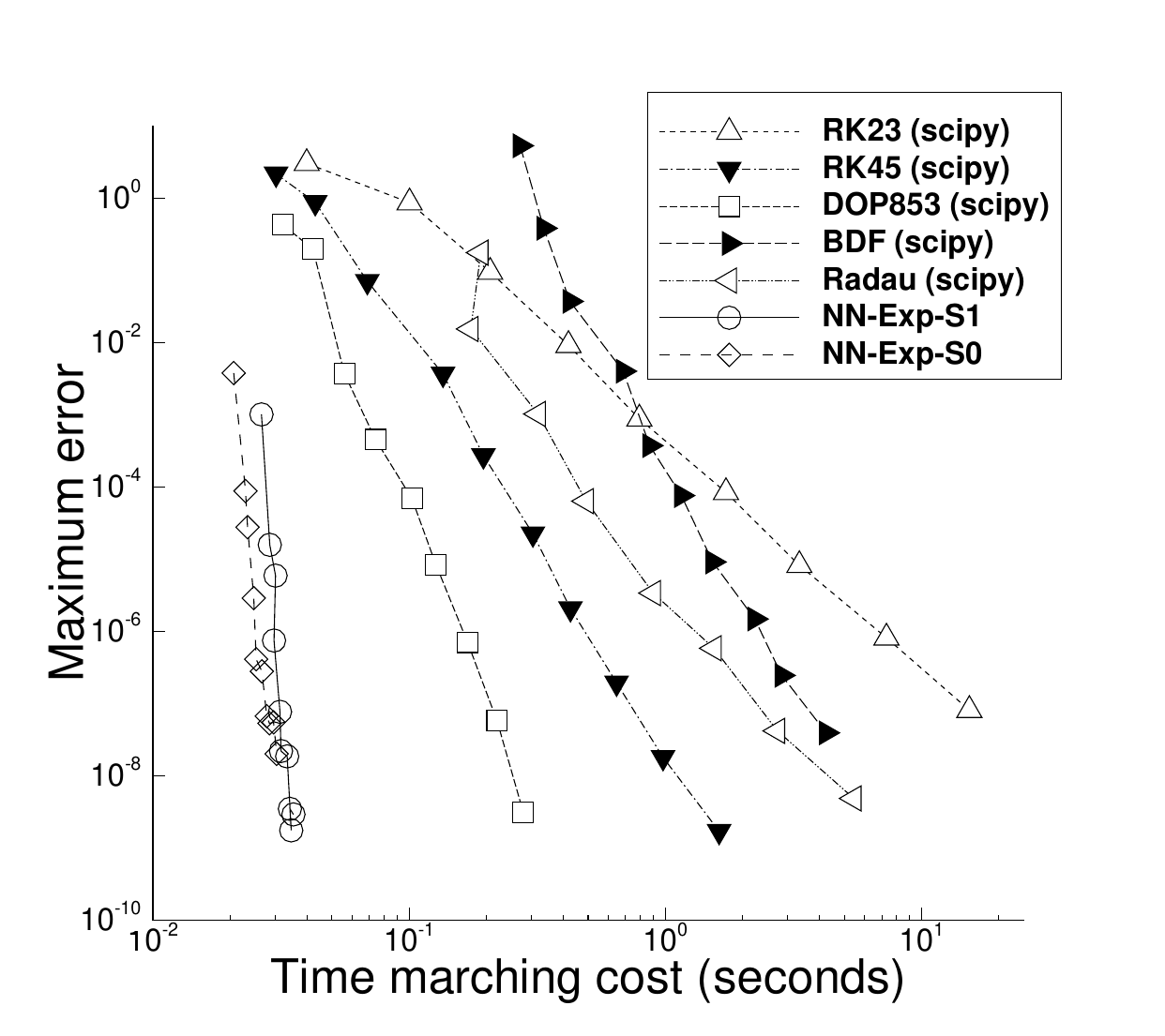}(a)
    \includegraphics[width=2in]{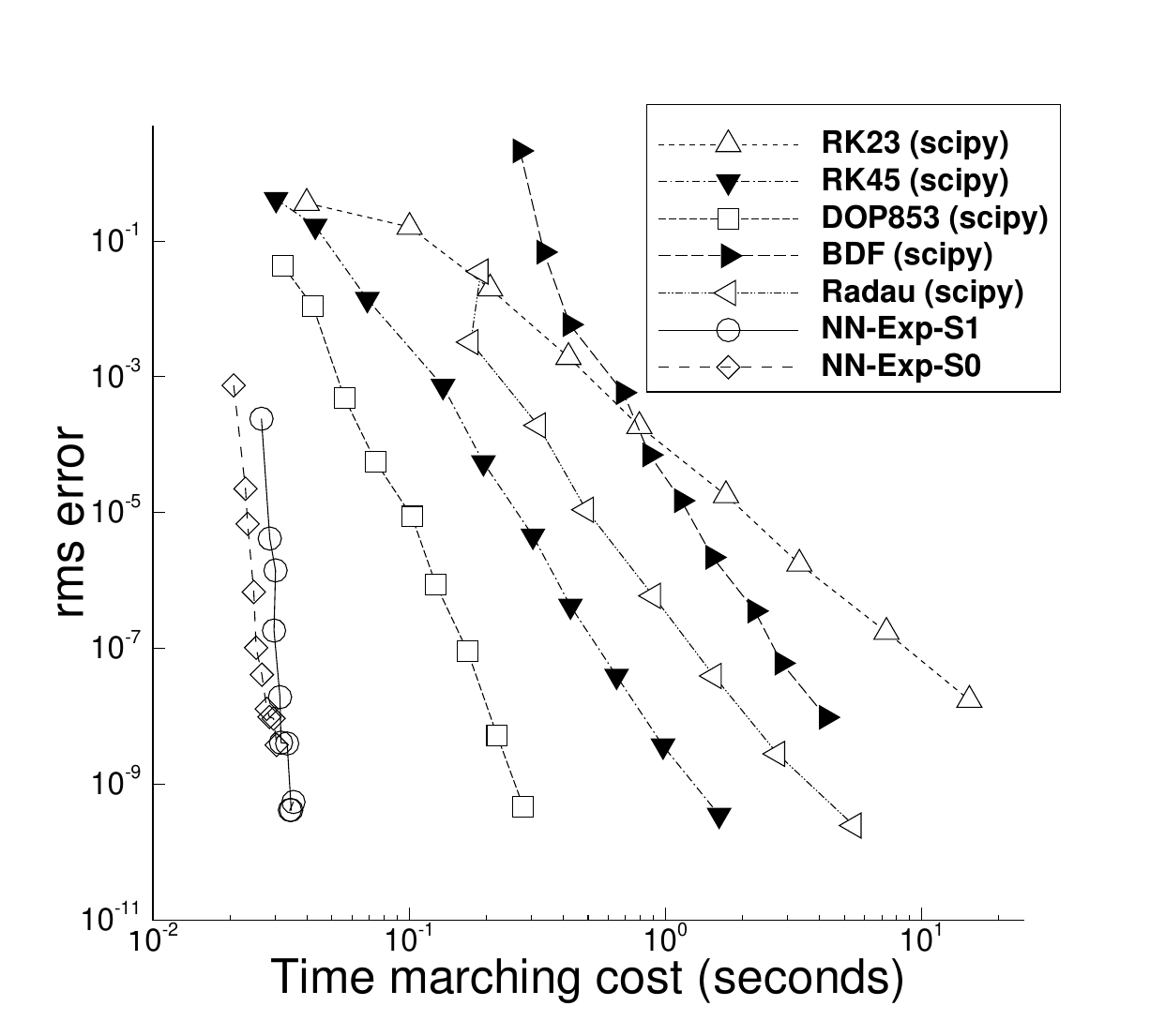}(b)
  }
  \caption{Free pendulum:
    Comparison of (a) the maximum and (b) the rms time-marching errors versus the time marching
    cost (wall time) between the NN algorithms
    (NN-Exp-S0, NN-Exp-S1) and the scipy methods.
    Data for the NN algorithms correspond to
    those of NN-Exp-S0 and NN-Exp-S1 in Figure~\ref{fg_10}.
    Scipy methods: absolute tolerance $10^{-16}$, relative tolerance is varied,
    dense output on time instants
    corresponding to $\Delta t=0.2$ for $t\in[0,200]$.
  }
  \label{fg_11}
\end{figure}

Figure~\ref{fg_11} is a comparison of the computational performance
(accuracy versus cost) between the
NN algorithms (NN-Exp-S0 and NN-Exp-S1) and the scipy methods for the free pendulum
problem. It shows the maximum and rms errors of the NN and scipy solutions
versus their time-marching cost. The parameters and the settings for
the NN algorithms correspond to those in Figure~\ref{fg_10}, in which
the hidden-layer width of the $\varphi$-subnet is varied.
The scipy solutions are obtained by varying the relative tolerance, with
a fixed absolute tolerance $10^{-16}$ for different methods.
Among the scipy methods, DOP853 shows the best performance,
followed by RK45 and Radau. NN-Exp-S0 exhibits a slightly better performance than NN-Exp-S1.
Both NN-Exp-S0 and
NN-Exp-S1 significantly outperform the scipy methods.


\begin{figure}
  \centerline{
    \includegraphics[width=2in]{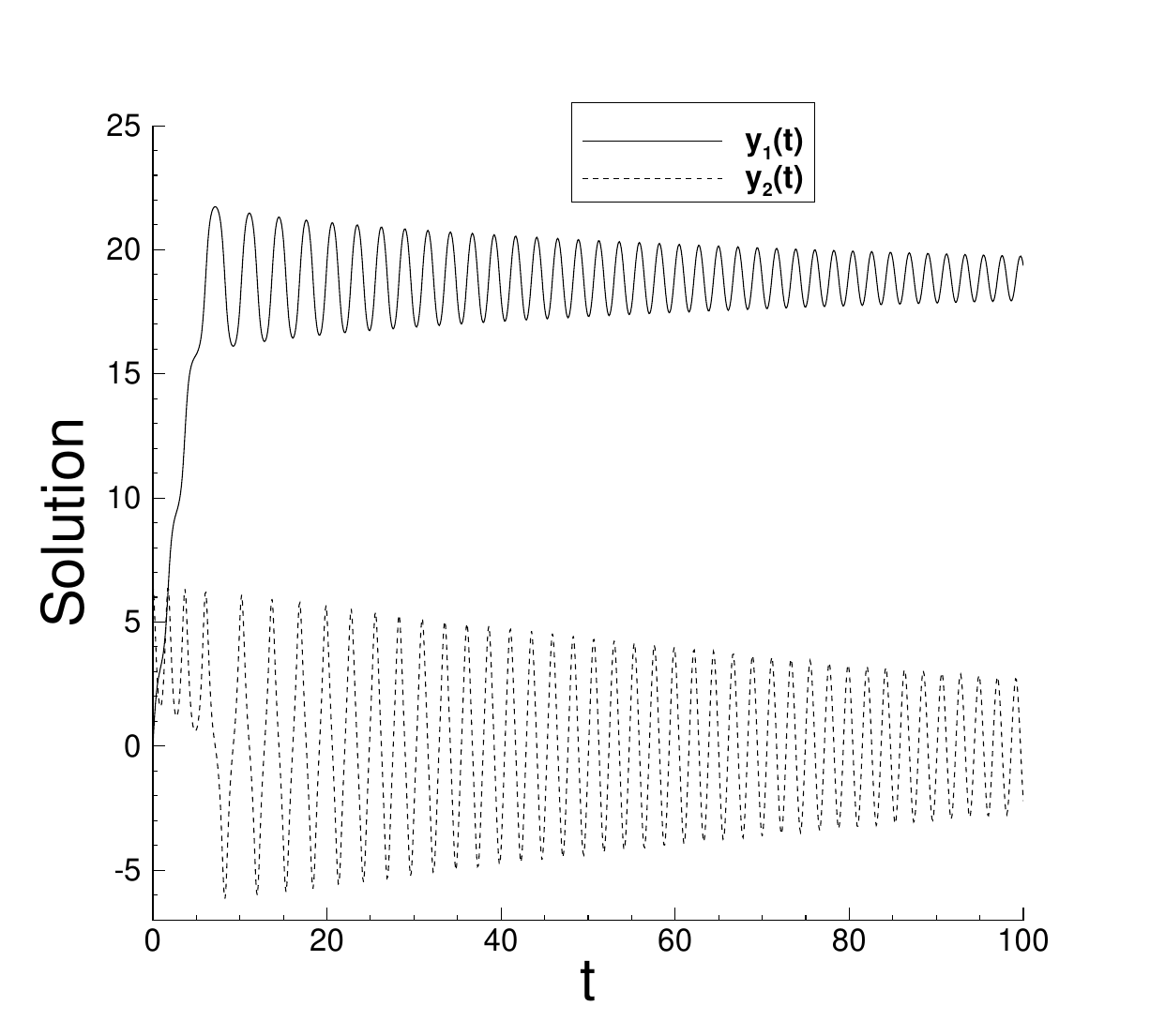}(a)
    \includegraphics[width=2in]{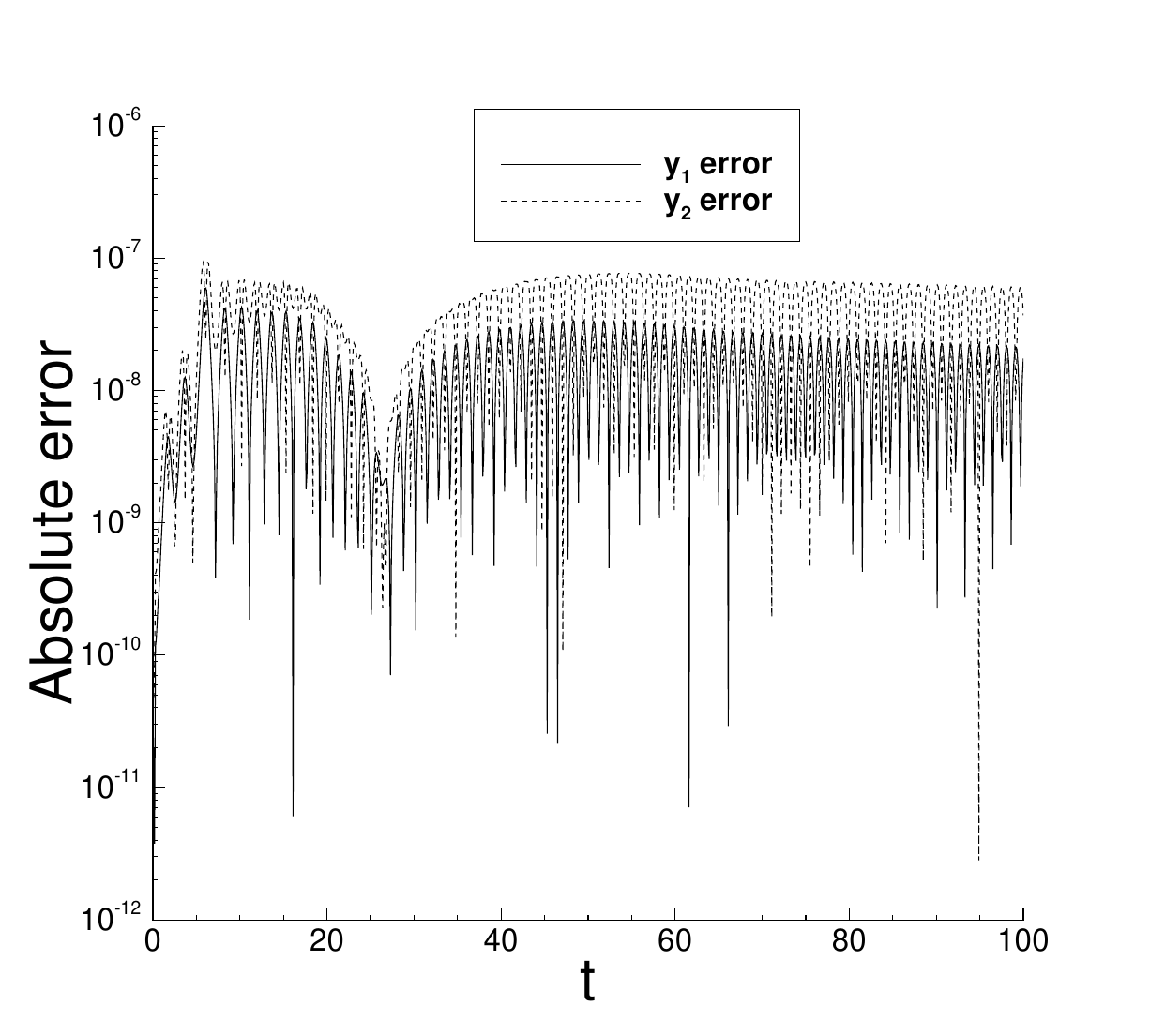}(b)
    }
  \centerline{
    \includegraphics[width=2in]{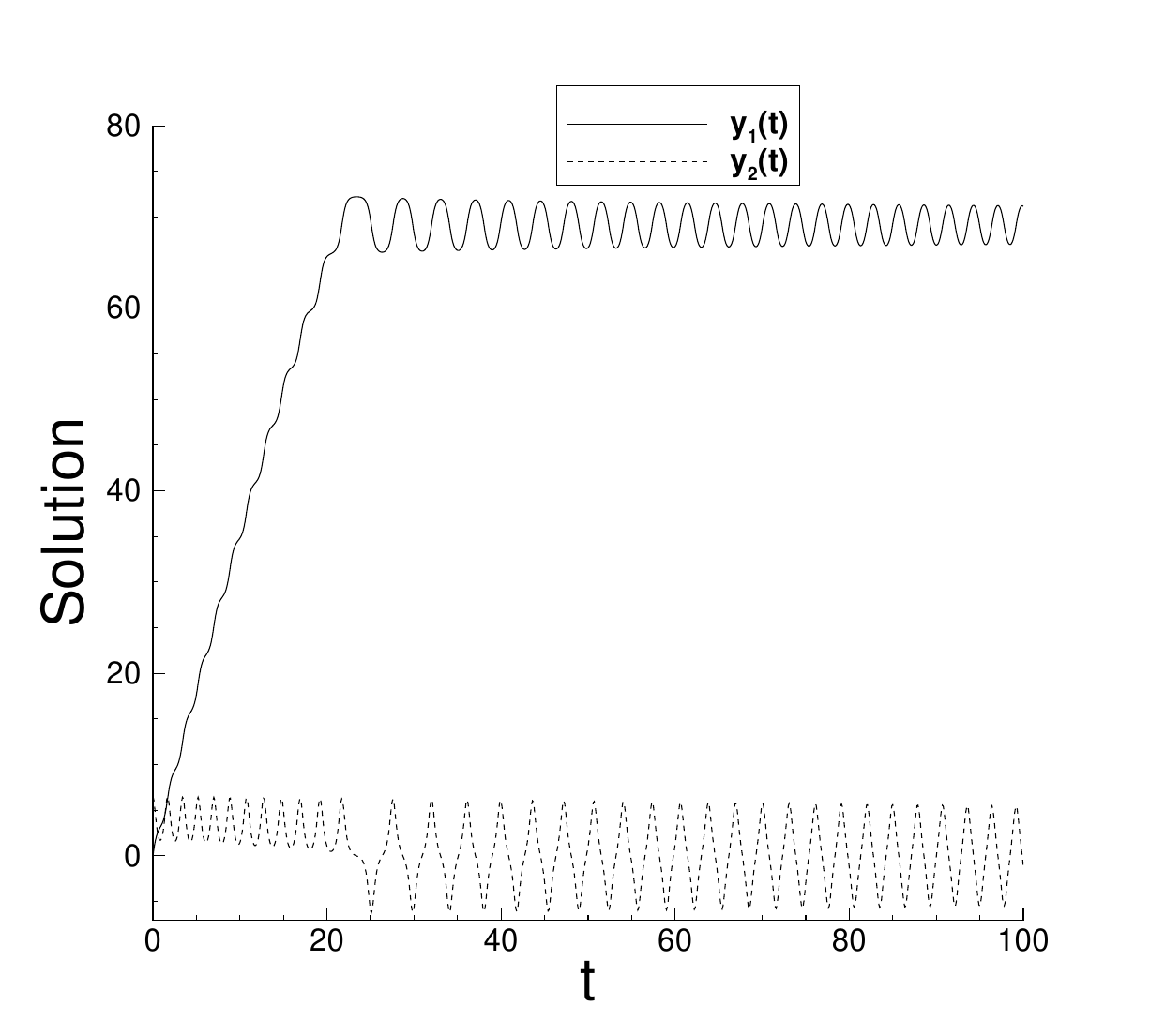}(c)
    \includegraphics[width=2in]{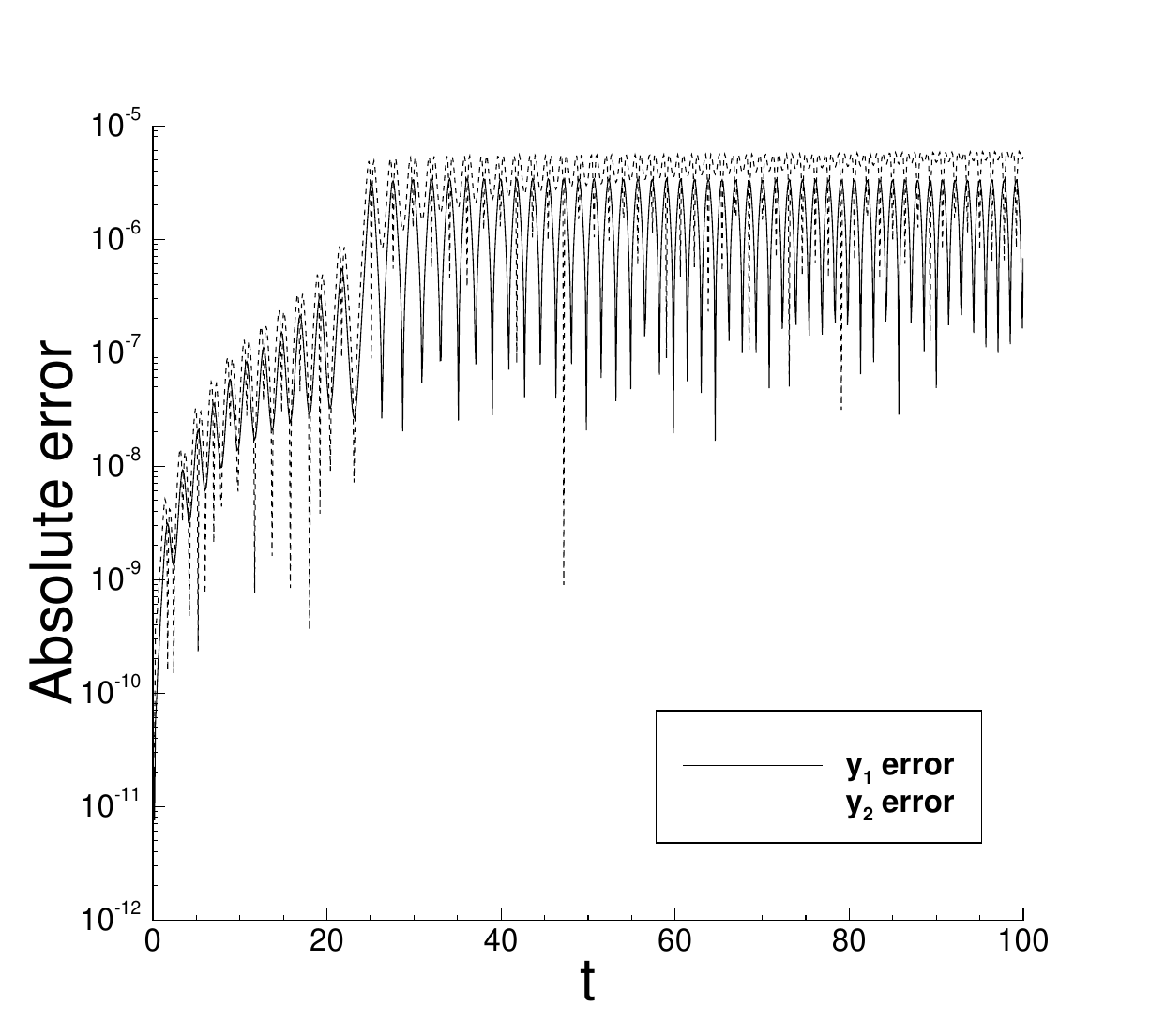}(d)
  }
  \caption{Free pendulum (exploiting periodicity of $f(y,t)$):
    NN-Exp-S1 solutions for $y_1$ and $y_2$ (left column) and their absolute errors (right column)
    for the problem~\eqref{eq_56} with the initial condition $(y_{01},y_{02})=(0,6.5)$ and
    the parameters $(\alpha,\beta)=(0.02,9.8)$ (top row)
    or $(\alpha,\beta)=(0.005,9.8)$ (bottom row).
    Training domain: $(y_{01},y_{02},\xi)\in[-3.2,3.2]\times[-6.7,6.7]\times[0,0.15]$,
    $3$ uniform sub-domains along $y_{01}$ with an enlargement factor $r=0.1$;
    On each sub-domain, NN: $[3,800,2]$,
    Gaussian activation function, $R_m=0.6$, $Q=1500$,
    $\Delta t=0.1$ in time marching for $t\in[0,100]$.
    Error computed relative to a reference solution (by scipy DOP853).
  }
  \label{fg_12}
\end{figure}

Since the RHS of system~\eqref{eq_56} is periodic in $y_1$ with a period $2\pi$,
i.e.~$f(y_1+2\pi,y_2,t)=f(y_1,y_2,t)$, the domain for the NN training
only needs to cover one period along the $y_{01}$ direction, in light of
Theorem~\ref{thm_a2}.
This will be crucial for computing the NN solutions corresponding to large
initial angular velocity ($y_{02}$) values, in which the pendulum angle ($y_1$)
can assume large magnitudes.
Figure~\ref{fg_12} shows the NN-Exp-S1 solutions for two such cases,
which exploit the periodicity of $f(y_1,y_2,t)$ during time-marching as
discussed in Remark~\ref{rem_210}.
These results correspond to an initial condition $(y_{01},y_{02})=(0,6.5)$ at $t_0=0$,
and the problem parameters $(\alpha,\beta)=(0.02,9.8)$ and $(\alpha,\beta)=(0.005,9.8)$,
respectively. Because the initial angular velocity $y_{02}=6.5$ is sufficiently
large, the pendulum spins around the axis for several periods before
settling down to oscillate around the equilibrium position (at an elevated angle $y_1$).
The simulation parameters for this group of tests are specified in the figure caption.
In particular, the NN is trained on the domain
$(y_{01},y_{02},\xi)\in[-3.2,3.2]\times[-6.7,6.7]\times[0,0.15]$, which
is slightly larger than one period $[-\pi,\pi]$ along $y_{01}$.
For time-marching using the learned $\psi(y_{01},y_{02},\xi)$, we exploit
the periodicity of $f(y_1,y_2,t)$ as follows (see also Remark~\ref{rem_210}).
Given $(y_{1k},y_{2k},t_k)$ and the step size $\Delta t$,
we compute $y_{1k}^*=\text{mod}(y_{1k}+\pi,2\pi)-\pi$,
$q=\lfloor\frac{y_{1k}+\pi}{2\pi} \rfloor$, and
$(y_{1,k+1},y_{2,k+1})=\psi(y_{1k}^*,y_{2k},\Delta t)+(2\pi q, 0)$.
The results in Figure~\ref{fg_12} confirm the effectiveness of
the method described in Remark~\ref{rem_210}
for taking advantage of the RHS $f(y,t)$ that is periodic with respect to $y$,
and show that the NN algorithm is accurate for this type of problems.
Without exploiting this periodicity, one would need to use a domain having a
large dimension along $y_{01}$ to train the NN,
which would pose significant challenges to network training for such problems.

\subsubsection{Forced Pendulum Oscillation}
\label{sec_322}

\begin{table}[tb]
  \centering
  \begin{tabular}{l|l}
    \hline
    domain: $(y_{01},y_{02},t_0,\xi)\in [-2,2]\times[-4,4]\times[0,2.01]\times[0,0.11]$
    & sub-domains: 3, along $y_{01}$, uniform\\
    NN ($\varphi$-subnet) architecture: $[4, M, 2]$ ($M$ varied) & $r$: $0.1$ \\
    activation function: Gaussian & $\delta_m$: $5.0$ or $7.0$  \\
    $Q$: $2000$ or $1500$, random & $R_m$: to be specified \\
    $\Delta t$: $0.1$ (for time marching) & time: $t\in[0,200]$ \\
    \hline
  \end{tabular}
  \caption{NN simulation parameters for the forced
    pendulum problem (Section~\ref{sec_322}).
  }
  \label{tab_2}
\end{table}

We next consider a non-autonomous system with the forced pendulum equation,
\begin{subequations}\label{eq_a56}
  \begin{align}
    &
    \frac{dy_1}{dt} = y_2, \quad
    \frac{dy_2}{dt} = -\alpha y_2 - \beta \sin y_1 + \gamma\cos(\pi t), \label{eq_a56a} \\
    &
    y_1(t_0) = y_{01}, \quad y_2(t_0) = y_{02},
  \end{align}
\end{subequations}
where $\gamma$ is an additional constant and the other variables are
the same as in Section~\ref{sec_321}.
We employ $(\alpha,\beta,\gamma)=(0.1,9.8,0.2)$, 
$y_0=(y_{01},y_{02})=(1.0,-1.0)$, and $t_0=0$ for this problem.

We train the NN using the parameter values from Table~\ref{tab_2}
to learn the algorithmic function $\psi(y_{01},y_{02},t_0,\xi)$.
Since the forcing term is periodic, the RHS of~\eqref{eq_a56a} satisfies
$f(y_1,y_2,t+2)=f(y_1,y_2,t)$.
Therefore, we only need to learn $\psi(y_{01},y_{02},t_0,\xi)$
on a domain that is not smaller than $2.0$ along $t_0$
(see Theorem~\ref{thm_a1}).
We will use $[0,2.01]$ as the training domain along $t_0$.

\begin{figure}
  \centerline{
    \includegraphics[width=2in]{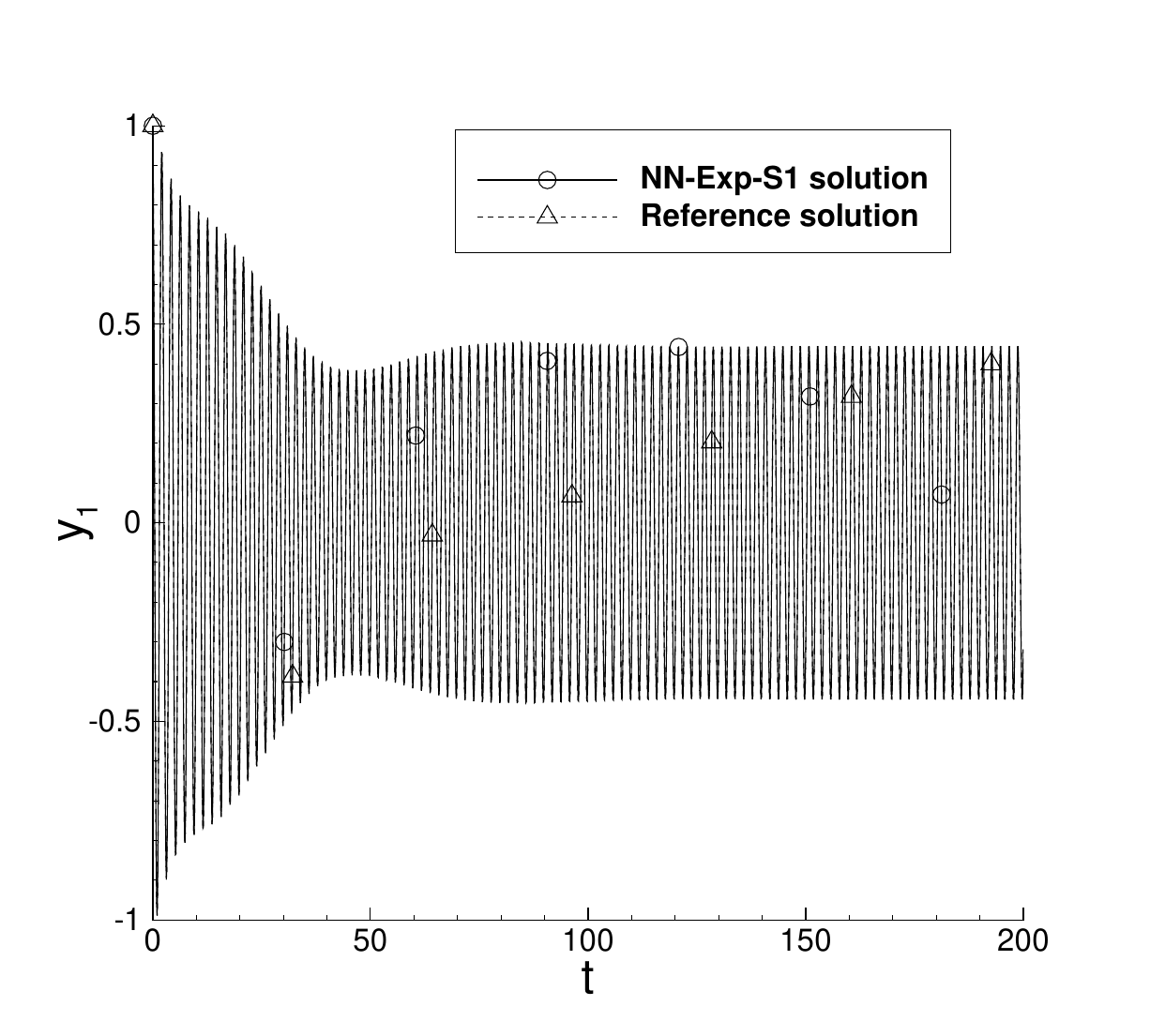}(a)
    \includegraphics[width=2in]{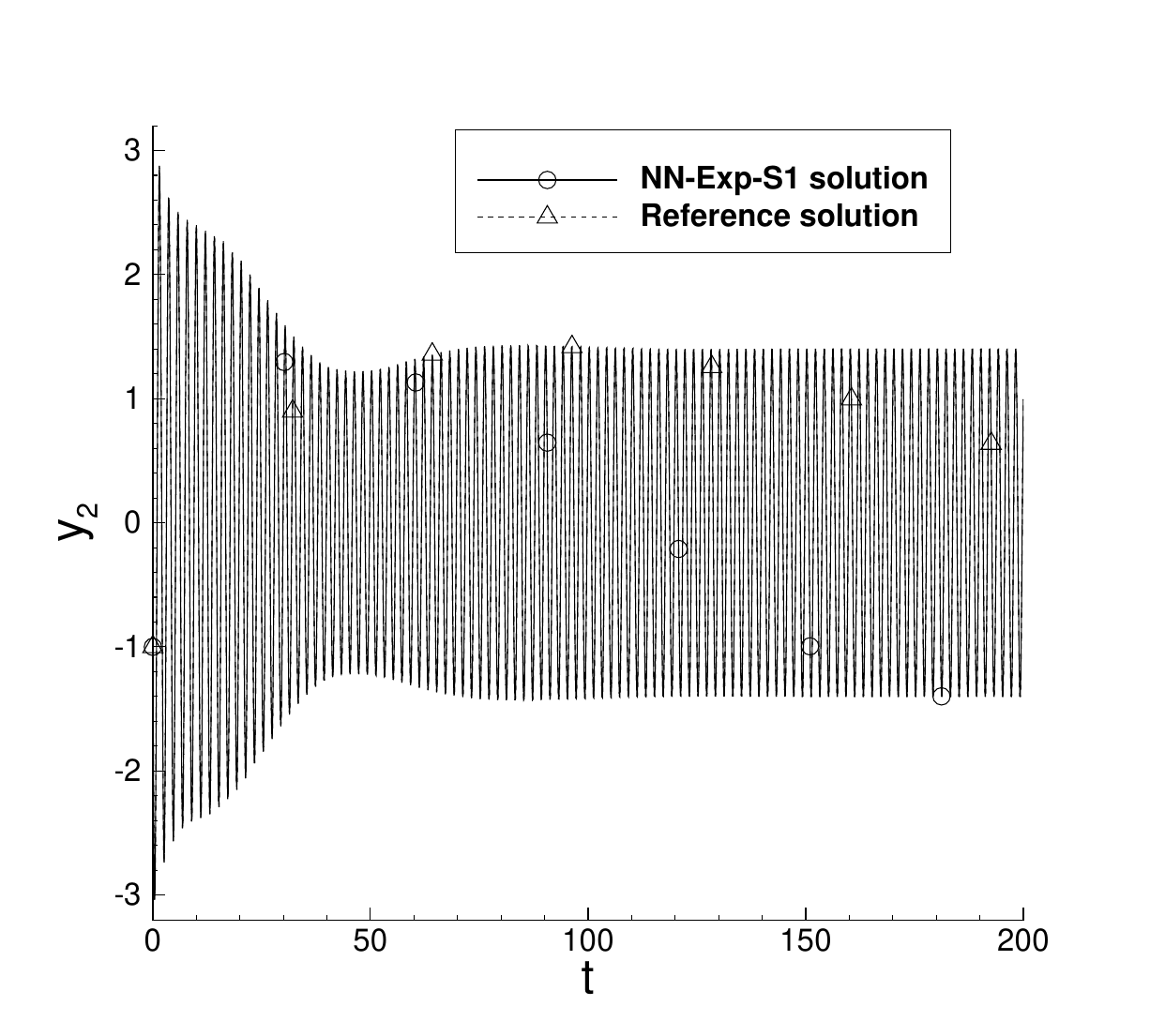}(b)
    \includegraphics[width=2in]{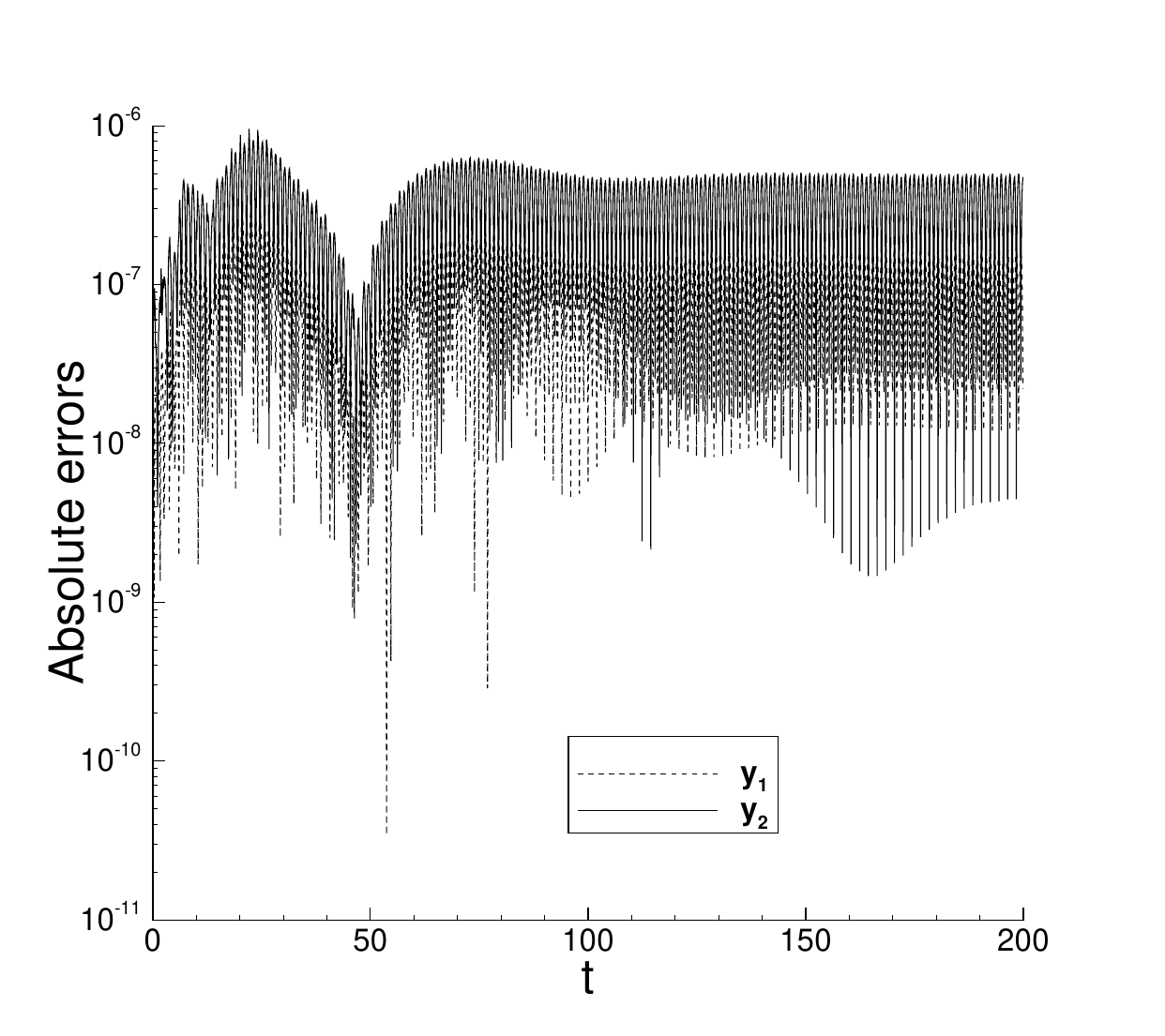}(c)
  }
  \caption{Forced pendulum: Comparison 
    of (a) $y_1(t)$ and (b) $y_2(t)$ between the NN-Exp-S1 solution and
    the reference solution. (c) Absolute-error histories of
    the NN-Exp-S1 solution for $y_1(t)$ and $y_2(t)$.
    NN-Exp-S1: $M=1200$, $Q=2000$, $\delta_m=5.0$,
    and $R_m=0.4$; Other parameter values are given in Table~\ref{tab_2}.
    Reference solution: obtained by the scipy DOP853 method with relative tolerance $10^{-13}$
    and absolute tolerance $10^{-16}$.
  }
  \label{fg_13}
\end{figure}

Figure~\ref{fg_13} illustrates the solutions for $y_1(t)$ and $y_2(t)$,
and their absolute errors, for $t\in[0,200]$ obtained by NN-Exp-S1.
The errors are computed against a reference solution attained by the scipy DOP853 method
with absolute tolerance $10^{-16}$ and  relative tolerance
$10^{-13}$. The reference solution is also shown in Figures~\ref{fg_13}(a,b) for comparison.
The parameter values  are provided in the figure
caption or in Table~\ref{tab_2}.
The NN solutions for $y_1(t)$ and $y_2(t)$  are highly accurate,
with a maximum error on the order
of $10^{-7}$ for a step size $\Delta t=0.1$.

\begin{figure}
  \centerline{
    \includegraphics[width=2in]{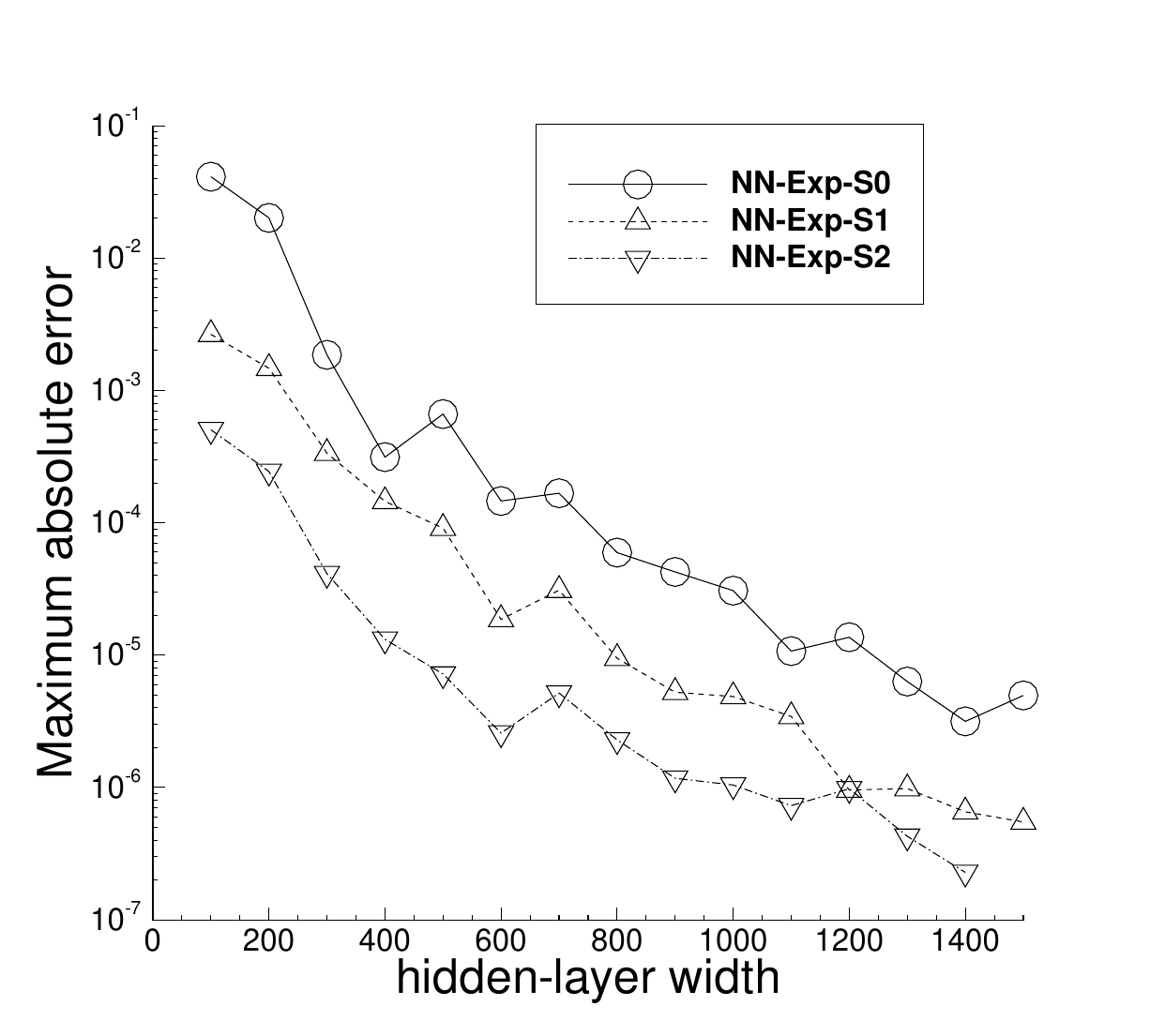}(a)
    \includegraphics[width=2in]{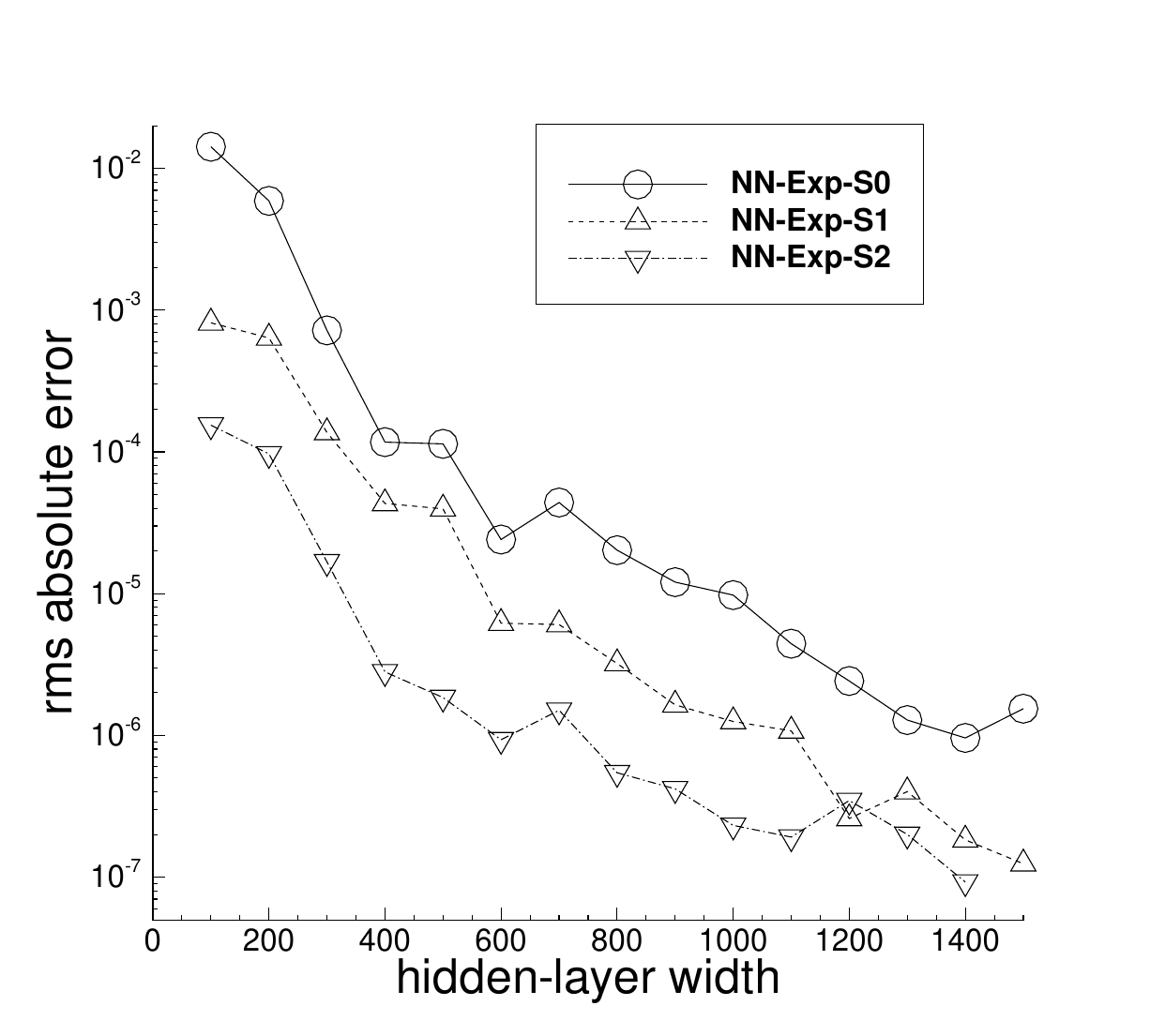}(b)
    \includegraphics[width=2in]{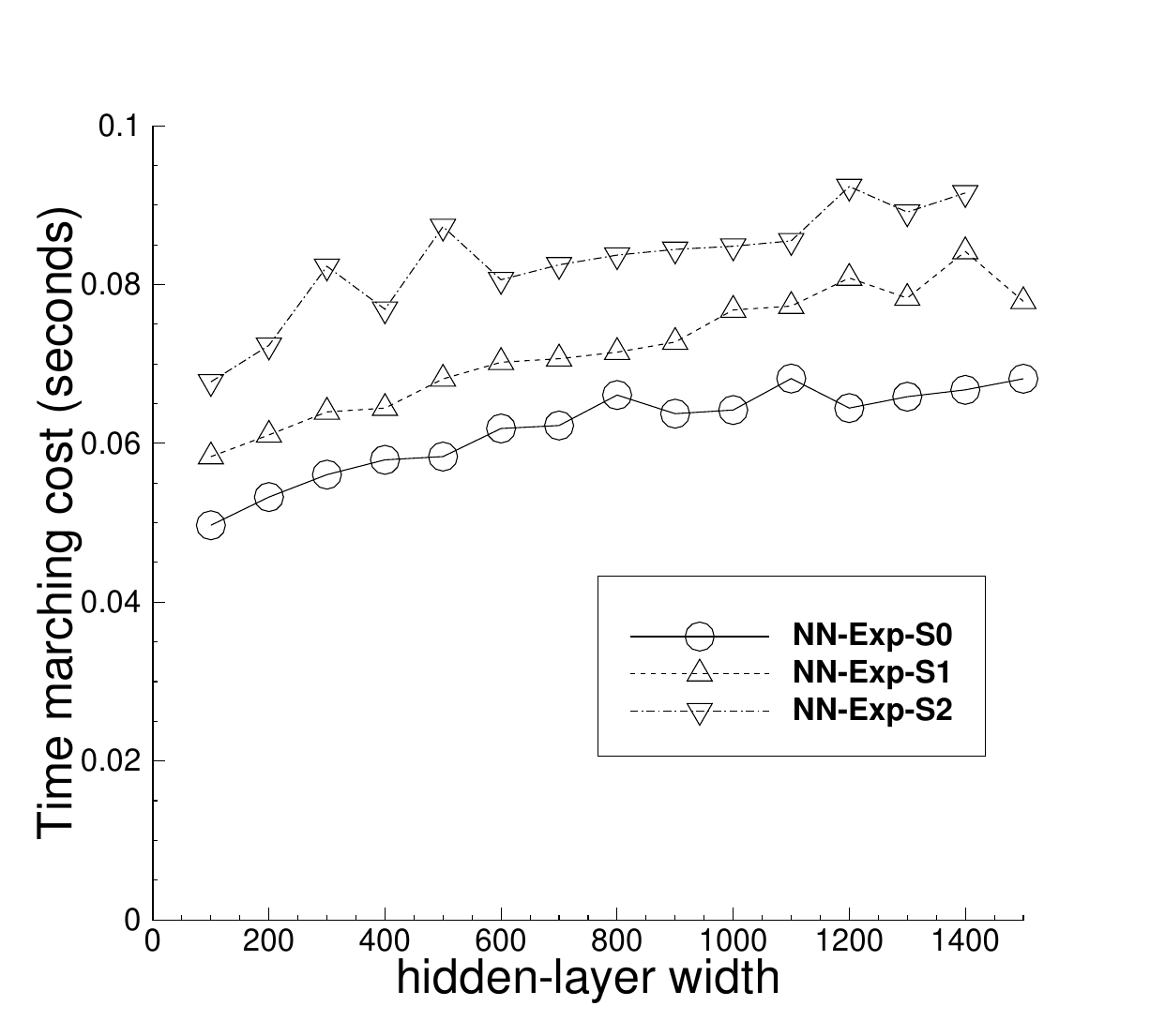}(c)
  }
  \caption{Forced pendulum:
    Comparison of (a) the maximum and (b) the rms solution errors,
    and (c) the time marching cost (wall
    time) versus the hidden-layer width $M$ for
    the NN algorithms (NN-Exp-S0, NN-Exp-S1, and NN-Exp-S2).
    NN-Exp-S0: $Q=2000$, $R_m=0.3$, $\delta_m=7.0$;
    NN-Exp-S1: $Q=2000$, $R_m=0.4$, $\delta_m=5.0$.
    NN-Exp-S2: $Q=1500$, $R_m=0.4$, $\delta_m=5.0$.
    Other parameter values are given in Table~\ref{tab_2}.
  }
  \label{fg_14}
\end{figure}

A comparison of the solution errors and the time-marching cost of different NN algorithms
is provided in Figure~\ref{fg_14}. Here we plot the maximum and rms solution errors
on $t\in[0,200]$, and the time-marching time,
versus the hidden-layer width $M$
obtained by NN-Exp-S0, NN-Exp-S1 and NN-Exp-S2.
The parameter values are provided in the figure caption or in
Table~\ref{tab_2}.
The solution errors decrease nearly exponentially
with increasing  $M$, while the time-marching cost grows
only quasi-linearly as $M$ increases.
NN-Exp-S2 is more accurate than NN-Exp-S1, which in turn
is more accurate than NN-Exp-S0.
In terms of the time-marching cost, NN-Exp-S2 is the most expensive
and NN-Exp-S0 is the least expensive among them.

\begin{figure}
  \centerline{
    \includegraphics[width=2in]{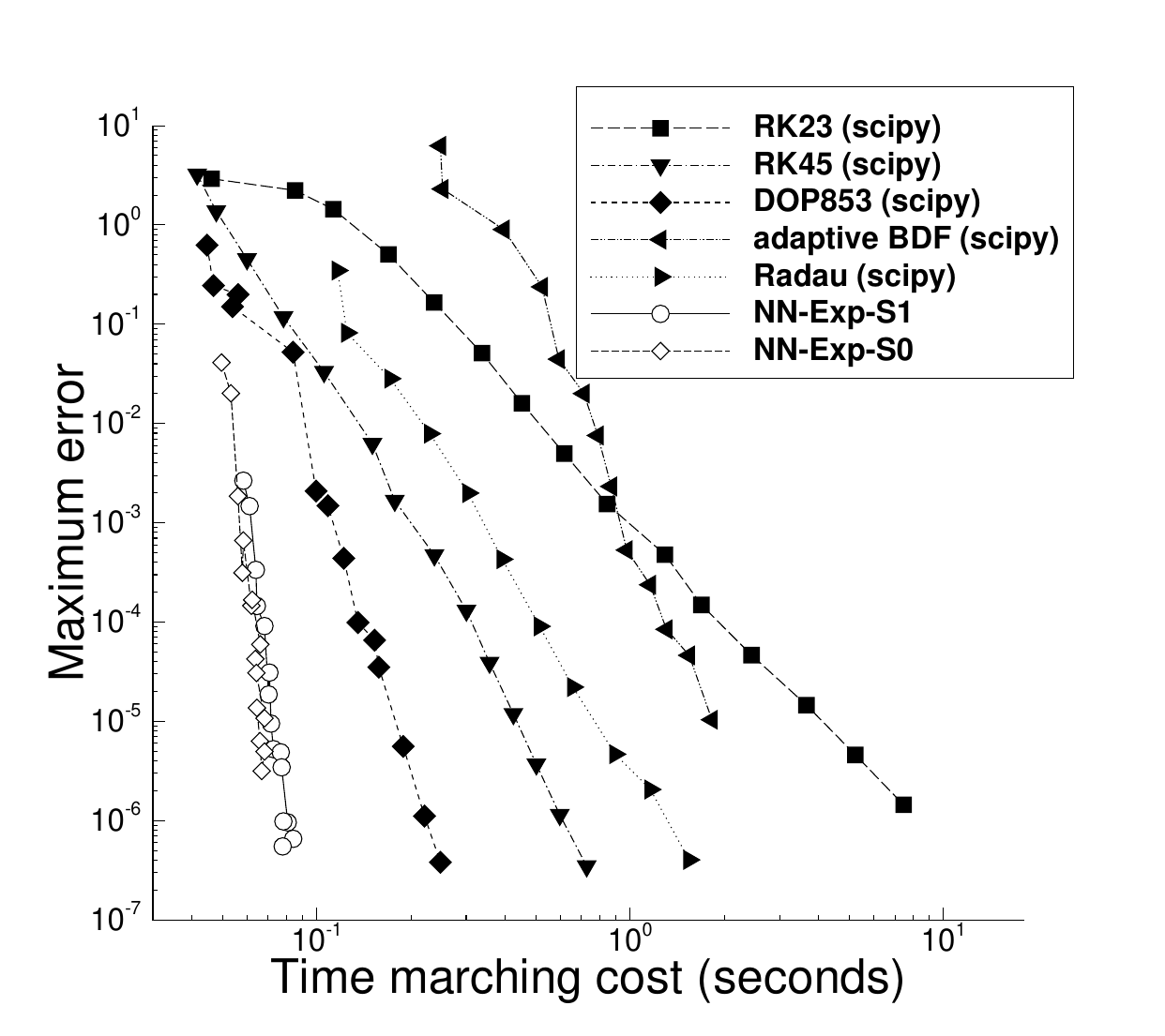}(a)
    \includegraphics[width=2in]{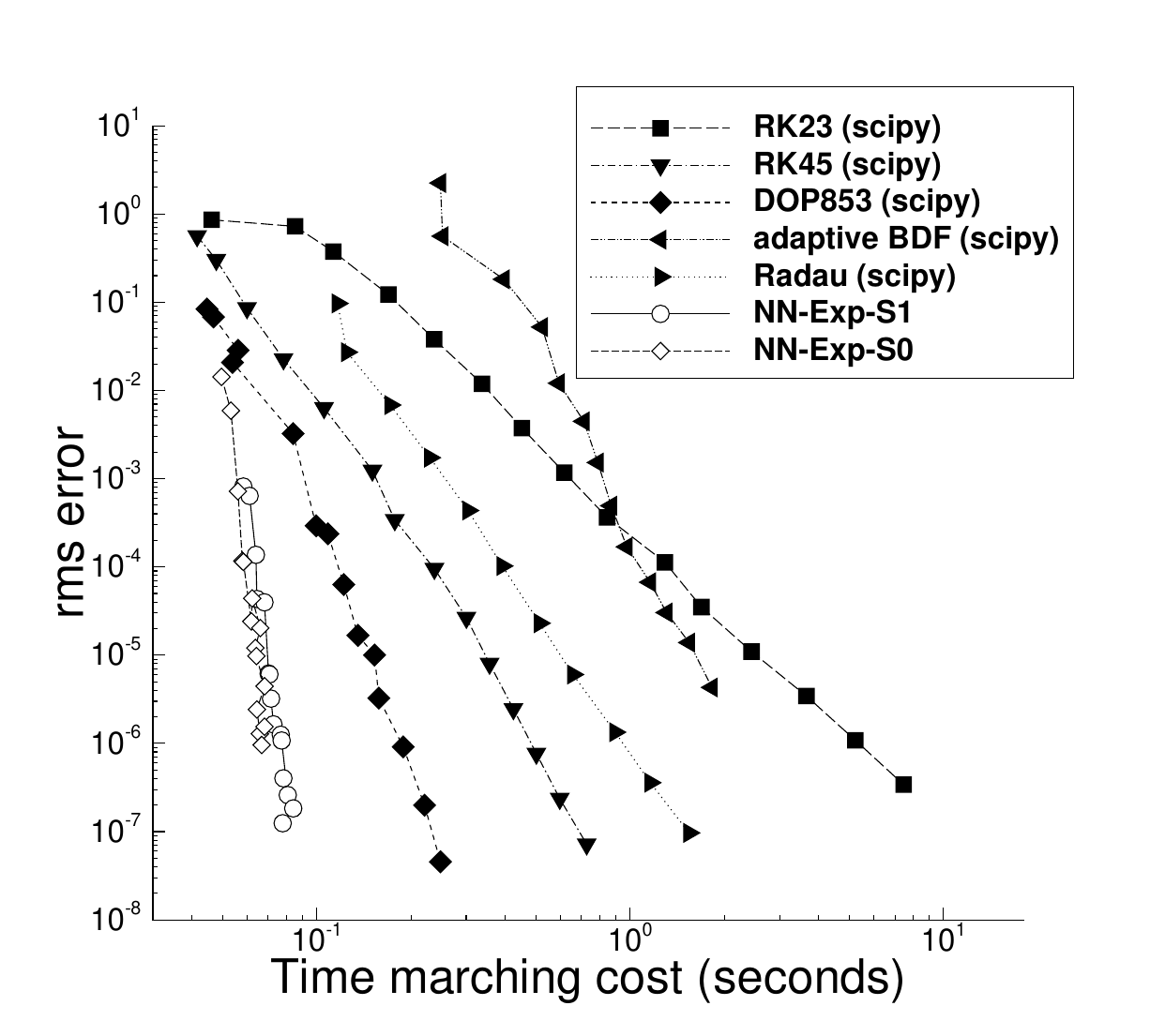}(b)
  }
  \caption{Forced pendulum:
    Comparison of (a) the maximum and (b) the rms solution errors versus the time-marching
    cost (wall time) between the NN algorithms (NN-Exp-S0 and NN-Exp-S1)
    and the scipy methods.
    Data for NN-Exp-S0 and NN-Exp-S1 correspond to those from Figure~\ref{fg_14}.
    Scipy methods: absolute tolerance $10^{-16}$,
    data points corresponding to different relative tolerance values,
    dense output on points corresponding to $\Delta t=0.1$
    for $t\in[0,200]$.
  }
  \label{fg_15}
\end{figure}

A comparison of the computational performance (accuracy versus cost)
between the NN algorithms and the scipy methods is shown in Figure~\ref{fg_15}.
Here we plot the maximum and rms solution errors as a function of
the time-marching time (for $t\in[0,200]$) obtained by the NN-Exp-S0 and NN-Exp-S1
algorithms and the scipy methods.
The two NN algorithms exhibit essentially the same performance.
Among the scipy methods, DOP853 shows the best performance,  followed
by RK45 and the other methods. Both NN-Exp-S0 and NN-Exp-S1 markedly outperform DOP853 and the
other scipy methods.


\subsection{Van der Pol Oscillator}
\label{sec_vdp}

\begin{table}[tb]
  \centering
  \begin{tabular}{l|l|l}
    \hline
   $\mu=5$ & domain: $(y_{01},y_{02},\xi)\in [-2.05,2.05]\times[-8,8]\times[0,0.035]$
    & NN ($\varphi$-subnet): $[3, M, 2]$ \\
    & sub-domains: 3, along $y_{01}$, uniform & activation function: Gaussian \\
   & $r$: $0.1$  & $\delta_m$: $1$  \\
   & $Q$: $1500$, random & $R_m$: to be specified \\
    & $\Delta t$: $0.03$ (time-marching) & time: $t\in[0,120]$  \\
    \hline
    $\mu=100$ & domain: $(y_{01},y_{02},\xi)\in [-2.05,2.05]\times[-140,140]\times[0,h_{\max}]$
    & NN ($\varphi$-subnet): $[3, M, 2]$ \\
    & sub-domains: 3 or 5 along $y_{01}$ (uniform); & activation function: Gaussian \\
    & $\quad$ 5 along $y_{02}$ (non-uniform), sub-domain & $\delta_m$: $1$  \\
    & $\quad$ boundaries: $[-140, -0.5, -0.03, 0.03, 0.5, 140]$ & $R_m$: to be specified \\
    & $\quad$ $h_{\max}$ varied for different sub-domains & $\Delta t$: quasi-adaptive  \\
    & $r$: $0.1$ along $y_{01}$ and $0.05$ along $y_{02}$, or $r=0$ in $y_{01}$ and $y_{02}$
    & time: $t\in[0,300]$  \\
   & $Q$: $1000$ or $1400$, random &  \\
    \hline
  \end{tabular}
  \caption{NN simulation parameters for the van der Pol oscillator
    (Section~\ref{sec_vdp}).
  }
  \label{tab_3}
\end{table}

In this test we evaluate the learned NN algorithms
using the van der Pol oscillator problem:
\begin{subequations}\label{eq_57}
  \begin{align}
    &
    \frac{dy_1}{dt} = y_2, \quad
    \frac{dy_2}{dt} = \mu (1-y_1^2)y_2 - y_1, \\
    &
    y_1(t_0) = y_{01}, \quad y_2(t_0) = y_{02},
  \end{align}
\end{subequations}
where $y(t)=(y_1(t),y_2(t))$ are the unknowns,
$\mu>0$ is a constant parameter, $t_0$ is the initial time, and
$y_0=(y_{01},y_{02})$ are the initial data. We employ $t_0=0$
and  $(y_{01},y_{02})=(2,0)$  for time integration.
This problem becomes stiff when $\mu$ is large. We use
$\mu=5$ for a non-stiff case and $\mu=100$ for a stiff case to
test the performance of the learned NN algorithms.

Let's first consider the problem~\eqref{eq_57} with $\mu=5$.
The simulation parameters related to the NN algorithm
are listed in Table~\ref{tab_3}.
In particular, we partition the training domain into uniform sub-domains
along the $y_{01}$ direction, and employ an ELM network with
architecture $[3,M,2]$ and the Gaussian activation function
for the $\varphi$-subnet
to learn $\psi(y_{01},y_{02},\xi)$,
where the hidden-layer width $M$ is varied. In this table $R_m$ again
denotes the scope of random hidden-layer coefficients.

\begin{figure}
  \centerline{
    \includegraphics[width=2in]{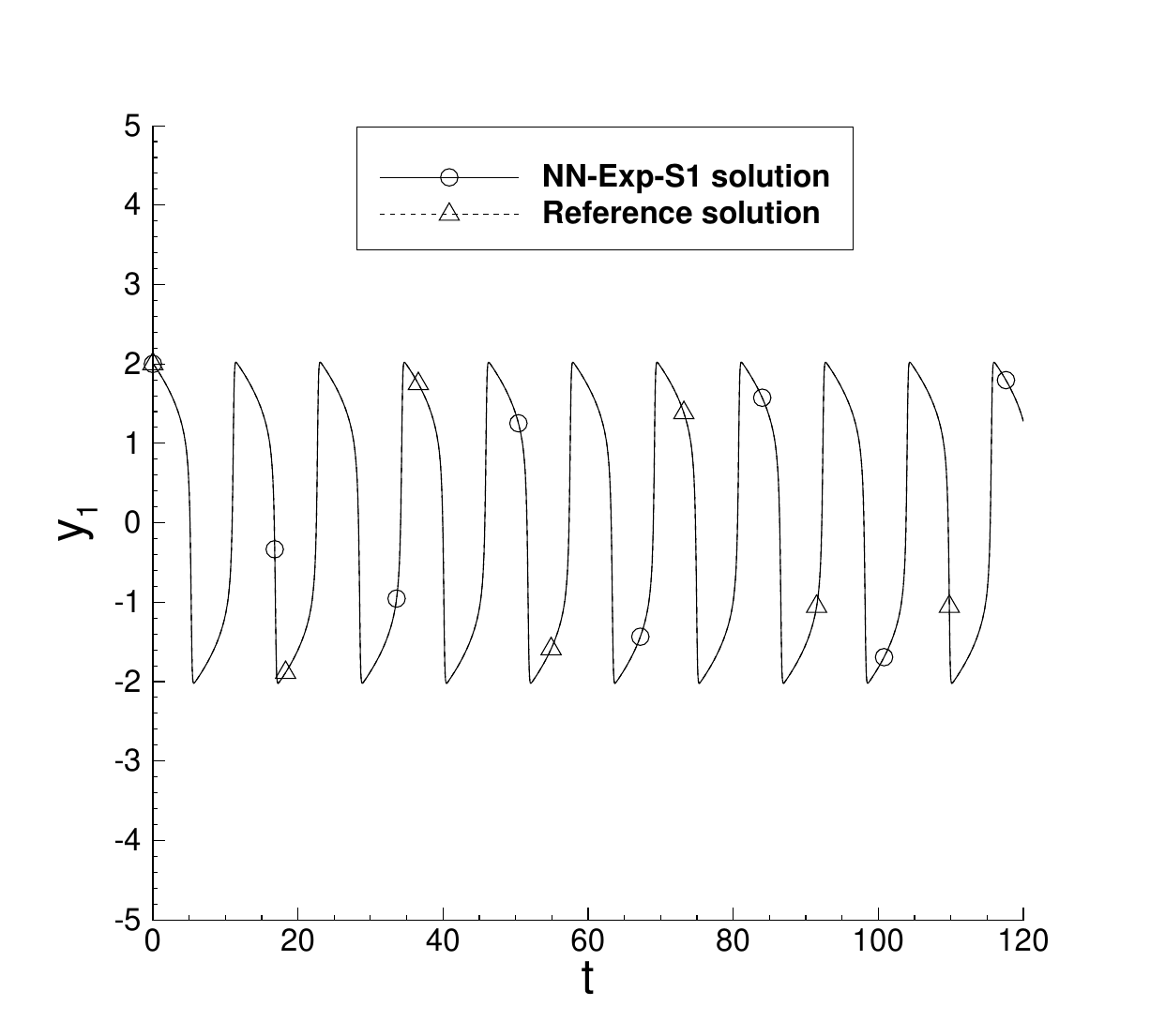}(a)
    \includegraphics[width=2in]{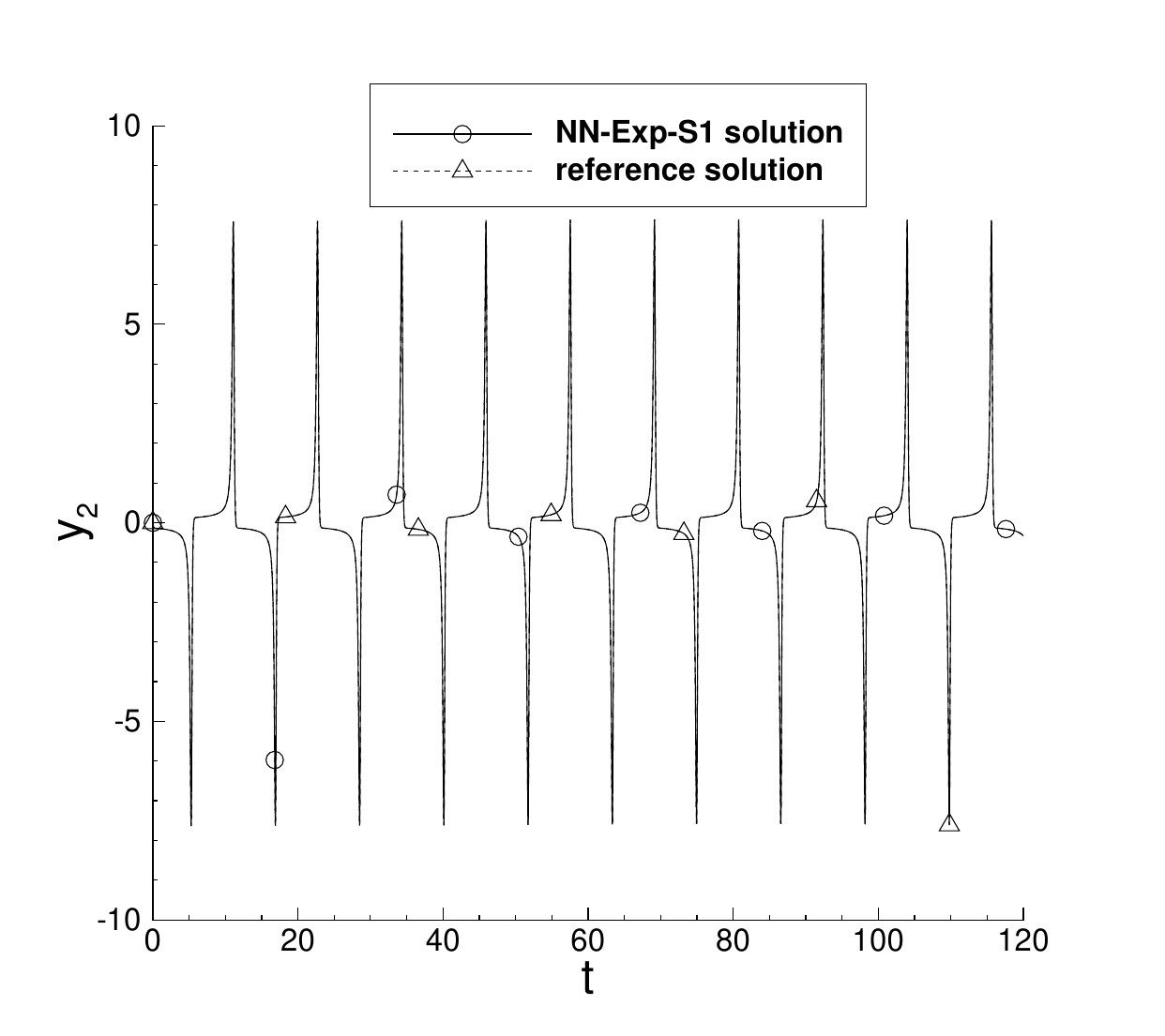}(b)
    \includegraphics[width=2in]{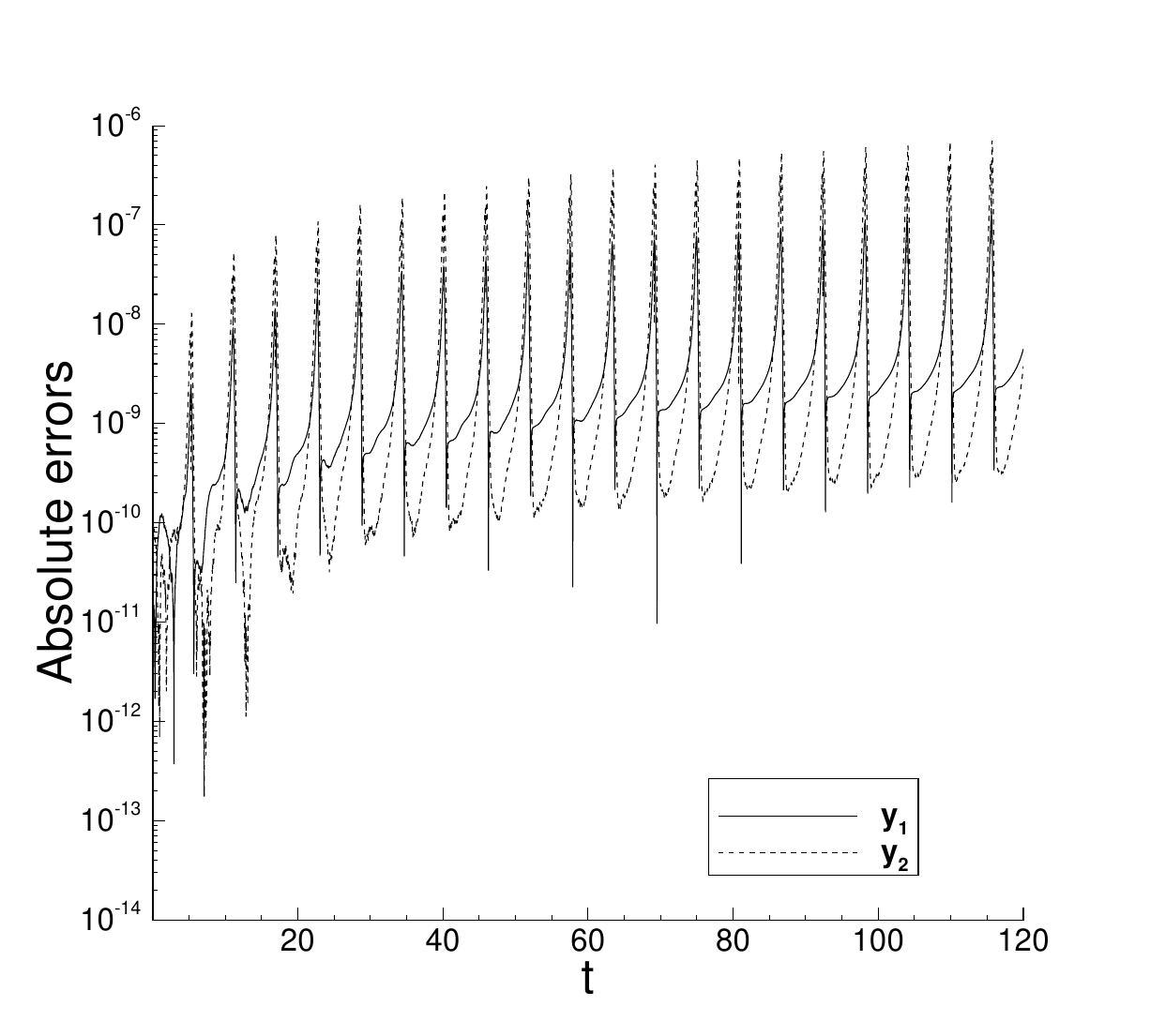}(c)
  }
  \caption{Van der Pol oscillator ($\mu=5$):
    Comparison of  (a) $y_1(t)$ and (b) $y_2(t)$ between the
    NN-Exp-S1 solution and the reference solution. (c) Absolute-error histories of
    the NN-Exp-S1 solution for $y_1(t)$ and $y_2(t)$.
    NN-Exp-S1: $M=1100$, $Q=1500$,
    $R_m=0.5$; Other parameter values are
    given in Table~\ref{tab_3}.
    Reference solution: obtained by the scipy DOP853 method with absolute tolerance $10^{-16}$ and
    relative tolerance $10^{-13}$.
  }
  \label{fg_16}
\end{figure}

Figure~\ref{fg_16} illustrates characteristics of the NN solutions
and their accuracy by comparing the  $y_1(t)$ and $y_2(t)$
obtained by NN-Exp-S1 and a reference solution computed using the scipy DOP853 method
with a sufficiently small tolerance. The absolute errors for $y_1(t)$ and
$y_2(t)$ are also shown. The parameter values for these results
are either specified in the caption or listed in Table~\ref{tab_3}.
The NN-Exp-S1 solution is observed to be highly accurate, with a maximum error
on the order of $10^{-7}$ for $t\in[0,120]$.

\begin{figure}
  \centerline{
    \includegraphics[width=2in]{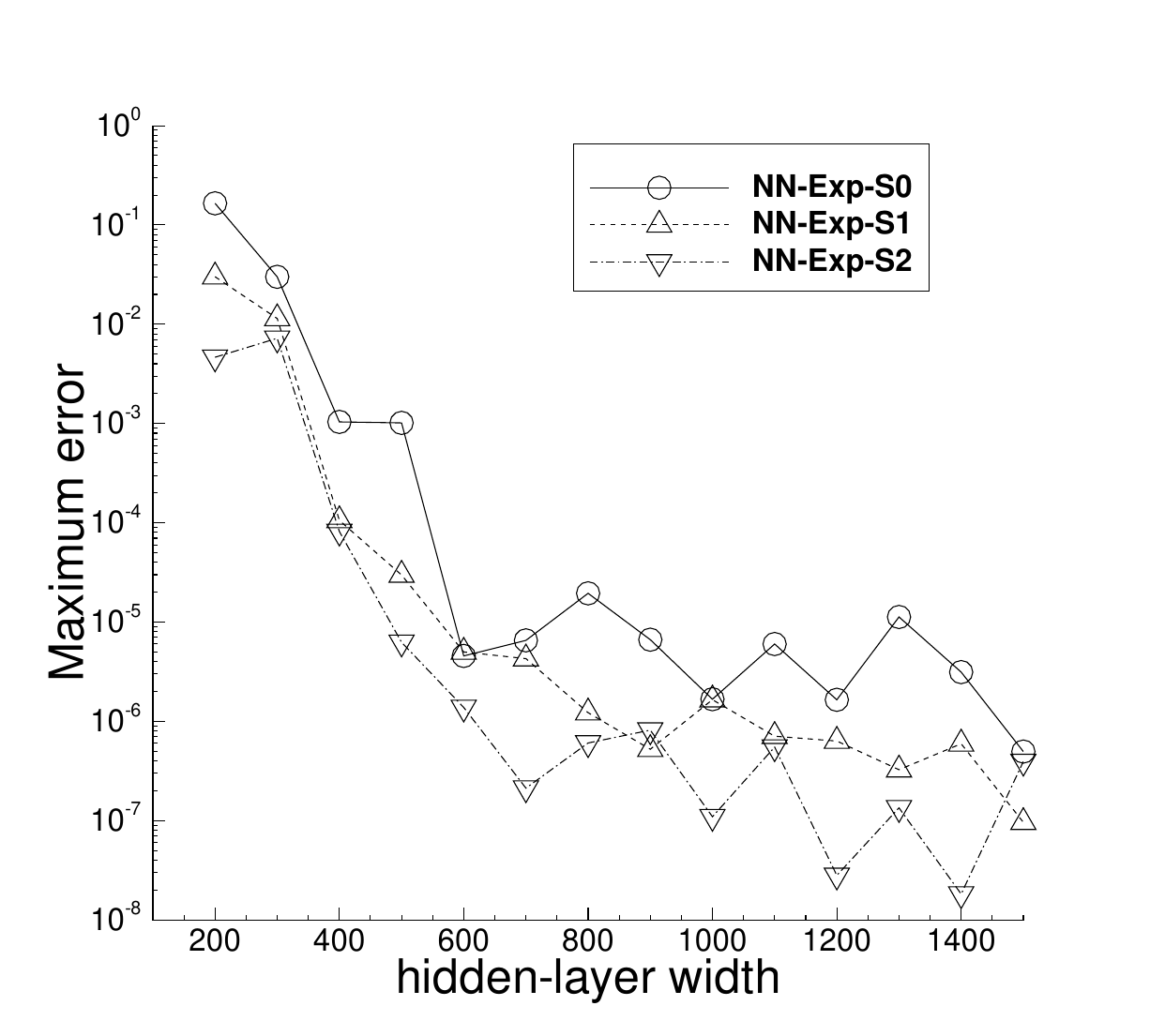}(a)
    \includegraphics[width=2in]{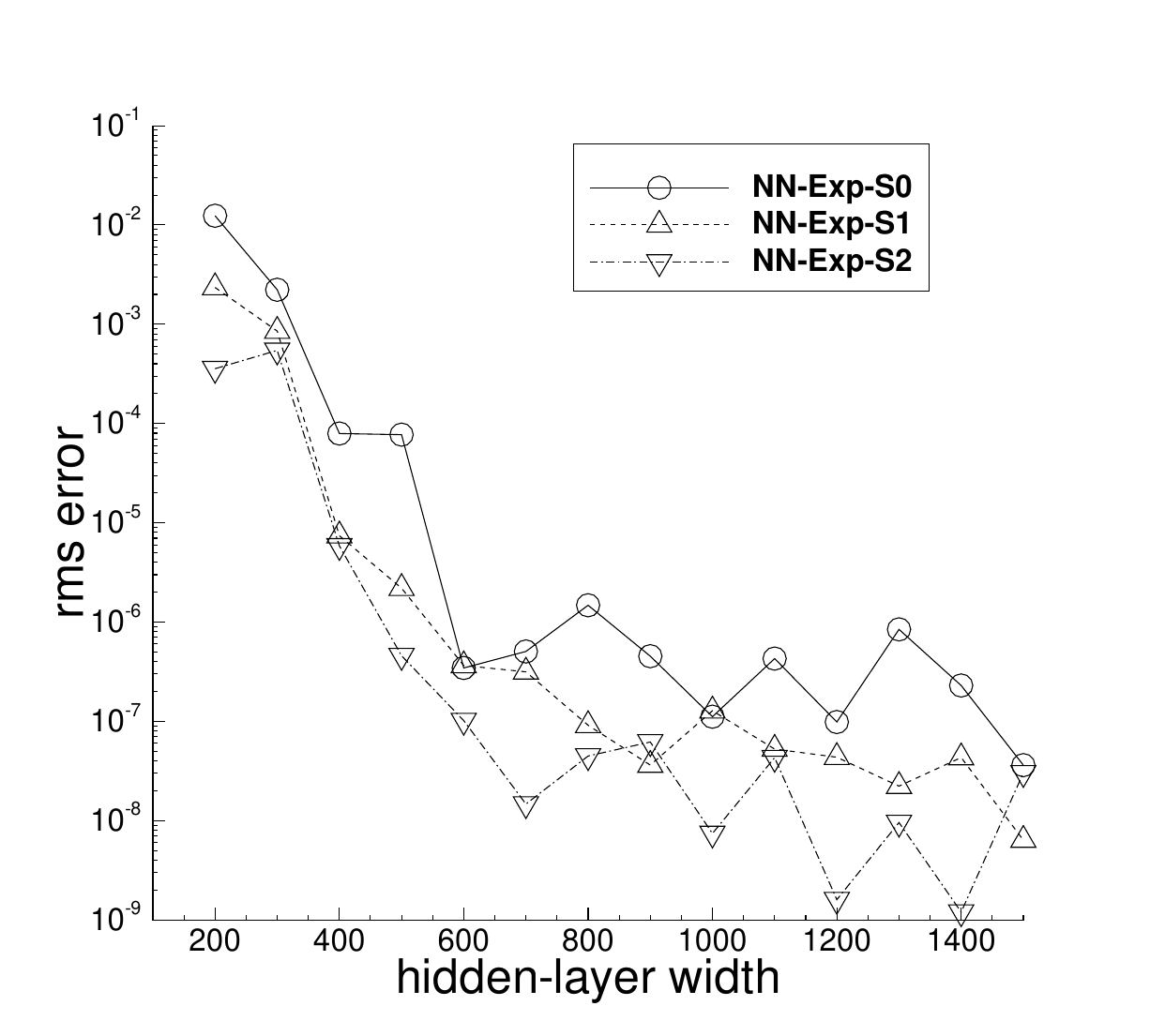}(b)
    \includegraphics[width=2in]{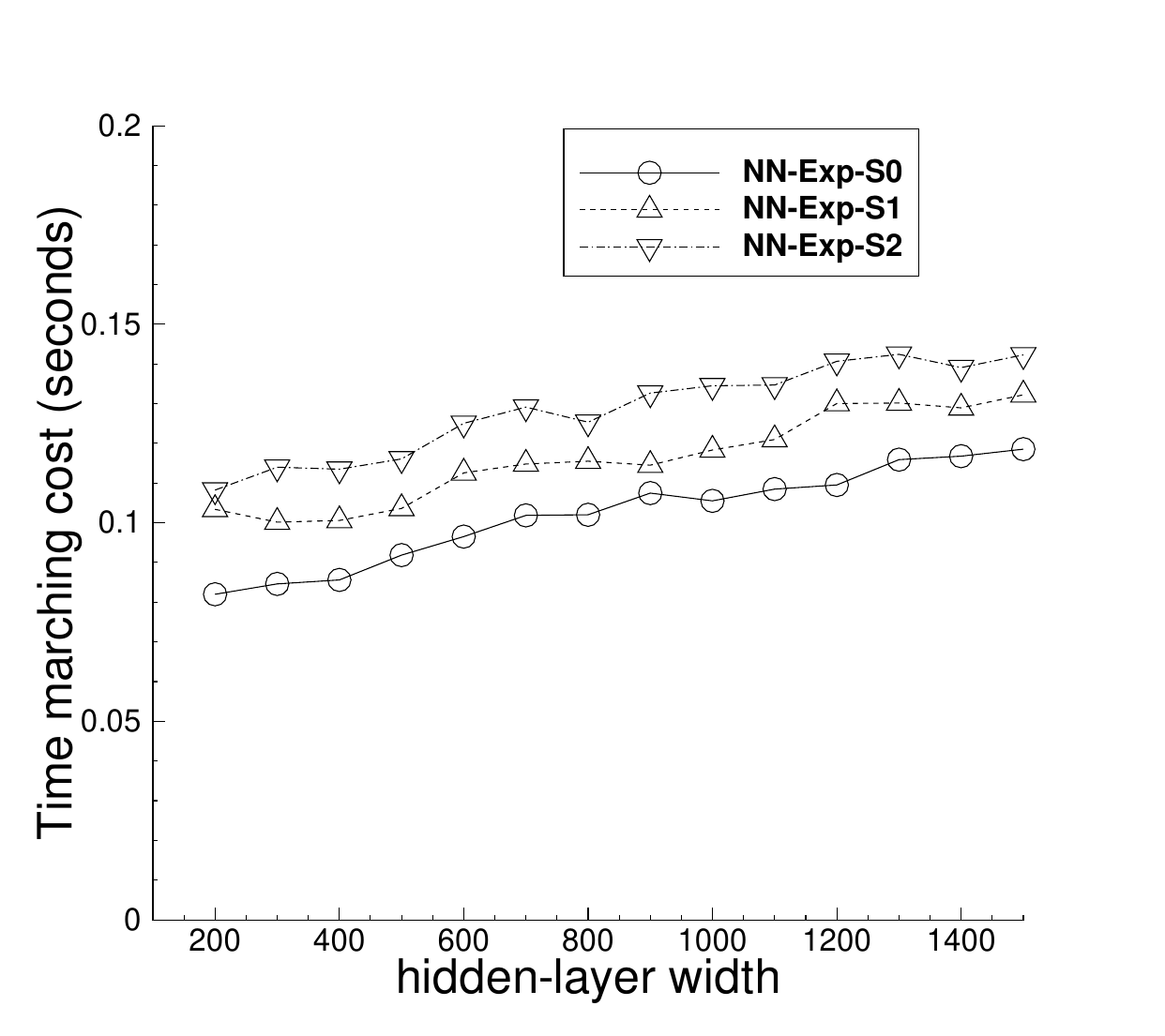}(c)
  }
  \caption{Van der pol oscillator ($\mu=5$):
    Comparison of (a) the maximum and (b) the rms time-marching errors, and (c) the
    time-marching cost (wall time) versus the hidden-layer width $M$
    for the NN algorithms (NN-Exp-S0, NN-Exp-S1, NN-Exp-S2).
    NN-Exp-S0: $R_m=0.45$;
    NN-Exp-S1: $R_m=0.5$; NN-Exp-S2: $R_m=0.55$.
    See Table~\ref{tab_3} for the other parameter values.
  }
  \label{fg_17}
\end{figure}

The convergence behavior of the NN algorithms is exemplified by Figure~\ref{fg_17}.
Here we plot the maximum and rms time-marching errors, as well as the time-marching cost,
of NN-Exp-S0, NN-Exp-S1 and NN-Exp-S2 for $t\in[0,120]$ as a function of
the hidden-layer width $M$ in the $\varphi$-subnet architecture.
The figure caption and Table~\ref{tab_3} provide all the parameter values
corresponding to these results.
It is evident that the  NN solution errors decrease nearly exponentially
with increasing $M$ (before saturation), and that
the time-marching cost grows only quasi-linearly
as $M$ increases.

\begin{figure}
  \centerline{
    \includegraphics[width=2in]{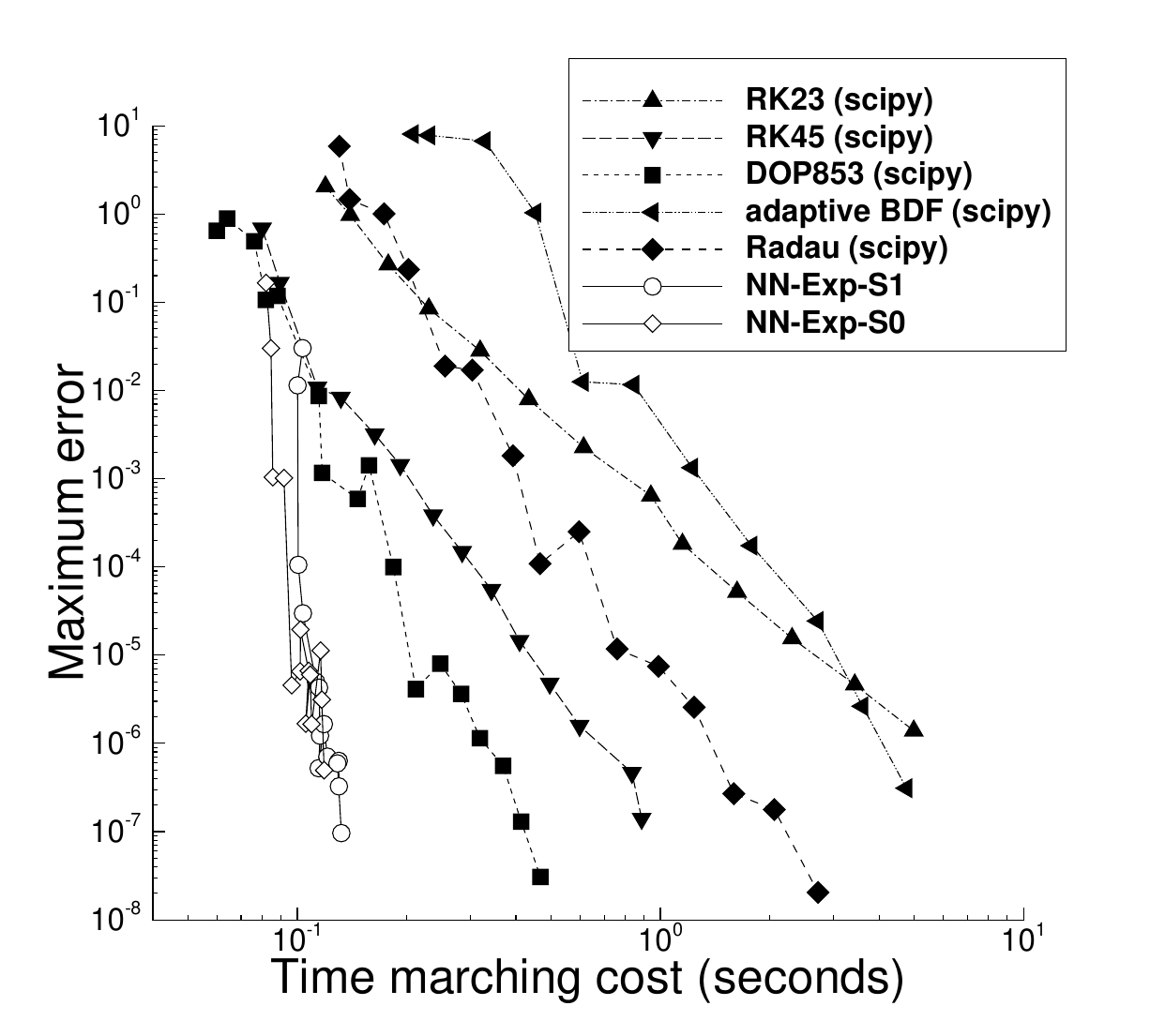}(a)
    \includegraphics[width=2in]{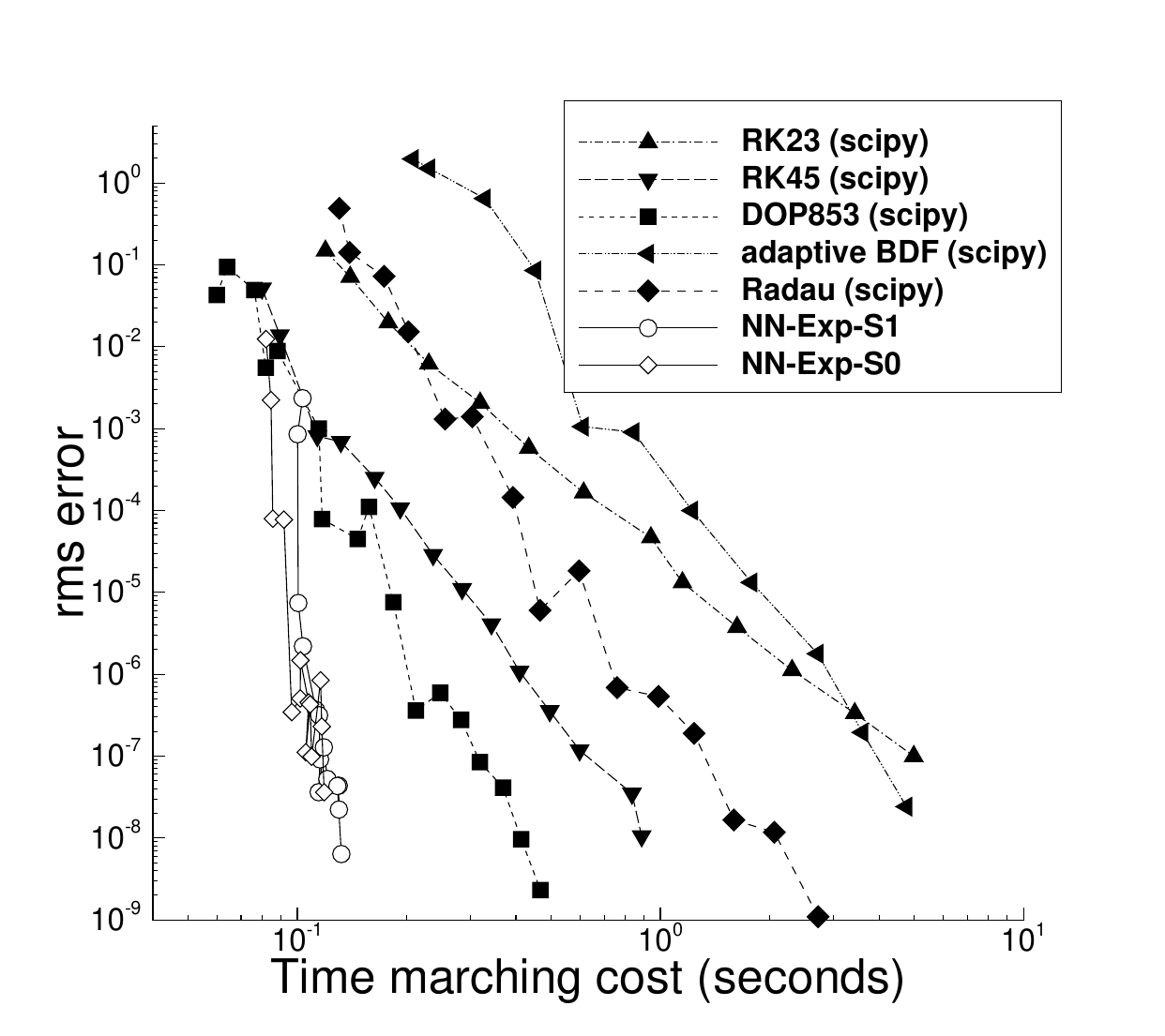}(b)
  }
  \caption{Van der Pol oscillator ($\mu=5$):
    Comparison of (a) the maximum error and (b) the rms error versus the time marching
    cost (wall time) between the NN algorithms (NN-Exp-S0 and NN-Exp-S1)
    and the scipy methods.
    Data for the NN algorithms correspond to those of NN-Exp-S0 and NN-Exp-S1 in
    Figure~\ref{fg_17}.
    Scipy methods: absolute tolerance $10^{-16}$,
    relative tolerance varied for different data points,
    dense output on points corresponding to $\Delta t=0.03$ for $t\in[0,120]$.
  }
  \label{fg_18}
\end{figure}

A performance comparison (accuracy vs. cost) between the current NN algorithms
and the scipy methods is provided in Figure~\ref{fg_18}.
Here we plot the maximum and rms time-marching errors (for $t\in[0,120]$)
as a function of the time-marching time obtained by 
NN-Exp-S0 and NN-Exp-S1, and by the scipy methods.
The data for NN-Exp-S0 and NN-Exp-S1 correspond to those in Figure~\ref{fg_17},
while those for the scipy methods are attained by varying the relative tolerance values.
The performance of NN-Exp-S0 and NN-Exp-S1 is close,
with the former appearing slightly better.
The DOP853 method exhibits the best performance among the scipy methods,
followed by RK45 and Radau.
The NN-Exp-S0 and NN-Exp-S1 algorithms perform notably better than
the scipy methods for this non-stiff case  with the van der Pol oscillator.

\begin{figure}
  \centerline{
    \includegraphics[width=2in]{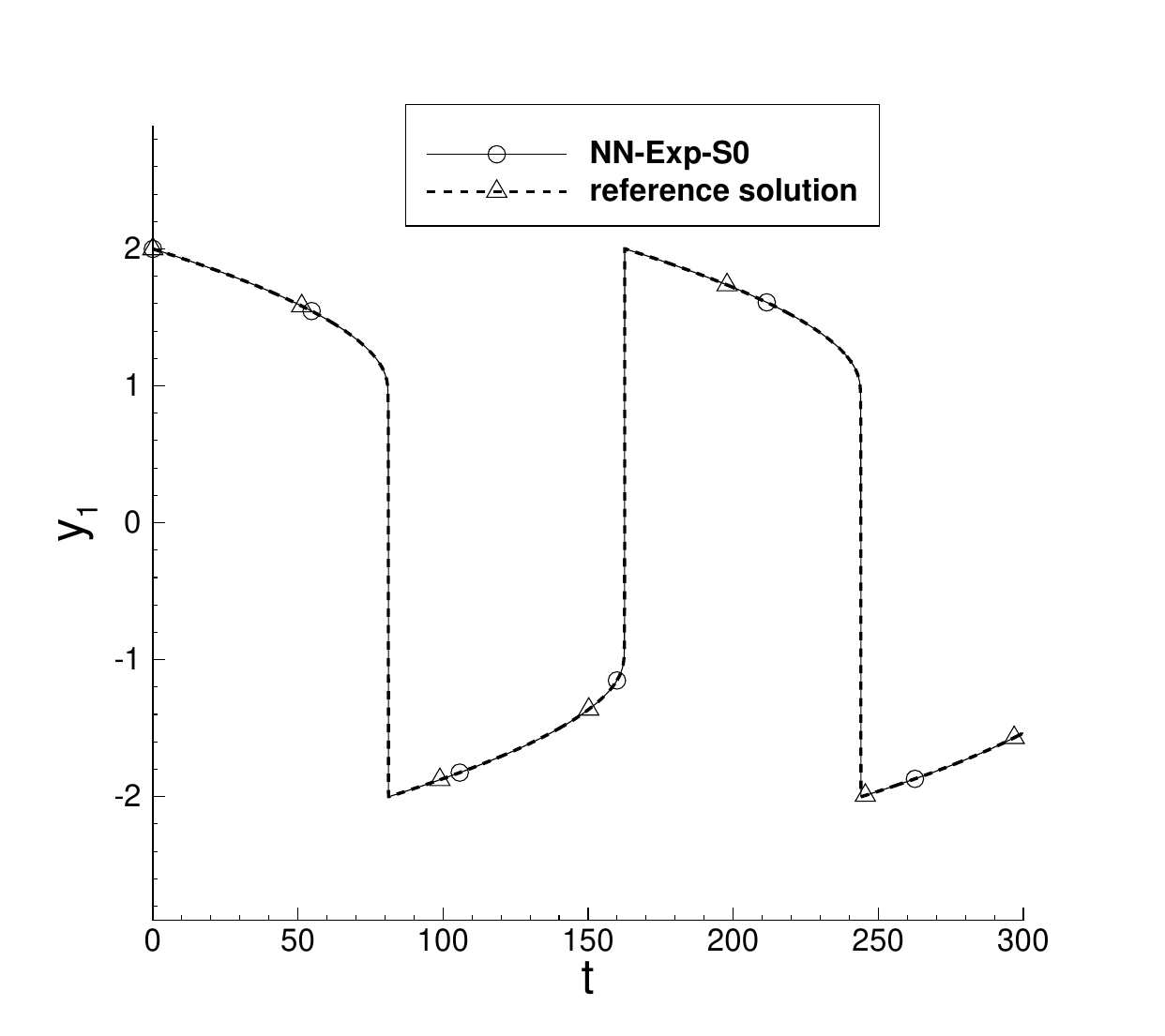}(a)
    \includegraphics[width=2in]{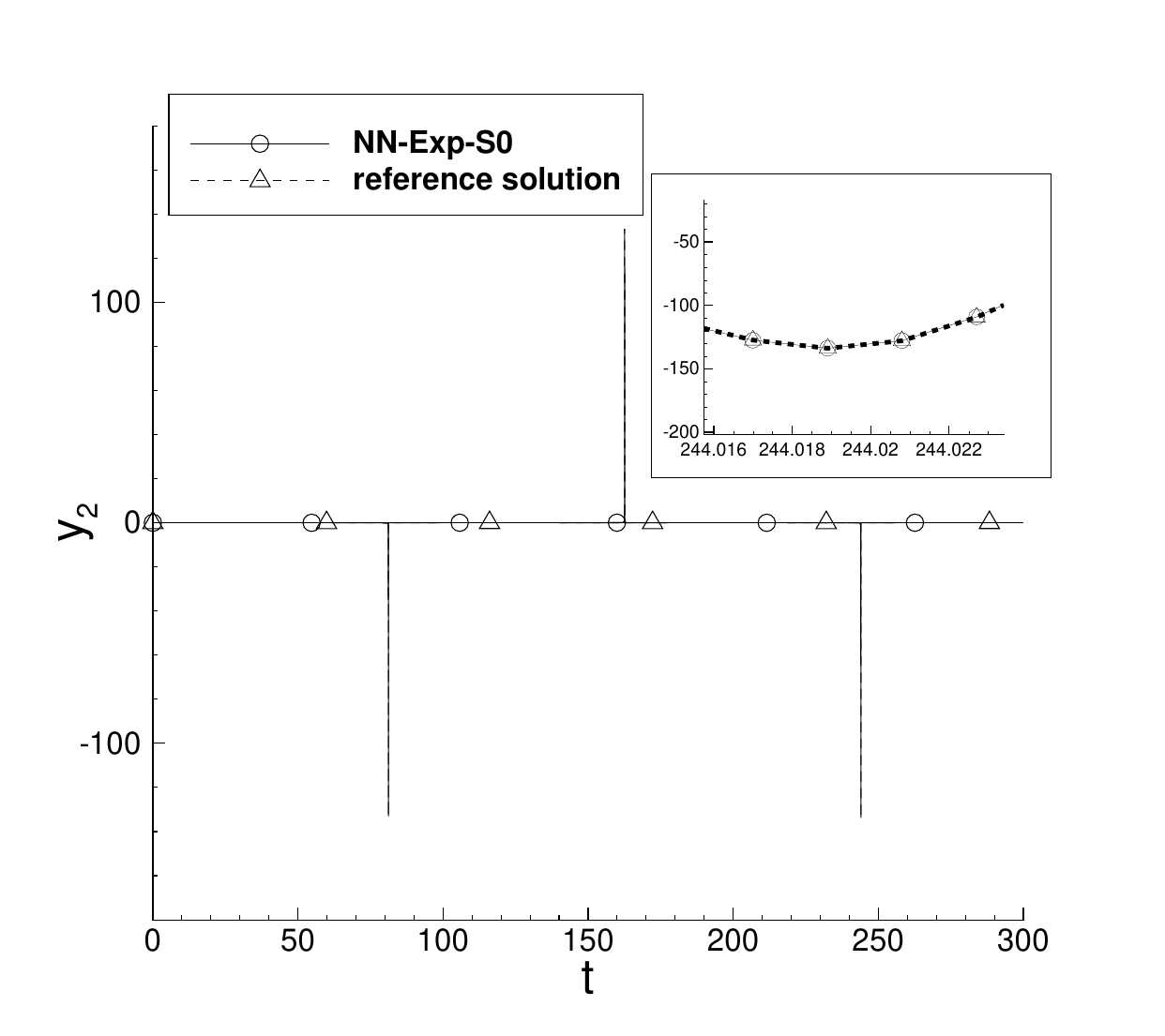}(b)
    \includegraphics[width=2in]{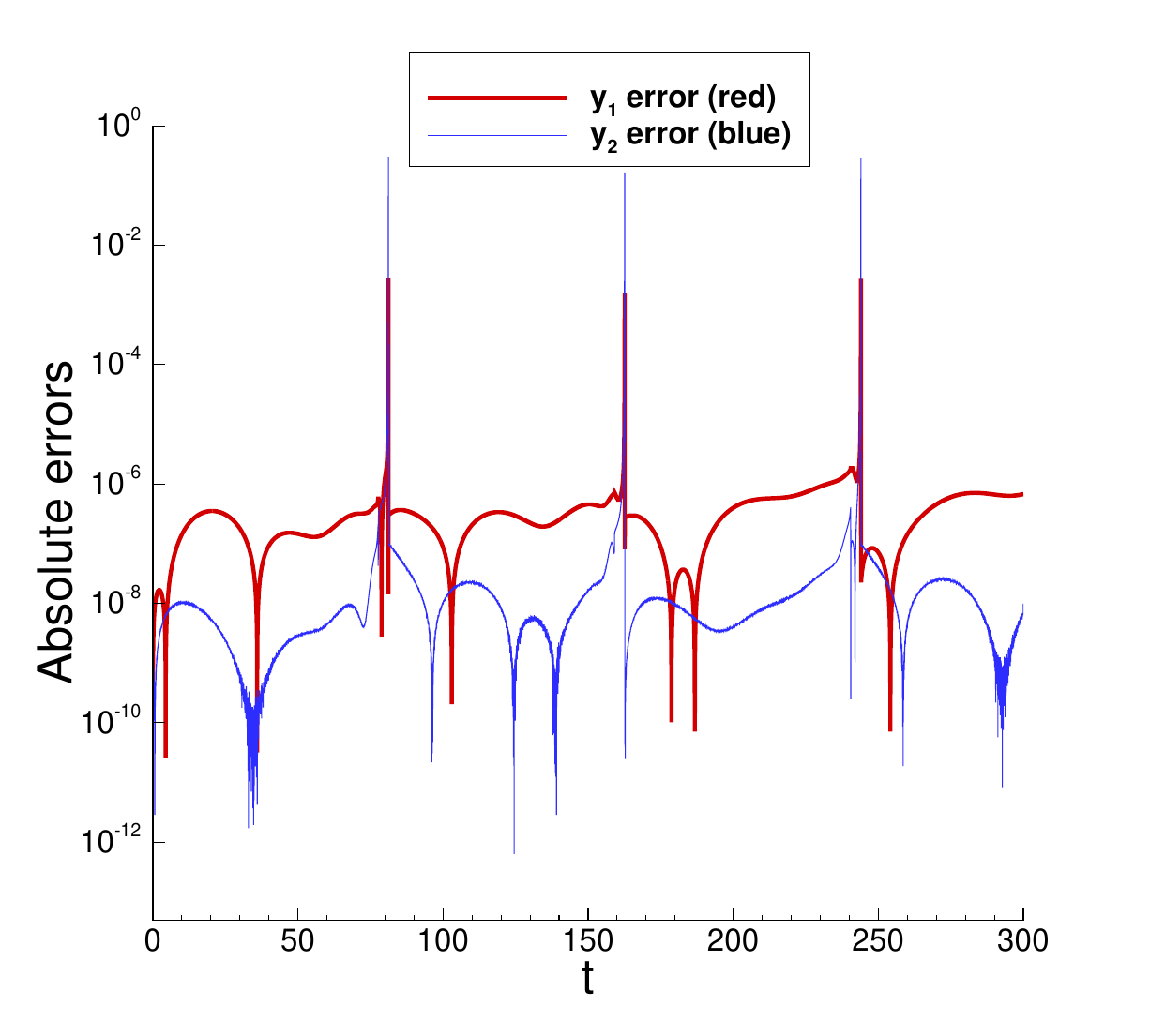}(c)
  }
  \caption{Van der Pol oscillator ($\mu=100$):
    Comparison of (a) $y_1(t)$ and (b) $y_2(t)$ between
    the NN-Exp-S0 solution and the reference solution.
    (c) Absolute-error histories of the NN-Exp-S0 solution
    for $y_1(t)$ and $y_2(t)$.
    NN-Exp-S0: 
    $3$ uniform sub-domains in $y_{01}$, with enlargement factor $r=0.1$;
    $5$ non-uniform
    sub-domains in $y_{02}$, with enlargement factor $r=0.05$,
    NN architecture: $[3, 800, 2]$
    with Gaussian activation;
    $R_m=0.75$, and $Q=1400$ on each sub-domain;
    See Table~\ref{tab_3} or the text for the other parameter values.
    Reference solution is computed by the scipy Radau method, with absolute tolerance
    $10^{-16}$ and relative tolerance $10^{-13}$.
    The inset of plot (b) shows a magnified view around $t=244.02$.
  }
  \label{fg_19}
\end{figure}

We next consider the case $\mu=100$, with which the problem~\eqref{eq_57} becomes stiff.
Since the velocity of the oscillator (i.e.~$y_2$) is very small in the
majority of time
but can increase to large magnitudes in bursts, using a constant step size
for time integration becomes less efficient and it is necessary to incorporate
some adaptive strategy into the NN algorithm.

We employ the following quasi-adaptive strategy for network training
and for time integration when solving this problem. We first choose a training domain
for $(y_{01},y_{02})\in[A_1,A_2]\times[B_1,B_2]$, and then partition this domain
into $m$ ($m\geqslant 1$) sub-domains along $y_{01}$ and
$n$ ($n\geqslant 1$) sub-domains along $y_{02}$.
Let $\Omega_{ij}=[a_i,a_{i+1}]\times[b_j,b_{j+1}]$ ($0\leqslant i\leqslant m-1$,
$0\leqslant j\leqslant n-1$) denote a sub-domain, where $A_1=a_0<a_1<\cdots<a_m=A_2$
and $B_1=b_0<b_1<\cdots<b_n=B_2$.
On $\Omega_{ij}$ we train a local NN  to learn $\psi(y_{01},y_{02},\xi)$
for  $(y_{01},y_{02},\xi)\in[a_i,a_{i+1}]\times[b_j,b_{j+1}]\times[0,h_{\max}^{(ij)}]$, where
$h_{\max}^{(ij)}$ depends on  $\Omega_{ij}$ and
can differ in different sub-domains. This provides the opportunity, during
time marching, to vary the step size based on the current value of $(y_1,y_2)$.
Specifically for this problem, we choose $h_{\max}$ to depend only on $j$,
i.e.~$(y_{01},y_{02},\xi)\in[a_i,a_{i+1}]\times[b_j,b_{j+1}]\times[0,h_{\max}^{(j)}]$
for $\Omega_{ij}$, and we partition the domain non-uniformly along
the $y_{02}$ direction.
%
For time integration, given $(y_{1k},y_{2k})$ at time $t_k$,
we compute the next approximation as follows:
\begin{enumerate}[(i)]
\item Determine the sub-domain $\Omega_{ij}$ such that $(y_{1k},y_{2k})\in\Omega_{ij}$,
  and let $h=0.95h_{\max}^{(j)}$.
\item Compute
  $y_{k+1}=(y_{1,k+1},y_{2,k+1})=\psi(y_{1k},y_{2k},h)$ and $t_{k+1}=t_k+h$.
\end{enumerate}

We use a training domain $(y_{01},y_{02})\in[-2.05,2.05]\times[-140,140]$ for this stiff case,
as shown in Table~\ref{tab_3}, and partition the domain into $5$ non-uniform
sub-domains along $y_{02}$, with the sub-domain boundaries
given by $y_{02}=[-140, -0.5, -0.03, 0.03, 0.5, 140]$.
Along the $y_{01}$ direction we partition the domain into $3$ uniform
sub-domains (for NN-Exp-S0), or $5$ uniform sub-domains (for NN-Exp-S1).
We use
\begin{equation*}
  \begin{split}
    &\mbs h_{\max}=(h_{\max}^{(0)},h_{\max}^{(1)},\dots,h_{\max}^{(4)})
    =(0.002,0.011,0.018,0.011,0.002),
    \quad \text{(for NN-Exp-S0)}; \\
    & \mbs h_{\max}=(h_{\max}^{(0)},h_{\max}^{(1)},\dots,h_{\max}^{(4)})
    =(0.002,0.011,0.025,0.011,0.002),
    \quad \text{(for NN-Exp-S1)}.
  \end{split}
\end{equation*}
The NN-Exp-S2 algorithm does not work for this stiff case (training fails to
converge), and
the implicit NN algorithms (NN-Imp-S1, NN-Imp-S2) are
not as competitive as the explicit NN-Exp-S0 and NN-Exp-S1 algorithms.

Figure~\ref{fg_19} is a comparison of $y_1(t)$ and $y_2(t)$
between the NN-Exp-S0 solution and a reference solution obtained by
the scipy Radau method, and the absolute errors
of the NN-Exp-S0 solution.
The values for the simulation parameters are provided in the figure
caption or Table~\ref{tab_3}.
The NN solution agrees well with the reference solution,
with the maximum absolute error on the order of $10^{-3}$ for $y_1$
and $10^{-1}$ for $y_2$. The inset of Figure~\ref{fg_19}(b) shows that
even at the bursts the NN method has captured the solution
accurately.

\begin{figure}
  \centerline{
    \includegraphics[width=2in]{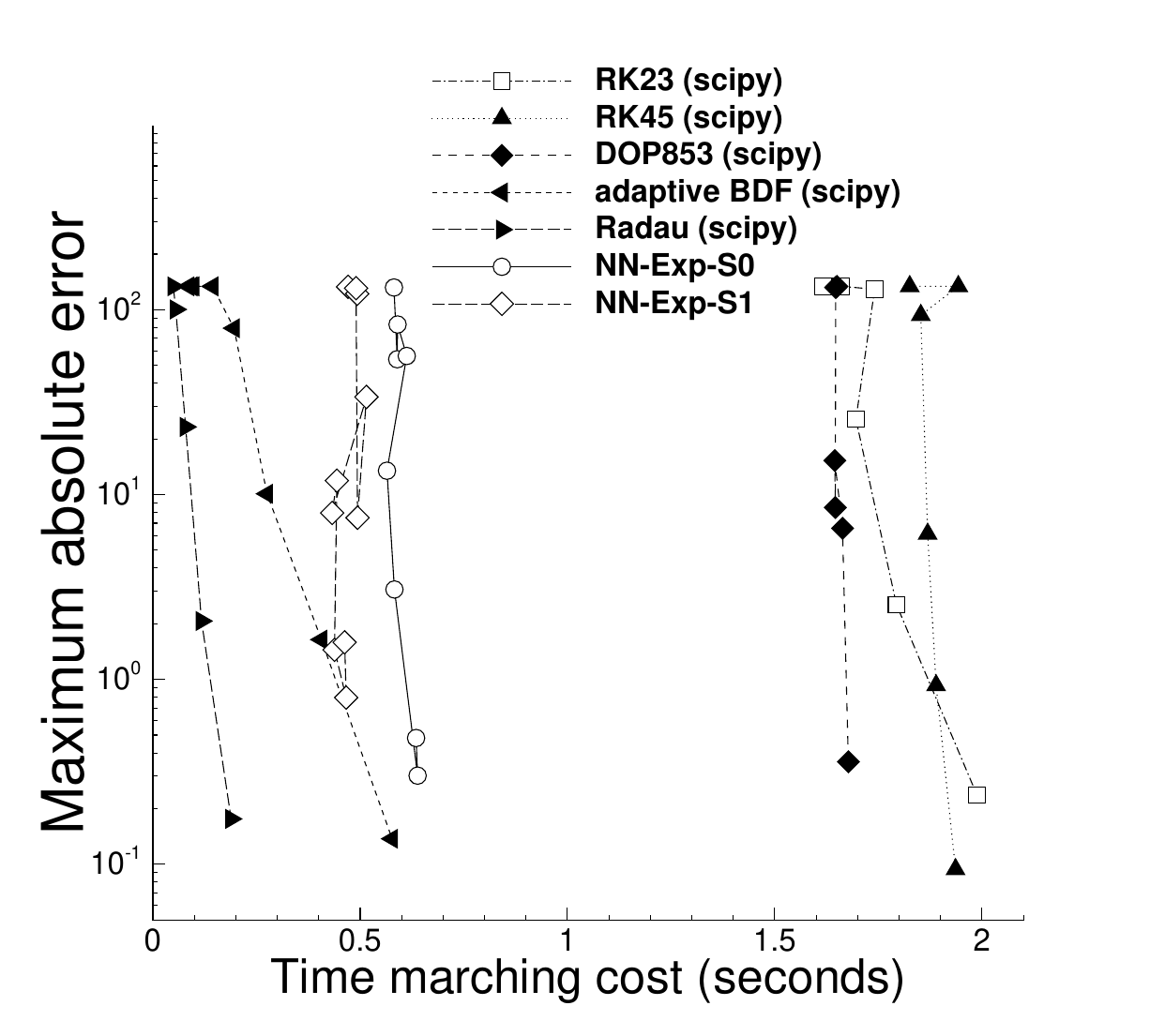}(a)
    \includegraphics[width=2in]{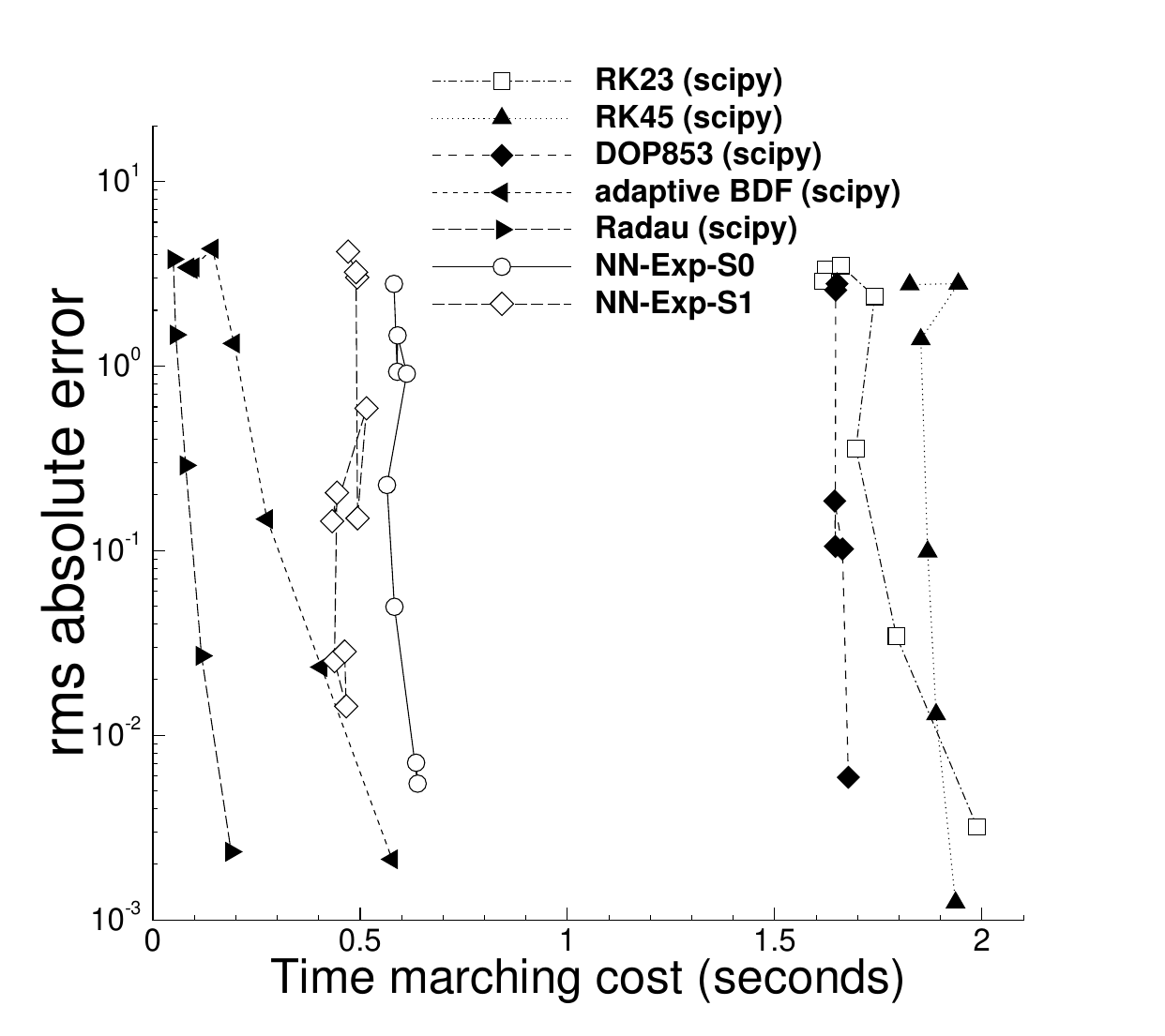}(a)
  }
  \caption{Van der Pol oscillator ($\mu=100$):
    Comparison of (a) the maximum and (b) the rms time-marching error 
    versus the time marching cost (wall time) between the NN algorithms (NN-Exp-S0,
    NN-Exp-S1) and the scipy methods for $t\in[0,300]$.
    Simulation parameters for NN-Exp-S0 follow those of Figure~\ref{fg_19}.
    NN-Exp-S1:
    5 uniform sub-domains in $y_{01}$ with enlargement factor $r=0$;
    $5$ non-uniform sub-domains in $y_{02}$ with enlargement factor $r=0$;
    $R_m=0.5$, and $Q=1000$ on each sub-domain;
    See Table~\ref{tab_3} and the text for the other parameter values.
    Scipy methods:
    absolute tolerance
    $10^{-16}$, data points corresponding to different relative tolerance values.
    Errors are computed with respect to a
    reference solution obtained by the scipy Radau method with absolute tolerance
    $10^{-16}$ and relative tolerance $10^{-13}$.
  }
  \label{fg_20}
\end{figure}

Figure~\ref{fg_20} compares the computational performance (accuracy versus cost)
of the current NN algorithms and the scipy methods for this stiff case.
It shows the maximum and the rms time-integration errors ($e_{\max}$ and $e_{\text{rms}}$)
on $t\in[0,300]$ as a function of the time-marching time for 
the NN-Exp-S0 and NN-Exp-S1 algorithms and the scipy methods.
The parameter values for the NN algorithms
are provided in the figure caption or in Table~\ref{tab_3}.
The scipy Radau method shows the best performance, followed by
the scipy BDF method. The NN-Exp-S0 and NN-Exp-S1 algorithms
are less competitive than the scipy Radau and BDF methods for this case.
The performance of the explicit scipy methods (DOP853, RK45, and RK23)
is significantly worse than the Radau/BDF and the NN-Exp-S0/NN-Exp-S1
algorithms.

\subsection{Lorenz63 Chaotic System}
\label{sec_lorenz63}

\begin{table}[tb]
  \centering
  \begin{tabular}{l|l}
    \hline
    domain: $(y_{01},y_{02},y_{03},\xi)\in [-20,20]\times[-25,25]\times[2,46]\times[0,0.012]$,
    & NN ($\varphi$-subnet): $[4, M, 3]$  \\
    $\quad$ or $(y_{01},y_{02},y_{03},\xi)\in [-20,15]\times[-27,16]\times[2,46]\times[0,0.012]$ & activation function: Gaussian  \\
    domain decomposition: none & $\delta_m$: $0.2$ \\
    $r$: $0.0$  & $R_m$: to be specified  \\
    $Q$: $1000$ or $1200$, random & time: $t\in[0,t_f]$, $t_f=200$, $17$  \\
     $\Delta t$: $0.01$ (time-marching) &   \\
    \hline
  \end{tabular}
  \caption{NN simulation parameters for the Lorenz63 model
    (Section~\ref{sec_lorenz63}).
  }
  \label{tab_4}
\end{table}

\begin{figure}
  \centerline{
    \includegraphics[width=1.9in]{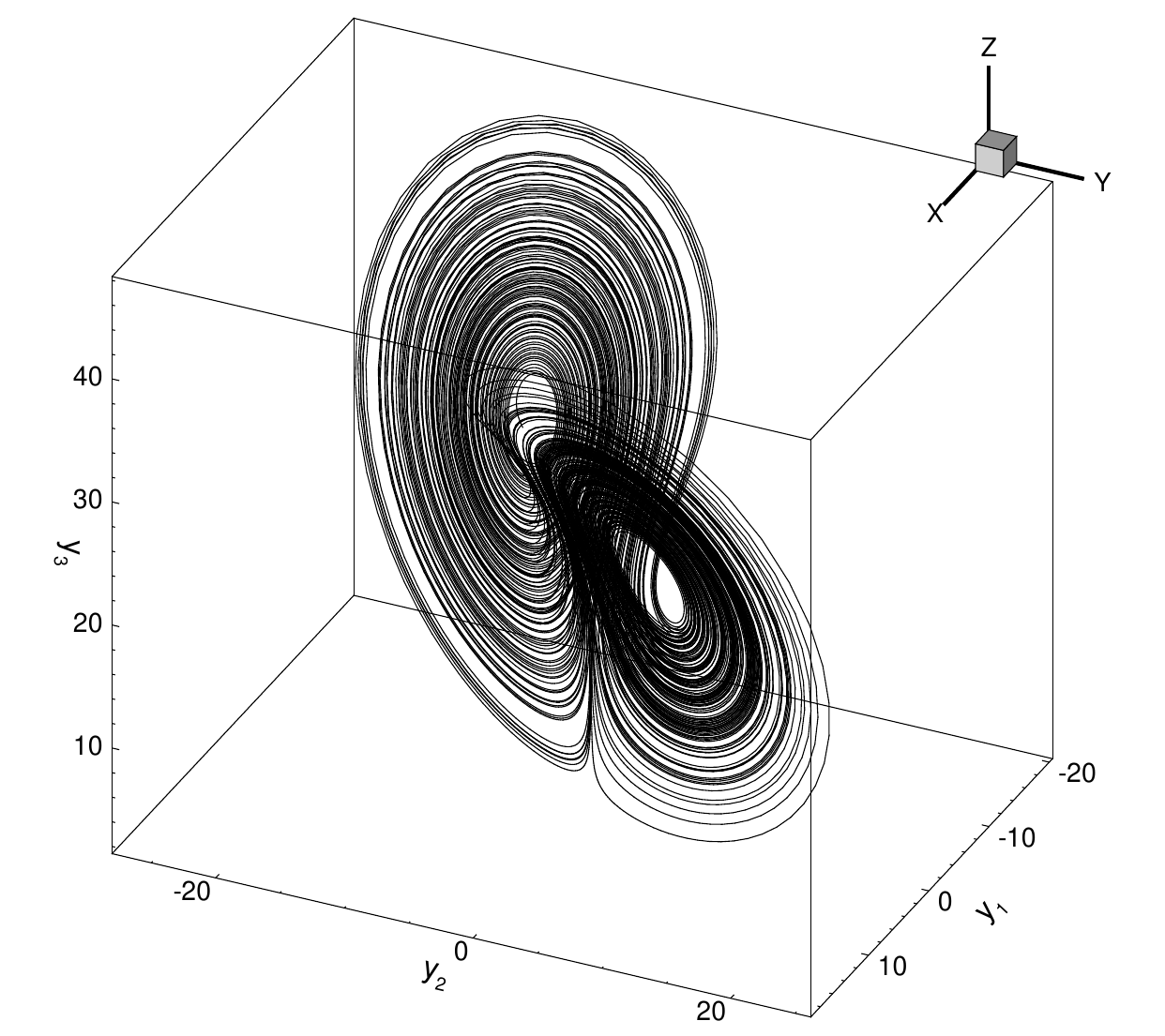}(a)
    \includegraphics[width=1.9in]{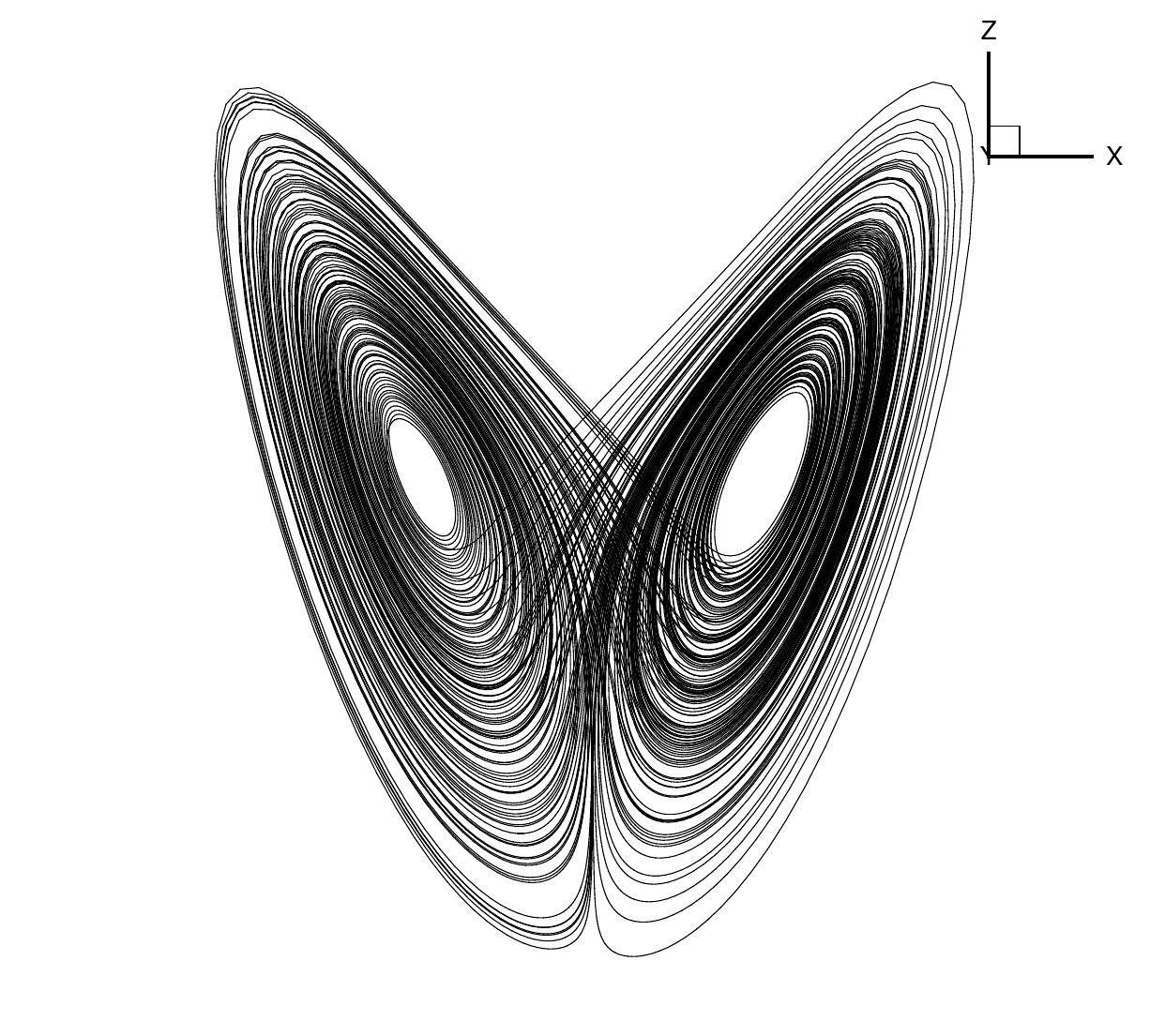}(b)
    \includegraphics[width=2.2in]{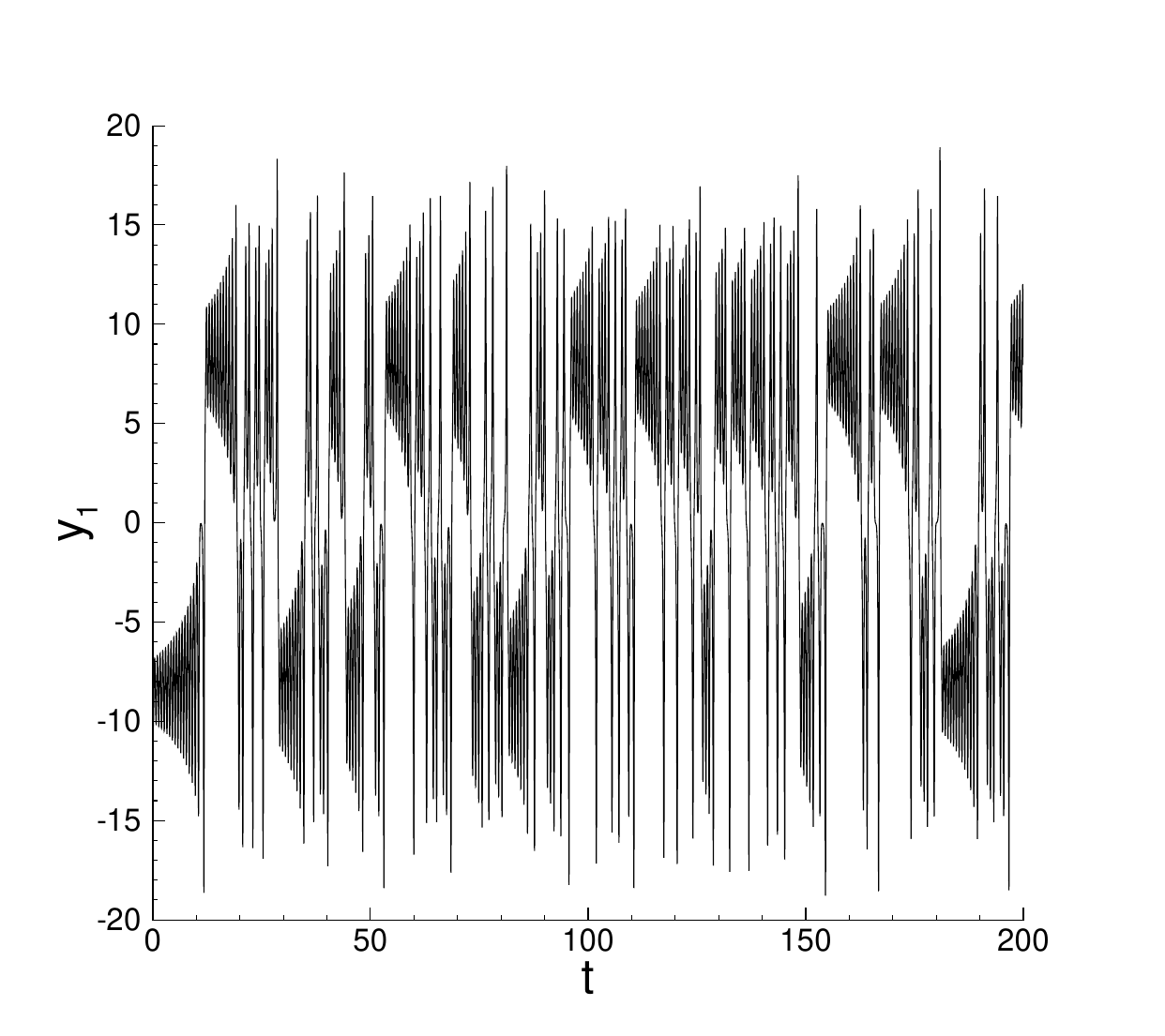}(c)
  }
  \caption{Lorenz63 system: (a) Perspective view, and (b) projection to the
    $y_1$-$y_3$ plane
    of the phase space trajectories obtained by NN-Exp-S0.
    (c) $y_1(t)$ history.
    NN-Exp-S0: training domain
    $(y_{01},y_{02},y_{03},\xi)\in[-20,20]\times[-25,25]\times[2,46]\times[0,0.012]$,
    NN: $[4,900,3]$, $Q=1200$, $R_m=0.12$, time integration for $t\in[0,200]$.
    See Table~\ref{tab_4} for the other parameter values.
  }
  \label{fg_21}
\end{figure}

In this subsection we test the NN algorithms using the Lorenz63 chaos model,
\begin{subequations}
\begin{align}
  & \frac{dy_1}{dt} = \sigma (y_2-y_1),\quad
   \frac{dy_2}{dt} = y_1(\rho-y_3) - y_2, \quad
   \frac{dy_3}{dt} = y_1y_2 - \beta y_3,\\
  & y_1(t_0)=y_{01}, \quad y_2(t_0) = y_{02}, \quad y_3(t_0) = y_{03},
\end{align}
\end{subequations}
where $y(t)=(y_1(t),y_2(t),y_3(t))$ are the unknowns,
$\sigma = 10$, $\rho=28$, $\beta=8/3$, and we employ the initial conditions
$y_0=(y_{01},y_{02},y_{03})=(-10,-10,25)$ with $t_0=0$ for time integration.

We learn the algorithmic function
$\psi(y_{01},y_{02},y_{03},\xi)$ using ELM with an
architecture $[4, M, 3]$ and Gaussian activation function
for the $\varphi$-subnet, where
$M$ is varied. Table~\ref{tab_4} provides values of the simulation
parameters related to the NN algorithms. No domain decomposition is used for
this problem.

Figure~\ref{fg_21} is an overview of the solution 
obtained by NN-Exp-S0. It shows the trajectory
in phase space with two different views, as well as the time history
of $y_1(t)$. The parameter values corresponding to these results
are listed in the figure caption or in Table~\ref{tab_4}.
The chaotic nature of the system is unmistakable.

\begin{figure}
  \centerline{
    \includegraphics[width=2in]{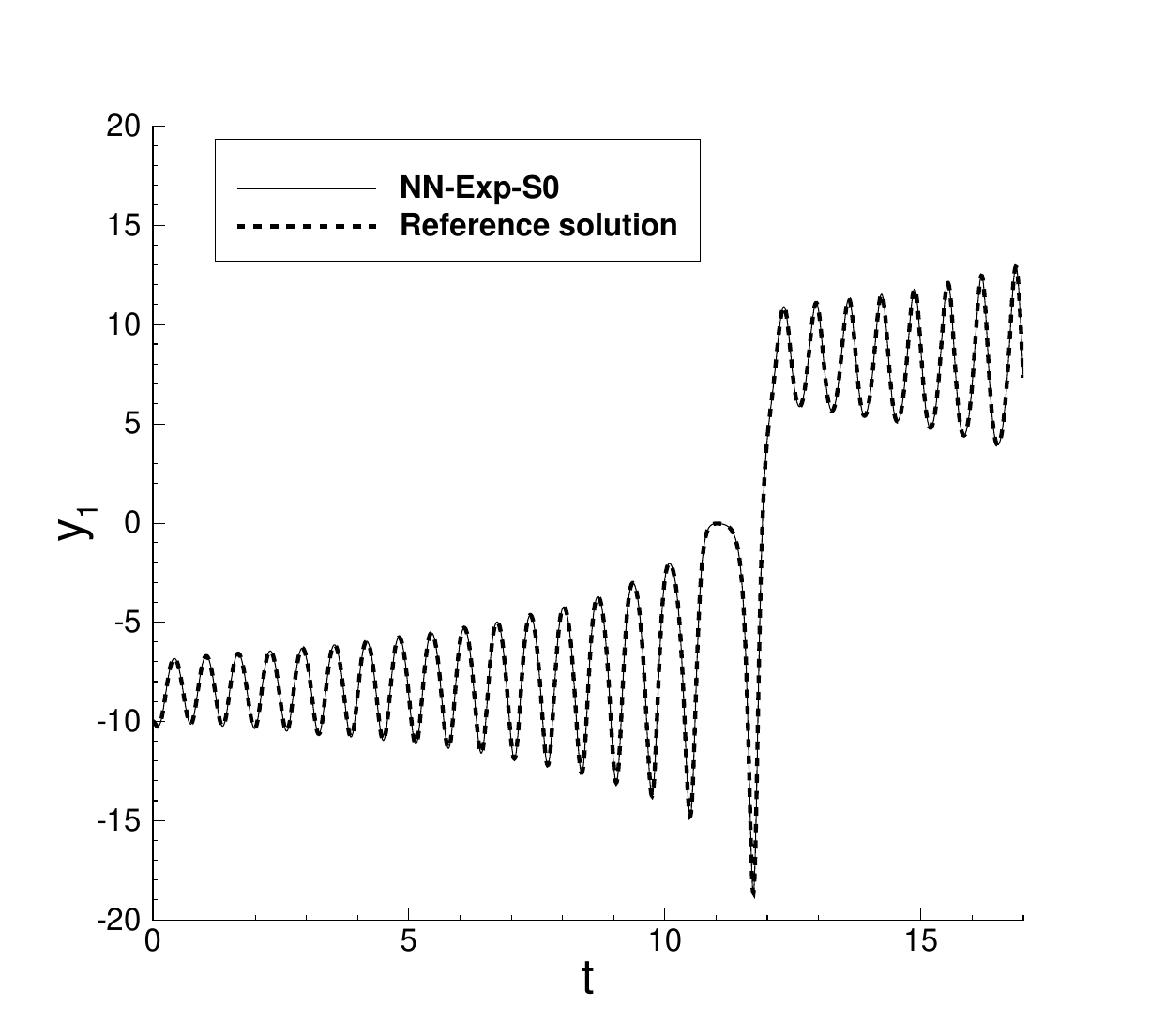}(a)
    \includegraphics[width=2in]{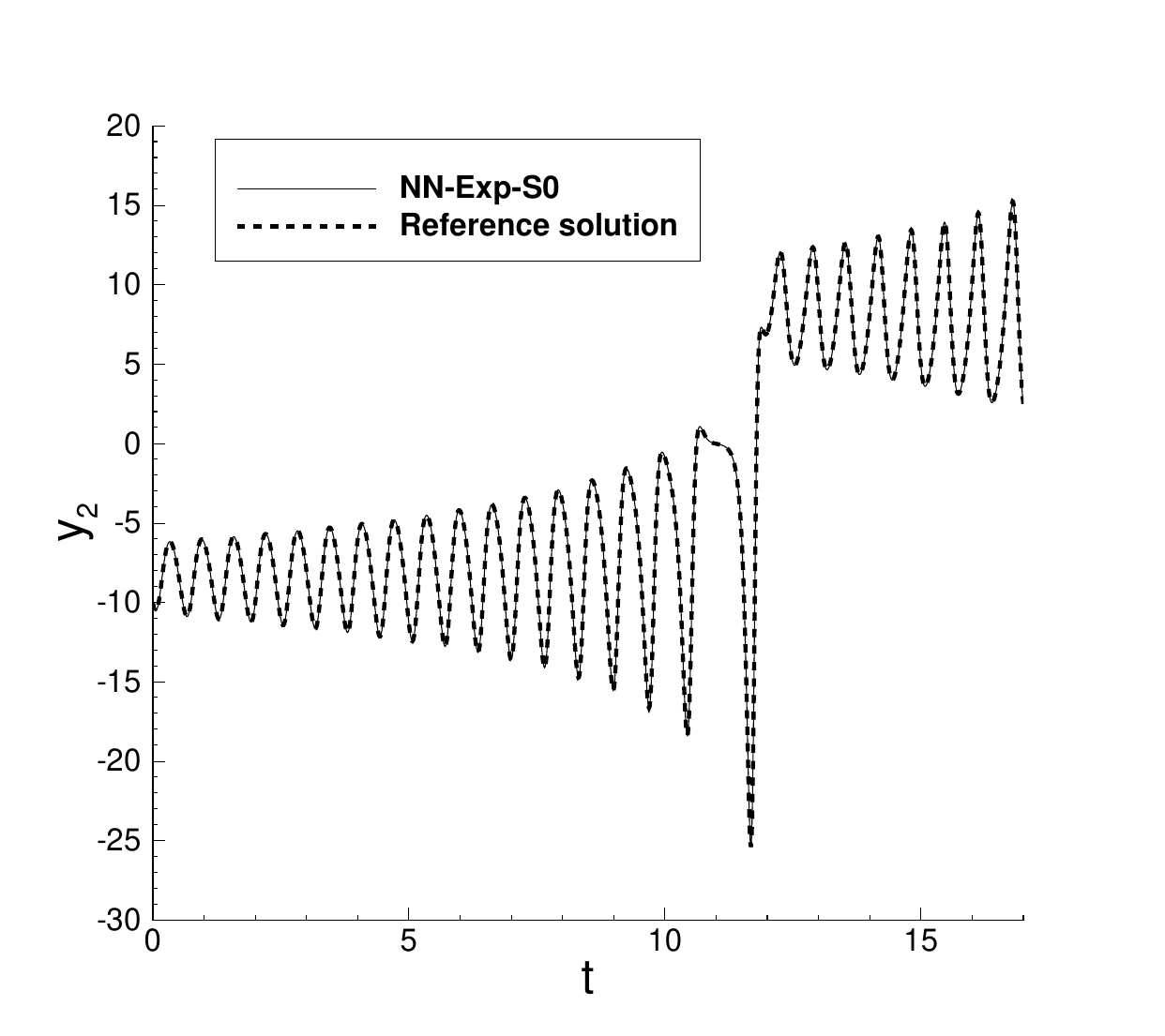}(b)
    \includegraphics[width=2in]{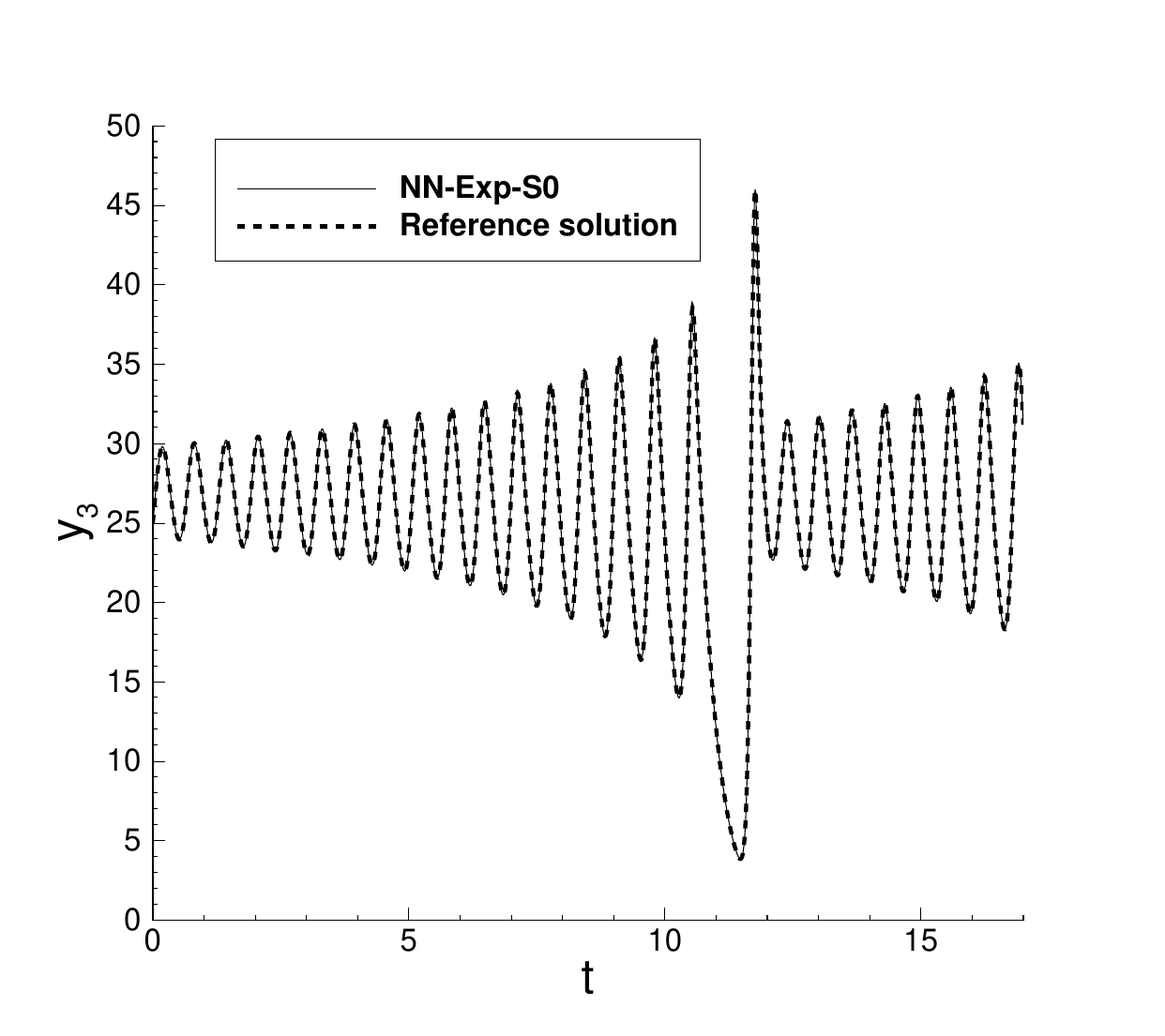}(c)
  }
  \centerline{
    \includegraphics[width=2in]{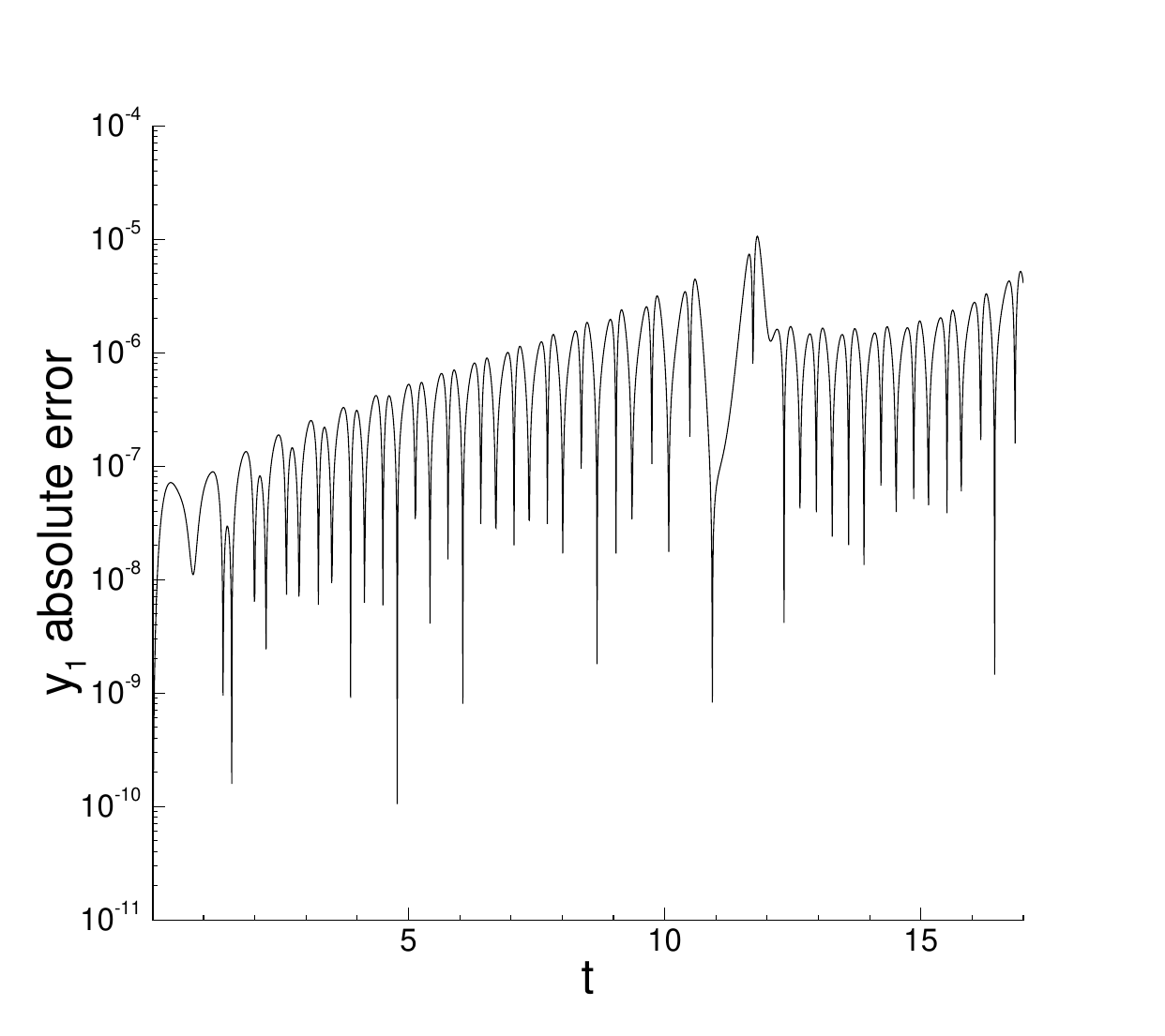}(d)
    \includegraphics[width=2in]{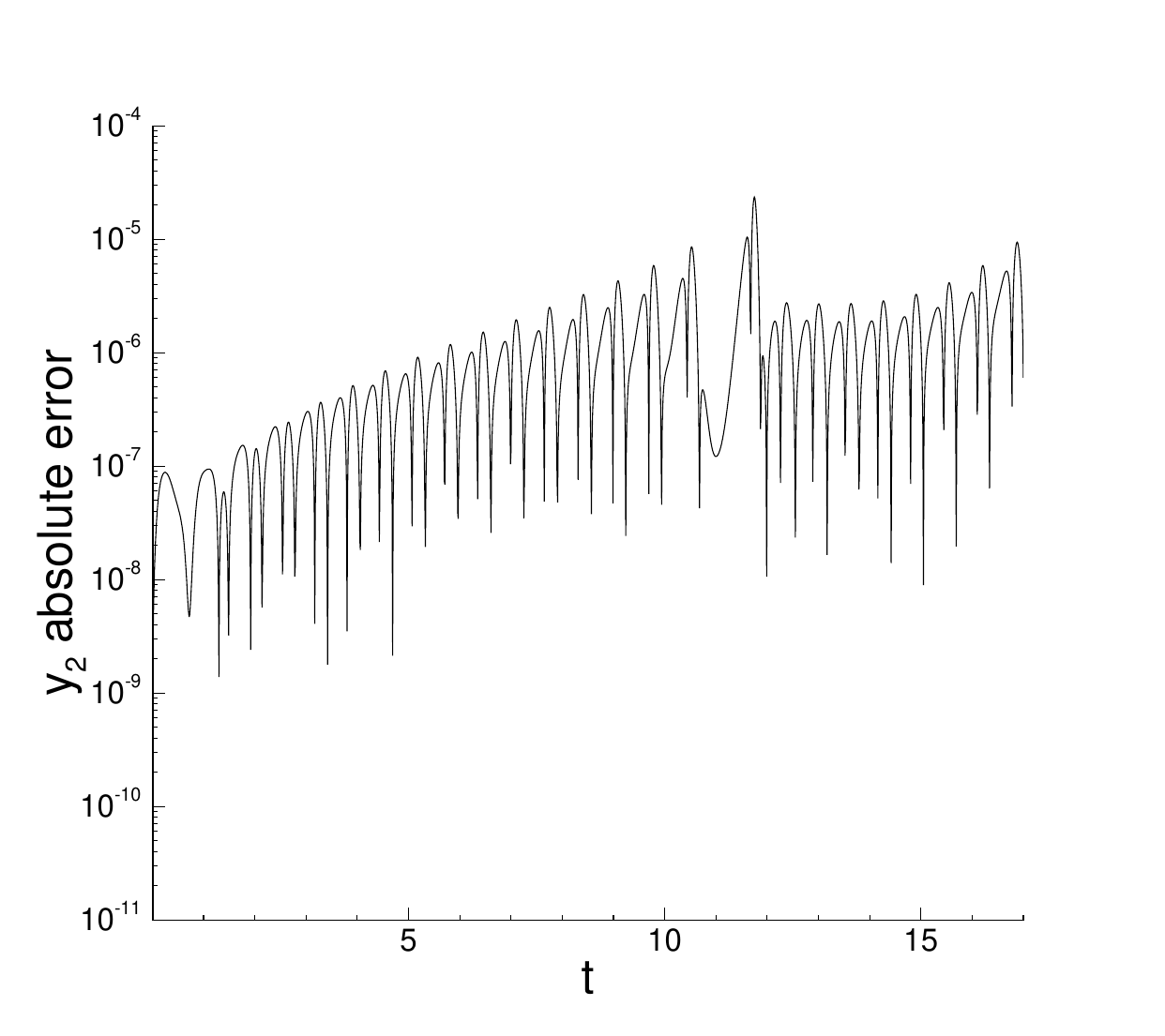}(e)
    \includegraphics[width=2in]{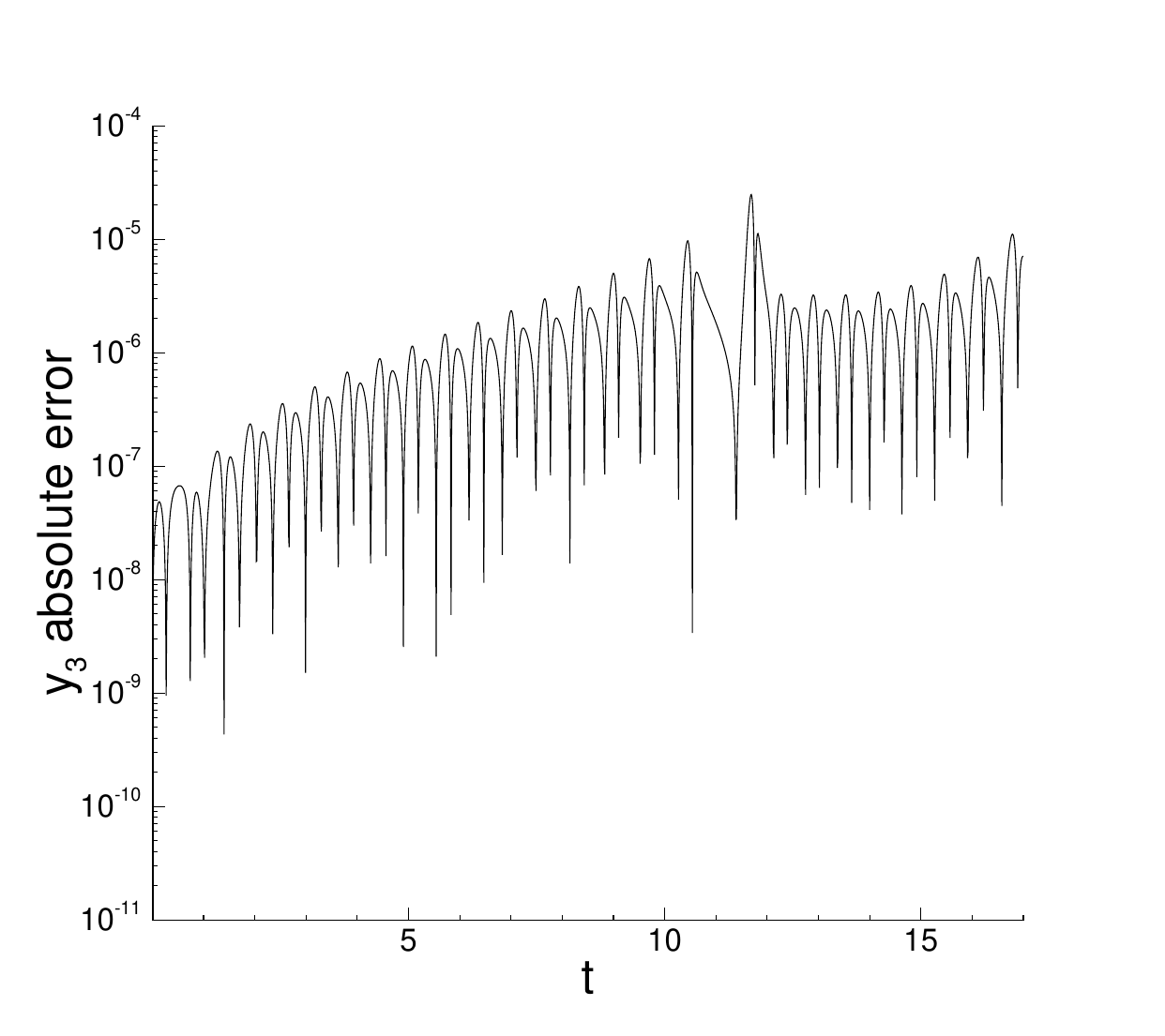}(f)
  }
  \caption{Lorenz63 system ($t_f=17$): Comparison between the NN-Exp-S0 solutions and
    the reference solutions (top row) for $y_1$, $y_2$ and $y_3$, and
    the absolute errors (bottom row) 
     of the NN-Exp-S0 solutions.
     Parameter values for NN-Exp-S0 follow those of Figure~\ref{fg_21},
     except here for $t\in[0,17]$.
     Reference solution is obtained by the scipy ``DOP853'' method
     with absolute tolerance $10^{-16}$
     and relative tolerance $10^{-13}$, with
     dense output on points corresponding to $\Delta t=0.01$.
  }
  \label{fg_22}
\end{figure}

Because the system is chaotic, comparison of different methods to compute their errors
for long time integration becomes impractical. However,
if the time horizon is not very long, comparing two solutions to compute
the error is still physically meaningful.
We next concentrate on a shorter time horizon $t\in[0,17]$ and compare
the NN algorithms with the scipy methods to evaluate their performance.

A comparison between the NN-Exp-S0 solution and a reference solution
obtained by the scipy DOP853 method is provided in Figure~\ref{fg_22}.
It shows the histories of $y_1$, $y_2$ and $y_3$ ($t\in[0,17]$)
obtained by NN-Exp-S0 and DOP853,  as well
as the absolute errors between these two solutions.
The parameter values corresponding to these results are
provided in the caption or in Table~\ref{tab_4}.
The NN solution agrees very well with the scipy solution,
with the maximum error on the order of $10^{-5}$ in the domain.

\begin{figure}
  \centerline{
    \includegraphics[width=2in]{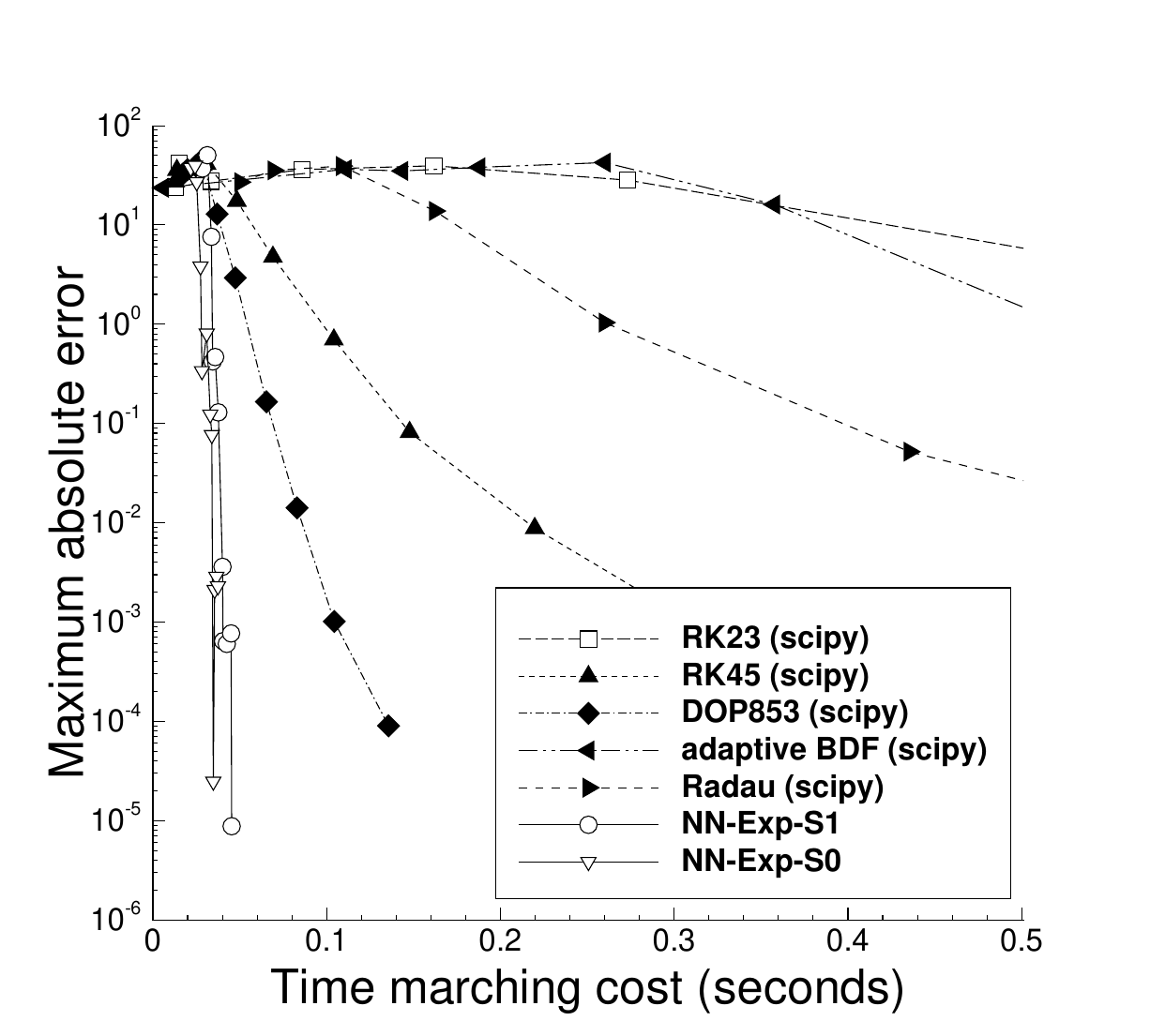}(a)
    \includegraphics[width=2in]{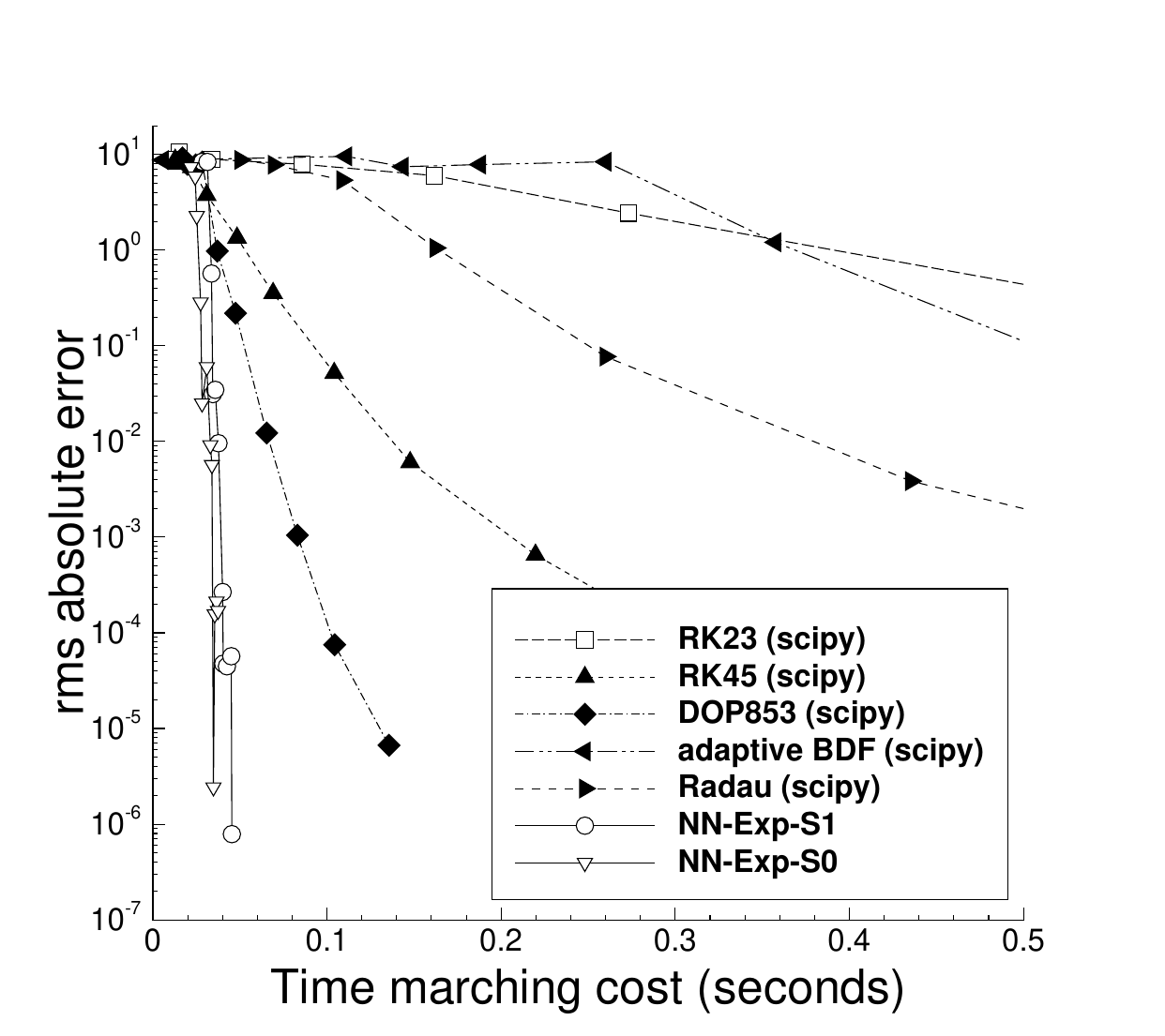}(b)
  }
  \caption{Lorenz63 system ($t_f=17$): Comparison of (a) the maximum and
    (b) the rms solution errors
    versus the time marching cost (wall time)
    between the NN algorithms (NN-Exp-S0, NN-Exp-S1) and the scipy methods.
    Parameter values for NN-Exp-S0 follow those of Figure~\ref{fg_22}.
    NN-Exp-S1: training domain $(y_{01},y_{02},y_{03},\xi)\in[-20,15]\times[-27,16]\times[2,46]\times[0,0.012]$,
    $R_m=0.1$, $Q=1000$; Other parameter values are given in Table~\ref{tab_4}.
    Scipy methods: 
    absolute tolerance $10^{-16}$, 
    relative tolerance varied,  dense output
    on points corresponding to $\Delta t=0.01$.
  }
  \label{fg_23}
\end{figure}

We compare the performance (accuracy versus cost) of the learned
NN algorithms and the scipy
methods  in Figure~\ref{fg_23}, which shows
the maximum and rms time-marching errors ($e_{\max}$, $r_{\text{rms}}$) obtained by
NN-Exp-S0/-S1 and the scipy methods as a function of the time-marching time
on $t\in[0,17]$. The figure caption and Table~\ref{tab_4} provide values
of the simulation parameters for these methods.
The data points for NN-Exp-S0 and NN-Exp-S1 correspond to different $M$
in the $\varphi$-subnet architecture $[4,M,3]$, and those data points for the scipy methods
correspond to different relative tolerance values.
The errors are computed against a reference solution attained by
the scipy DOP853 method with an absolute tolerance $10^{-16}$ and
a relative tolerance $10^{-13}$.
Among the scipy methods, DOP853 exhibits the best performance, followed by
RK45 and other methods. The NN-Exp-S0 and NN-Exp-S1 algorithms demonstrate
a similar performance, and both notably outperform the DOP853, RK45 and
the other scipy methods for this problem.

\subsection{Hindmarsh-Rose Neuronal Model}
\label{sec_hr}

\begin{table}[tb]
  \centering
  \begin{tabular}{l|l}
    \hline
    domain: $(y_{01},y_{02},y_{03},\xi)\in [-1.5,1.8]\times[-8,0.7]\times[2.7,3.3]\times[0,0.012]$,
    & NN ($\varphi$-subnet): $[4, M, 3]$  \\
    sub-domains: 4 along $y_{01}$ (non-uniform), sub-domain  & activation function: Gaussian  \\
    \quad\quad\quad boundaries: $[-1.5, -0.8, 0, 0.9, 1.8]$ & $\delta_m$: $1$ \\
    $r$: $0.0$  & $R_m$: to be specified  \\
    $Q$: $2000$, random & time: $t\in[0,t_f]$, $t_f=499.8$  \\
     $\Delta t$: $0.06$ (time-marching) &   \\
    \hline
  \end{tabular}
  \caption{NN simulation parameters for the Hindmarsh-Rose model
    (Section~\ref{sec_hr}).
  }
  \label{tab_5}
\end{table}

In the next example we test the NN algorithms using
the Hindmarsh-Rose model~\cite{HindmarshR1984},
which describes the spiking-bursting behavior of the membrane potential
in neurons. The model is given by,
\begin{subequations}
  \begin{align}
    & \frac{dy_1}{dt} = y_2 - y_1^3 + 3y_1 - y_3 + I, \quad
     \frac{dy_2}{dt} = 1-5y_1^2 - y_2, \quad
     \frac{dy_3}{dt} = 4\alpha\left(y_1 + \frac85 \right) - \alpha y_3, \\
    & y_1(t_0) = y_{01}, \quad y_2(t_0)=y_{02}, \quad y_3(t_0) = y_{03},
  \end{align}
\end{subequations}
where we employ $I=3.1$, $\alpha=0.006$, and the initial conditions
$(y_{01},y_{02},y_{03})=(-1,-3.5,3)$ with $t_0=0$.
Here $y_1$ represents the membrane potential, $y_2$ measures the
rate of ion transport,  and
$y_3$ is an adaption current that modulates the neuron firing rate.

Table~\ref{tab_5} lists the parameter values related to the
NN algorithms for learning $\psi(y_{01},y_{02},y_{03},\xi)$.
In particular, we employ $4$ non-uniform sub-domains
along the $y_{01}$ direction, with the sub-domain boundaries at
$y_{01}=-1.5$, $-0.8$, $0$, $0.9$, and $1.8$, with
an enlargement factor $r=0$ for network training.
The $\varphi$-subnet architecture is $[4,M,3]$ on each sub-domain, with
$M$ varied in the tests.

\begin{figure}
  \centerline{
    \includegraphics[width=2in]{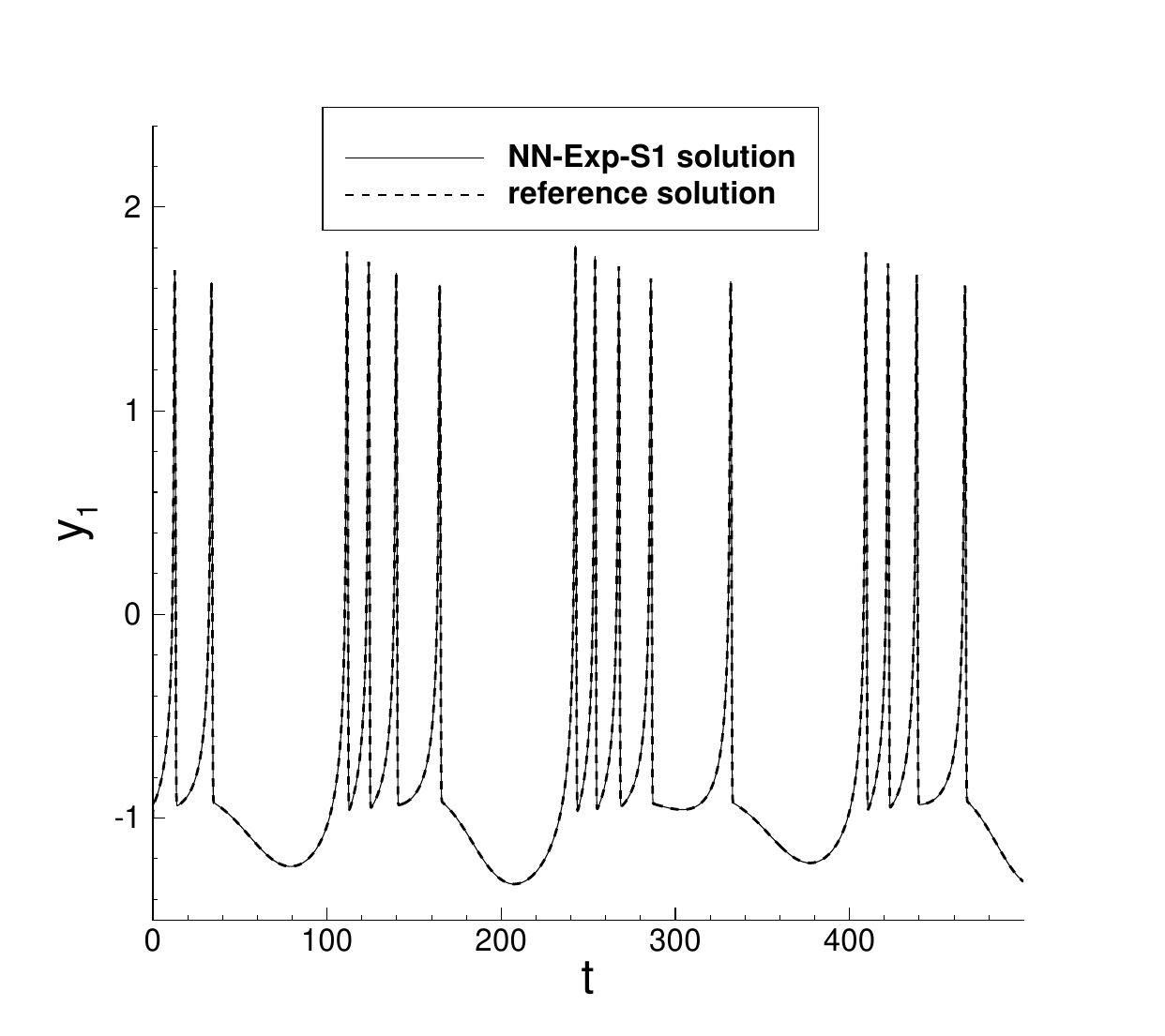}(a)
    \includegraphics[width=2in]{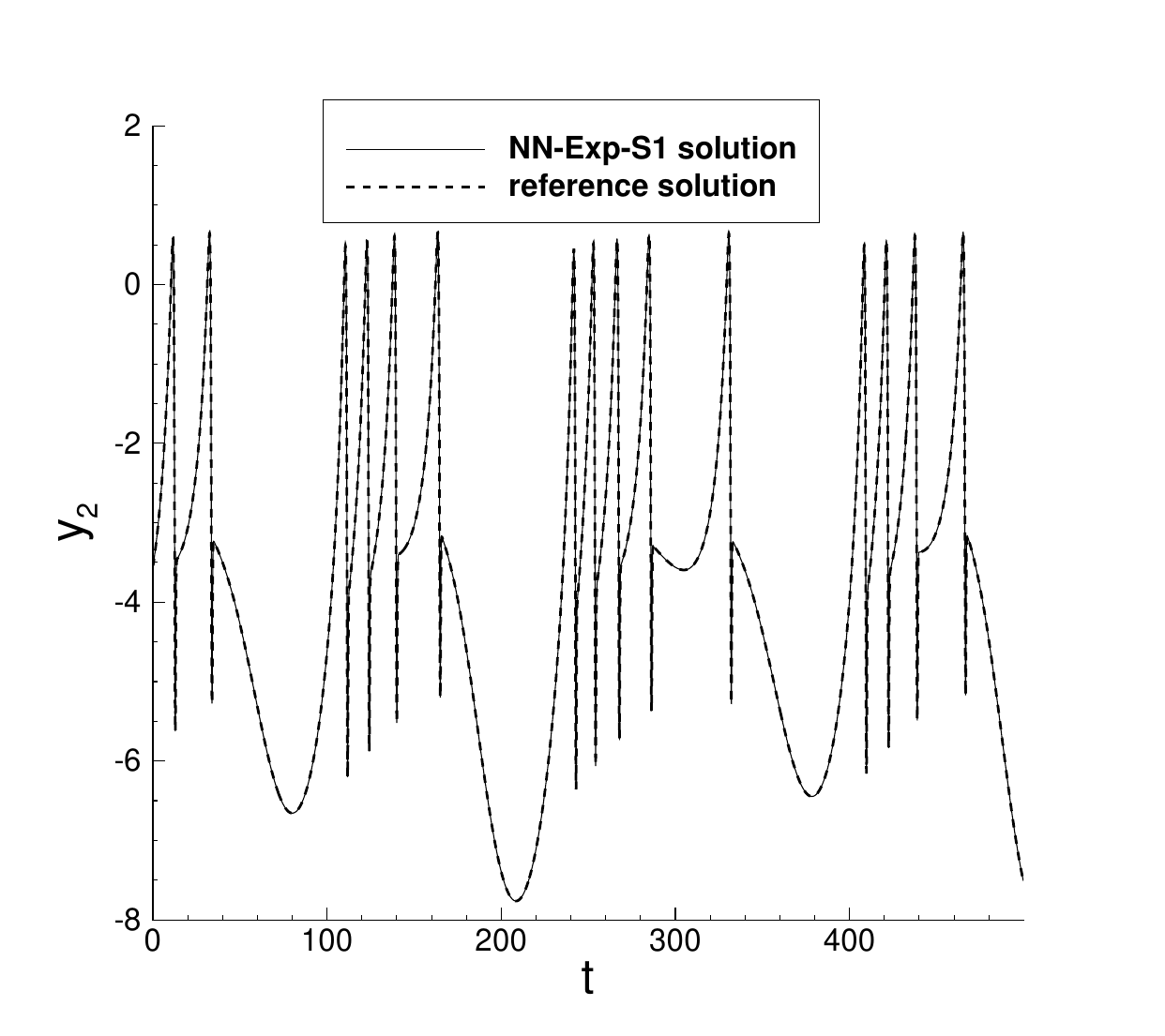}(b)
    \includegraphics[width=2in]{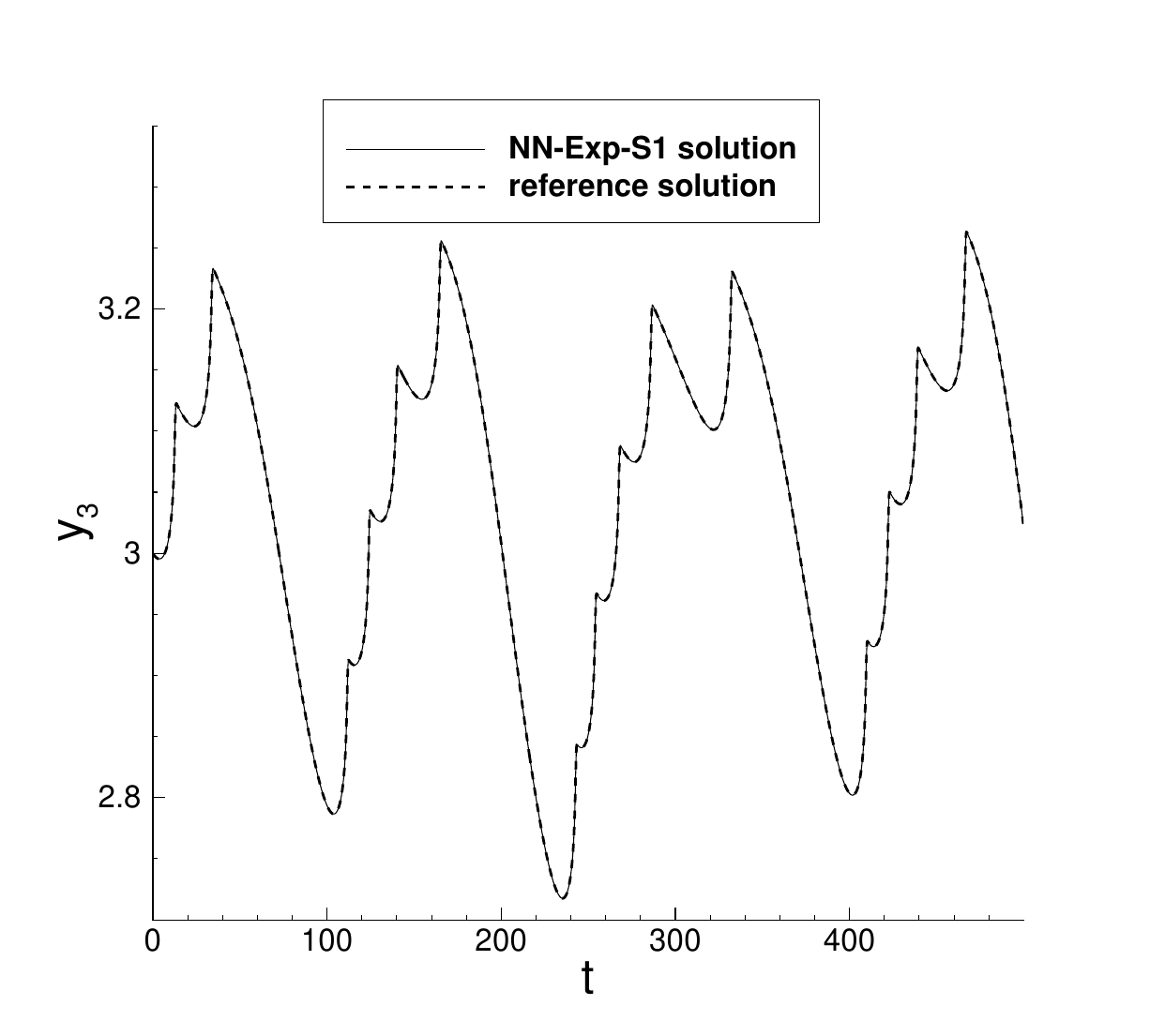}(c)
  }
  \centerline{
    \includegraphics[width=2in]{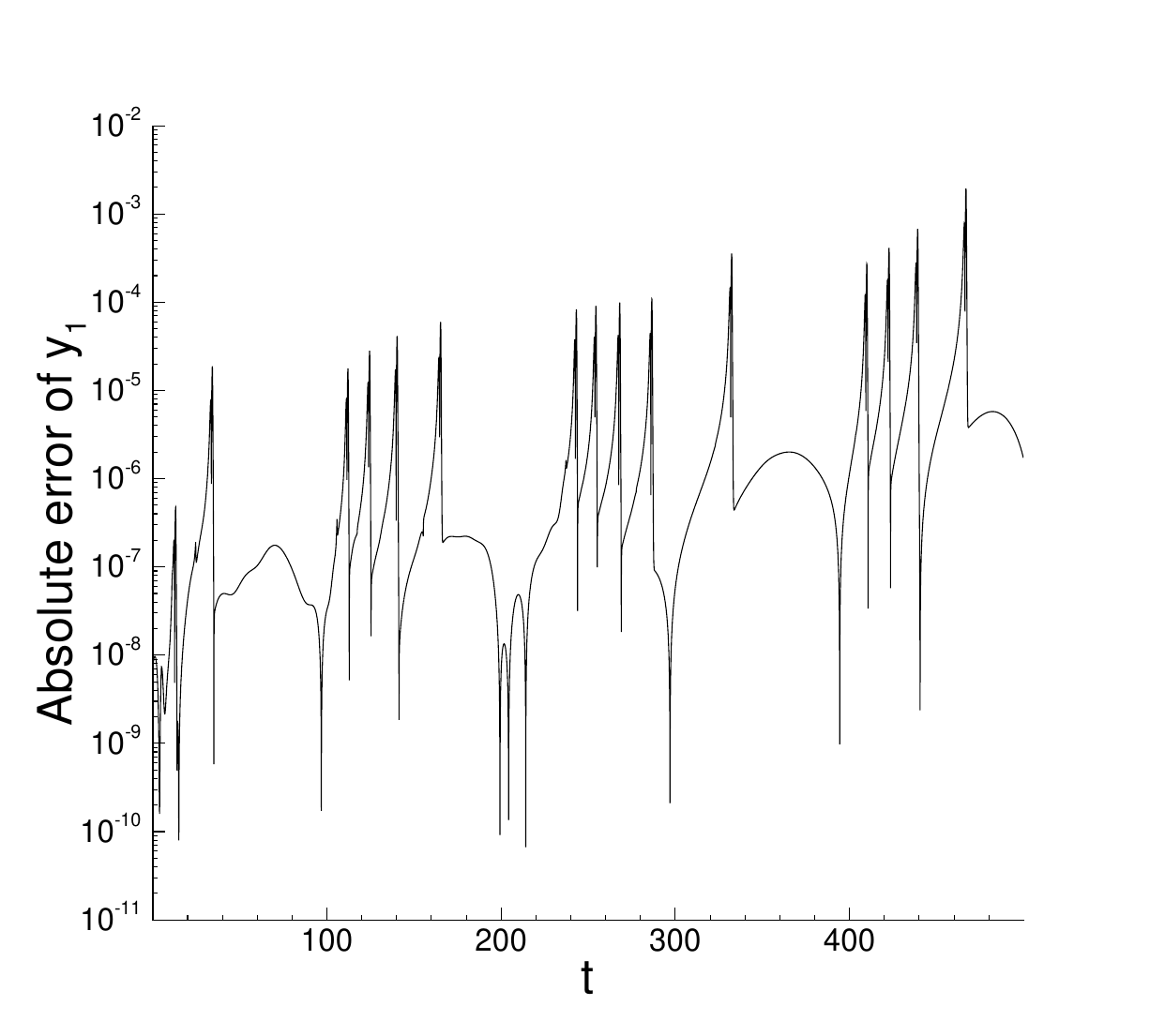}(d)
    \includegraphics[width=2in]{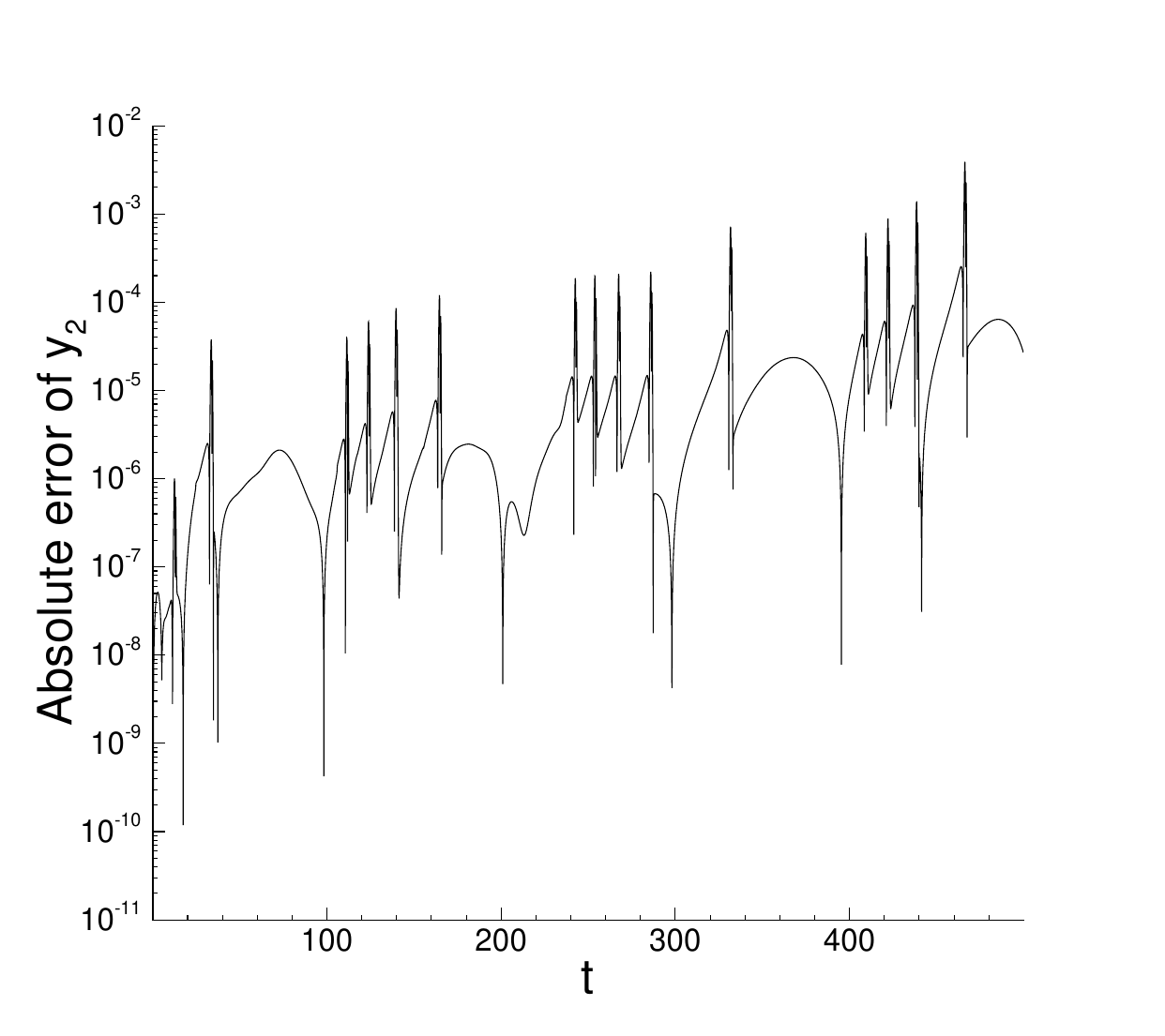}(e)
    \includegraphics[width=2in]{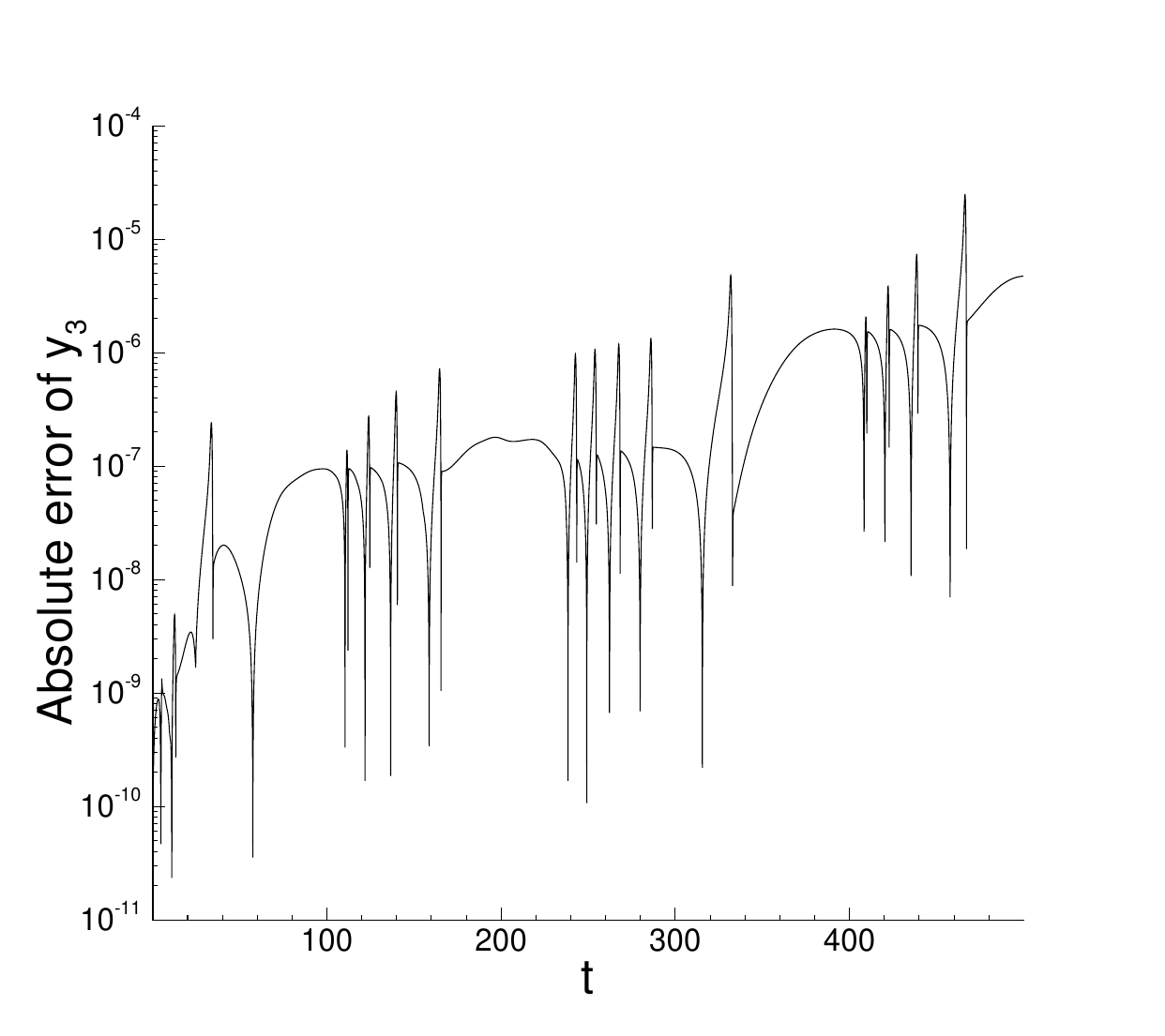}(f)
  }
  \caption{Hindmarsh-Rose neuron model:
    Comparison of  (a) $y_1(t)$, (b) $y_2(t)$, and (c) $y_3(t)$ between the NN-Exp-S1 solution
    and the reference solution. The absolute error for (d) $y_1(t)$, (e) $y_2(t)$,
    and (f) $y_3(t)$ of the NN-Exp-S1 solution.
    NN-Exp-S1:
    architecture [4,1200,3] on each sub-domain,
    $R_m=0.39$; Other simulation parameters are given in Table~\ref{tab_5}.
    Reference solution: computed by scipy DOP853, with absolute tolerance $10^{-16}$
    and relative tolerance $10^{-13}$, dense output on points corresponding to $\Delta t=0.06$.
  }
  \label{fg_24}
\end{figure}

Figure~\ref{fg_24} illustrates the solution characteristics for $y_1(t)$, $y_2(t)$
and $y_3(t)$ obtained by NN-Exp-S1. 
A reference solution for these variables obtained by the scipy DOP853 method
is also shown in Figures~\ref{fg_24}(a,b,c) for comparison. The
absolute errors  are shown in the bottom row of the plots.
The  parameter values for these results
are provided in the figure caption or in Table~\ref{tab_5}.
The bursting behavior of the solution is evident from the history plots.
The NN algorithm has captured the dynamics of the system accurately.

\begin{figure}
  \centerline{
    \includegraphics[width=2in]{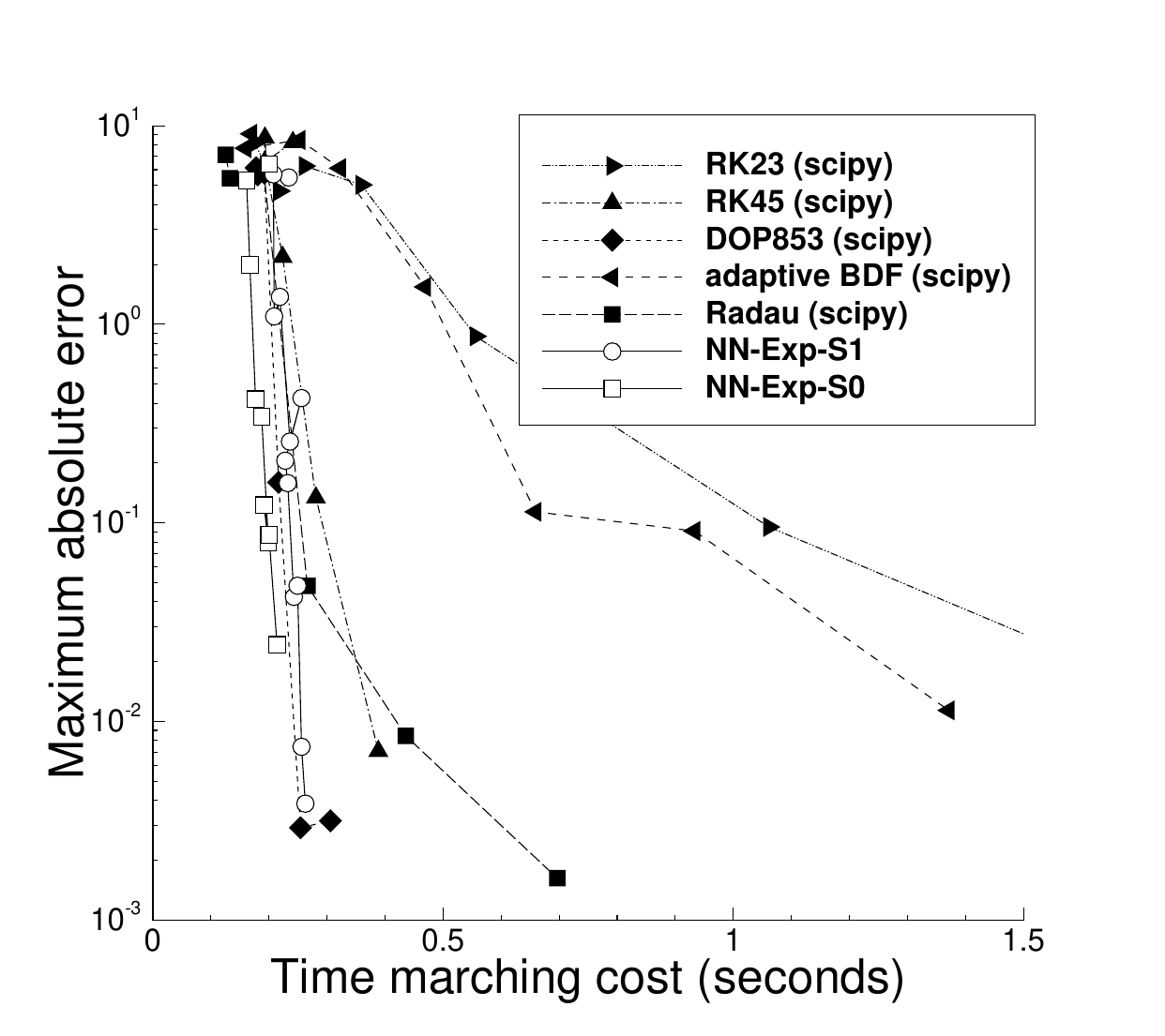}(a)
    \includegraphics[width=2in]{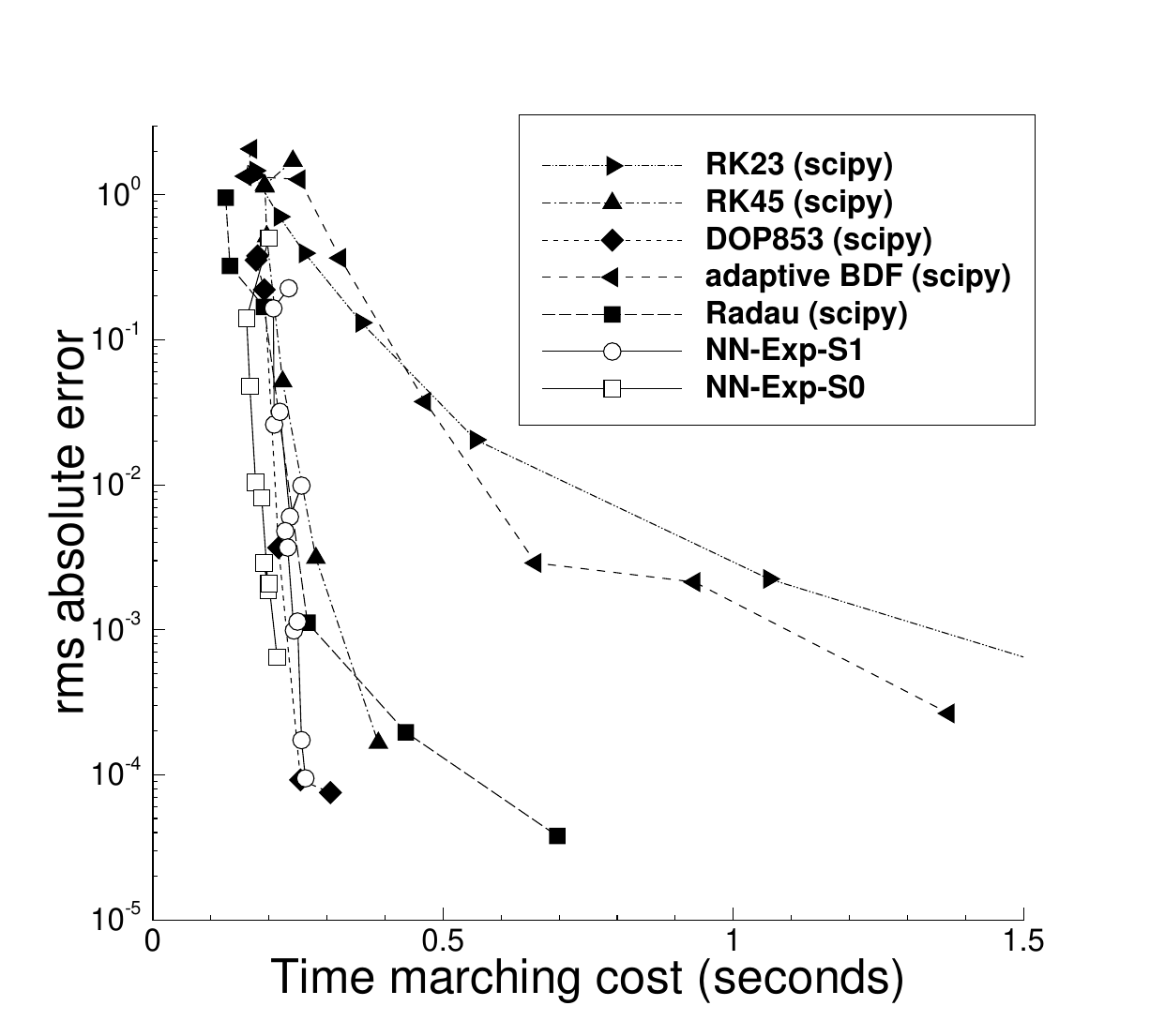}(b)
  }
  \caption{Hindmarsh-Rose neuron model:
    Comparison of (a) the maximum (b) the rms time-marching errors
    versus the time marching
    cost (wall time) between the NN algorithms (NN-Exp-S0, NN-Exp-S1)
    and the scipy methods.
    $R_m=0.12$ for NN-Exp-S0, $R_m=0.39$ for NN-Exp-S1; See Table~\ref{tab_5} for
    the other parameter values of the NN algorithms.
    Scipy methods:
    absolute tolerance $10^{-16}$,
    data points corresponding to different relative tolerance values.
  }
  \label{fg_25}
\end{figure}

In Figure~\ref{fg_25} we compare the computational performance (accuracy versus cost)
between the NN algorithms and the scipy methods for the Hindmarsh-Rose model.
It shows the maximum and rms time-marching errors
as a function of the time-marching cost (for $t\in[0,499.8]$) obtained
by the NN-Exp-S0/NN-Exp-S1 algorithms and the scipy methods.
The performance of NN-Exp-S0, NN-Exp-S1 and DOP853 is close to one another,
with NN-Exp-S0 being slightly better.
These three methods show a better performance than the Radau and RK45 methods.
All these methods perform significantly better than RK23 and BDF.

\subsection{Lorenz96 Chaotic System}
\label{sec_lorenz96}

\begin{table}[tb]
  \centering
  \begin{tabular}{l|l}
    \hline
    domain: $(y_{01},y_{02},y_{03},y_{04},y_{05},\xi)\in [-5,10]^5\times[0,0.011]$,
    & NN ($\varphi$-subnet): $[6, M, 5]$ \\
    \quad or $[-9,12]\times[-7,10]\times[-5.5,11]\times[-6,7]\times[-4,12]\times[0,0.011]$ &
    activation function: Gaussian\\
    sub-domains: 2 or 4 along $y_{01}$, uniform  & $\delta_m$: $1$  \\
    $r$: $0.0$ & $R_m$: to be specified   \\
    $Q$: $1500$ or $2000$, random  & time: $t\in[0,t_f]$, $t_f=5$ or $50$ \\
     $\Delta t$: $0.01$ (time-marching) &   \\
    \hline
  \end{tabular}
  \caption{NN simulation parameters for the Lorenz96 model
    (Section~\ref{sec_lorenz96}).
  }
  \label{tab_6}
\end{table}

In the last example we employ the Lorenz96 chaotic system~\cite{Lorenz1996} to
evaluate the performance of the  NN  algorithms.
This model is given by (in $5$ dimensions),
\begin{subequations}
  \begin{align}
    & \frac{dy_1}{dt} = (y_2 - y_4)y_5 - y_1 + F, \quad
     \frac{dy_2}{dt} = (y_3 - y_5)y_1 - y_2 + F, \quad
     \frac{dy_3}{dt} = (y_4 - y_1)y_2 - y_3 + F, \\
    & \frac{dy_4}{dt} = (y_5 - y_2)y_3 - y_4 + F, \quad
     \frac{dy_5}{dt} = (y_1 - y_3)y_4 - y_5 + F, \\
    & y_1(t_0)=y_{01},\ y_2(t_0)=y_{02},\ y_3(t_0)=y_{03},\ y_4(t_0)=y_{04},\ y_5(t_0)=y_{05},
  \end{align}
\end{subequations}
where $F=8$, and the initial conditions are
$y_0=(y_{01},y_{02},y_{03},y_{04},y_{05})=(-0.99,-1,-1,-1,-1)$ with $t_0=0$.

We learn the algorithmic function
function $\psi(y_{01},y_{02},y_{03},y_{04},y_{05},\xi)$ using
an ELM network with the architecture $[6,M,5]$ and the Gaussian activation
function for the $\varphi$-subnet, where $M$ is varied. The simulation parameters
for the NN algorithms are summarized in Table~\ref{tab_6}.

\begin{figure}
  \centerline{
    \includegraphics[height=2.2in]{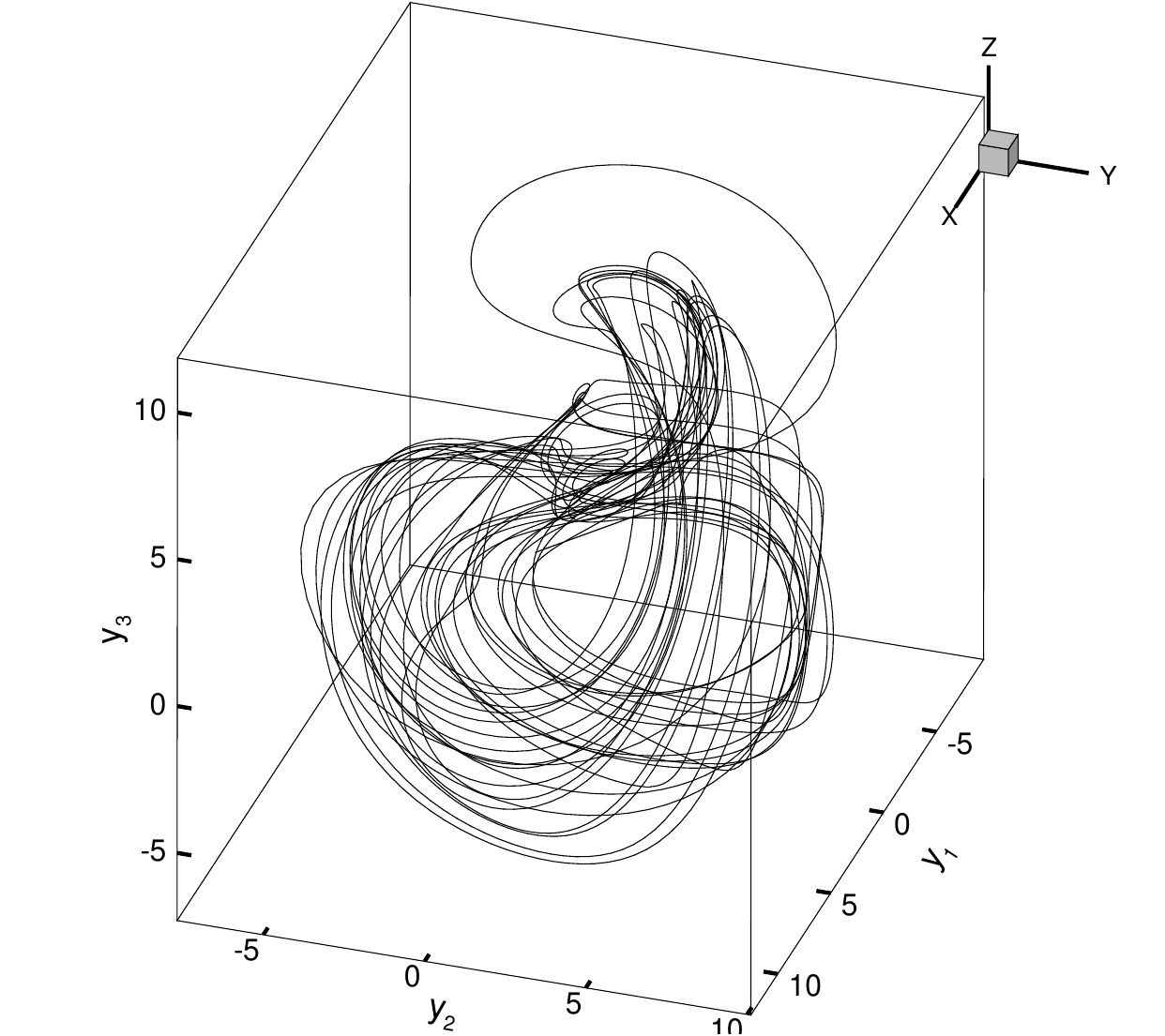}(a)
    \includegraphics[height=2.5in]{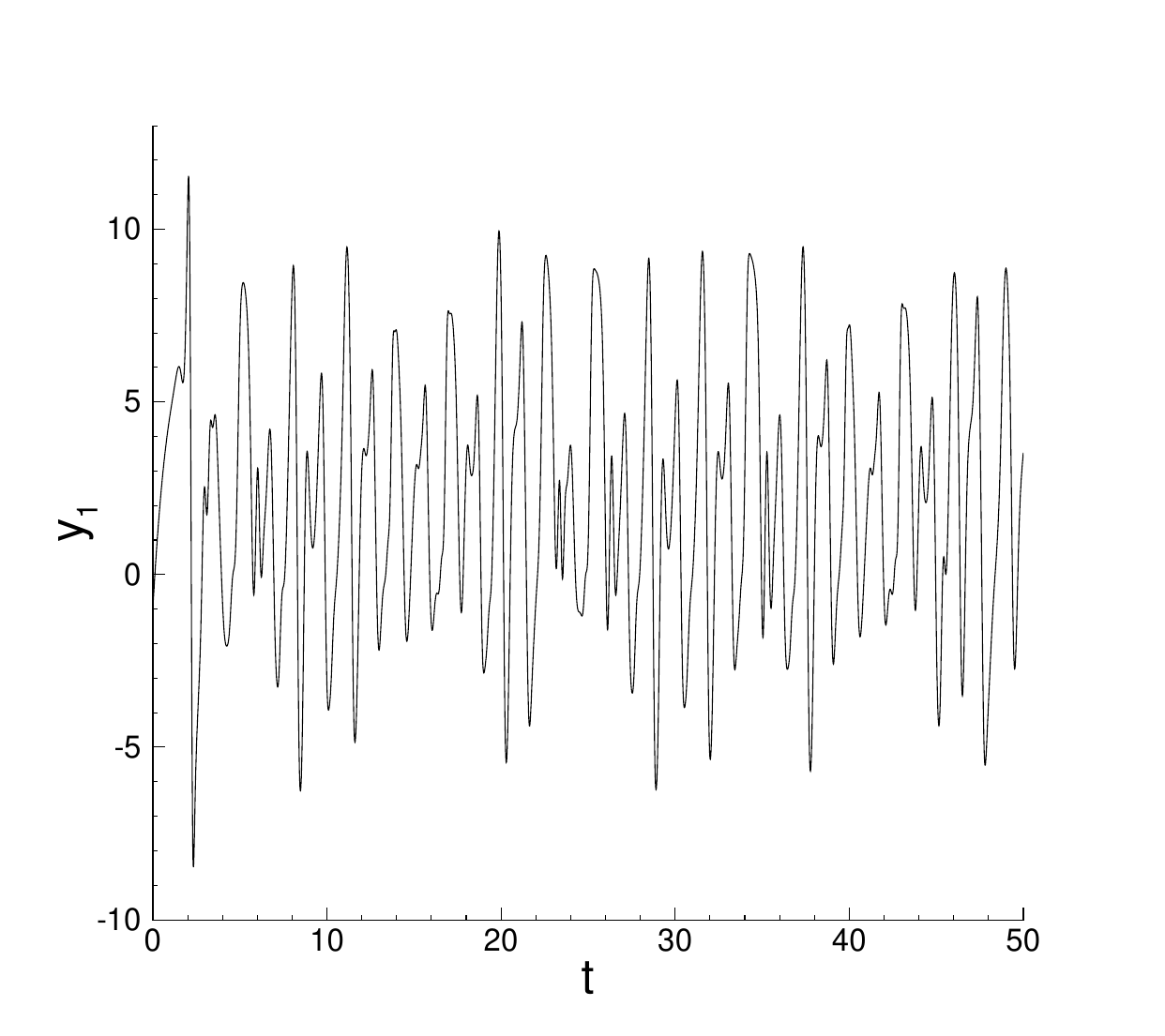}(b)
  }
  \caption{Lorenz96 system: (a) Trajectory in the $(y_1,y_2,y_3)$-phase space, and (b)
    the $y_1$ history, obtained by the NN-Exp-S1 algorithm.
    Training domain: $(y_{01},y_{02},y_{03},y_{04},y_{05},\xi)\in[-5,10]^5\times[0,0.011]$,
    with 2 uniform sub-domains along $y_{01}$;
    NN: $[6,800,5]$ on each sub-domain;
    $R_m=0.15$, $Q=1500$,  time integration for $t\in[0,50]$;
    The Other parameter values
    are given in Table~\ref{tab_6}.
  }
  \label{fg_26}
\end{figure}

Figure~\ref{fg_26} provides an overview of the solution characteristics.
It shows the system trajectory  ($t\in[0,50]$) in
the $(y_1,y_2,y_3)$ phase space and the time history of $y_1$ obtained by
the NN-Exp-S1 algorithm. The figure caption and Table~\ref{tab_6} list
the parameter values for this set of tests.

\begin{figure}
  \centerline{
    \includegraphics[width=1.8in]{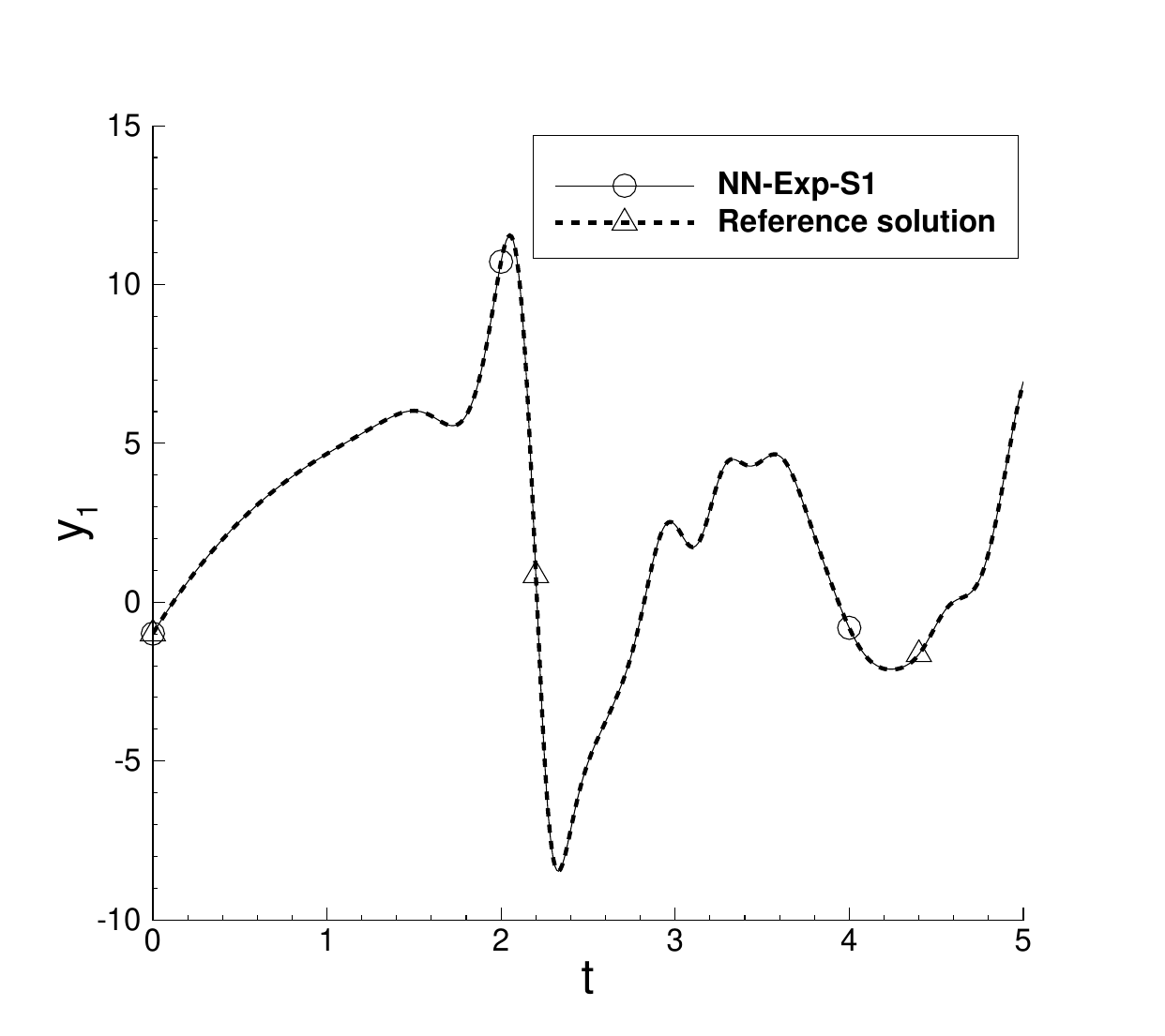}(a)
    \includegraphics[width=1.8in]{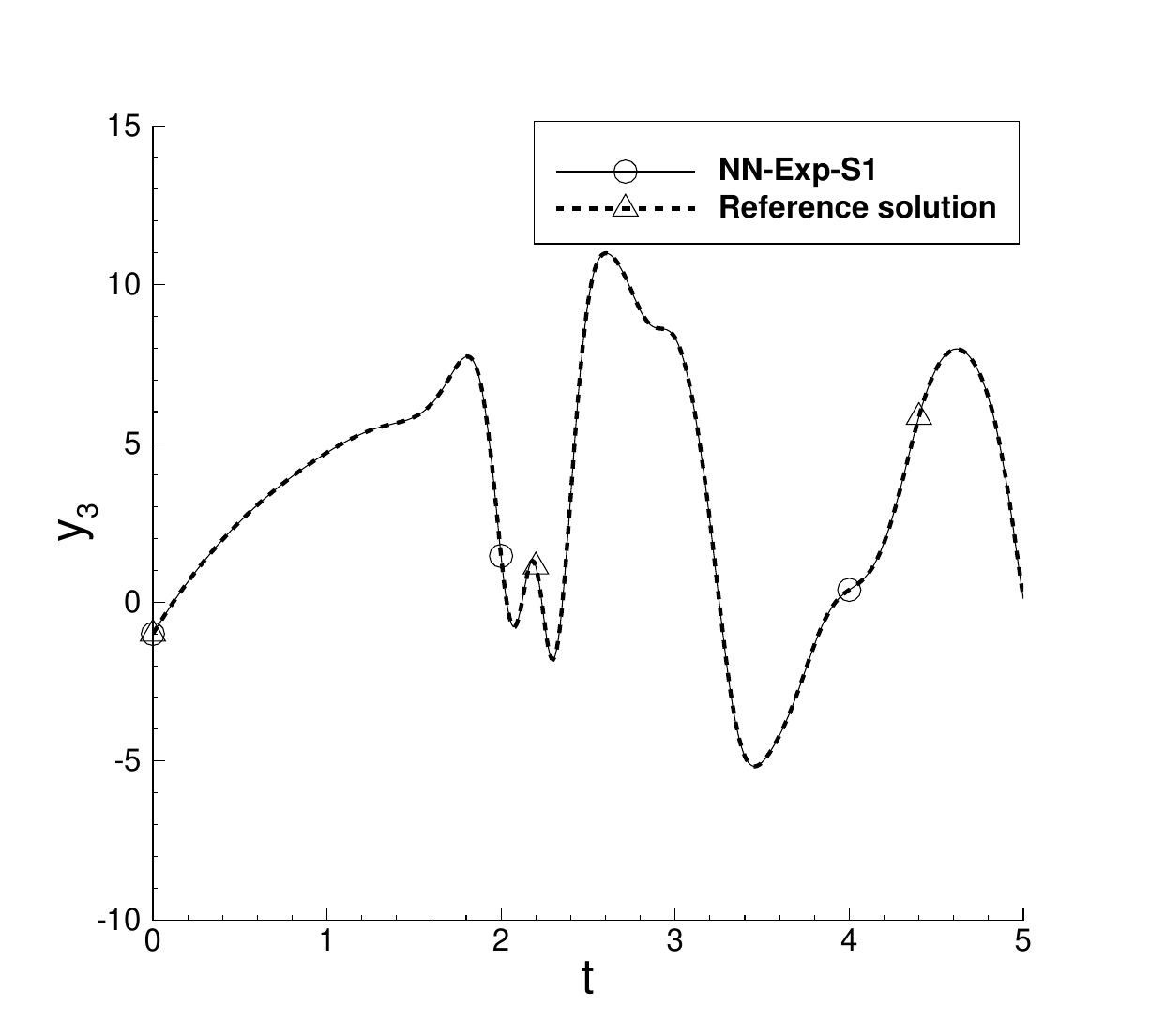}(b)
    \includegraphics[width=1.8in]{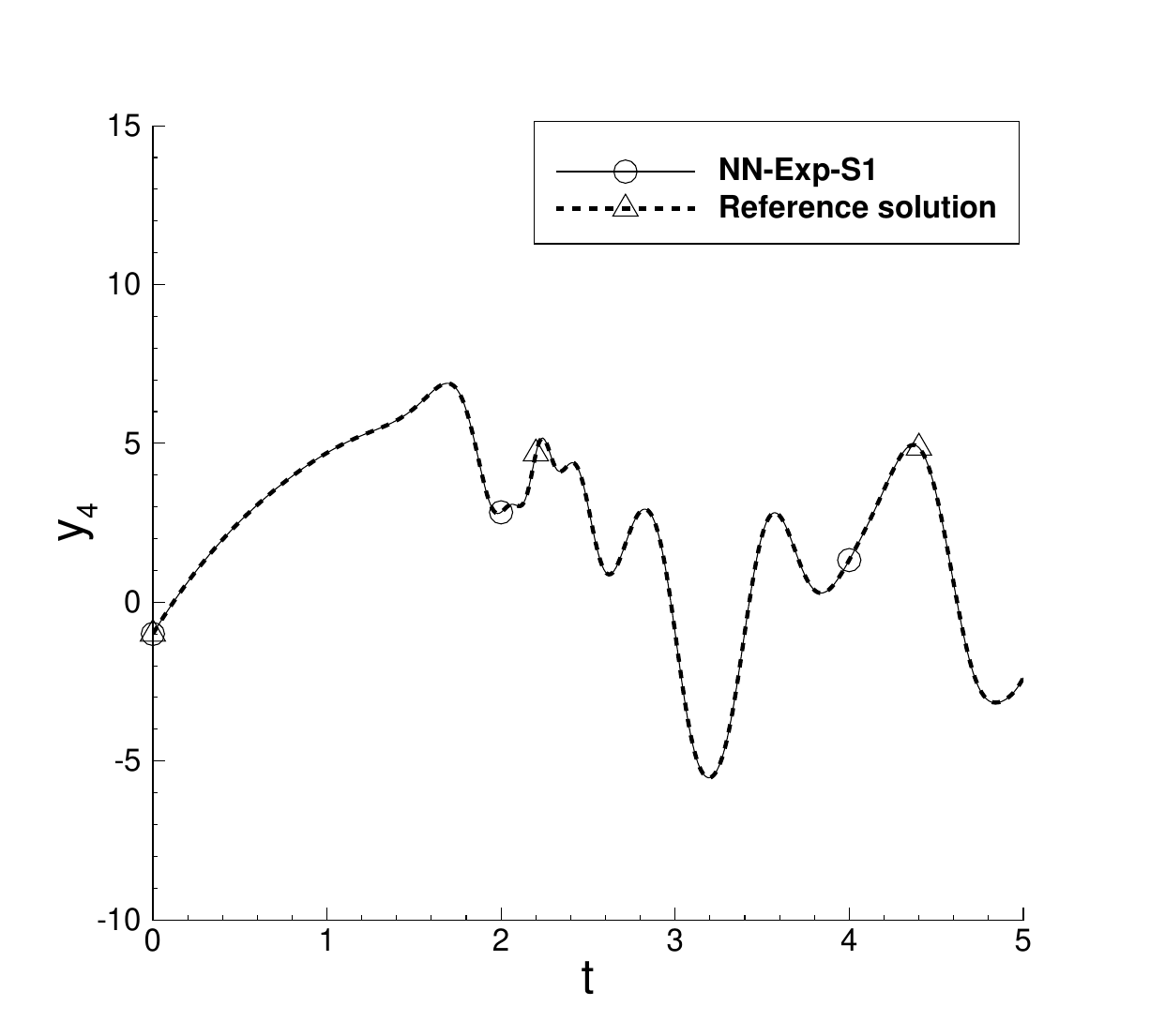}(c)
  }
  \centerline{
    \includegraphics[width=1.8in]{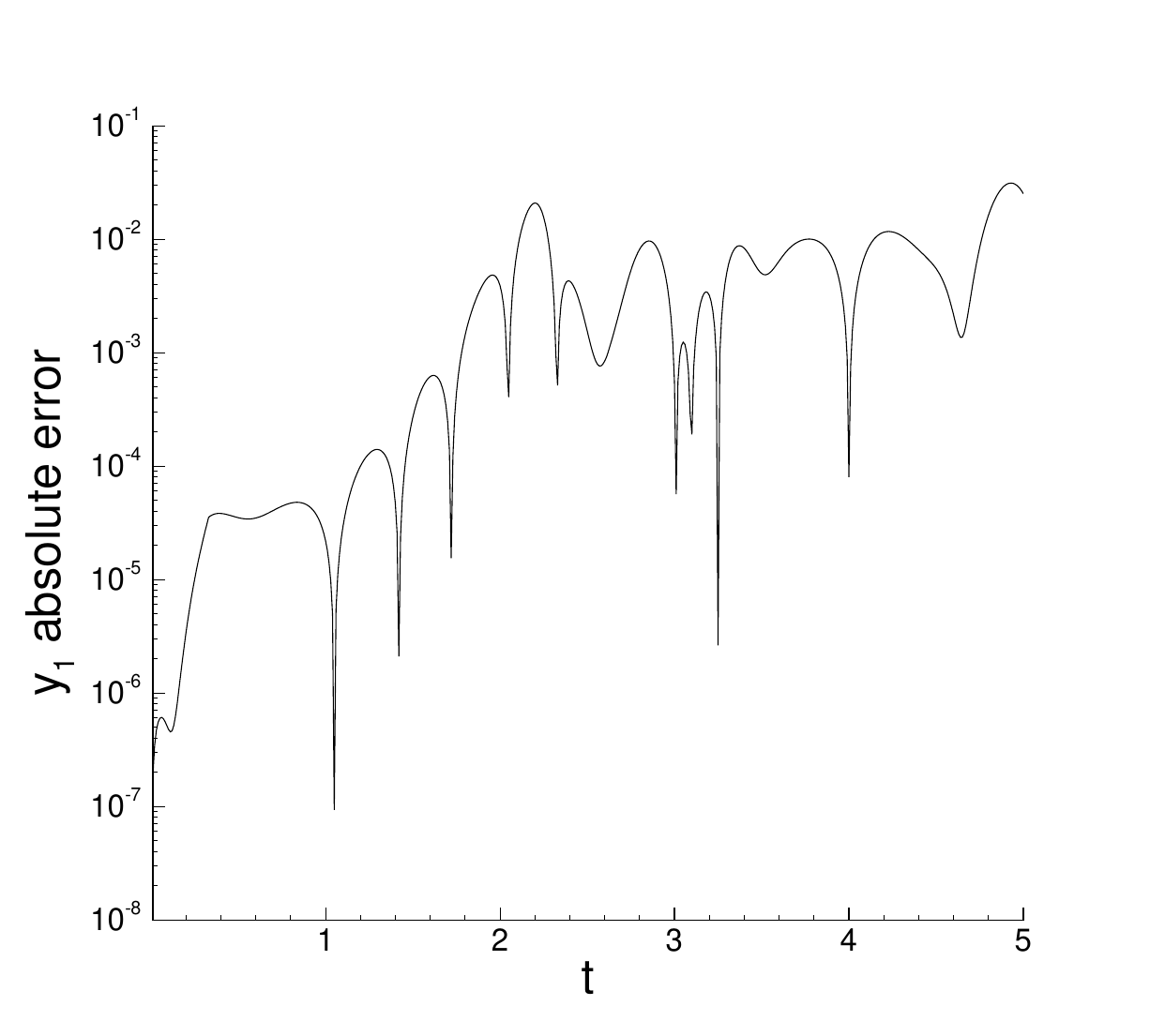}(d)
    \includegraphics[width=1.8in]{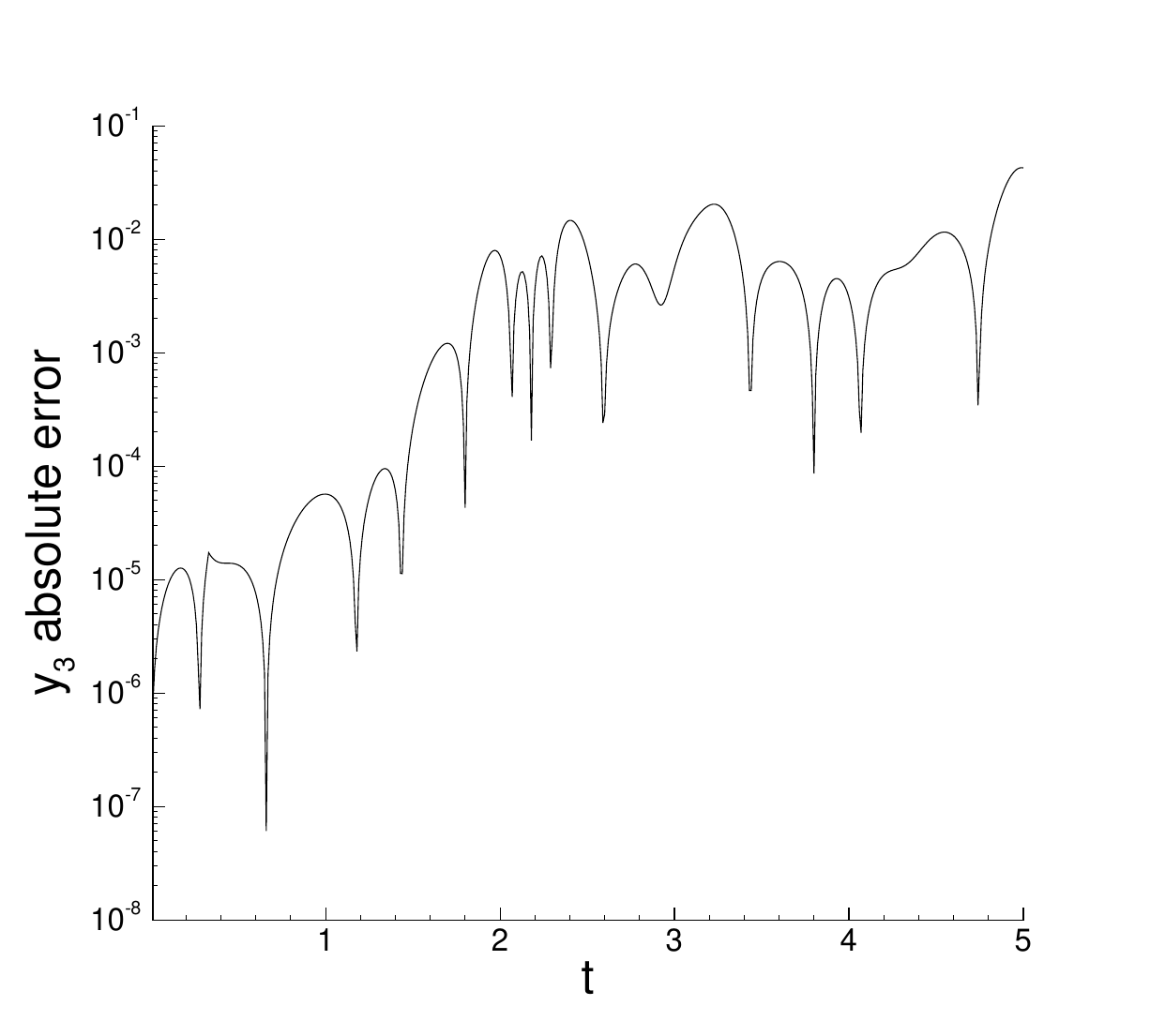}(e)
    \includegraphics[width=1.8in]{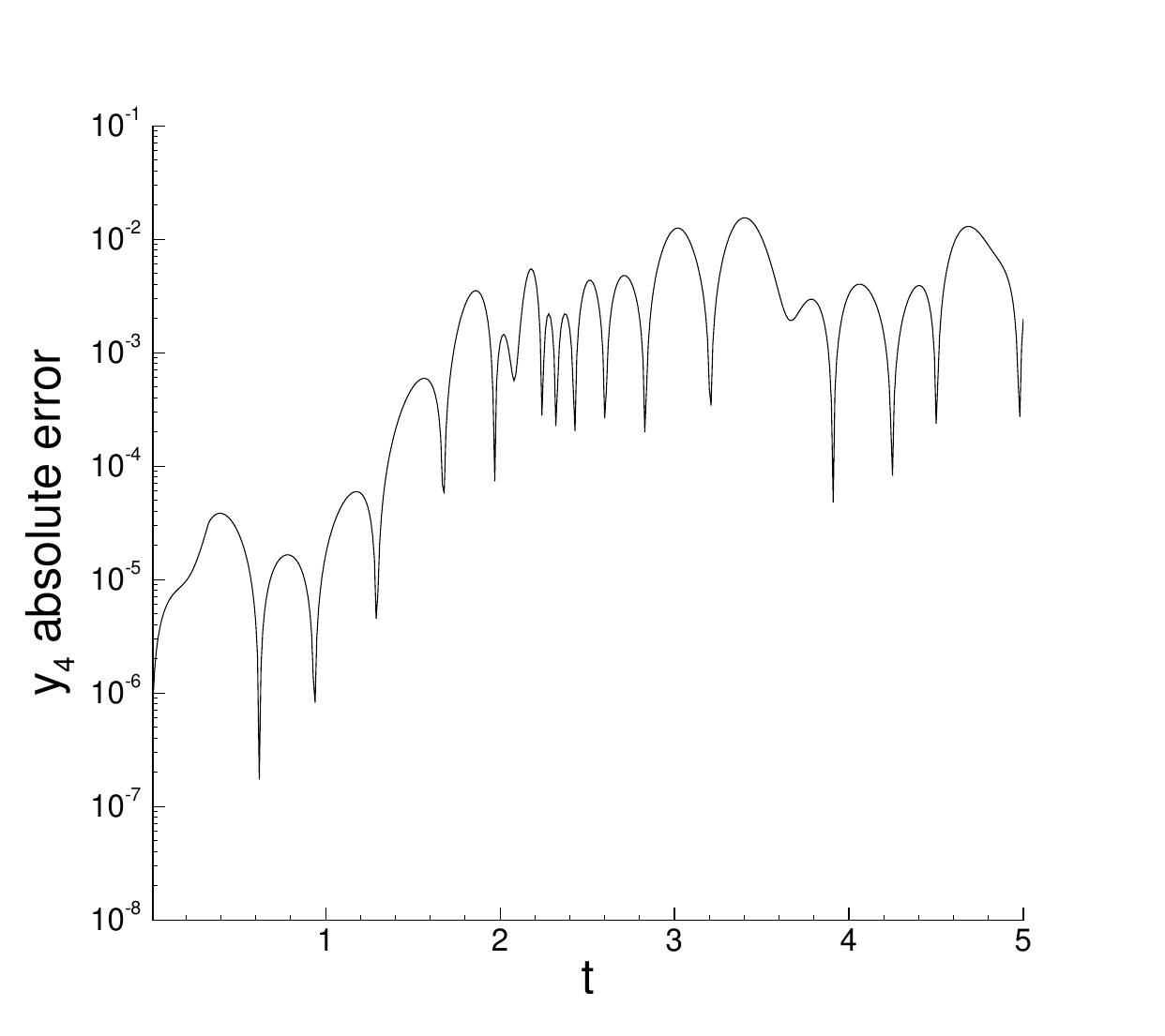}(f)
  }
  \caption{Lorenz96 system ($t_f=5$): Comparison 
    of (a) $y_1$, (b) $y_3$, and (c) $y_4$ between the NN-Exp-S1 solution and
    the reference solution. Absolute errors for
    $y_1(t)$, $y_3(t)$ and $y_4(t)$ of the NN-Exp-S1 solution.
    NN-Exp-S1: training domain $(y_{01},y_{02},y_{03},y_{04},y_{05},\xi)\in[-9,12]\times[-7,10]\times[-5.5,11]\times[-6,7]\times[-4,12]\times[0,0.011]$, $4$ uniform sub-domains
    along $y_{01}$;
    NN: $[6,1000,5]$ on each sub-domain,
    $R_m=0.15$, $Q=1500$, time integration for $t\in[0,5]$; See Table~\ref{tab_6}
    for the other parameter values.
    Reference solution: obtained by scipy DOP853, with
    absolute tolerance $10^{-16}$ and relative tolerance $10^{-13}$,
    dense output on points corresponding to $\Delta t=0.01$.
  }
  \label{fg_27}
\end{figure}

Due to the chaotic nature of the system, computing the error by comparing
the solution histories in long-term evolution becomes less meaningful.
We therefore consider a shorter history ($t\in[0,5]$) and compare
the current NN algorithms with the scipy methods.
Figure~\ref{fg_27} compares the solution histories for $y_1$, $y_3$ and $y_4$
from NN-Exp-S1 and  a reference solution computed by the scipy
DOP853 method (with a sufficiently small tolerance). It also shows
the absolute errors for these variables from the NN-Exp-S1 solution.
The NN algorithm has  captured the solution accurately.

\begin{figure}
  \centerline{
    \includegraphics[width=2in]{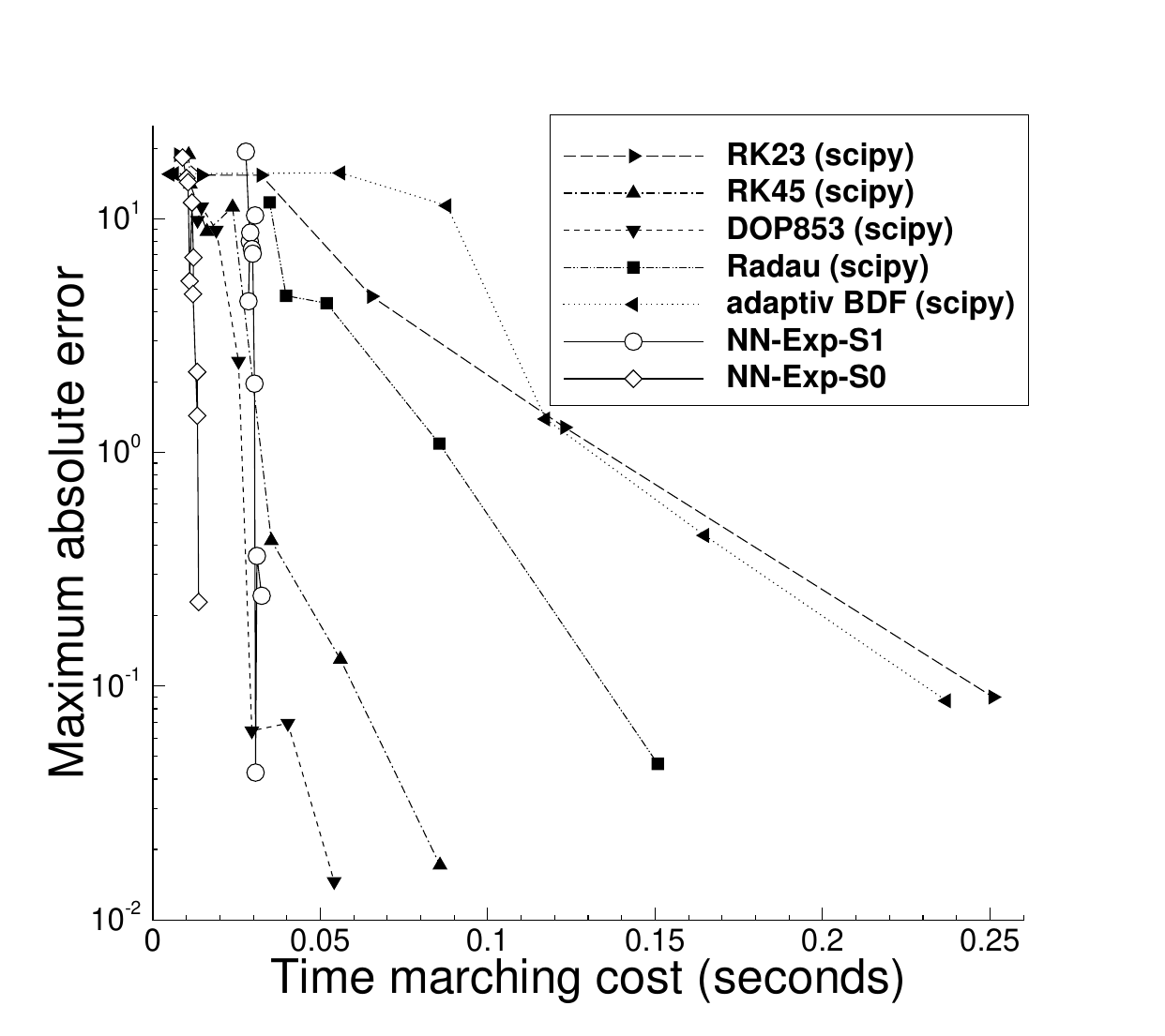}(a)
    \includegraphics[width=2in]{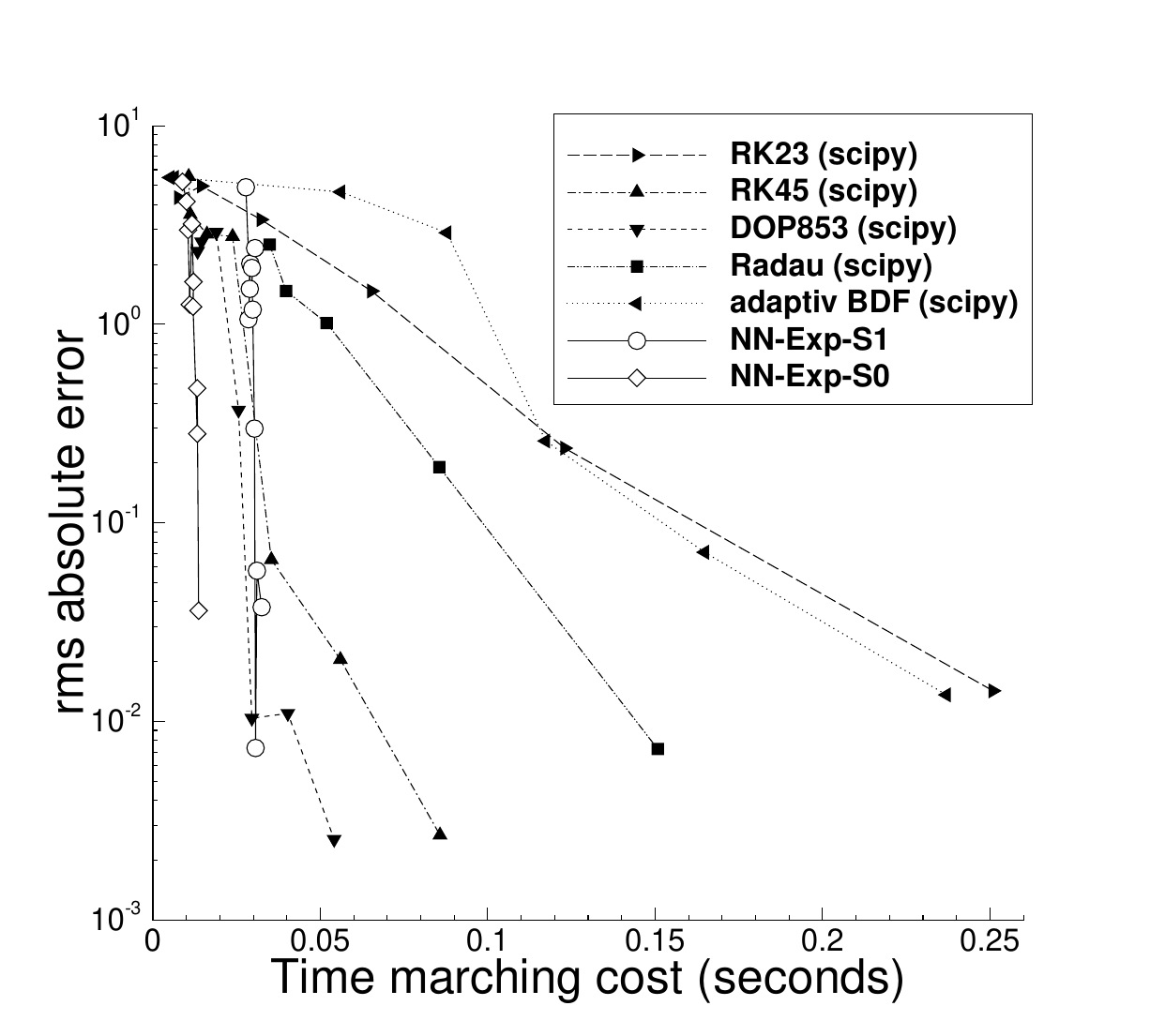}(b)
  }
  \caption{Lorenz96 system ($t_f=5$): Comparison of (a) the maximum
    and (b) the rms time-marching errors versus the time marching cost (wall time)
    between the NN algorithms (NN-Exp-S0, NN-Exp-S1) and the scipy methods.
    NN algorithms: training domain $(y_{01},y_{02},y_{03},y_{04},y_{05},\xi)\in[-9,12]\times[-7,10]\times[-5.5,11]\times[-6,7]\times[-4,12]\times[0,0.011]$, $4$ uniform sub-domains
    along $y_{01}$;
    NN: $[6,M,5]$ ($M$ varied);
    $R_m=0.05$ and $Q=2000$ for NN-Exp-S0, and $R_m=0.15$ and $Q=1500$ for NN-Exp-S1.
    See Table~\ref{tab_6} for the other parameter values.
    Scipy methods: 
    absolute tolerance $10^{-16}$, data points corresponding to different relative tolerance
    values, dense output on points corresponding to
    $\Delta t=0.01$ for $t\in[0,5]$.
  }
  \label{fg_28}
\end{figure}

We next compare the computational performance (accuracy versus cost)
between the NN algorithms and
the scipy methods for the Lorenz96 system.
Figure~\ref{fg_28} shows the maximum and rms time-marching errors
($t\in[0,5]$) as a function of the time-marching time
obtained by the current NN-Exp-S0 and NN-Exp-S1 algorithms
and the scipy methods. The errors are computed with respect to a reference
solution obtained by the scipy DOP853 method with a sufficiently small
tolerance. The parameter values for these tests are given in
the figure caption or in Table~\ref{tab_6}.
The data points for NN-Exp-S0 and NN-Exp-S1 correspond to different
$M$ in the $\varphi$-subnet architecture $[6,M,5]$, while those for
the scipy methods correspond to different tolerance values.
Among the scipy methods DOP853 shows the best performance, followed by RK45,
Radau and the other methods.  The performance of NN-Exp-S1 and scipy
DOP853 is close to each other. NN-Exp-S0 exhibits a better performance
than NN-Exp-S1, DOP853, and the other scipy methods for this problem.

\section{Concluding Remarks}
\label{sec_summary}


We have presented a method for learning the exact time integration
algorithm for non-autonomous and autonomous systems  based on ELM-type randomized
neural networks.
The trained NN  serves as a time-marching algorithm for solving
the initial value problems, for arbitrary initial data or step sizes
from some domain.

For a given non-autonomous (including autonomous) system,
the exact time integration algorithm can be represented
by an algorithmic function in higher dimensions,
which satisfies an associated system of
partial differential equations  with corresponding boundary conditions.
We learn this high-dimensional algorithmic function based on a physics informed approach,
by solving the associated PDE system using ELM-type randomized NNs.
We have presented several explicit and implicit formulations for
the algorithmic function, which accordingly lead to explicit and implicit
time integration algorithms, and discussed how to train the NN  by
the nonlinear least squares method.
For more challenging problems,  we find that it is
crucial to incorporate 
domain decomposition into the algorithm learning. 
In this case the algorithmic function is represented by a collection
of local randomized NNs, one for each sub-domain.
Importantly, the local NNs  for different sub-domains
are not coupled, due to the nature of the system of equations
for the algorithmic function, and thus they can be trained individually or in parallel.

Choosing the domain for the input variables  for
training the neural network is crucial to the accuracy of
the resultant time integration algorithm and to the efficiency of
network training. The training domain should be sufficiently large,
in order to adequately cover the region of interest in the
phase space. An overly large training domain,
on the other hand, increases the difficulty in network training.
When the right hand side of the non-autonomous system exhibits a periodicity
with respect to any of its arguments, while the solution itself to
the problem is not periodic, we show that in this case the algorithmic function
is either periodic, or when it is not, shows a well-defined relation
for different periods.
This property about the algorithmic function
can greatly simply the choice of the training domain and
the algorithm learning.


We have performed extensive numerical experiments with  benchmark problems
(non-stiff, stiff, or chaotic systems) to
evaluate the performance of the learned NN algorithms, and in particular
to compare them with the leading traditional algorithms
from the scipy library (DOP853, RK45, RK23, Radau, BDF).
We have made the following observations:
\begin{itemize}

\item
  The learned NN  algorithms produce highly accurate
  solutions with a high computational efficiency.
  Their solution accuracy  increases nearly exponentially,
  while their time-marching cost
  grows only quasi-linearly, as the
  number of degrees of freedom
  (training collocation points,  training parameters)
  in the ELM network increases.

\item
  The implicit NN algorithms are not as competitive as the explicit ones,
  in terms of both the accuracy and the time-marching cost.

\item
  Among the explicit NN algorithms, the accuracy increases from NN-Exp-S0
  to NN-Exp-S1, and to NN-Exp-S2, under comparable simulation
  parameters (training data points and training parameters).
  The reverse is true in terms of the time-marching cost.

\item
  The NN-Exp-S0 algorithm is observed to work well for both non-stiff and stiff
  problems. NN-Exp-S1 appears also able to work with both types of problems, but
  not as well as NN-Exp-S0 for stiff ones. On the other hand,  NN-Exp-S2 encounters
  difficulties for stiff problems (e.g.~failure to converge in NN training).

\item
  The computational performance of the learned NN algorithms is highly competitive
  compared with the traditional algorithms from scipy. The NN-Exp-S0 algorithm
  outperforms the best of the scipy methods in most of the tests,
  achieving superior accuracy under a comparable time-marching cost or the same
  accuracy under a lower time-marching cost.

\end{itemize}
The numerical results demonstrate that the learned NN  algorithms are
accurate, efficient, and computationally competitive.

%


An outstanding problem with the learned NN algorithms is the current lack of an
effective adaptive-stepping strategy, which would be crucial for
efficiently solving stiff problems.
The strategy described in Section~\ref{sec_vdp} for
the stiff van der Pol oscillator problem takes advantage of
domain decomposition and employs different $h_{\max}$ on different sub-domains
for training the ELM  network. During time integration,
the time step size is adapted depending on the sub-domain that
the current solution approximation falls in. This is a useful strategy,
but it is only quasi-adaptive. That the NN algorithms fail to
outperform the scipy Radau (or BDF) method  for the stiff van der Pol
oscillator case is largely due to the lack of an effective adaptive-stepping
strategy. This is an important research problem and will be pursued
in a future endeavor.


\section*{Acknowledgement}
This work was partially supported by
NSF (DMS-2012415). 

\section*{Appendix: Proof of Theorems from Section~\ref{sec_a221}}

\noindent\underline{\bf Proof of Lemma~\ref{lem_1}:}\\
Let $z(t)=y(t+T)=y(\eta)$, where $\eta=t+T$, for $t\in\mbb R$. Then
\begin{align}\label{eq_64}
  &
  \frac{dz}{dt} = \frac{dy}{d\eta} = f(y(\eta),\eta) = f(z(t),t+T) = f(z,t),
\end{align}
where we have used equations~\eqref{eq_1} and~\eqref{eq_a9}.
If $y(t_0)=y(t_0+T))$, then
\begin{align}\label{eq_65}
  z(t_0) = y(t_0+T) = y(t_0)=y_0.
\end{align}
Comparing the initial value problem consisting of~\eqref{eq_64}--\eqref{eq_65} and the problem~\eqref{eq_1},
and in light of the uniqueness under Assumption~\ref{ass_1},
we conclude that $z(t)=y(t)$ for all $t\in\mbb R$.
It follows that $y(t)$ is a periodic function with a period $T$.

\vspace{8pt}
\noindent\underline{\bf Proof of Theorem~\ref{thm_a1}:}\\
Let $\Psi(y_0,t_0,\xi)=\psi(y_0,t_0+T,\xi)=\psi(y_0,\eta,\xi)$,
where $\eta=t_0+T$, for
$(y_0,t_0,\xi)\in\mbb R^n\times\mbb R\times[0,h_{\max}]$.
Then
\begin{subequations}\label{eq_66}
  \begin{align}
    &
  \frac{\partial\Psi}{\partial\xi}
  = \left.\frac{\partial\psi}{\partial\xi}\right|_{(y_0,\eta,\xi)}
  = f(\psi(y_0,\eta,\xi),\eta+\xi)
  = f(\Psi(y_0,t_0,\xi),t_0+\xi+T) = f(\Psi,t_0+\xi), \\
  &
  \Psi(y_0,t_0,0) = \psi(y_0,t_0+T,0) = y_0,
\end{align}
\end{subequations}
where we have used equation~\eqref{eq_9} and the periodicity of $f(y,t)$.
Comparing the problems~\eqref{eq_66} and~\eqref{eq_9},
and in light of the uniqueness under Assumption~\ref{ass_1},
we conclude that $\Psi(y_0,t_0,\xi)=\psi(y_0,t_0,\xi)$
for all $(y_0,t_0,\xi)\in\mbb R^n\times\mbb R\times[0,h_{\max}]$.
It follows that $\psi(y_0,t_0,\xi)$ is periodic with respect to
$t_0$ with a period $T$
on $(y_0,t_0,\xi)\in\mbb R^n\times\mbb R\times[0,h_{\max}]$.

\vspace{8pt}
\noindent\underline{\bf Proof of Theorem~\ref{thm_a2}:}\\
Let $\Psi(y_0,t_0,\xi)=\psi(y_0+L_i\mbs e_i,t_0,\xi)=\psi(z_0,t_0,\xi)$,
where $z_0=y_0+L_i\mbs e_i\in\mbb R^n$, for
all $(y_0,t_0,\xi)\in\mbb R^n\times\mbb R\times[0,h_{\max}]$.
Then
\begin{subequations}\label{eq_67}
  \begin{align}
    &
    \frac{\partial\Psi}{\partial\xi}
    = \left.\frac{\partial\psi}{\partial\xi}\right|_{(z_0,t_0,\xi)}
    = f(\psi(z_0,t_0,\xi),t_0+\xi)
    = f(\Psi(y_0,t_0,\xi),t_0+\xi), \\
    &
    \Psi(y_0,t_0,0) = \psi(y_0+L_i\mbs e_i,t_0,0) = y_0+L_i\mbs e_i,
  \end{align}
\end{subequations}
where we have used equation~\eqref{eq_9}.

Let $\Phi(y_0,t_0,\xi) = \psi(y_0,t_0,\xi) + L_i\mbs e_i$,
for all $(y_0,t_0,\xi)\in\mbb R^n\times\mbb R\times[0,h_{\max}]$.
Then
\begin{subequations}\label{eq_68}
  \begin{align}
    &
    \frac{\partial\Phi}{\partial\xi}
    = \left.\frac{\partial\psi}{\partial\xi}\right|_{(y_0,t_0,\xi)}
    = f(\psi(y_0,t_0,\xi),t_0+\xi)
    = f(\psi(y_0,t_0,\xi)+L_i\mbs e_i,t_0+\xi) \notag \\
    &\qquad\qquad\qquad\qquad
    = f(\Phi(y_0,t_0,\xi),t_0+\xi), \\
    &
    \Phi(y_0,t_0,0) = \psi(y_0,t_0,0)+L_i\mbs e_i = y_0+L_i\mbs e_i,
  \end{align}
\end{subequations}
where we have used equation~\eqref{eq_9} and
the periodicity condition~\eqref{eq_a10} for $f(y,t)$.

Comparing the problem~\eqref{eq_67} with the problem~\eqref{eq_68},
and in light of the uniqueness under Assumption~\ref{ass_1},
we conclude that $\Psi(y_0,t_0,\xi)=\Phi(y_0,t_0,\xi)$
for all $(y_0,t_0,\xi)\in\mbb R^n\times\mbb R\times[0,h_{\max}]$.
This leads to equation~\eqref{eq_a11}.

\bibliographystyle{plain}
\bibliography{tstep,varpro,elm2,elm,mypub,dnn1,dnn}

\end{document}